
\def\service{T}
\catcode`\@=11
\def\unredoffs{\voffset=11mm \hoffset=0.5mm}

%
\newbox\leftpage \newdimen\fullhsize \newdimen\hstitle \newdimen\hsbody
\newdimen\hdim
\tolerance=400\pretolerance=800
%
%
\newif\ifsmall \smallfalse
\newif\ifdraft \draftfalse
\newif\iffrench \frenchfalse
\newif\ifeqnumerosimple \eqnumerosimplefalse
\nopagenumbers
\headline={\ifnum\pageno=1\hfill\else\hfil{\headrm\folio}\hfil\fi}
\def\draftstart{
\magnification=1200 \unredoffs\hsize=130mm\vsize=190mm
\hsbody=\hsize \hstitle=\hsize 
\nolabels
\iffrench
\dicof
\else
\dicoa
\fi
}

\font\elevrm=cmr9

\newdimen\chapskip
\font\twbf=cmssbx10 scaled 1200
\font\ssbx=cmssbx10

\font\caprm=cmr9
\font\capit=cmti9
\font\capbf=cmbx9
\font\capsl=cmsl9
\font\capmi=cmmi9
\font\capex=cmex9
\font\capsy=cmsy9
\chapskip=17.5mm
\def\makeheadline{\vbox to 0pt{\vskip-22.5pt
\line{\vbox to8.5pt{}\the\headline}\vss}\nointerlineskip}
\font\tbfi=cmmib10
\font\tenbi=cmmib7
\font\fivebi=cmmib5
\textfont4=\tbfi
\scriptfont4=\tenbi
\scriptscriptfont4=\fivebi
\font\headrm=cmr10

\font\eightrm=cmr6
\font\sixrm=cmr5
\font\eightmi=cmmi6
\font\sixmi=cmmi5
\font\eightsy=cmsy6
\font\sixsy=cmsy5
\font\eightbf=cmbx6
\font\sixbf=cmbx5
\skewchar\capmi='177 \skewchar\eightmi='177 \skewchar\sixmi='177
\skewchar\capsy='60 \skewchar\eightsy='60 \skewchar\sixsy='60

\def\elevenpoint{
\textfont0=\caprm \scriptfont0=\eightrm \scriptscriptfont0=\sixrm
\def\rm{\fam0\caprm}
\textfont1=\capmi \scriptfont1=\eightmi \scriptscriptfont1=\sixmi
\textfont2=\capsy \scriptfont2=\eightsy \scriptscriptfont2=\sixsy
\textfont3=\capex \scriptfont3=\capex \scriptscriptfont3=\capex
\textfont\itfam=\capit \def\it{\fam\itfam\capit} 
\textfont\slfam=\capsl  \def\sl{\fam\slfam\capsl} 
\textfont\bffam=\capbf \scriptfont\bffam=\eightbf
\scriptscriptfont\bffam=\sixbf
\def\bf{\fam\bffam\capbf} 
\textfont4=\tbfi \scriptfont4=\tenbi \scriptscriptfont4=\tenbi
\normalbaselineskip=13pt
\setbox\strutbox=\hbox{\vrule height9.5pt depth3.9pt width0pt}
\let\big=\elevenbig \normalbaselines \rm}

\catcode`\@=11

\font\tenmsa=msam10
\font\sevenmsa=msam7
\font\fivemsa=msam5
\font\tenmsb=msbm10
\font\sevenmsb=msbm7
\font\fivemsb=msbm5
\newfam\msafam
\newfam\msbfam
\textfont\msafam=\tenmsa  \scriptfont\msafam=\sevenmsa
  \scriptscriptfont\msafam=\fivemsa
\textfont\msbfam=\tenmsb  \scriptfont\msbfam=\sevenmsb
  \scriptscriptfont\msbfam=\fivemsb

\def\hexnumber@#1{\ifcase#1 0\or1\or2\or3\or4\or5\or6\or7\or8\or9\or
	A\or B\or C\or D\or E\or F\fi }

\font\teneuf=eufm10
\font\seveneuf=eufm7
\font\fiveeuf=eufm5
\newfam\euffam
\textfont\euffam=\teneuf
\scriptfont\euffam=\seveneuf
\scriptscriptfont\euffam=\fiveeuf
\def\frak{\ifmmode\let\next\frak@\else
 \def\next{\Err@{Use \string\frak\space only in math mode}}\fi\next}
\def\goth{\ifmmode\let\next\frak@\else
 \def\next{\Err@{Use \string\goth\space only in math mode}}\fi\next}
\def\frak@#1{{\frak@@{#1}}}
\def\frak@@#1{\fam\euffam#1}

\edef\msa@{\hexnumber@\msafam}
\edef\msb@{\hexnumber@\msbfam}

\def\Bbb{\ifmmode\let\next\Bbb@\else
 \def\next{\errmessage{Use \string\Bbb\space only in math mode}}\fi\next}
\def\Bbb@#1{{\Bbb@@{#1}}}
\def\Bbb@@#1{\fam\msbfam#1}

\catcode`\@=12
\def\sla#1{\mkern-1.5mu\raise0.4pt\hbox{$\not$}\mkern1.2mu #1\mkern 0.7mu}
\def\Dbar{\mkern-1.5mu\raise0.4pt\hbox{$\not$}\mkern-.1mu {\rm D}\mkern.1mu}
\def\Abar{\mkern1.mu\raise0.4pt\hbox{$\not$}\mkern-1.3mu A\mkern.1mu}
\def\dicof{
\gdef\Resume{RESUME}
\gdef\Toc{Table des mati\`eres}
\gdef\soumisa{Soumis \`a:}
}
\def\dicoa{
\gdef\Resume{ABSTRACT}
\gdef\Toc{Table of Contents}
\gdef\soumisa{Submitted to}
}

\def\uniset{\rlap{\elevrm 1}\kern.15em 1}
\def\bkR{{\rm I\kern-.17em R}}
\def\bkC{{\rm \kern.24em
            \vrule width.05em height1.4ex depth-.05ex
            \kern-.26em C}}

\def\frac#1#2{{\textstyle{#1\over#2}}}

\def\leaderfill{\leaders\hbox to 1em{\hss.\hss}\hfill}
\def\saclay{\if S\service \spec \else \spht \fi}
\def\spht{
\centerline{Service de Physique Th\'eorique, CEA-Saclay}
\centerline{F-91191 Gif-sur-Yvette Cedex, FRANCE}}
\def\spec{
\centerline{Service de Physique de l'Etat Condens\'e, CEA-Saclay}
\centerline{F-91191 Gif-sur-Yvette Cedex, FRANCE}}
\def\logo{
\if S\service 
\font\sstw=cmss10 scaled 1200
\font\ssx=cmss8
\vtop{\hsize 9cm
{\sstw {\twbf P}hysique de l'{\twbf E}tat {\twbf C}ondens\'e \par}
\ssx SPEC -- DRECAM -- DSM\par
\vskip 0.5mm
\sstw CEA -- Saclay \par
}
\else 
\vtop{\hsize 9cm
}
\fi }
\catcode`\@=11
\def\deqalignno#1{\displ@y\tabskip\centering \halign to
\displaywidth{\hfil$\displaystyle{##}$\tabskip0pt&$\displaystyle{{}##}$
\hfil\tabskip0pt &\quad
\hfil$\displaystyle{##}$\tabskip0pt&$\displaystyle{{}##}$
\hfil\tabskip\centering& \llap{$##$}\tabskip0pt \crcr #1 \crcr}}
\def\deqalign#1{\null\,\vcenter{\openup\jot\m@th\ialign{
\strut\hfil$\displaystyle{##}$&$\displaystyle{{}##}$\hfil
&&\quad\strut\hfil$\displaystyle{##}$&$\displaystyle{{}##}$
\hfil\crcr#1\crcr}}\,}
\openin 1=\jobname.sym
\ifeof 1\closein1\message{<< (\jobname.sym DOES NOT EXIST) >>}\else%
\input\jobname.sym\closein 1\fi
\newcount\nosection
\newcount\nosubsection
\newcount\neqno
\newcount\notenumber
\newcount\figno
\newcount\tabno
\def\content{\jobname.toc}
\def\symbols{\jobname.sym}
\newwrite\toc
\newwrite\sym
\def\authorname#1{\centerline{\bf #1}\smallskip}
\def\address#1{ #1\medskip}
\newdimen\hulp
\def\maketitle#1{
\edef\oneliner##1{\centerline{##1}}
\edef\twoliner##1{\vbox{\parindent=0pt\leftskip=0pt plus 1fill\rightskip=0pt
plus 1fill
                     \parfillskip=0pt\relax##1}}
\setbox0=\vbox{#1}\hulp=0.5\hsize
                 \ifdim\wd0<\hulp\oneliner{#1}\else
                 \twoliner{#1}\fi}

\def\submitted#1{{\it {\soumisa} #1}\par}
\def\title#1{\gdef\titlename{#1}
\maketitle{
\twbf
{\titlename}}
\vskip3truemm\vfill
\nosection=0
\neqno=0
\notenumber=0
\figno=1
\tabno=1
\def\prefix{}
\def\eqprefix{}
\mark{\the\nosection}
\message{#1}
\immediate\openout\sym=\symbols
}
\def\preprint#1{\vglue-10mm
\line{ \logo \hfill {#1} }\vglue 20mm\vfill}
\def\abstract{\vfill\centerline{\Resume} \smallskip \begingroup\narrower
\elevenpoint\baselineskip10pt}
\def\endabstract{\par\endgroup \bigskip}
\def\mktoc{\centerline{\bf \Toc} \medskip\caprm
\parindent=2em
\openin 1=\jobname.toc
\ifeof 1\closein1\message{<< (\jobname.toc DOES NOT EXIST. TeX again)>>}%
\else\input\jobname.toc\closein 1\fi
 \bigskip}
\def\section#1\par{\vskip0pt plus.1\vsize\penalty-100\vskip0pt plus-.1
\vsize\bigskip\vskip\parskip
\message{ #1}
\ifnum\nosection=0\immediate\openout\toc=\content%
\edef\ecrire{\write\toc{\par\noindent{\ssbx\ \titlename}
\string\leaderfill{\noexpand\number\pageno}}}\ecrire\fi
\advance\nosection by 1\nosubsection=0
\ifeqnumerosimple
\else \xdef\eqprefix{\prefix\the\nosection.}\neqno=0\fi
\vbox{\noindent\bf\prefix\the\nosection\ #1}
\mark{\the\nosection}\bigskip\noindent
\xdef\ecrire{\write\toc{\string\par\string\item{\prefix\the\nosection}
#1
\string\leaderfill {\noexpand\number\pageno}}}\ecrire}

\def\appendix#1#2\par{\bigbreak\nosection=0
\notenumber=0
\neqno=0
\def\prefix{A}
\mark{\the\nosection}
\message{\appendixname}
\leftline{\ssbx APPENDIX}
\leftline{\ssbx\uppercase\expandafter{#1}}
\leftline{\ssbx\uppercase\expandafter{#2}}
\bigskip\noindent\nonfrenchspacing
\edef\ecrire{\write\toc{\par\noindent{{\ssbx A}\
{\ssbx#1\ #2}}\string\leaderfill{\noexpand\number\pageno}}}\ecrire}%

\def\subsection#1\par {\vskip0pt plus.05\vsize\penalty-100\vskip0pt
plus-.05\vsize\bigskip\vskip\parskip\advance\nosubsection by 1
\vbox{\noindent\it\prefix\the\nosection.\the\nosubsection\
\it #1}\smallskip\noindent
\edef\ecrire{\write\toc{\string\par\string\itemitem
{\prefix\the\nosection.\the\nosubsection} {#1}
\string\leaderfill{\noexpand\number\pageno}}}\ecrire
}
\def\note #1{\advance\notenumber by 1
\footnote{$^{\the\notenumber}$}{\sevenrm #1}}

\def\nolabels{\def\wrlabel##1{}\def\eqlabel##1{}\def\reflabel##1{}}
\def\writelabels{\def\wrlabel##1{\leavevmode\vadjust{\rlap{\smash%
{\line{{\escapechar=` \hfill\rlap{\sevenrm\hskip.03in\string##1}}}}}}}%
\def\eqlabel##1{{\escapechar-1\rlap{\sevenrm\hskip.05in\string##1}}}%
\def\reflabel##1{\noexpand\llap{\noexpand\sevenrm\string\string\string##1}}}
\global\newcount\refno \global\refno=1
\newwrite\rfile
\def\ref{[\the\refno]\nref}
\def\nref#1{\xdef#1{[\the\refno]}\writedef{#1\leftbracket#1}%
\ifnum\refno=1\immediate\openout\rfile=\jobname.ref\fi
\global\advance\refno by1\chardef\wfile=\rfile\immediate
\write\rfile{\noexpand\item{#1\ }\reflabel{#1\hskip.31in}\pctsign}\findarg}
\def\findarg#1#{\begingroup\obeylines\newlinechar=`\^^M\pass@rg}
{\obeylines\gdef\pass@rg#1{\writ@line\relax #1^^M\hbox{}^^M}%
\gdef\writ@line#1^^M{\expandafter\toks0\expandafter{\striprel@x #1}%
\edef\next{\the\toks0}\ifx\next\em@rk\let\next=\endgroup\else\ifx\next\empty%
\else\immediate\write\wfile{\the\toks0}\fi\let\next=\writ@line\fi\next\relax}}
\def\striprel@x#1{}
\def\em@rk{\hbox{}}

\def\addref#1{\immediate\write\rfile{\noexpand\item{}#1}} 
\def\listrefs{
\ifnum\refno=1 \else
\immediate\closeout\rfile\writestoppt\baselineskip=14pt%
\vskip0pt plus.1\vsize\penalty-100\vskip0pt plus-.1
\vsize\bigskip\vskip\parskip\centerline{{\bf References}}\bigskip%
{\frenchspacing%
\parindent=20pt\escapechar=` \input \jobname.ref\vfill\eject}%
\nonfrenchspacing
\fi}
\def\startrefs#1{\immediate\openout\rfile=\jobname.ref\refno=#1}
\def\xref{\expandafter\xr@f}\def\xr@f[#1]{#1}
\def\refs#1{[\r@fs #1{\hbox{}}]}
\def\r@fs#1{\ifx\und@fined#1\message{reflabel \string#1 is undefined.}%
\xdef#1{(?.?)}\fi \edef\next{#1}\ifx\next\em@rk\def\next{}%
\else\ifx\next#1\xref#1\else#1\fi\let\next=\r@fs\fi\next}
%
\newwrite\lfile
{\escapechar-1\xdef\pctsign{\string\%}\xdef\leftbracket{\string\{}
\xdef\rightbracket{\string\}}}

\def\writestop{\def\writestoppt{\immediate\write\lfile{\string\pageno%
\the\pageno\string\startrefs\leftbracket\the\refno\rightbracket%
\string\def\string\secsym\leftbracket\secsym\rightbracket%
\string\secno\the\secno\string\meqno\the\meqno}\immediate\closeout\lfile}}
\def\writestoppt{}\def\writedef#1{}
\def\eqnn{\global\advance\neqno by 1 \ifinner\relax\else%
\eqno\fi(\eqprefix\the\neqno)}
%
\def\eqnd#1{\global\advance\neqno by 1 \ifinner\relax\else%
\eqno\fi(\eqprefix\the\neqno)\eqlabel#1
{\xdef#1{($\eqprefix\the\neqno$)}}
\edef\ewrite{\write\sym{\string\def\string#1{($\eqprefix%
\the\neqno$)}}%
}\ewrite%
}
%
\def\eqna#1{\wrlabel#1\global\advance\neqno by1
{\xdef #1##1{\hbox{$(\eqprefix\the\neqno##1)$}}}
\edef\ewrite{\write\sym{\string\def\string#1{($\eqprefix%
\the\neqno$)}}%
}\ewrite%
}
\def\em@rk{\hbox{}}
\def\xeqn{\expandafter\xe@n}\def\xe@n(#1){#1}
\def\xeqna#1{\expandafter\xe@na#1}\def\xe@na\hbox#1{\xe@nap #1}
\def\xe@nap$(#1)${\hbox{$#1$}}
\def\eqns#1{(\e@ns #1{\hbox{}})}
\def\e@ns#1{\ifx\und@fined#1\message{eqnlabel \string#1 is undefined.}%
\xdef#1{(?.?)}\fi \edef\next{#1}\ifx\next\em@rk\def\next{}%
\else\ifx\next#1\xeqn#1\else\def\n@xt{#1}\ifx\n@xt\next#1\else\xeqna#1\fi
\fi\let\next=\e@ns\fi\next}
\def\fig{fig.~\the\figno\nfig}
\def\nfig#1{\xdef#1{\the\figno}%
\immediate\write\sym{\string\def\string#1{\the\figno}}%
\global\advance\figno by1}%
\def\xfig{\expandafter\xf@g}\def\xf@g fig.\penalty\@M\ {}%
\def\figs#1{figs.~\f@gs #1{\hbox{}}}%
\def\f@gs#1{\edef\next{#1}\ifx\next\em@rk\def\next{}\else%
\ifx\next#1\xfig #1\else#1\fi\let\next=\f@gs\fi\next}%
\long\def\figure#1#2#3{\midinsert
#2\par
{\elevenpoint
\setbox1=\hbox{#3}
\ifdim\wd1=0pt\centerline{{\bf Figure\ #1}\hskip7.5mm}%
\else\setbox0=\hbox{{\bf Figure #1}\quad#3\hskip7mm}
\ifdim\wd0>\hsize{\narrower\noindent\unhbox0\par}\else\centerline{\box0}\fi
\fi}
\wrlabel#1\par
\endinsert}
\def\tab{table~\uppercase\expandafter{\romannumeral\the\tabno}\ntab}
\def\ntab#1{\xdef#1{\the\tabno}
\immediate\write\sym{\string\def\string#1{\the\tabno}}
\global\advance\tabno by1}
\long\def\table#1#2#3{\topinsert
#2\par
{\elevenpoint
\setbox1=\hbox{#3}
\ifdim\wd1=0pt\centerline{{\bf Table
\uppercase\expandafter{\romannumeral#1}}\hskip7.5mm}%
\else\setbox0=\hbox{{\bf Table
\uppercase\expandafter{\romannumeral#1}}\quad#3\hskip7mm}
\ifdim\wd0>\hsize{\narrower\noindent\unhbox0\par}\else\centerline{\box0}\fi
\fi}
\wrlabel#1\par
\endinsert}
\catcode`@=12
\def\draftend{\immediate\closeout\sym\immediate\closeout\toc
}
\draftstart
\preprint{T93/039}
\title{Connection between the harmonic analysis on the sphere and the
harmonic analysis on the one-sheeted hyperboloid: an analytic continuation
viewpoint}
\authorname{J. Bros}
\address{\saclay}
\authorname{G.A. Viano}
\address{\centerline{Istituto Nazionale Fisica Nucleare, Sezione di Genova}
\centerline{Dipartimento di Fisica dell'Universit\'a di Genova,}
\centerline{I-16146 Genova, ITALY}
}
\abstract
In a previous paper $\lbrack$B,V-1$\rbrack$, an algebra of
holomorphic
\lq\lq perikernels\rq\rq\ on a complexified hyperboloid $ X^{(c)}_{d-1} $ (in
$ \Bbb C^d) $ has been
introduced; each perikernel $ {\cal K} $ can be seen as the analytic
continuation of a kernel $ {\bf K} $ on the unit sphere $ \Bbb S^{d-1} $ in an
appropriate \lq\lq cut-domain\rq\rq , while the jump of $ {\cal K} $ across
the
corresponding \lq\lq cut\rq\rq\ defines a Volterra kernel $ K $ (in the sense
of
J. Faraut $\lbrack$Fa-1$\rbrack$) on the one-sheeted hyperboloid $ X_{d-1} $
(in $ \Bbb R^d). $ \par
In the present paper, we obtain results of harmonic
analysis for classes of perikernels which are invariant under the
group $ {\rm SO}(d,\Bbb C) $ and of moderate growth at infinity. For each
perikernel $ {\cal K} $ in such a class, the Fourier-Legendre coefficients of
the corresponding kernel $ {\bf K} $ on $ \Bbb S^{d-1} $ admit a carlsonian
analytic
interpolation $ \tilde  F(\lambda) $ in a half-plane, which is the \lq\lq
spherical
Laplace transform\rq\rq\ of the associated
Volterra kernel $ K $ on $ X_{d-1}. $ Moreover, the composition law $ {\cal K
}= {\cal K}_1\ast^{( c)}{\cal K}_2 $
for perikernels (interpreted in terms of convolutions for the
classes considered) is turned into the corresponding product of
the Laplace transforms $ (\tilde  F(\lambda)  = \tilde  F_1(\lambda)
\cdot\tilde  F_2(\lambda) ). $ \par
Our method for proving these results consists in showing
that, for any dimension $ d, $ $ d\geq 2, $ they can be reduced to special
properties of Fourier-Laplace analysis in $ \Bbb C: $ for the case $ d=2 $
(circle and hyperbola), this reduction is straightforward via an
appropriate parametrization of $ X^{(c)}_1; $ for the general case $ d>2, $ it
requires the use of a \lq\lq Radon-Abel transformation\rq\rq\ involving
integration on \lq\lq complex horocycles\rq\rq\ of $ X^{(c)}_{d-1}. $ \par
This work is closely connected with two studies of
mathematical physics: the first one (at the origin of the present
work) concerns the complex angular momentum analysis in quantum
field theory $\lbrack$Fa,V$\rbrack$, $\lbrack$B,V-2$\rbrack$; the second one
(developed as an
application of the present results) concerns the setting of
quantum field theory on a hyperboloid-shaped space-time manifold
(i.e. the de Sitter universe) $\lbrack$B$\rbrack$ $\lbrack$B,G,M$\rbrack$.
\endabstract
\vfill
\submitted{J. de Math. Pures et Appliqu\'ees}
\eject
\eject
\input definit.tex
\magnification=1200
\baselineskip=20.0truept
\pageno=2
\noindent \input mssymb
\overfullrule=0pt
{\bf 1. INTRODUCTION AND DESCRIPTION OF THE RESULTS} \par
\smallskip
The group theoretical generalizations of standard Fourier
analysis have led to two separate categories of results: while
the notion of Fourier series has found its natural extension in
the theory of compact groups (see e.g. $\lbrack$Ru$\rbrack$), an analogue of
the Fourier
integral transformation has been introduced for a certain class
of non-compact semi-simple Lie groups $\lbrack$Fu$\rbrack$. One can then think
as
an attractive idea that, for couples of compact and non-compact
groups admitting the same complexified group (e.g. $ {\rm SO}(d;\Bbb R) $
and $ {\rm SO}(1,d-1)), $ the two previous categories of results should be
connected in some sense by \lq\lq an analytic continuation procedure\rq\rq .
It is our purpose to show that, at least at the level of
homogeneous spaces and for suitable classes of analytic functions
on these spaces, such a connection can in fact be displayed. In
the present paper, we will only treat elementary examples of
associated homogeneous spaces namely the sphere $ \Bbb S^{d-1} $ and the
one-sheeted hyperboloid $ X_{d-1} $ (affiliated respectively with the
groups $ {\rm SO}(d;\Bbb R) $ and $ {\rm SO}(1,d-1)), $ which we will consider
as
real submanifolds of the same complexified hyperboloid $ X^{(c)}_{d-1} $ in $
\Bbb C^d. $
However, a similar analysis might be envisaged for other couples
$ (S,X) $ of (compact and non-compact) homogeneous spaces embedded in
a common complexified space $ X^{(c)}. $ \par
\smallskip
The basic feature which appears in this study is {\sl \lq\lq dual
analyticity\rq\rq\/} , namely the fact that the classes of analytic
functions $ {\cal F} $ considered on $ X^{(c)} $ are in correspondence with
classes
of analytic functions $ \tilde  F $ on a space $ \tilde  X $ of \lq\lq complex
characters\rq\rq ,
each function $ \tilde  F $ being associated with a function $ {\cal F} $ via
the
following {\sl double connection\/}: \par
i) $ \tilde  F $ is a Carlson-type interpolation (in a suitable
domain of $ \tilde  X) $ of the set of Fourier coefficients of $ {\cal
F}_{\vert S}. $ \par
ii) $ \tilde  F $ is a Fourier-Laplace-type transform of an
appropriate \lq\lq jump\rq\rq\ $ F = \Delta{\cal F} $ of $ {\cal F}, $ defined
on the non-compact real
space $ X. $ \par
\smallskip
As a matter of fact, the two simplest examples of dual
analyticity are provided respectively by the Taylor series $ \sum^{ \infty}_{
n=0}a(n)\zeta^ n $
and by the Legendre series $ \sum^{ \infty}_{ n=0}a(n)P_n(\zeta) . $ The close
parallel which
exists between these two cases was noticed by Stein and Wainger
who gave (in $\lbrack$S,W$\rbrack$) a first systematic study of the
correspondence between the analyticity properties of the sum of
the series $ {\cal S}(\zeta) $ in the $ \zeta $-plane and of the coefficient
function $ a(\lambda) $
in the $ \lambda $-plane. In particular, in both cases, the analyticity of
the function $ {\cal S}(\zeta) $ in a cut-plane, cut along the positive real
axis from 1 to $ \infty , $ was shown to be equivalent (under certain
additional $ \gq\gq L^2 $-assumptions\rq\rq ) to the analyticity (and
appropriate
$ \gq\gq L^2 $-property\rq\rq ) of the function $ a(\lambda) $ in the
half-plane $ \Bbb C^{(-1/2)}_+= \left\{ \lambda ; {\rm Re} \ \lambda  > - {1
\over 2} \right\} . $
For the Taylor series, the discovery of properties of this type
traces back to a result of Leroy $\lbrack$L$\rbrack$ in 1900, while for the
Legendre series, various partial results were previously obtained
in connection with problems of physics, in particular by Regge
$\lbrack$R$\rbrack$, Froissart $\lbrack$F$\rbrack$, Gribov
$\lbrack$G$\rbrack$, Martin $\lbrack$M$\rbrack$ and Khuri $\lbrack$K$\rbrack$,
after
the pioneering works of Poincar\'e $\lbrack$P$\rbrack$, Watson
$\lbrack$Wa$\rbrack$ and Sommerfeld
$\lbrack$So$\rbrack$. If the interplay of rotational invariance properties
(under
$ {\rm SO}(3;\Bbb R)) $ and of a certain background of analyticity
properties had clearly been at the origin of all these works on the
Legendre series, it had remained obscure, however, which geometrical
interpretation could be given to the following integral relation
expressing the coefficient function $ a(\lambda) $ in terms of the jump $
\underline{{\rm f}}  = \Delta{\cal S} $
of $ {\cal S}(\zeta) $ across the cut $ [1,+\infty[ $ (for all $ \lambda $ in
$ \Bbb C^{(-1/2)}_+): $
$$ a(\lambda)  = {\rm cst} \times \int^{ \infty}_ 1Q_\lambda( \zeta)
\underline{{\rm f}}(\zeta) {\rm d} \zeta \ , \eqno (1.1) $$
with the help of the Legendre function of the second kind $ Q_\lambda( \zeta)
. $
This, formula, introduced at first in $\lbrack$F$\rbrack$ and
$\lbrack$G$\rbrack$, and
independently in the more rigorous study of $\lbrack$S,W$\rbrack$, was
obtained in
these works by a direct use of the properties of the Legendre
functions of the first and second kinds and of their
interconnection. \par
It is only in later works by J. Faraut and G.A. Viano
($\lbrack$Fa-1$\rbrack$, $\lbrack$Fa,V$\rbrack$) about the
harmonic analysis on the one-sheeted hyperboloid that the
integral transformation expressed by Eq.(1.1) was given a
geometrical meaning; in fact, it was reobtained there from the
definition of the \lq\lq spherical Laplace transform\rq\rq\ of an $ {\rm
SO}(1,2)
$-invariant \lq\lq Volterra kernel\rq\rq\ on the one-sheeted hyperboloid. \par
Here, we shall introduce a relationship based on analytic
continuation between the harmonic analysis on the one-sheeted
hyperboloid $ X_{d-1} = \left\{ z \in  \Bbb R^{d^{\ }}; \right. $ $ z=
\left(z^{(0)},...,z^{(d-1)} \right); $ $ z^2\equiv z^{(0)^2}-z^{(1)^2}-...
\allowbreak \left.-z^{(d-1)^2}=-1 \right\} $ and the harmonic analysis on the
sphere $ \Bbb S^{d-1}, $ identified to the submanifold $ \left( {\rm Re} \
z^{(0)}=0, \right. $ $ \left. {\rm Im} \ z^{(1)}=...= {\rm Im} \ z^{(d-1)}=0
\right) $
of the corresponding complexified hyperboloid $ X^{(c)}_{d-1}. $ In the
special cases $ d=2 $ and $ d=3, $ our results will provide a geometrical
interpretation of properties of dual analyticity analogous to
those obtained in $\lbrack$S,W$\rbrack$ concerning respectively the Taylor
series
and the Legendre series; in the general case, similar properties
of dual analyticity for series of the form $ \sum^{ }_ na(n)P^{(d)}_n $ (the $
P^{(d)}_n $
being \lq\lq generalized Legendre polynomials\rq\rq ) will be derived and
interpreted in terms of the appropriate Laplace-type
transformation on $ X_{d-1}. $ \par
A Volterra kernel $ K \left(z,z^{\prime} \right) $ on $ X_{d-1} $ is (in the
sense of J.
Faraut) a kernel on $ X_{d-1} $ whose support is contained in the region
$ \Gamma^ +_{d-1}= \left\{ \left(z,z^{\prime} \right) \in  X_{d-1}\times
X_{d-1}; \right. $ $ \left.z \grsim  z^{\prime} \right\} , $ where the
ordering relation $ z \grsim  z^{\prime} $
in $ X_{d-1} $ is induced by the corresponding minkowskian ordering
relation in $ \Bbb R^d, $ namely $ z \grsim  0 $ if $ z^2 \equiv
z^{(0)^2}-z^{(1)^2}-...-z^{(d-1)^2}\geq 0, $
with $ z^{(0)}\geq 0. $ In a previous paper $\lbrack$B,V-1$\rbrack$, we
introduced the
following type of relationship between Volterra kernels on $ X_{d-1} $
and kernels on the sphere $ \Bbb S^{d-1}, $ considered as embedded in $
X^{(c)}_{d-1}; $
we considered triplets $ ({\cal K}, {\bf K} ,K) $ on the complexified
hyperboloid $ X^{(c)}_{d-1} $
such that: \par
\smallskip
i) $ {\cal K} \left(z,z^{\prime} \right) $ is an analytic function, called
\lq\lq perikernel\rq\rq ,
whose domain is $ X^{(c)}_{d-1}\times X^{(c)}_{d-1} \backslash  \Sigma_ \mu ,
$ where $ \Sigma_ \mu $ is the following subset
of real codimension 1 (called a \lq\lq cut\rq\rq ):
$$ \Sigma_ \mu  = \left\{ \left(z,z^{\prime} \right)\in X^{(c)}_{d-1}\times
X^{(c)}_{d-1}\ ;\ \ \left(z-z^{\prime} \right)^2=\rho \ ;\ \rho  \geq  2(
{\rm cosh} \ \mu -1) \right\} $$
\par
ii) $ {\bf K} $ is a kernel on $ \Bbb S^{d-1} $ defined by taking the
restriction of $ {\cal K}, $ i.e. $ {\bf K}  = {\cal K}_{ \left\vert \Bbb
S^{d-1}\times \Bbb S^{d-1} \right.} $ \par
iii) $ K $ is a Volterra kernel on $ X_{d-1} $ defined via the
jump $ \Delta{\cal K} $ of $ {\cal K} $ across $ \Sigma_ \mu , $ i.e. $
K=\Delta{\cal K }_{ \left\vert \Gamma^ +_{d-1} \right.}. $ \par
\smallskip
In the present paper, we shall only consider the class of
{\sl \lq\lq invariant triplets\rq\rq\/}\ $ ({\cal K}, {\bf K} ,K), $ such that
$ {\cal K}, $ $ {\bf K} , $ $ K $ are respectively
invariant under the (diagonal) action of the groups $ {\rm
SO}(1,d-1)^{(c)}\approx $ $ {\rm SO}(d;\Bbb C), $
$ {\rm SO}(d;\Bbb R), $ $ {\rm SO}(1,d-1) $ on the corresponding spaces $
X^{(c)}_{d-1}\times X^{(c)}_{d-1}, $ $ \Bbb S^{d-1}\times \Bbb S^{d-1}, $
$ X_{d-1}\times X_{d-1}. $ \par
By introducing the \lq\lq base point\rq\rq\ $ z_0=(0,...,0,1) $ in $
X^{(c)}_{d-1}, $
and the subgroups $ H^{(c)}, $ $ {\bf H} , $ $ H $ which are the stabilizers
of $ z_0 $
respectively in $ {\rm SO}(1,d-1)^{(c)}, $ $ {\rm SO}(d;\Bbb R) $ and $ {\rm
SO}(1,d-1), $ we can
identify an invariant triplet $ ({\cal K}, {\bf K} ,K) $ with a triplet of
functions
$ ({\cal F}, {\bf F} ,F), $ such that: \par
\smallskip
i) $ {\cal F}(z) = {\cal K} \left(z,z_0 \right) $ is an analytic function
defined in the
cut-domain $ D_\mu  = X^{(c)}_{d-1}\backslash X^\mu $ of $ X^{(c)}_{d-1}, $
where $ X^\mu = \left\{ z\in X^{(c)}_{d-1};\ z^{(d-1)}= {\rm cosh} \ v\ ;\
v\geq \mu \right\} ; $
moreover $ {\cal F}(z) $ is invariant under $ H^{(c)} $ and therefore only
depends
on the single complex variable $ z^{(d-1)} = {\rm cos} \ \theta , $ $ \theta $
denoting the
complex angle between $ z $ and $ z_0; $ we will then put: $ {\cal F}(z)
\equiv  \underline{f} \left(z^{(d-1)} \right), $
with $ \underline{f} $ analytic in the cut-plane $ \underline{D}_\mu  = \Bbb C
\backslash  [ {\rm cosh} \ \mu ,+\infty[ ; $ if $ \bar  \omega $
denotes the orthogonal projection of $ X^{(c)}_{d-1} $ onto the complex
$ z^{(d-1)} $-axis, we then have: $ \underline{D}_\mu  = \bar  \omega
\left(D_\mu \right) $ and $ {\cal F }\equiv  \bar  \omega^ \ast(
\underline{f}). $ \par
ii) $ {\bf F}(z) $ is an $ {\bf H} $-invariant function on $ \Bbb S^{(d-1)}, $
defined by $ {\bf F}  = {\cal F}_{ \left\vert \Bbb S^{(d-1)} \right.}; $ since
$ \bar  \omega \left(\Bbb S^{(d-1)} \right) = [-1,+1], $ one
then has: $ {\bf F}(z) = \underline{{\bf f}} \left(z^{(d-1)} \right) $ with $
\underline{{\bf f}}  = \underline{f}_{\vert[ -1,+1]}. $ \par
iii) $ F $ is an $ H $-invariant function on $ X_{d-1} $ with support
contained in the set $ X^\mu_ += \left\{ z\in X_{d-1}; \right. $ $ z \grsim
z_0, $ $ \left.z^{(d-1)}\geq  {\rm cosh}  \mu \right\} , $ such
that
$ F_{ \left\vert X^\mu_ + \right.} = \Delta{\cal F}_{ \left\vert X^\mu_ +
\right.}, $ where $ \Delta{\cal F} $ is the jump of $ {\cal F} $ across $
X^\mu ; $ since $ \bar  \omega \left(X^\mu_ + \right)=[ {\rm cosh} \ \mu ,\
+\infty[ , $
one then has: $ F(z) = \underline{{\rm f}} \left(z^{(d-1)} \right), $ where $
\underline{{\rm f}} $ denotes the jump $ \Delta \underline{f} $ of $
\underline{f} $
across $ \Bbb R, $ whose support is contained in the \lq\lq cut\rq\rq\ $ [
{\rm cosh} \ \mu ,\ +\infty[ . $ \par
\smallskip
We will then say that $ ({\cal F}, {\bf F} ,F) $ is an invariant
triplet with domain $ D_\mu $ on $ X^{(c)}_{d-1}. $ \par
\smallskip
For $ d=3, $ the invariant triplets $ ({\cal F}, {\bf F} ,F) $ provide a
natural
geometrical interpretation for the Legendre series considered in
$\lbrack$S,W$\rbrack$; in fact, each function $ {\cal S}( {\rm cos} \ \theta)
= \sum^{ \infty}_{ \ell =0}a(\ell) P_{\ell}( {\rm cos} \ \theta) $
which admits an analytic continuation $ \underline{f}( {\rm cos} \ \theta) $
in the cut-plane $ \underline{D}_\mu $ $ (\mu \geq 0) $
can be identified with an $ {\bf H} $-invariant function $ {\bf F}(z) $ on the
sphere
$ \Bbb S^2 $ admitting the analytic continuation $ {\cal F }= \bar  \omega^
\ast( \underline{f}) $ in the
domain $ D_\mu $ of $ X^{(c)}_2; $ the jump $ \underline{{\rm f}}  = \Delta
\underline{f} $ of $ \underline{f} $ across $ \Bbb R $ is then
identified with an $ H $-invariant Volterra function $ F $ on $ X_2. $ \par
\smallskip
In order to make the harmonic analysis on $ X_{d-1} $ possible, we
shall always consider invariant triplets $ ({\cal F}, {\bf F} ,F) $ on $
X^{(c)}_{d-1} $ {\sl of
moderate growth at infinity\/}, which means that $ {\cal F} $ and $ F $ are
supposed to be \lq\lq essentially\rq\rq\ bounded by $ {\rm Cst} \left\vert
z^{(d-1)} \right\vert^ m, $ for a fixed
number $ m $ such that $ m > -1, $ called the order of moderate growth.
We will then give a unified proof (valid for all $ d, $ $ d \geq  2) $ of
the following property, which (in view of the pioneering works
$\lbrack$F$\rbrack$ and $\lbrack$G$\rbrack$) we call: \par
\medskip
{\bf Theorem F.G.:} \par
\smallskip
(I) For every invariant triplet $ ({\cal F}, {\bf F} ,F) $ with domain $
D_\mu $ on $ X^{(c)}_{d-1} $
and moderate
growth of order $ m $ $ (m > -1), $ there exists a function $ \tilde
F(\lambda) $ analytic
in the half-plane $ \Bbb C^{(m)}_+=\{ \lambda \in \Bbb C; $ $ {\rm Re} \
\lambda  > m\} $ enjoying the following properties: \par
a) $ \tilde  F(\lambda) $ is a Laplace-type transform of the $ (H
$-invariant) Volterra function $ F $ on $ X_{d-1}; $ it is expressed in terms
of $ \underline{{\rm f}} \left(z^{(d-1)} \right) = F(z) $ by a relation of the
following form (similar
to (1.1)):
$$ \tilde  F(\lambda)  = {\rm Cst} \int^{ \infty}_{ {\rm cosh} \ \mu}
Q^{(d)}_\lambda( \zeta) \underline{{\rm f}}(\zeta) \left(\zeta^ 2-1
\right)^{{d-3 \over 2}} {\rm d} \zeta \eqno (1.2) $$
where $ Q^{(d)}_\lambda $ is a generalized Legendre function of the second
kind. \par
b) The sequence of coefficients $ [ {\bf f}]_{\ell} $ of the Fourier-Legendre
expansion of the function $ {\bf F}(z) \left(\equiv \underline{{\bf f}}
\left(z^{(d-1)} \right) \right) $ on $ \Bbb S^{(d-1)}, $
namely $ \left\{[ {\bf f}]_{\ell} = \int^{ +1}_{-1}P^{(d)}_{\ell}( \zeta)
\underline{{\bf f}}(\zeta) \left(1-\zeta^ 2 \right)^{{d-3 \over 2}} {\rm d}
\zeta \right\} $ (where the $ P^{(d)}_{\ell} $ are the
corresponding generalized Legendre polynomials) is such that:
$$ \forall \ell ,\ \ {\rm with} \ \ \ell  \geq  m,\ \ \ \ [ {\bf f}]_{\ell}  =
\tilde  F(\ell) \eqno (1.3) $$
Moreover, $ \tilde  F(\lambda) $ is a carlsonian analytic interpolation of the
sequence $ \left\{[ {\bf f}]_{\ell} ;\ \ell  \geq  m \right\} $ in $ \Bbb
C^{(m)}_+, $ which is uniformly bounded
by $ {\rm Cst} \ {\rm e}^{-( {\rm Re} \ \lambda -m)\mu} . $ \par
\smallskip
(II) Conversely, if $ {\bf F} $ is an $ {\bf H} $-invariant function on $ \Bbb
S^{(d-1)} $
whose set of Fourier-Legendre coefficients $ [ {\bf f}]_{\ell} $ admit an
analytic
interpolation $ \tilde  F(\lambda) $ in $ \Bbb C^{(m)}_+ $ \lq\lq
essentially\rq\rq\ bounded by $ {\rm Cst\ e}^{-( {\rm Re} \ \lambda -m)\mu} ,
$
then $ {\bf F} $ admits an analytic continuation $ {\cal F} $ in $ D_\mu $
with moderate
growth of order $ m, $ and $ \tilde  F(\lambda) $ is related to the
corresponding
triplet $ ({\cal F}, {\bf F} ,F) $ by the transformation (1.2). \par
\medskip
More precise statements concerning specific classes of
triplets $ ({\cal F}, {\bf F} ,F) $ and the corresponding classes of functions
$ \tilde  F $
will be proved in the course of this work
and the converse part (II) will also be completed by inversion formulae
allowing to reconstruct the triplet $ ({\cal F}, {\bf F} ,F) $
in terms of a given function $ \tilde  F: $ these
various aspects of the theorem F.G. will be developed in theorems
1, 3, 3', 4, 5 (see \S 2.3 for the case $ d=2, $ and section 4 for the
general case $ d\geq 3). $ \par
\medskip
Our way of stating the previous properties of dual
analyticity might seem somewhat pedantic if one wanted to
consider these properties (as in $\lbrack$S,W$\rbrack$)
merely as theorems of complex analysis in one variable involving
generalized Legendre functions. However, we wish to show that the
geometrical (and thereby multidimensional) viewpoint adopted
throughout this work presents other features which go
substantially beyond the one-variable formulation, and are directly
related to the group-theoretical considerations made at the
beginning. \par
A first important feature which complements the theorem
F.G. concerns convolution products. In $\lbrack$B,V-1$\rbrack$, we introduced
a
composition product $ \ast^{( c)} $ on the class of perikernels on $
X^{(c)}_{d-1}, $
enjoying the following property. If $ {\cal K }= {\cal K}_1\ast^{( c)}{\cal
K}_2 $ and if $ \left({\cal K}_i, {\bf K}_i,K_i \right) $
$ (i=1,2), $ $ ({\cal K}, {\bf K} ,K) $ denote the corresponding triplets,
then one has: \par
a) $ {\bf K}  = {\bf K}_1\ast {\bf K}_2, $ where $ \ast $ denotes the
composition product of
kernels on $ \Bbb S^{(d-1)}. $ \par
b) $ K = K_1\diamond K_2, $ where $ \diamond $ denotes the composition product
of
Volterra kernels on $ X_{d-1} $ (in the sense of $\lbrack$Fa-1$\rbrack$ and
$\lbrack$Fa-V$\rbrack$). \par
In the case of invariant triplets of moderate growth
considered here, the composition products $ K_1\diamond K_2 $ and $ {\bf
K}_1\ast {\bf K}_2 $ can now
be interpreted as {\sl convolution products\/} respectively on $ X_{d-1} $ and
$ \Bbb S^{(d-1)}, $
involving the action of the respective groups $ {\rm SO}(1,d-1) $ and $ {\rm
SO}(d;\Bbb R). $
By extending a result of $\lbrack$Fa,V$\rbrack$, we will then show that the
relation $ K = K_1\diamond K_2 $ implies the following relation for the
corresponding Laplace-type transforms: $ \tilde  F(\lambda)  = \tilde
F_1(\lambda) \cdot\tilde  F_2(\lambda) . $
It follows that the corresponding property for convolution
products on $ \Bbb S^{(d-1)}, $ namely $ {\bf K}  = {\bf K}_1\ast {\bf K}_2
\Longrightarrow  \forall \ell , $ $ [ {\bf f}]_{\ell}  = \left[ {\bf f}_1
\right]_{\ell} \cdot \left[ {\bf f}_2 \right]_{\ell} , $
is then obtained from the previous one (for every $ \ell  \geq  m) $ as a
by-product of the theorem F.G. (property (I-b)): our theorem 2 (in
\S 2.3 for $ d=2 $ and in \S 4.3 for $ d\geq 3) $ will be devoted to this
convolution property. \par
\smallskip
Another feature which will in fact govern the
whole presentation of this work concerns our method for proving the
results announced above. The main idea of this method is that,
{\sl for every dimension\/} $ d $ $ (d \geq  2), $ the theorem F.G. (completed
by the
convolution properties) can be reduced to a simple property of
Fourier-Laplace analysis in $ \Bbb C, $ which we call \lq\lq property $ (
{\rm F.G} )_{\circ} \dq\dq : $
the latter concerns classes of $ (2\pi -) $ periodic functions $ f(\theta) $
analytic in a certain cut-plane $ \dot {\cal J}^{(\mu)} $ and satisfying an
additional
symmetry condition $ \left(S_d \right) $ corresponding to the $ d
$-dimensional problem
considered; the associated analytic function $ \tilde  F(\lambda) , $ which we
call
$ \tilde  F={\cal L}(f), $ coincides with the
usual Laplace transform $ L( {\rm f}) $ of the jump $ {\rm f} $ of $ f $
across the cut
of $ \dot {\cal J}^{(\mu)} $ (along $ i \Bbb R^+). $ Our section 2 is devoted
to this property (\S 2.2) and
to the theorem F.G. for the case $ d=2. $ The latter, which concerns harmonic
analysis on the circle and on the hyperbola (\S 2.3), can be
directly reduced to the property $ {\rm (F.G)}_{\circ} $ by using
the parametrization $ z^{(0)} = -i\ {\rm sin} \ \theta , $ $ z^{(1)} = {\rm
cos} \ \theta $ of $ X^{(c)}_1. $ One
will also notice that the results on Taylor series with analytic
coefficients $ \sum^{ }_ na(n) \zeta^ n $ of the type obtained in
$\lbrack$L$\rbrack$ or $\lbrack$S.W$\rbrack$
reduce to the primary version of the property $ ( {\rm F.G} )_{\circ} $
(namely without the symmetry
condition $ \left(S_d \right), $ as presented in \S 2.1) via the change of
variable $ \zeta  = {\rm e}^{-i\theta} . $ \par
In the general case $ d \geq  3, $ the reduction of the theorem
F.G to the property $ {\rm (F.G} )_{\circ} $ necessitates the use of an
integral
transformation $ {\cal R}^{(c)}_d $ which we call \lq\lq Radon-Abel
transformation\rq\rq\ and
whose properties are studied in our section 3. After this
preparatory work, the function $ \tilde  F(\lambda) $ associated with a given
triplet $ ({\cal F}, {\bf F} ,F) $ will be introduced in section 4 (\S 4.1)
via a
direct $ d $-dimensional extension of the \lq\lq spherical Laplace
transformation\rq\rq\ on $ X_2 $ of $\lbrack$Fa-1$\rbrack$ and
$\lbrack$Fa,V$\rbrack$, from which we specially
retain the geometrical meaning. We shall put:
$$ \tilde  F = L_d(F) = \left[L\circ{\cal R}_d \right](F)\ , \eqno (1.4) $$
where $ {\cal R}_d $ denotes an integral transformation involving the
integration of the invariant Volterra function $ F $ on appropriate
\lq\lq horocycles\rq\rq\ of $ X_{d-1} $ (with compact support); $ {\cal
R}_d(F) $ is then a
function with support contained in $ \Bbb R^+, $ on which the one-dimensional
Laplace
transformation $ L $ can be applied. \par
The transformation $ {\cal R}^{(c)}_d $ operates (in a way similar to $ {\cal
R}_d) $
by integration on \lq\lq complex horocycles\rq\rq\ which are $ (d-2) $-cycles
on
$ X^{(c)}_{d-1} $  obtained by continuous distortion of the horocycles of $
X_{d-1}. $
By acting on the analytic function $ {\cal F} $ of the given triplet, it
yields an analytic function $ \hat  f = {\cal R}^{(c)}_d({\cal F}) $ of one
complex variable
which satisfies the conditions of the property $ {\rm (F.G)}_{\circ} , $ and
whose
jump $ \hat {\rm f} ${\bf\ }across the cut of $ \dot {\cal J}^{(\mu)} $ (along
$ i\Bbb R^+) $ coincides with $ {\cal R}_d(F) $
(see \S 3.6). \par
The reduction of the theorem F.G to the property $ {\rm (F.G} )_{\circ} $
will then appear as a consequence of formula (1.4), of its
\lq\lq complex analogue\rq\rq :
$$ \tilde  F = {\cal L}_d({\cal F}) = \left[{\cal L}\circ{\cal R}^{(c)}_d
\right]({\cal F}) \eqno (1.5) $$
(expressing the fact that $ \tilde  F={\cal L} \left(\hat  f \right) = L
\left(\hat {\rm f} \right)), $ and of
the treatment of
the harmonic analysis on the sphere $ \Bbb S^{(d-1)} $ by means of the
transformation $ {\cal R}^{(c)}_d, $ given in \S 4.2; the complete arguments
for this
reduction are presented in \S 4.3 and \S 4.4, devoted respectively to
the direct and converse parts (I and II) of the theorem F.G.. \par
\medskip
Let us now mention the following complementary features
and results which will also be found in sections 3 and 4. \par
\smallskip
{\sl a) Volterra convolutions and algebras of perikernels\/} \par
As a first result, it is proved (in \S 3.3, proposition 9
i)) that the invariant Volterra kernels of moderate growth of a
given order $ m $ $ (m>-1) $ form a subalgebra of the Volterra algebra on
$ X_{d-1}. $ A similar result (quoted in \S 3.6, proposition 12,
but proved elsewhere) holds for the invariant perikernels of
moderate growth on $ X^{(c)}_{d-1}. $ \par
Secondly, this structure of convolution algebra is shown
to be preserved by the transformation $ {\cal R}_d $ in the following sense:
the transform $ {\cal R}_d \left(F_1\diamond F_2 \right) $ of the Volterra
convolution product of $ F_1 $ and
$ F_2 $ on $ X_{d-1} $ is equal to the one-dimensional convolution product
(on $ \Bbb R^+) $ of the corresponding transforms $ {\cal R}_d \left(F_1
\right) $ and $ {\cal R}_d \left(F_2 \right). $
The corresponding property for the Laplace
transforms (i.e. $ \tilde  F(\lambda)  = \tilde  F_1(\lambda) \cdot\tilde
F_2(\lambda) ) $ will then follow directly
from Eq.(1.4). This argument, which relies directly on the use of
horocycles and on the interplay of Iwasawa-type and Cartan-type
decompositions of the group $ {\rm SO}(1,d-1), $ replaces the computation of
$\lbrack$Fa,V$\rbrack$ based on a \lq\lq product formula\rq\rq\ for the
Legendre
functions of the second kind. \par
\smallskip
{\sl b) Inversion formulae\/} \par
Integral formulae expressing the triplet $ ({\cal F}, {\bf F} ,F) $ in terms
of $ \tilde  F $ are obtained by inverting formulae (1.4) and (1.5) and by
taking into account the explicit expressions for $ L^{-1}, $ $ {\cal L}^{-1} $
and for
$ {\cal R}^{-1}_d, $ $ {\cal R}^{(c)-1}_d $ derived respectively in \S 2.2 and
\S 3.7. These integral
formulae, presented in \S 4.4 (theorem 4), involve a family of
\lq\lq elementary perikernels\rq\rq , represented by analytic functions $
\Psi^{( d)}_\lambda( {\rm cos} \ \theta) $
whose properties are studied in \S 3.8. An equivalent system of
inversion formulae involving Legendre functions of the first kind
$ P^{(d)}_\lambda( {\rm cos} \ \theta) $ as integral kernels is also obtained
in theorem 4. A
geometrical interpretation of the latter is given in our final
subsection \S 4.5: the Legendre functions $ P^{(d)}_\lambda $ can also be
characterized as a second family of \lq\lq elementary perikernels\rq\rq ,
thanks to the introduction of a modified complex Radon-Abel
transformation $ {\cal R}^{\ast( c)}_d, $ allowing the following variant of
Eq.(1.5)
to be justified:
$$ \tilde  F = \left[{\cal L}\circ{\cal R}^{\ast( c)}_d \right]({\cal F})\ .
\eqno (1.5^\prime ) $$
\par
Algebraic relations between the functions $ \Psi^{( d)}_\lambda $
and $ P^{(d)}_\lambda $ (used in our proof of theorem 4 for showing the
equivalence of the two systems of inversion formulae) are derived
in \S 3.8 with the help of useful identities on Abel-type
transforms, established in the Appendix. \par
\smallskip
{\sl c) Functional spaces\/} \par
A more technical aspect of the theorem F.G
concerns the functional spaces which the various functions $ {\cal F}, $ $
\tilde  F $
etc... are supposed to belong to; some useful norm estimates
concerning the transformations $ {\cal R}^{(c)}_d $ and $ {\cal R}^{(c)-1}_d $
will be given in \S 3.6
and 3.7. Our choice of $ L^1 $-type or $ L^{\infty}
$-type norms (convenient under several respects, in particular for
the study of convolution products) will not allow us
to obtain a {\sl strict\/} equivalence between classes of triplets $ ({\cal
F}, {\bf F} ,F) $
and classes of functions $ \tilde  F, $ but theorems 3, 3' and 5 will however
appear as
a satisfactory \lq\lq converse of theorem 1\rq\rq . \par
\medskip
In a recent work $\lbrack$Fa-2$\rbrack$, J. Faraut has generalized his
\lq\lq spherical Laplace transformation\rq\rq\ to a large class of real
symmetric spaces on which the notion of invariant Volterra
kernel is still meaningful. One can express the hope that
the present analysis might also be extended (under conditions
which remain to be defined) to the corresponding complexified
spaces. \par
\bigskip
\noindent{\bf 2. THE PROPERTY }$ {\bf (F.G.)}_{\circ} ; ${\bf\ FOURIER-LAPLACE
TRANSFORMATION ON THE
COMPLEX HYPERBOLA }$ {\bf X}^{( {\bf c})}_ {\bf 1} $ \par
\medskip
{\bf 2.1 Connection between the Fourier series and
the Laplace transformation} \par
\smallskip
In the complex plane $ \Bbb C $ of the variable $ \theta  = u+iv, $ we
consider, for each $ \mu >0, $ the \lq\lq cut-domain\rq\rq\ $ {\cal
J}^{(\mu)}_ +={\cal J}_+\backslash \Xi( \mu) , $ where
$ {\cal J}_+ = \{ \theta \in \Bbb C;\ {\rm Im} \ \theta  > 0\} , $ and $ \Xi(
\mu) =\{ \theta \in \Bbb C; $ $ \theta =2k\pi +iv; $ $ k\in \Bbb Z, $ $ v
\geq  \mu\} . $ We
will use the notation $ \dot  A = A/2\pi  \Bbb Z $ for every subset $ A $ of $
\Bbb C $
which is invariant under the translation group $ 2\pi  \Bbb Z $ (e.g. $ \dot
{\cal J}^{(\mu)}_ +, $
$ \dot {\cal J}_+, $ $ \dot  \Xi( \mu) , $ $ \dot {\Bbb R}, $ $ \dot {\Bbb C}
$ etc...). \par
We shall call $ \Gamma $ the class of all cycles $ \gamma $ of $ \dot {\cal
J}^{(\mu)}_ + $ such
that: \par
\smallskip
a) $ \gamma $ is homologous in $ H^F_1 \left(\dot {\cal J}^{(\mu)}_ + \right)
$ to any cycle $ \gamma_ a $ with
support $ ]-a+i\infty ,-a] \cup  [-a,+a]\cup[ a,a+i\infty[ , $ $ 0 < a < 2\pi
, $ oriented from $ -a+i\infty $ to
$ a+i\infty . $ \par
\smallskip
b) $ {\rm supp\ } \gamma $ contains linear infinite branches of the form $
\left]-a_-+i\infty ,-a_-+iv_\gamma \right], $ \par
\noindent$ \left[a_++iv_\gamma ,a_++i\infty \right[ $ (see Fig.1a)). \par
\medskip
As in $\lbrack$B.V-1$\rbrack$, we consider the space $ {\cal O}^{(0)}
\left(\dot {\cal J}^{(\mu)}_ + \right) $ of
functions $ f(\theta) $ which satisfy the following conditions: \par
i) $ \forall k\in \Bbb Z, $ $ f(\theta)  = f(\theta +2k\pi) . $ \par
ii) $ f $ is holomorphic in $ \dot {\cal J}^{(\mu)}_ + $ and admits a
continuous
boundary value on $ \dot {\Bbb R}, $ denoted by $ {\bf f}(u). $ \par
iii) The limits $ f_\varepsilon( v) = \dlowlim{ {\rm lim}}{
\matrix{\displaystyle u \longrightarrow 0 \cr\displaystyle u>0 \cr}} \
f(\varepsilon u+iv), $ $ \varepsilon  = + $ or $ -, $ exist
in $ {\cal C}^{(0)} \left(\Bbb R^+ \right) $ and one puts: $ {\rm f}(v) = i
\left[f_+(v)-f_-(v) \right] = \Delta( f)(v); $
this jump $ {\rm f} $ of $ f $ across the cut $ \dot  \Xi( \mu) $ belongs to
the space $ {\cal C}^{(0)}_\mu \left(\Bbb R^+ \right) $ of continuous
functions on $ \Bbb R^+ $ whose support is contained in $ [\mu ,+\infty[ . $
\par
For such functions $ f, $ we put for any given $ m $ in $ \Bbb R, $
$$ g_m[f](v) = {\rm e}^{-mv} \dlowlim{ {\rm sup}}{0\leq u\leq 2\pi}  \vert
f(u+iv)\vert , \eqno (2.1) $$
which implies in particular that:
$$ {\rm e}^{-mv} \left\vert f_{\pm}( v) \right\vert  \leq  g_m[f](v)\ \ \ \
{\rm and} \ \ \ {\rm e}^{-mv}\vert {\rm f}(v)\vert  \leq  2\ g_m[f](v) \eqno
(2.1^\prime ) $$
\par
We now introduce the subspace $ {\cal O}_m \left(\dot {\cal J}^{(\mu)}_ +
\right) $ of functions $ f(\theta) $
satisfying conditions i), ii), iii) and the following one: \par
\smallskip
iv) $ g_m[f](v) \in  L^1 \left(\Bbb R^+ \right); $ we then put
$$ \left\Vert g_m[f] \right\Vert_ 1 = \int^{ \infty}_ 0g_m[f](v) {\rm d} v
\eqno (2.1^{\prime\prime} ) $$
\par
We also denote by $ \underline{{\cal O}}_m \left(\dot {\cal J}^{(\mu)}_ +
\right) $ (resp. $ \underline{{\cal O}}^\ast_ m \left(\dot {\cal J}^{(\mu)}_ +
\right)) $ the subspace of functions $ f(\theta) $
obtained by replacing condition iv) by the following condition v)
(resp. $ {\rm v}^\ast )): $ \par
v) $ g_m[f](v) \in  L^{\infty} \left(\Bbb R^+ \right); $ we then put:
$$ \left\Vert g_m[f] \right\Vert_{ \infty}  = \dlowlim{ {\rm sup}}{v;v\geq 0}
\left\vert g_m[f](v) \right\vert \ . \eqno (2.1^{\prime\prime} ^\prime ) $$
\par
$ {\rm v}^\ast ) $ $ g^\ast_ m[f](v) = (1+v) g_m[f](v) \in  L^{\infty}
\left(\Bbb R^+ \right). $ \par
\smallskip
With each function $ f $ in $ {\cal O}_m \left(\dot {\cal J}^{(\mu)}_ +
\right) $ or $ \underline{{\cal O}}_m \left(\dot {\cal J}^{(\mu)}_ + \right),
$ we can
associate: \par
\smallskip
a) The sequence of Fourier coefficients of $ {\bf f}(u): $
$$ \forall \ell  \in  \Bbb Z,\ \ [ {\bf f}]_{\ell}  = \int^ \pi_{ -\pi} {\rm
e}^{i\ell u} {\bf f}(u) {\rm d} u \eqno (2.2) $$
\par
b) The image $ \tilde  F $ of the jump $ {\rm f} =\Delta( f) $ of $ f $ by the
Laplace transformation $ L: $
$$ [L( {\rm f})](\lambda)  = \tilde  F(\lambda)  = \int^{ +\infty}_ \mu {\rm
e}^{-\lambda v} {\rm f}(v) {\rm d} v \eqno (2.3) $$
\par
We then have: \par
\smallskip
{\bf Proposition 1}: The following properties hold for any
function $ f $ in $ {\cal O}_m \left(\dot {\cal J}^{(\mu)}_ + \right) $ $
(m\in  \Bbb R, $ $ \mu >0): $ \par
i) the corresponding transform $ \tilde  F(\lambda) $ $ (\tilde  F=L(\Delta
f)) $ is analytic in the
half-plane $ \Bbb C^{(m)}_+ = \{ \lambda =\sigma +i\nu \in  \Bbb C;\ \sigma  >
m\} , $ continuous in $ \bar {\Bbb C}^{(m)}_+ $
and satisfies the uniform majorization:
$$ \left\vert\tilde  F(\sigma +i\nu) \right\vert  \leq  2 \left\Vert g_m[f]
\right\Vert_ 1 {\rm e}^{-(\sigma -m)\mu} \ , \eqno (2.4) $$
\par
ii) $ \tilde  F $ can be defined directly in terms of $ f $ by the
following Fourier-Laplace transformation $ {\cal L}: $
$$ \forall \lambda \in  \bar {\Bbb C}^{(m)}_+\ ,\ \ \tilde  F(\lambda)  =
[{\cal L}(f)](\lambda)  = \int^{ }_ \gamma {\rm e}^{i\lambda \theta}
f(\theta) {\rm d} \theta \ , \eqno (2.5) $$
where $ \gamma $ is any cycle in $ \Gamma . $ \par
\smallskip
iii) the corresponding Fourier coefficients $ [ {\bf f}]_{\ell} $ (of $ {\bf
f}(u)) $
and the transform $ \tilde  F(\lambda) $ are related together by the following
set
of equations:
$$ \forall \ell  \in  \Bbb Z\ ,\ \ {\rm with} \ \ \ell  \geq  m\ ,\ \ [ {\bf
f}]_{\ell}  = \tilde  F(\ell) \eqno (2.6) $$
\vfill\eject
{\bf Proof:} \par
i) The properties described in i) follow from the second
inequality (2.1') which implies the majorization of the r.h.s. of
Eq.(2.3) by $ 2 \int^{ +\infty}_ \mu {\rm e}^{-(\sigma -m)v}g_m[f](v) {\rm d}
v, $ for all $ \lambda  = \sigma +i\nu $ in $ \bar {\Bbb C}^{(m)}_+. $ \par
\smallskip
ii) Since $ g_m[f] \in  L^1 \left(\Bbb R^+ \right), $ each integral
$$ I_\gamma( \lambda)  = \int^{ }_ \gamma {\rm e}^{i\lambda \theta} f(\theta)
{\rm d} \theta \eqno  $$
is well-defined for $ \lambda \in  \bar {\Bbb C}^{(m)}_+ $ provided $ \gamma $
belongs to the
class $ \Gamma ; $ in order to show that $ I_\gamma( \lambda) $ is independent
of the
representative $ \gamma , $ let us consider two cycles $ \gamma , $ $
\gamma^{ \prime} $ in $ \Gamma , $ which we
respectively approximate by sequences of paths $ \left\{ \gamma^{( n)}
\right\} $ and $ \left\{ \gamma^{ \prime( n)} \right\} $
obtained from $ \gamma $ and $ \gamma^{ \prime} $ by chopping the infinite
parts of the
latter contained in the set $ \{ \theta =u+iv; $ $ \left.v\geq v_n \right\} $
(Fig.1b); we assume
that $ \left\{ v_n \right\} $ is an increasing sequence such that $ v_n
\longrightarrow +\infty $ and $ g_m[f] \left(v_n \right) \longrightarrow 0 $
for $ n \longrightarrow \infty $ (which is legitimate since $ g_m[f]\in L_1
\left(\Bbb R^+ \right)). $ \par
We then have: $ \int^{ }_{ \gamma^{( n)}} {\rm e}^{i\lambda \theta} f(\theta)
{\rm d} \theta \longrightarrow I_\gamma( \lambda) , $ $ \int^{ }_{ \gamma^{
\prime( n)}} {\rm e}^{i\lambda \theta} f(\theta) {\rm d} \theta
\longrightarrow I_{\gamma^{ \prime}}( \lambda) , $
and in view of the Cauchy theorem and of Eq.(2.1):
$$ \left\vert \int^{ }_{ \gamma^{( n)}-\gamma^{ \prime( n)}} {\rm
e}^{i\lambda \theta} f(\theta) {\rm d} \theta \right\vert  \leq  4\pi \ {\rm
e}^{-(\sigma -m)v_n}g_m[f] \left(v_n \right) $$
(since $ \gamma^{( n)}-\gamma^{ \prime( n)} $ is homologous to a sum of two
paths contained
respectively in $ \left[-2\pi +iv_n, \right. $ $ \left.iv_n \right] $ and $
\left[iv_n, \right. $ $ \left.2\pi +iv_n \right], $ as shown by
Fig.1b). It follows from the latter inequality that $ \dlowlim{ {\rm lim}}{n
\longrightarrow \infty}  \left\vert \int^{ }_{ \gamma^{( n)}-\gamma^{ \prime(
n)}} {\rm e}^{i\lambda \theta} f(\theta) {\rm d} \theta \right\vert =0, $
which entails that $ I_\gamma( \lambda)  = I_{\gamma^{ \prime}}( \lambda) , $
from now on denoted by $ [{\cal L}(f)](\lambda) . $ \par
By choosing for $ \gamma $ any representative $ \gamma_ a $ with $ a
\longrightarrow 0 $ (see our
definition of $ \Gamma , $ condition a)), we can then write:
$$ \matrix{\displaystyle[{\cal L}(f)](\lambda) & \displaystyle = \dlowlim{
{\rm lim}}{a \longrightarrow 0} \left\{ -i \int^{ +\infty}_ 0 {\rm
e}^{i\lambda( iv-a)}f(-a+iv) {\rm d} v... \right. \hfill \cr\displaystyle  &
\displaystyle \left.+ \int^ a_{-a} {\rm e}^{i\lambda u} {\bf f}(u) {\rm d} u +
i \int^{ +\infty}_ 0 {\rm e}^{i\lambda( iv+a)}f(a+iv) {\rm d} v \right\}
\hfill \cr\displaystyle  & \displaystyle = \int^{ +\infty}_ \mu {\rm
e}^{-\lambda v} {\rm f}(v) {\rm d} v = \tilde  F(\lambda) \ . \hfill \cr} $$
\par
\smallskip
iii) By choosing the representative $ \gamma =\gamma_ \pi $ (Fig.1a) and
making use of the $ 2\pi $-periodicity of the integrand $ {\rm e}^{i\lambda
\theta} f(\theta) , $ when $ \lambda =\ell $
is an integer, we can write (for $ \ell  \geq  m): $
$$ \tilde  F(\ell)  = [{\cal L}(f)](\ell)  = \int^{ }_{ \gamma_ \pi} {\rm
e}^{i\ell \theta} f(\theta) {\rm d} \theta  = \int^{ +\pi}_{ -\pi} {\rm
e}^{i\ell u} {\bf f}(u) {\rm d} u = [ {\bf f}]_{\ell} $$
\par
\smallskip
{\bf Remarks}: \par
i) As a corollary, the majorization (2.4) of $ \tilde  F(\lambda) $ yields
correspondingly: $ \forall \ell  \geq  m, $ $ \left\vert[ {\bf f}]_{\ell}
\right\vert  \leq  {\rm Cst.} \cdot {\rm e}^{-(\ell -m)\mu} , $ the latter
being
already implied by the analyticity of $ f $ in the strip $ \Theta_ \mu  = \{
\theta =u+iv\ ; $
$ 0 < v < \mu\} $ (provided $ f $ is continuous in $ \bar  \Theta_ \mu ). $
\par
ii) In view of Carlson's theorem (see e.g. $\lbrack$Bo$\rbrack$ p.153) the
function $ \tilde  F(\lambda) $ represents the unique analytic interpolation
of the
sequence $ \left\{[ {\bf f}]_{\ell} ;\ell \geq m \right\} $ in $ \Bbb
C^{(m)}_+ $ satisfying the majorization (2.4).
As a matter of fact, this uniqueness property will appear as a
direct by-product of the converse of proposition 1 which we shall
now establish. \par
\smallskip
{\bf Proposition 2} : Let $ {\bf f}(u) $ be a continuous function on $ \dot
{\Bbb R} $
whose set of Fourier coefficients $ [ {\bf f}]_{\ell} $ (defined by Eq.(2.2))
admit an
analytic interpolation $ \tilde  F(\lambda) $ satisfying the following
properties: \par
i) $ \tilde  F(\lambda) $ is analytic in the half-plane $ \Bbb C^{(m)}_+ $ and
continuous in $ \bar {\Bbb C}^{(m)}_+, $ $ m $ being a given number in
$ \Bbb R\backslash \Bbb Z. $ \par
ii) $ \forall \ell  \in  \Bbb Z, $ with $ \ell >m: $ $ \tilde  F(\ell)  = [
{\bf f}]_{\ell} . $ \par
iii) In $ \bar {\Bbb C}^{(m)}_+, $ $ \tilde  F $ satisfies a uniform
majorization of the following form:
$$ \dlowlim{ {\rm sup}}{\sigma \geq m} {\rm e}^{(\sigma -m)\mu}
\left\vert\tilde  F(\sigma +i\nu) \right\vert  = G^{(\mu)}_ m \left[\tilde  F
\right](\nu) \ , \eqno (2.7) $$
with $ G^{(\mu)}_ m \left[\tilde  F \right] \in  L^1(\Bbb R). $ \par
\smallskip
Then there exists a (unique) analytic function $ f $ in $ \underline{{\cal
O}}_m \left(\dot {\cal J}^{(\mu)}_ + \right) $
whose boundary value on $ \dot {\Bbb R} $ is $ {\bf f} $ and whose jump $
{\rm f} =\Delta( f) $
across the cut $ \dot  \Xi( \mu) $ admits $ \tilde  F $ as its Laplace
transform
$ (\tilde  F=L( {\rm f})={\cal L}(f)). $ Moreover the function $ f $ satisfies
the
following uniform majorization in its domain $ \dot {\cal J}^{(\mu)}_ +: $
$$ \dlowlim{ {\rm sup}}{u} \vert f(u+iv)\vert  \leq  {1 \over \pi\vert {\rm
sin} \ \pi m\vert}  \left\Vert G^{(\mu)}_ m \left[\tilde  F \right]
\right\Vert_ 1 {\rm e}^{mv} + \dlowlim{ {\rm sup}}{u}\vert {\bf f}(u)\vert
{\rm e}^{E[m]v} \eqno (2.8) $$
\par
\smallskip
{\bf Proof:} \par
We first introduce the function $ {\rm f}(v) $ by the formula:
$$ {\rm e}^{-mv} {\rm f}(v) = {1 \over 2\pi}  \int^{ +\infty}_{ -\infty} {\rm
e}^{i\nu v}\tilde  F(m+i\nu) {\rm d} \nu \ ; \eqno (2.9) $$
in view of condition iii), $ {\rm f}(v) $ is continuous and satisfies the
majorization:
$$ \left\vert {\rm e}^{-mv} {\rm f}(v) \right\vert  \leq  {1 \over 2\pi}
\left\Vert G^{(\mu)}_ m \left[\tilde  F \right] \right\Vert_ 1 $$
\par
By calling $ {\rm L}_\sigma , $ for $ \sigma  \geq  m, $ the line $ \{
\lambda =\sigma +i\nu ;\ \nu \in  \Bbb R\} $ of
the $ \lambda $-plane (oriented from $ \nu =-\infty $ to $ \nu =+\infty ), $
we can rewrite the
definition of $ {\rm f} $ given by Eq.(2.9) as follows:
$$ {\rm f}(v) = {1 \over 2\pi i} \int^{ }_{ {\rm L}_m} {\rm e}^{\lambda
v}\tilde  F(\lambda) {\rm d} \lambda \ ; \eqno (2.9^\prime ) $$
then, in view of condition iii), we can justify (as in the
contour distortion argument of proposition 1) the shifting of the
integration cycle in Eq.(2.9') from its initial situation (i.e. $ {\rm L}_m) $
to an arbitrary line $ {\rm L}_\sigma $ (with $ \sigma  > m), $ which yields
(in view of
(2.7)) the following majorization:
$$ \forall \sigma ,\ \ \sigma  \geq  m\ ,\ \ \vert {\rm f}(v)\vert  \leq  {
\left\Vert G^{(\mu)}_ m \left[\tilde  F \right] \right\Vert_ 1 \over 2\pi}
{\rm e}^{m\mu} \cdot {\rm e}^{-\sigma( \mu -v)}\ ; $$
the latter immediately implies that the support of $ {\rm f} $ is contained
in the half-line $ [\mu ,+\infty[ . $ \par
We will now construct the requisite function $ f(\theta) $ as a sum
$ f=h+r, $ where $ h(\theta)  = {1 \over 2\pi}  \sum^{ }_{ \ell >m}[ {\bf
f}]_{\ell} {\rm e}^{-i\ell \theta} $ will be shown to be analytic in
the domain $ \dot {\Bbb C}\backslash\dot  \Xi( \mu) $ and such that $ \Delta(
h) = {\rm f} ; $ the
remainder $ r(\theta)  = {1 \over 2\pi}  \sum^{ }_{ \ell <m}[ {\bf f}]_{ {\bf
\ell}} {\rm e}^{-i\ell \theta} $ will be analytic in $ \dot {\cal J}_+, $
continuous
in $ \dbar {\cal J}_+, $ and therefore will not contribute to the jump of $ f
$ across $ \dot  \Xi( \mu) $
(i.e. $ \Delta( r) = 0). $ \par
Let us first study the function $ h $ and show that $ h_{ \left\vert\dot
{\cal J}^{(\mu)}_ + \right.} $
belongs to the space $ \underline{{\cal O}}_m \left(\dot {\cal J}^{(\mu)}_ +
\right). $ \par
The definition and properties of $ h $ will be obtained by an
adaptation of Watson's resummation method $\lbrack$Wa$\rbrack$, namely, we
introduce the following integrals:
$$ h_\varepsilon( u) = {i \over 4\pi}  \int^{ }_ C{\tilde  F(\lambda) {\rm
e}^{-i\lambda( u-\varepsilon \pi)} \over {\rm sin} \ \pi \lambda}  {\rm d}
\lambda \ ,\ \ \ {\rm with} \ \ \ \varepsilon  = + {\rm or}  -\ \ ; \eqno
(2.10) $$
in the latter, the contour $ C $ encircles the half-line $ [E(m)+1,+\infty[ $
and belongs to $ \Bbb C^{(m)}_+, $ as shown on Fig.2. Since by assumption $
\tilde  F $
is analytic in $ \Bbb C^{(m)}_+ $ and satisfies a majorization of the form
(2.7), the integrals (2.10) are convergent and moreover remain
unchanged when the contour $ C $ is distorted (in $ \bar {\Bbb
C}^{(m)}_+\backslash\{ \lambda \in \Bbb R; $
$ \lambda  \geq  E(m)+1\} ) $ and replaced by the line $ {\rm L}_m= \{
\lambda =m+i\nu ; $ $ \nu \in  \Bbb R\} , $ provided
the real variable $ u $ is kept respectively in $ [0,2\pi] $ for $ h_+(u) $
and
in $ [-2\pi ,0] $ for $ h_-(u). $ \par
On the other hand, by applying the residue theorem to the
integrals (2.10) and by taking into account our assumption
$ \tilde  F(\ell)  = [ {\bf f}]_{\ell} , $ for $ \ell  > m, $ we finally
obtain the following relations: \par
\vfill\eject
for\nobreak\ \nobreak\  $ u \in  [0,2\pi] \ , $
$$ h_+(u) = - {1 \over 4\pi}  \int^{ +\infty}_{ -\infty}{\tilde  F(m+i\nu)
{\rm e}^{-i(m+i\nu)( u-\pi)} \over {\rm sin} \ \pi( m+i\nu)}  {\rm d} \nu = {1
\over 2\pi}  \sum^{ }_{ \ell >m}[ {\bf f}]_{\ell} {\rm e}^{-i\ell u} \eqno
(2.11) $$
\par
for\nobreak\ \nobreak\  $ u \in  [-2\pi ,0]\ , $
$$ h_-(u) = - {1 \over 4\pi}  \int^{ +\infty}_{ -\infty}{\tilde  F(m+i\nu)
{\rm e}^{-i(m+i\nu)( u+\pi)} \over {\rm sin} \ \pi( m+i\nu)}  {\rm d} \nu = {1
\over 2\pi}  \sum^{ }_{ \ell >m}[ {\bf f}]_{\ell} {\rm e}^{-i\ell u} \eqno
(2.12) $$
(i.e. $ h_+(u) = h_-(u-2\pi) , $ for $ u\in[ 0,2\pi] ). $ \par
If we now substitute the complex variable $ \theta  = u+iv $ to $ u $ in
the integral of Eq.(2.11) and take into account the following majorization
(where $ Y $ denotes the Heaviside function):
$$ \left\vert {\rm e}^{(-im+\nu)( u-\pi)} \right\vert \leq  2\ {\rm cosh} \
\pi \nu \ \left[Y(\nu) {\rm e}^{\nu( u-2\pi)} +Y(-\nu)  {\rm e}^{\nu u}
\right]\ , \eqno (2.12^\prime ) $$
we see that this integral converges uniformly, and therefore defines an
analytic continuation of $ h_+ $ in the strip $ \{ \theta =u+iv\in \Bbb C; $ $
0 < u < 2\pi\} , $
continuous in the closure of the latter. In fact, we have:
$$ \matrix{\displaystyle \forall \theta  = u+iv\ \ ,\ \ {\rm with} \ 0 \leq  u
\leq  2\pi \ \ ,\ \ v \in  \Bbb R\ , \cr\displaystyle {\rm e}^{-mv}h_+(u+iv) =
{1 \over 2\pi}  \int^{ \infty}_{ -\infty} {\rm e}^{i\nu v}\tilde
H^{(u)}_m(\nu) {\rm d} \nu \ , \cr} \eqno (2.13) $$
with
$$ \tilde  H^{(u)}_m(\nu)  = - {\tilde  F(m+i\nu) {\rm e}^{-i(m+i\nu)(
u-\pi)} \over 2\ {\rm sin} \ \pi( m+i\nu)} \ . \eqno (2.13^\prime ) $$
{}From (2.7) and (2.12'), it then follows that $ \left\vert\tilde
H^{(u)}_m(\nu) \right\vert  \leq  {G^{(\mu)}_ m \left[\tilde  F \right](\nu)
\over\vert {\rm sin} \ \pi m\vert} $
which thus yields the majorization:
$$ {\rm e}^{-mv} \left\vert h_+(u+iv) \right\vert  \leq  {1 \over 2\pi\vert
{\rm sin} \ \pi m\vert}  \left\Vert G^{(\mu)}_ m \left[\tilde  F \right]
\right\Vert_ 1\ . \eqno (2.14) $$
\par
The function $ h_-(\theta)  = h_+(\theta +2\pi) $ being similarly defined in
the strip $ \{ \theta  = u+iv; $ $ -2\pi \leq u\leq 0\} , $ we can compute the
discontinuity $ h_+(iv)-h_-(iv) $
by replacing $ u $ by $ iv $ in the integrals (2.11) and (2.12) and by
subtracting Eq.(2.12) from Eq.(2.11) side by side; we then obtain
(in view of Eq.(2.9')):
$$ h_+(iv)-h_-(iv) = {-i \over 2\pi}  \int^{ +\infty}_{ -\infty}\tilde
F(m+i\nu) {\rm e}^{(m+i\nu) v} {\rm d} \nu  = -i {\rm f}(v) \eqno (2.15) $$
It then follows that $ h_+ $ and $ h_- $ have a common analytic continuation
in $ \dot {\cal J}^{(\mu)}_ +\cup\dbar {\cal J}_- $ which (in view of formulae
(2.14) and (2.15)) is the
requisite analytic
function $ h(\theta) $ satisfying $ \Delta( h) = {\rm f} $ and $ h_{
\left\vert\dot {\cal J}^{(\mu)}_ + \right.} \in  \underline{{\cal O}}_m
\left(\dot {\cal J}^{(\mu)}_ + \right). $ \par
The function $ r(\theta) $ can now be studied as follows. Since $ {\bf f}(u) $
and $ {\bf h}(u) $ $ (= h_{\vert \dot {\Bbb R}}) $ are continuous, $ {\bf r}
= {\bf f} - {\bf h} $ is a
continuous periodic function represented by the Fourier series $ {1 \over
2\pi}  \sum^{ }_{ \ell <m}[ {\bf f}]_{\ell} {\rm e}^{-i\ell u}. $
It is therefore convenient to introduce the function $ \hat  r \left( {\rm
e}^{i\theta} \right) = {\rm e}^{+iE(m)\theta} r(\theta) ={1 \over 2\pi}
\sum^{ }_{ \ell^{ \prime} \geq 0}[ {\bf f}]_{E[m]-\ell^{ \prime}} {\rm
e}^{i\ell^{ \prime} \theta} $
and to consider the sequence of Cesaro polynomials $ \left\{{\cal C}_{\ell}
\left( {\rm e}^{i\theta} \right) \right\} $
associated with this Fourier series; in view of Fejer's theorem,
this sequence converges uniformly {\sl on the unit circle\/} to the
continuous function $ \hat  r(z)_{\vert\vert z\vert =1}, $ and therefore it
converges
uniformly in the {\sl closed\/} unit disk to the analytic function $ \hat
r(z). $
It follows that $ r(\theta) $ is analytic in $ \dot {\cal J}_+, $ continuous
in $ \dbar {\cal J}_+ $ and
satisfies the following majorization:
$$ \forall v,\ \ v\geq 0,\ \ {\rm e}^{-E[m]v} \dlowlim{ {\rm sup}}{u} \vert
r(u+iv)\vert  \leq  \dlowlim{ {\rm sup}}{u} \vert {\bf r}(u)\vert  \leq
\dlowlim{ {\rm sup}}{u}\vert {\bf f}(u)\vert + \dlowlim{ {\rm sup}}{u}\vert
{\bf h}(u)\vert \ , $$
which also implies, in view of the inequality (2.14):
$$ {\rm e}^{-E[m]v} \dlowlim{ {\rm sup}}{u}\vert r(u+iv)\vert  \leq
\dlowlim{ {\rm sup}}{u}\vert {\bf f}(u)\vert  + {1 \over 2\pi\vert {\rm sin} \
\pi m\vert}  \left\Vert G^{(\mu)}_ m \left[\tilde  F \right] \right\Vert_ 1
\eqno (2.16) $$
\par
The last statement of proposition 2 (inequality (2.8))
follows directly from the inequalitites (2.14) and (2.16), which
thereby achieves the proof of the fact that $ f(\theta) $ belongs to $
\underline{{\cal O}}_m \left(\dot {\cal J}^{(\mu)}_ + \right), $
q.e.d. \par
We now give an analogue of proposition 2 for the case when
$ m $ is an integer, which necessitates a slightly more refined
treatment. \par
\smallskip
{\bf Proposition 2':} If the assumptions of proposition 2 are
satisfied with $ m \in  \Bbb Z $ and if \par
\smallskip
$\alpha$) $ G^{(\mu)}_ m \left[{ {\rm d}\tilde  F \over {\rm d} \lambda}
\right] \in  L^1(\Bbb R)\ , $ \par
$\beta$) the function $ \tilde  F_m(\nu)  = \tilde  F(m+i\nu) $
is, in a given interval $ \left\{ \nu ;\vert \nu\vert  \leq  \nu_ 0 \right\} $
the integral of a function
of bounded variation, \par
$\gamma$) $ \tilde  F(m)=[ {\bf f}]_m, $ \par
\noindent then the conclusions of proposition 2 are still valid, up to the
following change: the function $ f $ belongs to $ \underline{{\cal O}}^\ast_ m
\left(\dot {\cal J}^{(\mu)}_ + \right) $ and
satisfies the following majorization in $ \dot {\cal J}^{(\mu)}_ +: $
$$ \dlowlim{ {\rm sup}}{u}\vert f(u+iv)\vert  \leq  C \left[ \left\Vert
G^{(\mu)}_ m \left[\tilde  F \right] \right\Vert_ 1+ \left\Vert G^{(\mu)}_ m
\left[{ {\rm d}\tilde  F \over {\rm d} \lambda} \right] \right\Vert_
1+V_{\nu_ 0} \left({ {\rm d}\tilde  F_m \over {\rm d} \nu} \right) \right] {
{\rm e}^{mv} \over 1+v} + \dlowlim{ {\rm sup}}{u} \vert {\bf f}(u)\vert  {\rm
e}^{(m-1)v}, \eqno (2.8^\prime ) $$
where $ V_{\nu_ 0} \left({ {\rm d}\tilde  F_m \over {\rm d} \nu} \right) $
denotes the sum of the total variations of the
functions $ {\rm Re} \left({ {\rm d}\tilde  F_m \over {\rm d} \nu} \right) $
and $ {\rm Im} \left({ {\rm d}\tilde  F_m \over {\rm d} \nu} \right) $ in the
interval $ \left[-\nu_ 0,+\nu_ 0 \right], $ and $ C $ is
a suitable constant. \par
\smallskip
\vfill\eject
{\bf Proof:} \par
The argument is similar to the one given in the proof of
proposition 2, the function $ f(\theta) $ being now expressed by $ f=r+h+{1
\over 2\pi}  [ {\bf f}]_m {\rm e}^{-im\theta} , $
with $ r $ and $ h $ defined as before. While $ r(u+iv) $ behaves at infinity
as $ {\rm Cst.\ e}^{(m-1)v}, $ it is now the function $ h_1(\theta)  =
h(\theta)  + {1 \over 2\pi}  [ {\bf f}]_m\ {\rm e}^{-im\theta} $
(instead of $ h) $ which can be checked to belong to $ \underline{{\cal O}}_m
\left(\dot {\cal J}^{(\mu)}_ + \right) $ by the
following argument. Let $ \chi( \nu) $ be a $ C^{\infty} $ function such that
$ 0<\chi <1, $ $ \chi( \nu)  = 0 $
for $ \vert \nu\vert  > \nu_ 0, $ $ \chi( \nu)  = 1 $ for $ \vert \nu\vert  <
\nu_ 1 $ $ \left(\nu_ 1<\nu_ 0 \right), $ and let $ \tilde  \psi_ u = \chi \
\tilde  H^{(u)}_m, $
$ \tilde  \varphi_ u = (1-\chi)  \tilde  H^{(u)}_m; $ the expression (2.13')
of $ \tilde  H^{(u)}_m(\nu) $ must now be
replaced by the following one (in the sense of distributions):
$$ \tilde  H^{(u)}_m(\nu)  = \dlowlim{ {\rm lim}}{ \matrix{\displaystyle
\varepsilon \longrightarrow 0 \cr\displaystyle \varepsilon >0 \cr}}  -
{\tilde  F(m+i\nu) {\rm e}^{-i(m+i\nu)( u-\pi)} \over 2\ {\rm sin} \ \pi(
m+\varepsilon +i\nu)} \ . \eqno (2.17) $$
\par
In view of Eq.(2.13), we can now write: $ {\rm e}^{-mv}h_+(u+iv) = \psi_
u(v)+\varphi_ u(v), $ where $ \psi_ u $ and $ \varphi_ u $
are respectively the inverse Fourier-transforms of $ \tilde  \psi_ u(\nu) $
and $ \tilde  \varphi_ u(\nu) . $
The function $ \tilde  \varphi_ u(\nu) $ is continuous on $ \Bbb R $ and, in
view of (2.7)
and (2.17), it is majorized by $ A\times G^{(\mu)}_ m \left[\tilde  F
\right](\nu) $ (with e.g. $ A= {\rm e}^{\pi \nu_ 1}/ {\rm sinh} \ \pi \nu_ 1);
$
it follows that $ \left\vert \varphi_ u(v) \right\vert \leq A \left\Vert
G^{(\mu)}_ m \left[\tilde  F \right] \right\Vert_ 1. $ \par
By applying the same argument to the function $ { {\rm d} \over {\rm d} \nu}
\tilde  \varphi_ u(\nu) $ and
by taking into account our assumption $\alpha$), one would justify
similarly an inequality of the following form:
$$ (1+v) \left\vert \varphi_ u(v) \right\vert  \leq  A \left\Vert G^{(\mu)}_ m
\left[{ {\rm d}\tilde  F \over {\rm d} \lambda} \right] \right\Vert_ 1 +
A^{\prime} \left\Vert G^{(\mu)}_ m \left[\tilde  F \right] \right\Vert_{ 1\ ,}
\eqno (2.17^\prime ) $$
\par
We now estimate $ \psi_ u(v) $ by making use of our
assumption $\beta$): in view of
Eq.(2.17), we can assert that the function $ \tilde  \alpha_ u(\nu)  = \pi
\nu\tilde  \psi_ u(\nu) $ $ =\chi( \nu) \times \left[\pi \nu \ \tilde
H^{(u)}_m(\nu) \right] $
has the same regularity properties as the function $ \tilde  F_m(\nu) $ in the
interval $ \left[-\nu_ 0,\nu_ 0 \right]; $ therefore, the second derivative $
\tilde  \alpha^{ \prime\prime}_ u(\nu) $ of $ \tilde  \alpha_ u(\nu) $
(in the sense of distributions)
is a measure with compact support; it follows that the
inverse Fourier transform $ \alpha_ u(v) $ of $ \tilde  \alpha_ u $ is
majorized by $ {\rm Cst.}(1+\vert v\vert)^{ -2}, $
where the constant
$ {\rm Cst.} $ can be seen to be proportional to the \lq\lq local
semi-norm\rq\rq\ $ V_{\nu_ 0} \left({ {\rm d}\tilde  F_m \over {\rm d} \nu}
\right) $
of $ \tilde  F_m. $ By taking Eq.(2.17) into account, we can then write:
\smallskip
$$ \psi_ u(v) = {1 \over 2\pi}  \int^{ +\infty}_{ -\infty} {\rm e}^{i\nu v}
{\tilde  \alpha_ u(\nu) \over \pi( \nu -i\varepsilon)}  {\rm d} \nu  = {i
\over \pi}  \int^ v_{-\infty} \alpha_ u \left(v^{\prime} \right) {\rm d}
v^{\prime} \ ; \eqno  $$
it follows that $ \psi_ u(+\infty)  = {i \over \pi}  \tilde  \alpha_ u(0) = -
{1 \over 2\pi}  {\rm e}^{-imu}\tilde  F(m) $ and that one
has:
$$ \left\vert \psi_ u(v)-\psi_ u(+\infty) \right\vert_{ \left\vert \Bbb R^+
\right.} \leq  B\ V_{\nu_ 0} \left({ {\rm d}\tilde  F_m \over {\rm d} \nu}
\right)(1+v)^{-1}, \eqno (2.17^{\prime\prime} ) $$
for a suitable constant $ B. $ \par
Therefore, in view of the assumption $ \tilde  F(m) = [ {\bf f}]_m, $ one
concludes that:
$$ \matrix{\displaystyle {\rm e}^{-mv}h_1(u+iv) & \displaystyle = \varphi_
u(v)+\psi_ u(v)+ {1 \over 2\pi}  [ {\bf f}]_m {\rm e}^{-imu} \hfill
\cr\displaystyle  & \displaystyle = \varphi_ u(v) + \left[\psi_ u(v) - \psi_
u(+\infty) \right], \hfill \cr} \eqno (2.18) $$
so that (in view of the inequalities (2.17') and (2.17")):
$$ (1+v) {\rm e}^{-mv} \left\vert h_1(u+iv) \right\vert  \leq  A^{\prime}
\left\Vert G^{(\mu)}_ m \left[\tilde  F \right] \right\Vert_ 1 + A \left\Vert
G^{(\mu)}_ m \left[{ {\rm d}\tilde  F \over {\rm d} \lambda} \right]
\right\Vert_ 1+B\ V_{\nu_ 0} \left({ {\rm d}\tilde  F_m \over {\rm d} \nu}
\right) \eqno (2.18^\prime ) $$
\par
A majorization of $ r(\theta) $ similar to that of proposition 2
yields:
$$ \left\vert {\rm e}^{-(m-1)v}r(u+iv) \right\vert  \leq  \dlowlim{ {\rm
sup}}{u} \vert {\bf f}(u)\vert  + \dlowlim{ {\rm sup}}{u} \left\vert {\bf
h}_1(u) \right\vert \eqno (2.19) $$
Inequalities (2.18') and (2.19) then imply the inequality (2.8'),
which achieves the proof of proposition 2'. \par
{\bf Remark:} \par
The function $ {\rm e}^{-mv}h_1(u+iv) $ given by the above expression
(2.18) also appears (via the previous analysis) as the inverse
Fourier transform of the distribution:
$$ \tilde  H^{(u)}_{m,1}(\nu)  = \dlowlim{ {\rm lim}}{ \matrix{\displaystyle
\varepsilon \longrightarrow 0 \hfill \cr\displaystyle \varepsilon >0 \hfill
\cr}}  - {1 \over 2} {\tilde  F(m+i\nu) {\rm e}^{-i(m+i\nu)( u-\pi)} \over
{\rm sin} \ \pi( m-\varepsilon +i\nu)} \eqno (2.19^\prime ) $$
(to be compared with the distribution $ \tilde  H^{(u)}_m(\nu) $ given by
Eq.(2.17)). \par
\medskip
{\bf Convolution products} \par
In proposition 1 of $\lbrack$B.V-1$\rbrack$, we introduced a convolution
product, denoted by $ f = f_1\ast^{( c)}f_2, $ such that if $ f_i\in{\cal
O}^{(0)} \left(\dot {\cal J}^{ \left(\mu_ i \right)}_+ \right), $
$ i = 1,2, $ $ f \in  {\cal O}^{(0)} \left(\dot {\cal J}^{ \left(\mu_ 1+\mu_ 2
\right)}_+ \right); $ this product enjoys the following two
properties: \par
\smallskip
i) $ {\bf f}  = {\bf f}_1\ast {\bf f}_2, $ where $ \ast $ denotes the usual
convolution
product on $ \Bbb R/2\pi  \Bbb Z $ (i.e. $ \int^{ 2\pi}_ 0 {\bf f}_1
\left(u-u^{\prime} \right) $ $ \times {\bf f}_2 \left(u^{\prime} \right) {\rm
d} u^{\prime} ), $ \par
ii) $ \Delta f = \Delta f_1\diamond \Delta f_2, $ where $ \diamond $ denotes
the convolution product
on $ \Bbb R^+ $ (i.e. $ \int^ v_0 {\rm f}_1 \left(v-v^{\prime} \right) $ $
\times {\rm f}_2 \left(v^{\prime} \right) {\rm d} v^{\prime} ). $ \par
\smallskip
We will now show: \par
\smallskip
{\bf Proposition 3:} : If $ f_1\in{\cal O}_m \left(\dot {\cal J}^{ \left(\mu_
1 \right)}_+ \right) $ and $ f_2\in{\cal O}_m \left(\dot {\cal J}^{
\left(\mu_ 2 \right)}_+ \right), $ then $ f=f_1\ast^{( c)}f_2\in{\cal O}_m
\left(\dot {\cal J}^{ \left(\mu_ 1+\mu_ 2 \right)}_+ \right). $
The corresponding transforms $ \tilde  F_1, $ $ \tilde  F_2, $ $ \tilde  F $
of $ f_1, $ $ f_2, $ $ f $ satisfy the
relation:
$$ \tilde  F(\lambda)  = \tilde  F_1(\lambda) \cdot\tilde  F_2(\lambda) \eqno
(2.20) $$
in their common domain $ \bar {\Bbb C}^{(m)}_+. $ This relation interpolates
(for $ \ell  \geq  m) $ the corresponding set of equations $ [ {\bf
f}]_{\ell}  = \left[ {\bf f}_1 \right]_{\ell} \cdot \left[ {\bf f}_2
\right]_{\ell} $ for
the Fourier-coefficients of $ {\bf f}_1, $ $ {\bf f}_2, $ $ {\bf f}  = {\bf
f}_1\ast {\bf f}_2. $ \par
\smallskip
{\bf Proof:} \par
The first statement $ (f\in{\cal O}_m \left(\dot {\cal J}^{ \left(\mu_ 1+\mu_
2 \right)}_+ \right)) $ is obtained by adding
to the proof of proposition 1 of $\lbrack$B.V-1$\rbrack$ a majorization of the
r.h.s. of Eq.(3) of this article, based on Eq.(2.1); this yields:
$$ \vert f(u+iv)\vert  \leq  {\rm e}^{mv} \int^{ }_{\dot  \gamma( \theta)} g_m
\left[f_1 \right] \left(v-v^{\prime} \right)g_m \left[f_2 \right]
\left(v^{\prime} \right) \left\vert {\rm d} \theta^{ \prime} \right\vert \ ,
$$
and by choosing a convenient representative of $ \dot  \gamma( \theta) $ (see
Fig.1 of
$\lbrack$B.V-1$\rbrack$), we obtain:
$$ g_m[f](v) = \dlowlim{ {\rm sup}}{u} \left\vert {\rm e}^{-mv}f(u+iv)
\right\vert \leq  {\rm Cst} \left\{ Y \left(\mu_ 1-v \right)+Y \left(v-\mu_ 1
\right) \left[g_m \left[f_1 \right]\diamond g_m \left[f_2 \right] \right](v)
\right\} $$
Since $ g_m \left[f_i \right] \in  L^1 \left(\Bbb R^+ \right) $ $ (i=1,2), $ $
g_m \left[f_1 \right] \diamond  g_m \left[f_2 \right] $ also belongs
to $ L^1 \left(\Bbb R^+ \right) $ and therefore $ g_m[f] \in L^1 \left(\Bbb
R^+ \right), $ which proves that
$ f \in  {\cal O}_m \left(\dot {\cal J}^{ \left(\mu_ 1+\mu_ 2 \right)}_+
\right). $ \par
On the other hand, Eq.(2.20) is given by the standard
theorem for the Laplace transform of a convolution product on $ \Bbb R^+ $
(applied here to $ {\rm f}  = {\rm f}_1 \diamond  {\rm f}_2). $ The last
statement of proposition 3
is then a direct consequence of proposition 1, iii) (Eqs.(2.6)) of
the present paper. \par
\smallskip
{\bf Remark:} \par
All the previous results remain true for the case $ \mu =0, $ up
to obvious additional specifications: in the definition of the
spaces $ {\cal O}_m \left(\dot {\cal J}^{(0)}_+ \right) $ etc..., the
continuity of the functions $ f(\theta) $ is
required on the boundary $ \dot {\Bbb R} \cup  \dot  \Xi( 0) $ of $ \dot
{\cal J}^{(0)}_+ $ from both
sides of $ \dot  \Xi( 0). $ Proposition 1 holds as it is written (the
exponential factor in Eq.(2.4) being now equal to 1), and
propositions 2, 2' remain true provided the assumption (2.7) (now
written for $ \mu =0) $ is enhanced by a condition of the following form:
$$ \dlowlim{ {\rm sup}}{\nu} \left\vert\tilde  F(\sigma +i\nu) \right\vert  =
M \left[\tilde  F \right](\sigma) ,\ \ {\rm with} \ \ \ M \left[\tilde  F
\right] \in  L^1([m,+\infty[) \ . \eqno (2.20^\prime ) $$
It is also easy to adapt the proof of proposition 1 of $\lbrack$B.V-1$\rbrack$
to
the case $ \mu =0 $ and to justify correspondingly the previous
proposition 3. \par
\medskip
{\bf Application: a theorem of dual analyticity for power series} \par
Let us specialize the statements of propositions 1, 2, 2',
3 to classes of functions $ f(\theta) $ $ (f\in{\cal O}_m \left(\dot {\cal
J}^{(\mu)}_ + \right)) $ for which: $ \forall \ell  < 0, $ $ [ {\bf
f}]_{\ell} =0. $
In this case, the decomposition $ f=h+r $ (see the proof of
proposition 2) is such that $ r $ is an entire function $ (r\not= 0 $ if $
m>0), $
and therefore $ f(\theta) $ is (like $ h(\theta) ) $ analytic in the
cut-domain $ \dot {\Bbb C}\backslash\dot  \Xi( \mu) . $
By putting $ {\rm e}^{-i\theta} =\zeta , $ $ {\rm e}^\mu =\alpha , $ and
defining the analytic function $ {\cal S}(\zeta) $ by $ {\cal S} \left( {\rm
e}^{-i\theta} \right) = f(\theta) , $
we can restate the results of propositions 1, 2, 2', 3 as
properties of the corresponding power series $ {\cal S}(\zeta)  = \sum^{ }_{
\ell \geq 0}a_{\ell} \zeta^{ \ell} , $ $ (\forall \ell , $
$ \ell \geq 0, $$ $ $ a_{\ell}  = {1 \over 2\pi}  [ {\bf f}]_{\ell} ). $ \par
\smallskip
{\bf Theorem}: (i) The following two properties of a power
series $ {\cal S}(\zeta)  = \sum^{ }_{ \ell \geq 0}a_{\ell} \zeta^{ \ell} $
are \lq\lq essentially\rq\rq\ equivalent (in the sense
specified below in (iii)): \par
\smallskip
a) $ {\cal S}(\zeta) $ is analytic in a cut-plane $ \Bbb C\backslash[ \alpha
,+\infty[ $ $ (\alpha \geq 1) $ and
uniformly bounded by $ {\rm Cst}\vert \zeta\vert^ m. $ \par
b) there exists an analytic function $ a(\lambda) $ in the
half-plane $ \Bbb C^{(m)}_+, $ uniformly bounded by $ {\rm Cst} \ \times
\alpha^{ - {\rm Re} \ \lambda} , $ such
that: $ \forall \ell  \geq m, $ $ a(\ell)  = a_{\ell} . $ \par
\smallskip
(ii) The function $ a(\lambda) $ is obtained as the Mellin transform
of the jump: \par
\noindent$ s(\zeta)  = {1 \over i} \dlowlim{ {\rm lim}}{
\matrix{\displaystyle \varepsilon \longrightarrow 0 \cr\displaystyle
\varepsilon >0 \cr}}  [{\cal S}(\zeta +i\varepsilon) -{\cal S}(\zeta
-i\varepsilon)] $ across the cut $ [\alpha ,+\infty[ , $
namely (for $ {\rm Re} \ \lambda  > m): $
$$ a(\lambda)  = {1 \over 2\pi}  \int^{ +\infty}_ \alpha s(\zeta)  \zeta^{
-\lambda -1} {\rm d} \zeta $$
\par
Conversely, one has (for all $ \zeta $ in $ \Bbb C\backslash[ \alpha
,+\infty[ ): $
$$ {\cal S}(\zeta)  = \sum^{ }_{ 0\leq \ell <m}a_{\ell} \zeta^{ \ell}  - {1
\over 2} \int^{ +\infty}_{ -\infty}{ a(m+i\nu) \times( -\zeta)^{ m+i\nu}
\over {\rm sin} \ \pi( m+i\nu)}  {\rm d} \nu $$
If $ m $ is integer, the integral is taken in the sense of
distributions with $ \dlowlim{ {\rm lim}}{ \matrix{\displaystyle \varepsilon
\longrightarrow 0 \cr\displaystyle \varepsilon >0 \cr}}  {1 \over {\rm sin} \
\pi( m-\varepsilon +i\nu)} . $ \par
\smallskip
(iii) The mapping $ {\cal S}(\zeta)  \longrightarrow  a(\lambda) $ satisfies a
continuity
inequality of the form
$$ N^{(m,\alpha)}_{ \infty}( a) \leq  {\rm Cst} \ \Vert{\cal S}\Vert^{( m)}_1\
, $$
for the following choice of norms:
$$ \Vert{\cal S}\Vert^{( m)}_1 = \int^{ \infty}_ 0 \dlowlim{ {\rm sup}}{\vert
\zeta\vert =\rho}  \left[\vert \zeta\vert^{ -m}\vert{\cal S}(\zeta)\vert
\right] { {\rm d} \rho \over \rho} \ ,N^{(m,\alpha)}_{ \infty}( a) =
\dlowlim{ {\rm sup}}{\nu}  \left[ \dlowlim{ {\rm sup}}{\sigma \geq m}
\alpha^{( \sigma -m)}\vert a(\sigma +i\nu)\vert \right]\ . $$
The inverse mapping $ a(\lambda)  \longrightarrow  {\cal S}(\zeta) $ satisfies
a continuity inequality
of the form:
$$ \Vert{\cal S}\Vert^{( m)}_{\infty}  \leq  {\rm Cst} \ N^{(m)}_1(a) $$
for the following choice of norms (with $ a_m(\nu) =a(m+i\nu) ) $:
$$ \Vert{\cal S}\Vert^ m_{\infty}  = \dlowlim{ {\rm sup}}{\zeta \in \Bbb C}
\left[\vert \zeta\vert^{ -m}\vert{\cal S}(\zeta)\vert \right]\ ,\
N^{(m)}_1(a) = \int^{ +\infty}_{ -\infty} \left\vert a_m(\nu) \right\vert
{\rm d} \nu \ \ \ {\rm if} \ \ m \in  \Bbb R-\Bbb Z $$
$$ {\rm (resp.}\Vert{\cal S}\Vert^{( m)}_{\infty} = \dlowlim{ {\rm
sup}}{\zeta \in \Bbb C} \left[(1+ {\rm \ell n}(1+\vert \zeta\vert))\vert
\zeta\vert^{ -m}\vert{\cal S}(\zeta)\vert \right]\ , $$
$$ N^{(m)}_1(a)= {\rm max} \left[ \int^{ +\infty}_{ -\infty} \left\vert
a_m(\nu) \right\vert {\rm d} \nu \ ,\ \int^{ +\infty}_{ -\infty} \left\vert{
{\rm d} a_m \over {\rm d} \nu}( \nu) \right\vert {\rm d} \nu \ ,\ V_{\nu_ 0}
\left({ {\rm d} a_m \over {\rm d} \nu} \right) \right]\ , $$
if $ m \in  \Bbb Z, $ with $ V_{\nu_ 0} $ defined as in proposition 2'). \par
\smallskip
(iv) \lq\lq Multiplicative convolution\rq\rq : \par
Let $ {\cal S}^{(i)}(\zeta)  = \sum^{ }_{ \ell \geq 0}a^{(i)}_{\ell} \zeta^{
\ell} , $ $ i=1,2, $ satisfy (i)a) and let $ {\cal S}(\zeta)  = \sum^{ }_{
\ell \geq 0}a_{\ell} \zeta^{ \ell} , $
with $ \forall \ell , $ $ a_{\ell}  = a^{(1)}_{\ell} \cdot a^{(2)}_{\ell} . $
Then one has:
$$ {\cal S}(\zeta)  = {1 \over 2i\pi}  \int^{ }_ \gamma{\cal S}^{(1)}
\left(\zeta /\zeta^{ \prime} \right){\cal S}^{(2)} \left(\zeta^{ \prime}
\right) { {\rm d} \zeta^{ \prime} \over \zeta^{ \prime}} \ , $$
where $ \gamma $ can be chosen as being the circle $ \left\{ \zeta^{ \prime} \
; \left\vert \zeta^{ \prime} \right\vert  = \vert \zeta\vert \right\} , $ for
$ \vert \zeta\vert  < \alpha , $
and
$$ \forall \zeta >0,\ \ \ \ {\cal S}(\zeta)  = {1 \over 2i\pi}  \int^{ \zeta
/\alpha_ 1}_{\alpha_ 2}{\cal S}^{(1)} \left(\zeta /\zeta^{ \prime}
\right){\cal S}^{(2)} \left(\zeta^{ \prime} \right) { {\rm d} \zeta^{ \prime}
\over \zeta^{ \prime}} \ . $$
Moreover, the corresponding transforms $ a^{(i)}(\lambda) $ $ (i=1,2), $ $
a(\lambda) $
satisfy in $ \Bbb C^{(m)}_+ $ the relation $ a(\lambda)  = a^{(1)}(\lambda)
\cdot a^{(2)}(\lambda) . $ \par
\smallskip
The theorem on power series presented in $\lbrack$S,W$\rbrack$ (theorem 3)
is a statement of strict equivalence for $ {\cal S}(\zeta) $ and $ a(\lambda)
$ (obtained
thanks to the use of $ L^2 $-norms) which corresponds to the case $ m=-1/2 $
and $ \alpha =1 $ of the previous theorem. \par
\bigskip
\vfill\eject
{\bf 2.2 Fourier-Laplace transformation for analytic functions
in classes }$ {\cal O}^{(d)}_m {\bf (}\dot {\cal J}^{(\mu)} {\bf )} : $ {\bf
the property }$ {\bf (F.G.)}_ {\bf o} $ \par
\smallskip
Let $ {\cal J}^{(\mu)} = \{ \theta \in \Bbb C; $ $ \theta \not= iv+2k\pi , $ $
\vert v\vert  \geq  \mu , $ $ k \in  \Bbb Z\} , $ where
$ \mu \geq 0. $ For
each integer $ d, $ with $ d\geq 2, $ we then define the following space $
{\cal O}^{(d)}_m \left(\dot {\cal J}^{(\mu)} \right) $
of analytic functions in $ \dot {\cal J}^{(\mu)} : $ each function in $ {\cal
O}^{(d)}_m \left(\dot {\cal J}^{(\mu)} \right) $ is
defined as the analytic continuation (in $ \dot {\cal J}^{(\mu)} ) $ of a
function $ f(\theta) $
in $ {\cal O}_m \left(\dot {\cal J}^{(\mu)}_ + \right) $ whose boundary value
on $ \dot {\Bbb R} $
satisfies the additional symmetry condition:
$$ \left( {\bf S}_d \right)\ \ \ \ \ \ \ \ \ \ {\bf f}(u) = (-1)^d {\rm
e}^{i(d-2)u} {\bf f}(-u) \eqno (2.21) $$
This analytic continuation (obtained by the \lq\lq Schwarz symmetry
principle\rq\rq\ or \lq\lq one-dimensional edge-of-the-wedge property\rq\rq\
$\lbrack$Pa$\rbrack$) is
still denoted by $ f(\theta) $ and satisfies in its domain $ \dot {\cal
J}^{(\mu)} $ the
functional relation:
$$ \left(S_d \right)\ \ \ \ \ \ \ \ \ \ f(\theta)  = (-1)^d {\rm
e}^{i(d-2)\theta} f(-\theta) \eqno (2.22) $$
\par
We define similarly the spaces $ \underline{{\cal O}}^{(d)}_m \left(\dot
{\cal J}^{(\mu)} \right) $ (resp. $ \underline{{\cal O}}^{(d)\ast}_ m
\left(\dot {\cal J}^{(\mu)} \right)) $
by taking the
functions $ f $ in $ \underline{{\cal O}}_m \left(\dot {\cal J}^{(\mu)}_ +
\right) $ (resp. $ \underline{{\cal O}}^\ast_ m \left(\dot {\cal J}^{(\mu)}_ +
\right)) $ which satisfy the condition $ \left(S_d \right). $ \par
For all these functions $ f, $ we still denote by $ {\rm f} =\Delta( f) $ the
jump of $ f $ across the cut $ \dot  \Xi( \mu) $ defined as in \S 2.1 (the
jump of $ f $
across the opposite cut $ \{ \theta  \in  \dot {\Bbb C}; $ $ -\theta \in\dot
\Xi( \mu)\} , $ determined
by Eq.(2.22), will not be used in the following). \par
In our treatment of the Fourier-Laplace transformation on
the $ (d-1) $-dimensional complex hyperboloid (see section 4), we will
apply the Fourier-Laplace transformation $ {\cal L} $ of section 2.1 to the
(corresponding) subspace $ {\cal O}^{(d)}_m \left(\dot {\cal J}^{(\mu)}
\right) $ of $ {\cal O}_m \left(\dot {\cal J}^{(\mu)}_ + \right). $ We will
keep the
same notations:
$$ \forall f \in  {\cal O}^{(d)}_m \left(\dot {\cal J}^{(\mu)}_ + \right)\ \
,\ \ \ {\cal L}(f) = L( {\rm f}) = \tilde  F\ , $$
with $ L $ and $ {\cal L} $ defined by Eqs.(2.3) and (2.5) and, by adapting
the
results of the previous study, we will obtain the following
statements. The inclusion of the special case $ \mu =0 $ in all these
statements will be obtained by applying the prescriptions given in
our remark after proposition 3 (in particular concerning the
additional condition (2.20') when it is needed, i.e. in
propositions 5 and 5'). \par
\medskip
{\bf Proposition 4} : For every function $ f $ in $ {\cal O}^{(d)}_m
\left(\dot {\cal J}^{(\mu)} \right), $ the
Fourier coefficients $ [ {\bf f}]_{\ell} $ of $ {\bf f}  = f_{\vert \dot {\Bbb
R}} $ satisfy the set of
relations
$$ \forall \ell \in \Bbb Z\ ,\ \ \ [ {\bf f}]_{\ell}  = (-1)^d[ {\bf
f}]_{-(\ell +d-2)} \eqno (2.23) $$
so that one can write
$$ {\bf f}(u)=  {\rm e}^{+i \left({d-2 \over 2} \right)(u-\pi )} \times  {1
\over 2\pi}  \sum^{ }_{ \ell \in \Bbb Z}[ {\bf f}]_{\ell} {\rm cos} \left[
\left(\ell  + {d-2 \over 2} \right)u- \left({d-2 \over 2} \right)\pi \right]
\eqno (2.24) $$
\par
Moreover, they are related to the Fourier-Laplace transform
$ \tilde  F = {\cal L}(f) $ of $ f $ by the following set of equations:
$$ \forall \ell  \in  \Bbb Z\ ,\ \ {\rm with} \ \ \ell  \geq  m\ ,\ \ [ {\bf
f}]_{\ell}  = \tilde  F(\ell) \eqno (2.25) $$
\par
{\bf Proof} : Eqs.(2.23) and (2.24) are a direct consequence of
Eq.(2.21). Eq.(2.25) has been proved in proposition 1
(Eq.(2.6)). \par
We shall now give a reciprocal statement similar to
proposition 2. \par
\medskip
{\bf Proposition 5} : Let $ {\bf f} $ be a continuous function on $ \dot {\Bbb
R} $
satisfying the symmetry condition $ \left( {\bf S}_d \right) $ (Eq.(2.21)) and
such that
its set of Fourier coefficients $ [ {\bf f}]_{\ell} $ admit an
analytic interpolation $ \tilde  F $ in $ \Bbb C^{(m)}_+ $ which
fulfils the conditions of proposition 2, $ m $ being in $ \Bbb R \backslash
\Bbb Z. $
Then, there exists a
unique analytic function $ f $ in $ \underline{{\cal O}}^{(d)}_m \left(\dot
{\cal J}^{(\mu)} \right) $ whose restriction to $ \dot {\Bbb R} $
is $ {\bf f} $ and such that $ \tilde  F={\cal L}(f). $ \par
Moreover, if $ m > E \left[- {d-2 \over 2} \right], $ the following inversion
formulae express $ {\rm f} , $ $ {\bf f} $
and $ f $ in terms of $ \tilde  F: $ \par
\smallskip
a) $ \forall d, $ $ d \geq  2: $ for $ v \in  \Bbb R^+, $
$$ \eqalignno{ {\rm f}(v) & = {1 \over 2\pi} \int^{ +\infty}_{ -\infty}\tilde
 F(m+i\nu) {\rm e}^{(m+i\nu) v} {\rm d} \nu   & (2.26) \cr  &  = {i {\rm e}^{-
\left({d-2 \over 2} \right)v} \over \pi}  \int^{ +\infty}_{ -\infty}\tilde
F(m+i\nu) {\rm sin} \left\{ \left[\nu -i \left(m+{d-2 \over 2} \right)
\right]v \right\} {\rm d} \nu &  (2.27) \cr  &  = { {\rm e}^{- \left({d-2
\over 2} \right)v} \over \pi}  \int^{ +\infty}_{ -\infty}\tilde  F(m+i\nu)
{\rm cos} \left\{ \left[\nu -i \left(m+{d-2 \over 2} \right) \right]v
\right\} {\rm d} \nu &  (2.28) \cr} $$
\par
b) $ \forall u\ , $ $ u\in[ 0,2\pi] , $
$$ {\bf f}(u) = {\rm e}^{i \left({d-2 \over 2} \right)(u-\pi)} \left[{i \over
2\pi}  \int^{ }_{{\cal C}}{\tilde  F(\lambda) {\rm cos} \left[ \left(\lambda
+{d-2 \over 2} \right)(u-\pi) \right] \over {\rm sin} \ \pi \lambda}  {\rm d}
\lambda  + {\cal P}^{(d)}_m[ {\bf f}](u) \right]; \eqno (2.29) $$
\par
c) $ \forall \theta , $ $ {\rm Re} \ \theta  \in  [0,2\pi] , $
$$ f(\theta)  = {\rm e}^{i \left({d-2 \over 2} \right)(\theta -\pi)} \left[{-1
\over 2\pi}  \int^{ +\infty}_{ -\infty}{\tilde  F(m+i\nu) {\rm cos} \left[
\left(m+{d-2 \over 2}+i\nu \right)(\theta -\pi) \right] \over {\rm sin} \
\pi( m+i\nu)}  {\rm d} \nu  + {\cal P}^{(d)}_m[ {\bf f}](\theta) \right];
\eqno (2.30) $$
In these formulae, $ {\cal C} $ denotes the contour represented on Fig.2 and $
{\cal P}^{(d)}_m[ {\bf f}](\theta) $ denotes the
trigonometric polynomial
$$ {\cal P}^{(d)}_m[ {\bf f}](\theta)  = {1 \over 2\pi} \sum^{ }_{
-m-d+2<\ell <m}[ {\bf f}]_{\ell} {\rm cos} \left[ \left(\ell  + {_{d-2} \over
2} \right)\theta - \left({d-2 \over 2} \right)\pi \right] \eqno (2.31) $$
\par
{\bf Proof} : Proposition 2 ensures the existence of the analytic function $ f
$ in $ \underline{{\cal O}}_m \left(\dot {\cal J}^{(\mu)}_ + \right) $
and therefore also in $ \underline{{\cal O}}^{(d)}_m \left(\dot {\cal
J}^{(\mu)} \right), $ since condition $ {\bf S}_d $ is satisfy by $ {\bf f} .
$
Now, if $ m > E \left[- {d-2 \over 2} \right], $ the decomposition $ f = h+r,
$ with $ h(\theta)  = {1 \over 2\pi}  \sum^{ }_{ \ell >m}[ {\bf f}]_{\ell}
{\rm e}^{-i\ell \theta} , $
introduced in the proof of proposition 2, can be written as follows, in view
of
Eq.(2.23):
$$ f(\theta)  = h(\theta)  + (-1)^d {\rm e}^{i(d-2)\theta} h(-\theta)  + {\rm
e}^{i \left({d-2 \over 2} \right)(\theta -\pi)}{\cal P}^{(d)}_m[ {\bf
f}](\theta) \ , \eqno (2.32) $$
(which shows explicitly that $ f $ satisfies Eq.(2.22)). \par
By making use of the expression (2.10) of $ h(u) $ $ (h = h_+ $ for $ u \in
[0,2\pi] $ and $ h = h_- $
for $ u \in  [-2\pi ,0]), $ Eq.(2.32) readily yields Eq.(2.29).
Eq.(2.30) is obtained similarly by making use of Eqs.(2.11)
and (2.12) (with $ u $ replaced by $ \theta , $ and $ {\rm Re} \ \theta $
varying respectively in $ [0,2\pi] $ and in
$ [-2\pi ,0]). $ Finally, Eq.(2.26) can be obtained by taking the
discontinuity across $ \{ \theta =iv;\ v\geq 0\} $
of the r.h.s. of Eq.(2.32), namely of $ h(\theta)  + (-1)^d {\rm
e}^{i(d-2)\theta} h(-\theta) ; $ this discontinuity can be
directly derived from Eq.(2.15). (Note that the jump of $ h(-\theta) $ across
this half-line is
equal to zero, but its expression through a vanishing Cauchy integral
justifies the r.h.s.
of Eqs.(2.27) and (2.28)).
Eqs.(2.27) and (2.28) can also be checked to follow from Eq.(2.30) (by
computing $ i \left(f_+(v)-f_-(v) \right)) $ respectively for $ d $ even and $
d $ odd. \par
In the case when $ m $ is an integer, we now have: \par
\medskip
{\bf Proposition 5'} : If the function $ {\bf f} $ satisfies the assumptions
of proposition 5, with $ m \in  \Bbb Z, $
and if moreover the function $ \tilde  F_m(\nu) =\tilde  F(m+i\nu) $ satisfies
the special regularity assumptions
$\alpha$), $\beta$), $\gamma$) of proposition 2', then the
conclusions of proposition 5 are still valid with the following specifications
(for $ m>- {d-2 \over 2}): $ \par
i) the analytic function $ f $ belongs to $ \underline{{\cal O}}^{(d)\ast}_ m
\left(\dot {\cal J}^{(\mu)} \right) $ \par
ii) the
integral at the r.h.s. of formula (2.30) must be understood in the sense of
distributions, with $ {1 \over {\rm sin} \ \pi( m+i\nu)}  = \dlowlim{ {\rm
lim}}{ \matrix{\displaystyle \varepsilon \longrightarrow 0 \cr\displaystyle
\varepsilon >0 \cr}}  {1 \over {\rm sin} \ \pi( m-\varepsilon +i\nu)} . $ \par
The proof is similar to that of proposition 5, except that the results of
proposition 2' (instead of proposition 2) must now be applied; moreover,
formula (2.32)
is now replaced by the following one:
$$ f(\theta)  = h_1(\theta) +(-1)^d {\rm e}^{i(d-2)\theta} h_1(-\theta)  +
{\rm e}^{i \left({d-2 \over 2} \right)(\theta -\pi)}{\cal P}^{(d)}_m[ {\bf
f}](\theta) \ ; \eqno (2.32^\prime ) $$
then, by using the property of $ h_1 $ described in the remark after the proof
of proposition
2' (see Eq.(2.17')) one justifies the occurrence of the $ \gq\gq \varepsilon
$-prescription\rq\rq\
in formula (2.30). \par
\medskip
{\bf Introduction of appropriate subspaces of functions} \par
\smallskip
For our applications of the previous properties in section 4, we will be led
(by
results of section 3, namely propositions 12 and 15) to consider subspaces of
functions $ f $
in $ {\cal O}^{(d)}_m \left(\dot {\cal J}^{(\mu)} \right) $ (resp. $
\underline{{\cal O}}^{(d)}_m \left(\dot {\cal J}^{(\mu)} \right)) $ whose
{\sl derivatives\/} $ f^{(r)} $ up to a certain order $ p $
belong themselves to $ {\cal O}_m \left(\dot {\cal J}^{(\mu)}_ + \right) $
(resp. $ \underline{{\cal O}}_m \left(\dot {\cal J}^{(\mu)}_ + \right)). $ It
is then appropriate to equip these
subspaces with a normed space structure, expressed in
terms of the functions $ g_m \left[f^{(r)} \right] $ (see Eq.(2.1)). \par
\smallskip
a) {\sl Subspaces\/} $ {\cal O}^{(d)}_{m,p} \left(\dot {\cal J}^{(\mu)}
\right): $ \par
We put:
$$ \mid \mid \mid f\mid \mid \mid^{( 1)}_{m,p} = \dlowlim{ {\rm max}}{0\leq
r\leq p} \left\Vert g_m \left[f^{(r)} \right] \right\Vert_ 1\ ; \eqno (2.33)
$$
with this definition, $ {\cal O}^{(d)}_{m,0} \left(\dot {\cal J}^{(\mu)}
\right) $ is identical with $ {\cal O}^{(d)}_m \left(\dot {\cal J}^{(\mu)}
\right), $ equipped
with the norm:
$$ \left\Vert g_m[f] \right\Vert_ 1 = \mid \mid \mid f\mid \mid \mid^{( 1)}_m
\equiv  \mid \mid \mid f\mid \mid \mid^{( 1)}_{m,0} \eqno (2.33^\prime ) $$
\par
b) {\sl Subspaces\/} $ \underline{{\cal O}}^{(d)}_{m,p} \left(\dot {\cal
J}^{(\mu)} \right): $ \par
We put:
$$ \mid \mid \mid f\mid \mid \mid^{( \infty)}_{ m,p} = \dlowlim{ {\rm
max}}{0\leq r\leq p} \left\Vert g_m \left[f^{(r)} \right] \right\Vert_{
\infty} \ , \eqno (2.34) $$
so that in particular $ \underline{{\cal O}}^{(d)}_{m,0} \left(\dot {\cal
J}^{(\mu)} \right) \equiv  \underline{{\cal O}}^{(d)}_m \left(\dot {\cal
J}^{(\mu)} \right), $ equipped with the norm:
$$ \left\Vert g_m[f] \right\Vert_{ \infty}  = \mid \mid \mid f\mid \mid
\mid^{( \infty)}_ m \equiv  \mid \mid \mid f\mid \mid \mid^{( \infty)}_{ m,0}\
. \eqno (2.34^\prime ) $$
\par
We also introduce correspondingly the following Banach spaces of functions $
\tilde  F(\lambda) $
analytic in $ \Bbb C^{(m)}_+, $ continuous on $ \bar {\Bbb C}^{(m)}_+: $ \par
\smallskip
c) {\sl Subspaces\/} $ \tilde {\cal O}^{\mu ,p} \left(\Bbb C^{(m)}_+ \right):
$ \par
We put:
$$ \mid \mid \mid\tilde  F\mid \mid \mid^{ m\mu p}_{(1)} = \dlowlim{ {\rm
max}}{0\leq r\leq p} \left\Vert G^{(\mu)}_ m \left[\lambda^ r\cdot\tilde  F
\right] \right\Vert_ 1\ , \eqno (2.35) $$
where (as in Eq.(2.7)):
$$ G^{(\mu)}_ m \left[\lambda^ r\tilde  F \right](\nu)  = \dlowlim{ {\rm
sup}}{\sigma \geq m} {\rm e}^{(\sigma -m)\mu} \left\vert( \sigma +i\nu)^
r\tilde  F(\sigma +i\nu) \right\vert \ . \eqno (2.36) $$
\par
{\sl d) Subspaces\/} $ \underline{{\cal }\tilde  O \underline{}}^{\mu ,p}
\left(\Bbb C^{(m)}_+ \right): $ \par
We put:
$$ \mid \mid \mid\tilde  F\mid \mid \mid^{ m\mu p}_{(\infty)}  = \dlowlim{
{\rm max}}{0\leq r\leq p} \left\Vert G^{(\mu)}_ m \left[\lambda^ r\tilde  F
\right] \right\Vert_{ \infty} \ . \eqno (2.37) $$
\par
For the case $ m\in \Bbb Z, $ we shall make use of the following subspaces:
\par
\smallskip
{\sl e) Subspaces\/} $ \underline{{\cal O}}^{(d)\ast}_{ m,p} \left(\dot {\cal
J}^{(\mu)} \right): $ \par
We put
$$ \mid \mid \mid f\mid \mid \mid^{( \infty) \ast}_{ m,p} = \dlowlim{ {\rm
max}}{0\leq r\leq p} \left\Vert g^\ast_ m \left[f^{(r)} \right] \right\Vert_{
\infty} \ , \eqno (2.38) $$
\par
{\sl f) Subspaces\/} $ \tilde {\cal O}^{\mu ,p}_\ast \left(\Bbb C^{(m)}_+
\right): $ \par
We consider the subspaces of functions $ \tilde  F $ in $ \tilde {\cal
O}^{\mu ,p} \left(\Bbb C^{(m)}_+ \right) $ such that $ G^{(\mu)}_ m \left[{
{\rm d} \over {\rm d} \lambda} \left(\lambda^ r\tilde  F \right) \right]\in
L^1(\Bbb R), $
$ 0\leq r\leq p, $ and which moreover
satisfy the assumption $\beta$) of proposition 2', and we put:
$$ \mid \mid \mid\tilde  F\mid \mid \mid^{ m,\mu ,p}_{(1)\ast}  = \mid \mid
\mid\tilde  F\mid \mid \mid^{ m\mu p}_1+\mid \mid \mid{ {\rm d}\tilde  F
\over {\rm d} \lambda} \mid \mid \mid^{ m\mu p}_1 + V_{\nu_ 0} \left({ {\rm
d}\tilde  F_m \over {\rm d} \nu} \right)\ . \eqno (2.38^\prime ) $$
\par
We can now complete the results of propositions 4, 5 and 5' by the following
statements: \par
\medskip
{\bf Proposition 6:} The transformation $ {\cal L} $ defines a continuous
mapping from each
subspace $ {\cal O}^{(d)}_{m,p} \left(\dot {\cal J}^{(\mu)} \right) $ into the
corresponding subspace $ \underline{{\cal }\tilde  O \underline{}}^{\mu ,p}
\left(\Bbb C^{(m)}_+ \right). $ \par
\medskip
{\bf Proof:} \par
If $ f \in  {\cal O}^{(d)}_{m,p} \left(\dot {\cal J}^{(\mu)} \right), $ all
the functions $ f^{(r)} $ $ (0 \leq  r \leq  p) $ belong to $ {\cal O}_m
\left(\dot {\cal J}^{(\mu)}_ + \right) $ and
yield the respective transforms: $ L \left(\Delta f^{(r)} \right) = \lambda^
r\tilde  F(\lambda) ; $ by applying proposition 1 to
the latter, we can write the corresponding majorizations (2.4) as follows (in
view of Eq.(2.36)):
$$ {\rm for} \ \ 0 \leq r\leq p\ ,\ \ \ \left\Vert G^{(\mu)}_ m
\left[\lambda^ r\tilde  F \right] \right\Vert_{ \infty}  \leq  2 \left\Vert
g_m \left[f^{(r)} \right] \right\Vert_ 1\ , $$
which implies, in view of Eqs.(2.33) and (2.37):
$$ \mid \mid \mid\tilde  F\mid \mid \mid^{ m\mu p}_{(\infty)} \leq  2 \mid
\mid \mid f\mid \mid \mid^{( 1)}_{m,p} $$
\par
\smallskip
{\bf Remark :} \par
If $ f $ belongs to $ {\cal O}^{(d)}_{m,p} \left(\dot {\cal J}^{(\mu)}
\right)\cap  \underline{{\cal O}}^{(d)}_{m,p} \left(\dot {\cal J}^{(\mu)}
\right), $ all the functions $ f^{(r)}(0\leq r\leq p) $ are
such that $ {\rm e}^{-mv}\Delta f^{(r)}(v) $ belong to $ L^1 \left(\Bbb R^+
\right)\cap L^{\infty} \left(R^+ \right) $ and therefore to $ L^2 \left(\Bbb
R^+ \right). $ It
follows that all the transforms $ \lambda^ r\tilde  F(\lambda)  = L
\left(\Delta f^{(r)} \right) $ belong correspondingly to the Hardy
space $ {\bf H}^2 \left(\Bbb C^{(m)}_+ \right) = \left\{\tilde  h(\lambda)
\right. $ holomorphic in $ \Bbb C^{(m)}_+ $ such that: $ \left. \dlowlim{
{\rm sup}}{\sigma \geq m} \int^{ +\infty}_{ -\infty} \left\vert\tilde
h(\sigma +i\nu) \right\vert^ 2 {\rm d} \nu <\infty \right\} $
(see e.g. chapter 8 of $\lbrack$Ru$\rbrack$). \par
\medskip
{\bf The transformation} $ [{\cal L}]^{-1}_{(d,m)} $ \par
In the following statements, which concern the inversion of the transformation
$ {\cal L} $
in the subspaces $ {\cal O}^{(d)}_m \left(\dot {\cal J}^{(\mu)} \right) $ (or
$ \underline{{\cal O}}^{(d)}_m \left(\dot {\cal J}^{(\mu)} \right)), $ we will
assume (as in proposition 5, 5') that the real number $ m $ is such that $ m >
E \left[-{d-2 \over 2} \right]. $
This will be sufficient for the applications of these results in
section 4. \par
Since the Laplace transformation $ L $ is an injective mapping, the kernel $
{\cal N}^{(d)}_m $
(resp. $ \underline{{\cal N}}^{(d)}_m) $ of $ {\cal L} $ in $ {\cal O}^{(d)}_m
\left(\dot {\cal J}^{(\mu)} \right) $ (resp. $ \underline{{\cal O}}^{(d)}_m
\left(\dot {\cal J}^{(\mu)} \right)) $ is the set of functions
$ f $ in this subspace such that: $ {\rm f}  = \Delta f = 0. $ For such
functions $ f, $ the corresponding
function $ h(\theta) $ or $ h_1(\theta) $ (defined in terms of $ \tilde
F(\lambda) $ by Eqs.(2.13) (2.13'), (2.19')) is
equal to zero (since $ \tilde  F = {\cal L}(f) = 0), $ according to whether $
m\in \Bbb R\backslash \Bbb Z $ or $ m\in \Bbb Z. $
Then, in view of Eqs.(2.31), (2.32), (2.32'), we conclude that:
$$ \matrix{\displaystyle \hfill{\cal N}^{(d)}_m= \underline{{\cal N}}^{(d)}_m=
\left\{ f(\theta)  = {\rm e}^{i \left({d-2 \over 2} \right)(\theta -\pi)} P(
{\rm cos} \ \theta) ,\ \ {\rm with} \right. \cr\displaystyle \hfill P( {\rm
cos} \ \theta) = \left. \sum^{ }_{ -{d-2 \over 2}\leq \ell <m}a_{\ell} {\rm
cos} \left[ \left(\ell +{d-2 \over 2} \right)\theta - \left({d-2 \over 2}
\right)\pi \right] \right\} \cr} $$
\par
Let us now denote by $ {\cal M}^{(d)}_m \left(\dot {\cal J}^{(\mu)} \right) $
(resp. $ \underline{{\cal M}}^{(d)}_m \left(\dot {\cal J}^{(\mu)} \right)) $
the subspace of functions $ f $
in $ {\cal O}^{(d)}_m \left(\dot {\cal J}^{(\mu)} \right) $ (resp. $
\underline{{\cal O}}^{(d)}_m \left(\dot {\cal J}^{(\mu)} \right)) $ such that
$ {\cal P}^{(d)}_m[ {\bf f}]=0, $ i.e. such that $ [ {\bf f}]_{\ell} =0 $ for
$ -m-d+2<\ell <m. $ (We have $ {\cal M}^{(d)}_m\approx{\cal O}^{(d)}_m/{\cal
N}^{(d)}_m $ and $ {\cal M}^{(d)}_m={\cal O}^{(d)}_m $
if $ \left\vert m+{d-2 \over 2} \right\vert  < {1 \over 2}, $ with $ d $ odd).
Then, the transformation $ {\cal L}, $ considered as acting on $ {\cal
M}^{(d)}_m $
has an inverse, which we call $ [{\cal L}]^{-1}_{(d,m)}; $ the explicit
expression of the image $ [{\cal L}]^{-1}_{(d,m)} \left(\tilde  F \right) $ of
a
function $ \tilde  F $ is then given by Eq.(2.30), in which the last term $
{\cal P}^{(d)}_m[ {\bf f}](\theta) $ has
been dropped. We can now state: \par
\medskip
{\bf Proposition 7:} If $ m > E \left[- {d-2 \over 2} \right], $ $ m\in  \Bbb
R\backslash \Bbb Z, $ the transformation $ [{\cal L}]^{-1}_{(d,m)} $
defines a continuous mapping from each subspace $ \tilde {\cal O}^{\mu ,p}
\left(\Bbb C^{(m)}_+ \right) $ into the corresponding
subspace $ \underline{{\cal M}}^{(d)}_{m,p} \left(\dot {\cal J}^{(\mu)}
\right) = \underline{{\cal M}}^{(d)}_m \left(\dot {\cal J}^{(\mu)}
\right)\cap  \underline{{\cal O}}^{(d)}_{m,p} \left(\dot {\cal J}^{(\mu)}
\right). $ \par
\medskip
{\bf Proof}: \par
In view of the previous definition of $ [{\cal L}]^{-1}_{(d,m)}, $ we can
write for every $ \tilde  F \in  \tilde {\cal O}^{\mu ,p} \left(\Bbb C^{(m)}_+
\right) $
(according to Eq.(2.32)):
$$ f(\theta)  = [{\cal L}]^{-1}_{(d,m)} \left(\tilde  F \right)(\theta)  =
h(\theta)  + (-1)^d {\rm e}^{i(d-2)\theta} h(-\theta) \ . \eqno (2.39) $$
The majorization (2.14) of $ h $ implies (in view of Eq.(2.1)):
$$ \left\Vert g_m[h] \right\Vert_{ \infty}  \leq  {1 \over 2\pi \mid {\rm sin}
\ \pi m\mid} \left\Vert G^{(\mu)}_ m \left[\tilde  F \right] \right\Vert_ 1\ ;
\eqno (2.40) $$
on the other hand, we have:
$$ \forall v \geq  0,\ \ g_m \left[ {\rm e}^{i(d-2)\theta} h(-\theta)
\right](v) = {\rm e}^{-(m+d-2)v} \dlowlim{ {\rm sup}}{u}\vert h(u-iv)\vert ,
$$
and in view of the majorization (2.14) (now used for $ h \left(u+iv^{\prime}
\right), $ with $ v^{\prime} =-v<0): $
$$ g_m \left[ {\rm e}^{i(d-2)\theta} h(-\theta) \right](v) \leq  { {\rm
e}^{-(2m+d-2)v} \over 2\pi \mid {\rm sin} \ \pi m\mid}  \left\Vert G^{(\mu)}_
m \left[\tilde  F \right] \right\Vert_ 1 \eqno (2.41) $$
Let us first consider the case when $ m > - {d-2 \over 2}: $ in this case,
formulae (2.39), (2.40) and
(2.41) directly imply the following inequality:
$$ \left\Vert g_m[f] \right\Vert_{ \infty}  \leq  {1 \over \pi \mid {\rm sin}
\ \pi m\mid}  \left\Vert G^{(\mu)}_ m \left[\tilde  F \right] \right\Vert_ 1
\eqno (2.42) $$
\par
By applying a similar argument to the derivatives $ f^{(r)}(\theta)  = [{\cal
L}]^{-1}_{(d,m)} \left(\lambda^ r\tilde  F \right)(\theta) $ (for $ 1 \leq  r
\leq  p), $
we obtain (since $ f^{(r)}(\theta)  = h^{(r)}(\theta)  + {\rm
e}^{i(d-2)\theta}  \sum^{ }_{ 0\leq q\leq r}c_qh^{(q)}(-\theta) , $ with
appropriate
constants $ c_q): $
$$ \left\Vert g_m \left[f^{(r)} \right] \right\Vert_{ \infty}  \leq  { {\rm
Cst} \over \mid {\rm sin} \ \pi m\mid}  \dlowlim{ {\rm max}}{0\leq q\leq r}
\left\Vert G^{(\mu)}_ m \left[\lambda^ q\tilde  F \right] \right\Vert_ 1
\eqno (2.43) $$
\par
The set of inequalities (2.42) and (2.43) (for $ r \leq  p) $ then yield (in
view of
Eqs.(2.34) and (2.35)):
$$ \mid \mid \mid[{\cal L}]^{-1}_{(d,m)} \left(\tilde  F \right)\mid \mid
\mid^{( \infty)}_{ m,p} \leq  { {\rm Cst} \over \mid {\rm sin} \ \pi m\mid}
\mid \mid \mid\tilde  F\mid \mid \mid^{ m\mu p}_1 \eqno (2.44) $$
\par
We now consider the remaining case: $ E \left[-{d-2 \over 2} \right] < m \leq
- {d-2 \over 2} $ (with $ d $ odd). In this
case, the majorization (2.41) can be improved by considering the function $
{\rm e}^{-m^{\prime} v}h_+(u+iv), $
with $ m^{\prime}  = -m -d + 2 \geq  m $ (instead of $ {\rm e}^{-mv}h_+(u+iv))
$ and by applying to the latter the
corresponding inequality (2.14) (which amounts to replace $ {\cal C} $ by the
homotopous contour $ {\rm L}_{m^{\prime}} , $
instead of $ {\rm L}_m $ in the integral (2.10)). We then obtain in this way:
$$ \matrix{\displaystyle g_m \left[ {\rm e}^{i(d-2)\theta} h(-\theta)
\right](v) & \displaystyle \leq  { {\rm e}^{- \left(m+m^{\prime} +d-2
\right)v} \over 2\pi \left\vert {\rm sin} \ \pi m^{\prime} \right\vert}
\left\Vert G^{(\mu)}_ m \left[\tilde  F \right] \right\Vert_ 1 \hfill
\cr\displaystyle  & \displaystyle = {1 \over 2\pi\vert {\rm sin} \ \pi
m\vert}  \left\Vert G^{(\mu)}_ m \left[\tilde  F \right] \right\Vert_ 1
\hfill \cr} \eqno (2.41^\prime ) $$
Inequalities (2.40), (2.41') again imply the inequalities (2.42), (2.43)
and therefore (2.44). \par
\medskip
{\bf Proposition 7'}: If $ m > - {d-2 \over 2}, $ with $ m \in  \Bbb Z, $ the
transformation $ [{\cal L}]^{-1}_{(d,m)} $
defines a continuous mapping from each subspace $ \tilde {\cal O}^{\mu
,p}_\ast \left(\Bbb C^{(m)}_+ \right) $ into the corresponding
subspace $ \underline{{\cal M}}^{(d)\ast}_{ m,p} \left(\dot {\cal J}^{(\mu)}
\right)= \underline{{\cal M}}^{(d)}_m \left(\dot {\cal J}^{(\mu)} \right)\cap
\underline{{\cal O}}^{(d)\ast}_{ m,p} \left(\dot {\cal J}^{(\mu)} \right). $
\par
The argument is similar to that of proposition 7, since one now has $ [{\cal
L}]^{-1}_{(d,m)} \left(\tilde  F \right)(\theta)  = h_1(\theta) +(-1)^d {\rm
e}^{i(d-2)\theta} h_1(-\theta) $
(instead of Eq.(2.39)); the majorization (2.18') of $ h_1 $ now plays the same
role as the
majorization (2.14) of $ h $ throughout the proof (since
$ m > - {d-2 \over 2}), $ and one obtains (in view of Eq.(2.38')) a continuity
inequality of the form:
$$ \mid \mid \mid[{\cal L}]^{-1}_{(d,m)} \left(\tilde  F \right)\mid \mid
\mid^{ \infty \ast}_{ m,p} \leq  {\rm Cst} \ \mid \mid \mid\tilde  F\mid \mid
\mid^{ m\mu p}_{(1)\ast} \eqno (2.45) $$
(In the derivation of the latter, one must take into account the obvious
inequalities:
$ \forall r,\ \ \ V_{\nu_ 0} \left({ {\rm d} \left(\lambda^ r\tilde  F
\right)_m \over {\rm d} \nu} \right) \leq  {\rm Cst} \ V_{\nu_ 0} \left({
{\rm d}\tilde  F_m \over {\rm d} \nu} \right)). $ \par
\smallskip
According to our terminology of section 1, the set of results stated in the
present subsection (propositions 4,...,7') can be called \lq\lq property $ (
{\rm F.G.})_{\circ} $ with symmetry
condition $ \left(S_d \right)\dq\dq ; $ in section 4, they will be used for
proving the theorem F.G. in the
general $ d $-dimensional case, $ d\geq 3; $ in our next subsection (2.3) we
will treat in a
straightforward way the case $ d=2. $ \par
\medskip
\noindent{\bf 2.3 The case }$ d=2 ${\bf : harmonic analysis for invariant
perikernels on
the complex hyperbola} $ X^{(c)}_1 $ \par
\smallskip
In \S 1.3 of $\lbrack$B.V-1$\rbrack$, we had introduced the algebra of
perikernels $ {\cal K} \left(z,z^{\prime} \right) $ on the complex hyperbola:
$$ X^{(c)}_1 = \left\{ z= \left(z^{(0)},z^{(1)} \right)=z(\theta)  \in  \Bbb
C^2\ ;\ z^{(0)}=-i\ {\rm sin\ } \theta ,\ z^{(1)}= {\rm cos} \ \theta ,\
\theta \in \dot {\Bbb C} \right\} $$
\par
We consider $ X^{(c)}_1 $ as a homogeneous space of the group $ {\rm
SO}(1,1)^{(c)} $
and we are interested in the class of $ {\rm SO}(1,1)^{(c)} $-invariant
perikernels $ {\cal K} $ $ ({\cal K} \left(gz,gz^{\prime} \right) = {\cal K}
\left(z,z^{\prime} \right), $
for all $ g $ in $ {\rm SO}(1,1)^{(c)}). $ In \S 1.3 of
$\lbrack$B.V-1$\rbrack$, we only mentioned the case of
perikernels invariant under the connected component $ {\rm
SO}_0(1,1)^{(c)}=\{ g=g(\alpha) ; $ $ \alpha \in \dot {\Bbb C}=\Bbb C/2\pi
\Bbb Z\} $
of $ {\rm SO}(1,1)^{(c)} $ whose action on $ X^{(c)}_1 $ is defined by $
g(\alpha) z(\theta) =z(\theta +\alpha) ; $ by
fixing the \lq\lq base point\rq\rq\ $ z_0=z(0)=(0,1) $ of $ X^{(c)}_1, $ we
can identify each such perikernel $ {\cal K} $ with an analytic function $
{\cal F}(z)={\cal K} \left(z,z_0 \right), $
also represented in the $ \theta $-plane by $ f(\theta)  = {\cal
F}(z(\theta)) $ (in fact: $ {\cal K} \left(z(\theta) ,z \left(\theta^{
\prime} \right) \right) = f \left(\theta -\theta^{ \prime} \right)). $
If $ {\cal K} $ is invariant under the full group $ {\rm SO}(1,1)^{(c)}, $ it
is in particular invariant
under the symmetry $ J $ with respect to the $ z^{(1)} $-axis in $ \Bbb C^2, $
which implies
that $ {\cal F}(z) = {\cal F}(Jz), $ and therefore $ f(\theta)  = f(-\theta) .
$ \par
The analyticity domain $ X^{(c)}_1\times X^{(c)}_1\backslash \Sigma^{(
c)}_\mu $ of $ {\cal K} $ and the corresponding domain $ D_\mu
=X^{(c)}_1\backslash X^\mu $
of $ {\cal F} $ (see section 1) have the special property (for $ d=2) $ that
the cuts $ \Sigma^{( c)}_\mu $ and
$ X^\mu $ have two connected components, namely: $ \Sigma^{( c)}_\mu  =
\Sigma^{( c)+}_\mu  \cup  \Sigma^{( c)-}_\mu , $ where
$$ \Sigma^{( c)\pm}_ \mu  = \left\{ \left(z,z^{\prime} \right); z=z(\theta) ,\
z^{\prime} =z \left(\theta^{ \prime} \right);\ \theta -\theta^{ \prime}  = iv\
;\ \ \ \pm v\geq \mu \right\} \ ; $$
$ X^\mu  = X^\mu_ +\cup X^\mu_ -, $ where:
$$ X^\mu_{ \pm}  = \{ z;\ z=z(\theta) ;\ \theta =iv\ ,\ \ \pm v\geq \mu\}
\eqno (2.46) $$
We note that $ X^\mu_ + $ and $ X^\mu_ - $ are infinite arcs of the branch
containing $ z_0 $ of the real
hyperbola $ X_1 = X^{(c)}_1\cap  \Bbb R^2 $ (see Fig.3). \par
The corresponding function $ f(\theta)  = {\cal F}(z(\theta)) $ is then
analytic and even in
the domain $ \dot {\cal J}^{(\mu)} $ introduced in \S 2.2; it is also
convenient to introduce the
representation $ \underline{f} \left(z^{(1)} \right) = \underline{f}( {\rm
cos} \ \theta)  = f(\theta) $ of $ {\cal F}, $ $ \underline{f} $ being
analytic in the
cut-plane $ \underline{D}_\mu  = \Bbb C\backslash[ {\rm cosh} \ \mu ,\
+\infty[ . $ \par
Functions $ {\cal F} $ of moderate growth of order $ m $ in $ D_\mu $ will be
introduced as
follows: \par
\medskip
{\bf Definition 1}: For each $ J $-invariant function $ {\cal F}(z) $ analytic
in a given domain $ D_\mu , $ one
puts:
$$ \forall m,\ \ {\cal G}_m[{\cal F}](\rho)  = \rho^{ -m} \dlowlim{ {\rm
sup}}{ \left\{ z;z^{(1)}\in E_\rho \right\}}  \vert{\cal F}(z)\vert =\rho^{
-m} \dlowlim{ {\rm sup}}{z^{(1)}\in E_\rho} \left\vert \underline{f}_1
\left(z^{(1)} \right) \right\vert \eqno (2.47) $$
where $ \rho  \geq  1 $ and $ E_\rho $ is the ellipse $ \left\{ z^{(1)}= {\rm
cos}(u+iv); \right. $ $ {\rm cosh} \ v=\rho\} . $ \par
\smallskip
a) We call $ \left[{\cal V}_1 \right]^m_\mu \left(X^{(c)}_1 \right) $ the
space of such functions $ {\cal F} $ for which:
$$ \Vert{\cal F}\Vert^{( 1)}_m = \int^{ \infty}_ 1{\cal G}_m[{\cal F}](\rho)
{ {\rm d} \rho \over \sqrt{ \rho^ 2-1}} < \infty \eqno (2.48) $$
\par
b) We call $ \left[{\cal V}_{\infty} \right]^m_\mu \left(X^{(c)}_1 \right) $
the space of such functions $ {\cal F} $ for which:
$$ \Vert{\cal F}\Vert^{( \infty)}_ m = \dlowlim{ {\rm sup}}{\rho \geq 1}
{\cal G}_m[{\cal F}](\rho)  < \infty \ . \eqno (2.49) $$
\par
For $ m $ integer, we also use the space $ \left[{\cal V}^\ast_{ \infty}
\right]^m_\mu \left(X^{(c)}_1 \right) $ defined by the norm:
$$ \Vert{\cal F}\Vert^{( \infty) \ast}_ m = \dlowlim{ {\rm sup}}{\rho \geq 1}
(1+ {\rm \ell n} \ \rho)  {\cal G}_m[{\cal F}](\rho)  < \infty \ . \eqno
(2.50) $$
By using the representation $ f(\theta)  = {\cal F}(z(\theta)) , $ also called
in the following the
isomorphism $ i: $ $ {\cal F} \charlvmidup{ \ i\ }{ \rightarrowfill} f, $ we
will obtain homeomorphisms from the latter spaces into
the corresponding spaces $ {\cal O}^{(2)}_m \left(\dot {\cal J}^{(\mu)}
\right), $ $ \underline{{\cal O}}^{(2)}_m \left(\dot {\cal J}^{(\mu)} \right)
$ and $ \underline{{\cal O}}^{(2)\ast}_ m \left(\dot {\cal J}^{(\mu)} \right)
$ $ (\equiv \underline{{\cal O}}^{(2)\ast}_{ m,0} \left(\dot {\cal J}^{(\mu)}
\right), $
see Eq.(2.38)) introduced in \S 2.2, since the symmetry condition $ \left(S_2
\right), $ i.e. $ f(\theta) =f(-\theta) $
is satisfied in the present case. \par
\medskip
{\bf Lemma 1:} The following homeomorphisms hold:
$$ \left[{\cal V}_1 \right]^m_\mu \left(X^{(c)}_1 \right)\approx{\cal
O}^{(2)}_m \left(\dot {\cal J}^{(\mu)} \right),\ \left[{\cal V}_{\infty}
\right]^m_\mu \left(X^{(c)}_1 \right)\approx \underline{{\cal O}}^{(2)}_m
\left(\dot {\cal J}^{(\mu)} \right),\ \left[{\cal V}^\ast_{ \infty}
\right]^m_\mu \left(X^{(c)}_1 \right)\approx \underline{{\cal O}}^{(2)\ast}_ m
\left(\dot {\cal J}^{(\mu)} \right). $$
\par
\smallskip
{\bf Proof :} \par
The following inequalities, resulting from Eqs.(2.1) and (2.47):
$$ \forall v \geq 0,\ \ {\rm Cst} \ g_m[f](v) \leq  {\cal G}_m[{\cal F}](
{\rm cosh} \ v) \leq {\rm  Cst}  ^{\prime} g_m[f](v) $$
(and $ {\rm Cst}(1+v)g_m[f](v) \leq  [1+ {\rm \ell n}( {\rm cosh} \ v)]{\cal
G}_m[{\cal F}]( {\rm cosh} \ v) \leq  {\rm Cst}^{\prime}( 1+v)g_m[f](v)) $
imply
(in view of Eqs.(2.33'), (2.34'), (2.38), (2.48), (2.49), (2.50)) the
following
norm equivalences:
$$ \Vert{\cal F}\Vert^{( 1)}_m\sim \mid \mid \mid f\mid \mid \mid^{( 1)}_m,\ \
\Vert{\cal F}\Vert^{( \infty)}_ m \sim  \mid \mid \mid f\mid \mid \mid^{(
\infty)}_ m,\ \ \Vert{\cal F}\Vert^{( \infty) \ast}_ m\sim \mid \mid \mid
f\mid \mid \mid^{( \infty) \ast}_{ m,0}\ . $$
\par
\smallskip
{\sl J-invariant triplets of moderate growth\/} \par
\smallskip
With any function $ {\cal F} $ in a space $ \left[{\cal V}_1 \right]^m_\mu
\left(X^{(c)}_1 \right) $ or $ \left[{\cal V}_{\infty} \right]^m_\mu
\left(X^{(c)}_1 \right) $ is associated
a $ \gq\gq J $-invariant triplet of moderate growth\rq\rq\ $ ({\cal F}, {\bf
F} ,F) $ (representing the \lq\lq invariant triplet\rq\rq\ $ ({\cal K}, {\bf
K} ,K) $
such that $ {\cal F}(z) = {\cal K} \left(z,z_0 \right)), $ in which: \par
\smallskip
a) $ {\bf F}(z) = {\bf K} \left(z,z_0 \right) $ is the restriction of $ {\cal
F} $ to the circle
$$ S_1= \left\{ z\in X^{(c)}_1\ ;\ \ z=z(\theta) \ ,\ \ \theta  = u\in  \dot
{\Bbb R} \right\} \ , $$
namely $ \forall u, $ $ {\bf F}(z(u)) = {\bf f}(u) $ (with $ {\bf f}  =
f_{\vert \dot {\Bbb R}}). $ \par
\smallskip
b) $ F(z) = K \left(z,z_0 \right) = \Delta{\cal F}(z) $ is the jump of $
{\cal F} $ across the cut $ X^\mu_ + $ (see
Eq.(2.46)), namely:
$$ \forall v\geq 0\ ,\ F(z(iv)) = {\rm f}(v) = \Delta f(v) = i\ \dlowlim{
{\rm lim}}{ \matrix{\displaystyle \varepsilon \longrightarrow 0 \hfill
\cr\displaystyle \varepsilon >0 \hfill \cr}}  [{\cal F}(z(\varepsilon
+iv))-{\cal F}(z(-\varepsilon +iv))] $$
\par
\smallskip
{\bf The \lq\lq theorem F.G.\rq\rq\ for }$ d=2 ${\bf \nobreak\ :} \par
\smallskip
We are now in a position to transpose all the results of \S 2.2 (i.e.
\lq\lq property $ {\rm (F.G)}_{\circ} \dq\dq $ for $ d=2) $ by the isomorphism
$ ({\cal F}, {\bf F} ,F) \charlvmidup{ (i)}{ \rightarrowfill}  (f, {\bf f} ,
{\rm f}); $ we thus obtain
corresponding statements (theorems 1, 3, 3', 4, 5) which describe the various
features of what we called \lq\lq the theorem F.G\rq\rq\ for the case $ d=2 $
(see section 1)
supplemented by the convolution property (theorem 2). \par
\smallskip
We introduce the Fourier-Laplace transformation $ {\cal L}_2 $ for $ J
$-invariant
triplets $ ({\cal F}, {\bf F} ,F) $ of moderate growth on $ X^{(c)}_1 $ by the
following definition: $ {\cal L}_2 = {\cal L }\circ  i, $ where $ {\cal L} $
is the
transformation defined in proposition 1. \par
\smallskip
We then have: \par
\medskip
{\bf Theorem 1 }$ (d=2)\ : $ With every $ J $-invariant triplet $ ({\cal F},
{\bf F} ,F) $ on $ X^{(c)}_1, $ such
that $ {\cal F }\in  \left[{\cal V}_1 \right]^m_\mu \left(X^{(c)}_1 \right)  $
$ (m\in  \Bbb R), $ one can associate: \par
\smallskip
a) the transform $ \tilde  F = {\cal L}_2({\cal F}), $ which can be expressed
as follows in terms
of the function $ \underline{{\rm f}} \left(z^{(1)} \right) = F(z) $ $ (z \in
X^\mu_ +): $
$$ \tilde  F(\lambda)  = \int^{ +\infty}_{ {\rm cosh} \ \mu} \left(\zeta -
\sqrt{ \zeta^ 2-1} \right)^\lambda  \underline{{\rm f}}(\zeta)  { {\rm d}
\zeta \over \sqrt{ \zeta^ 2-1}} \eqno (2.51) $$
\par
b) the set of Fourier coefficients $ [ {\bf f}]_{\ell} $ of $ {\bf F} $ on $
S_1: $
$$ \forall \ell \in \Bbb Z,\ [ {\bf f}]_{\ell}  = [ {\bf f}]_{-\ell}  =
\int^{ +\pi}_{ -\pi} {\rm cos} \ \ell u\ \underline{{\bf f}}( {\rm cos} \ u)
{\rm d} u\ , \eqno (2.51^\prime ) $$
where $ \underline{{\bf f}} \left(z^{(1)} \right) = {\bf F}(z) $ $ (z \in
S_1). $ \par
\smallskip
Then the following properties hold: \par
\smallskip
i) $ \tilde  F \in  \tilde { \underline{{\cal O}} }^{\mu ,0} \left(\Bbb
C^{(m)}_+ \right) $ and the transformation $ {\cal L}_2 $ defines a continuous
mapping from $ \left[{\cal V}_1 \right]^m_\mu \left(X^{(c)}_1 \right) $ into $
\tilde { \underline{{\cal O}} }^{\mu ,0} \left(\Bbb C^{(m)}_+ \right). $ \par
\smallskip
ii) $ \forall \ell  \in  \Bbb Z,\ \ \ \ \ \ell  \geq  m,\ \ \ \ \ [ {\bf
f}]_{\ell}  = \tilde  F(\ell)
$ \par
{\bf Proof :} \par
Eq.(2.51) reduces to Eq.(2.3) (i.e. $ \tilde  F = L( {\rm f}) = {\cal L}(f) =
{\cal L }\circ  i({\cal F})), $
after the
identification $ \zeta = {\rm cosh} \ v, $ $ \underline{{\rm f}}( {\rm cosh} \
v) = {\rm f}(v), $ and Eq.(2.51') reduces to
Eq.(2.2) (with $ \underline{{\bf f}}( {\rm cos} \ u) = {\bf f}(u)). $ Then,
properties i) and ii) just express
the results of proposition 1 (or 3, with $ d=2), $ completed (as far as the
continuity of the mapping $ {\cal L}_2 $ is concerned) by the results of
proposition 6
(case $ d=2) $ and of lemma 1. \par
\medskip
{\bf Theorem 2} $ (d=2)\ : $ Let $ \left({\cal K}_1, {\bf K}_1,K_1 \right) $
and $ \left({\cal K}_2, {\bf K}_2,K_2 \right) $ be two invariant
triplets on $ X^{(c)}_1 $ such that $ {\cal F}_i(z) = {\cal K}_i \left(z,z_0
\right) \in  \left[{\cal V}_1 \right]^m_{\mu_ i} \left(X^{(c)}_1 \right), $ $
i=1,2. $
Then the invariant triplet $ ({\cal K}, {\bf K} ,K) $ obtained by convolution
of the latter,
namely:
$$ {\cal K }= {\cal K}_1\ast^{( c)}{\cal K}_2\ ,\ \ \ {\bf K} = {\bf K}_1\ast
{\bf K}_2\ ,\ \ \ K = K_1 \diamond  K_2\ , $$
is such that: $ {\cal F}(z) = {\cal K} \left(z,z_0 \right) \in  \left[{\cal
V}_1 \right]^m_{\mu_ 1+\mu_ 2} \left(X^{(c)}_1 \right)\ . $ \par
Moreover, the corresponding transforms $ \tilde  F_i = {\cal L}_2 \left({\cal
F}_i \right), $ $ \tilde  F = {\cal L}_2({\cal F}) $ satisfy the
relation $ \tilde  F(\lambda)  = \tilde  F_1(\lambda) \cdot\tilde
F_2(\lambda) , $ which interpolates in $ \bar {\Bbb C}^{(m)}_+, $ the
corresponding set of relations: $ [ {\bf f}]_{\ell} = \left[ {\bf f}_1
\right]_{\ell} \cdot \left[ {\bf f}_2 \right]_{\ell} $ for the Fourier
coefficients of
$ {\bf F} , $ $ {\bf F}_1 {\bf ,} ${\bf\ }$ {\bf F}_2 {\bf .} $ \par
\medskip
{\bf Proof\nobreak\ :} \par
This is a direct corollary of proposition 3 and of lemma 1; in fact, as
it is explained at the end of section 1 of $\lbrack$B.V-1$\rbrack$, the
composition product $ {\cal K}_1\ast^{( c)}{\cal K}_2 $
(see $\lbrack$B.V-1$\rbrack$, formula (18)) reduces in the case of invariant
perikernels to the
convolution product $ f = f_1\ast^{( c)}f_2, $ and therefore proposition 3
applies; moreover,
one immediately checks that, since $ {\bf f}_1 $ and $ {\bf f}_2 $ satisfy the
symmetry condition $ \left( {\bf S}_2 \right), $
the same holds for $ {\bf f}  = {\bf f}_1\ast {\bf f}_2, $ which implies that
$ f \in  {\cal O}^{(2)}_m \left(\dot {\cal J}^{ \left(\mu_ 1+\mu_ 2 \right)}
\right) $ and therefore
$ {\cal F }\in  \left[{\cal V}_1 \right]^m_{\mu_ 1+\mu_ 2} \left(X^{(c)}_1
\right). $ \par
\smallskip
{\bf Theorem 3} $ (d=2)\ : $ Let $ {\bf F}(z) $ be a function on $ S_1, $ such
that $ {\bf F}(z) = {\bf F}(Jz), $
whose set of Fourier coefficients $ [ {\bf f}]_{\ell} $ (Eq.(2.51')) admit an
analytic
interpolation $ \tilde  F(\lambda) $ in $ \Bbb C^{(m)}_+ $ such that: \par
\smallskip
a) $ m \in  \Bbb R\backslash  \Bbb Z $ and $ \tilde  F \in  \tilde {\cal
O}^{\mu ,0} \left(\Bbb C^{(m)}_+ \right) $ \par
\smallskip
b) $ \forall \ell , $ $ \ell >m, $ $ \tilde  F(\ell)  = [ {\bf f}]_{\ell} $
\par
\smallskip
Then $ {\bf F} $ admits an analytic continuation $ {\cal F} $ in the domain $
D_\mu $ such that: \par
\smallskip
i) $ {\cal F }\in  \left[{\cal V}_{\infty} \right]^m_\mu \left(X^{(c)}_1
\right) $ \par
\smallskip
ii) $ \tilde  F = {\cal L}_2({\cal F}) $ \par
\medskip
{\bf Proof\nobreak\ :} \par
This results from the first part of proposition 5 (case $ d=2) $ and from
lemma 1. \par
\medskip
Similarly proposition 5' (case $ d=2) $ and lemma 1 yield: \par
\medskip
{\bf Theorem 3'} $ (d=2)\ : $ If the assumptions of theorem 3 are satisfied
with the
following conditions (replacing a) and b)): \par
\smallskip
a') $ m\in \Bbb Z $ and $ \tilde  F\in\tilde {\cal O}^{\mu ,0}_\ast \left(\Bbb
C^{(m)}_+ \right) $ \par
b') $ \forall \ell , $ $ \ell  \geq m, $ $ \tilde  F(\ell)  = [ {\bf
f}]_{\ell} , $ \par
\smallskip
\noindent then the conclusions of theorem 3 remain valid, with the additional
specification: \par
\smallskip
i') $ {\cal F}\in \left[{\cal V}^\ast_{ \infty} \right]^m_\mu \left(X^{(c)}_1
\right). $ \par
\medskip
{\bf Theorem 4} $ (d=2)\ : $ Let $ ({\cal F}, {\bf F} ,F) $ be a $ J
$-invariant triplet whose transforms
$ \tilde  F={\cal L}_2({\cal F}) $ and $ \left\{[ {\bf f}]_{\ell} ;\ \ell \in
\Bbb Z \right\} $ satisfy the conditions a), b) of theorem 3, with $ m\in \Bbb
R^+\backslash \Bbb N. $
Then the following inversion formulae hold: \par
\smallskip
a) $ \forall z = z(u), $ $ z^{(1)}= {\rm cos} \ u\in[ -1,+1]\ , $
$$ {\bf F}(z(u)) = \underline{{\bf f}}( {\rm cos} \ u)={1 \over 2\pi} \sum^{
}_{ \ell \in \Bbb Z}[ {\bf f}]_{\ell} {\rm cos} \ \ell \theta \eqno (2.52) $$
\par
b) $ \forall z=z(iv), $ $ z^{(1)}= {\rm cosh} \ v \geq 1\ , $
$$ F(z(iv)) = \underline{{\rm f}}( {\rm cosh} \ v) = {1 \over \pi}  \int^{
+\infty}_{ -\infty}\tilde  F(m+i\nu) {\rm cos}[(\nu -im)v] {\rm d} \nu \eqno
(2.53) $$
\par
c) $ \forall z=z(\theta) , $ $ z^{(1)}= {\rm cos} \ \theta  \in
\underline{D}_\mu , $
$$ {\cal F}(z(\theta))  = \underline{f}( {\rm cos} \ \theta)  = -{1 \over
2\pi} \int^{ +\infty}_{ -\infty}{\tilde  F(m+i\nu) {\rm cos}[(m+i\nu)( \theta
-\pi)] \over {\rm sin} \ \pi( m+i\nu)} {\rm d} \nu +{\cal P}_m[ {\bf
F}](\theta) \eqno (2.54) $$
where
$$ {\cal P}_m[ {\bf F}](\theta)  = {1 \over 2\pi}  \sum^{ }_{\vert \ell\vert
<m}[ {\bf f}]_{\ell} {\rm cos} \ \ell \theta \ . \eqno (2.55) $$
\par
\smallskip
In the case when $ m \in  \Bbb N, $ formula (2.54) must be understood in the
sense of distributions with $ {1 \over {\rm sin} \ \pi( m+i\nu)} $ replaced $
\dlowlim{ {\rm lim}}{ \matrix{\displaystyle \varepsilon \longrightarrow 0
\hfill \cr\displaystyle \varepsilon >0 \hfill \cr}}  {1 \over {\rm sin} \
\pi( m-\varepsilon +i\nu)} , $ and the
special conditions a'), b') of theorem 3' must be assumed. \par
\medskip
{\bf Proof\nobreak\ :} \par
This is a direct corollary of formulae (2.28), (2.30), (2.31) of
proposition 5 (written for $ d=2), $ the last statement being justified by
proposition 5'. \par
\smallskip
Let us now introduce, for each $ m $ in $ \Bbb R^+, $ the transformation $
\left[{\cal L}_2 \right]^{-1}_m=i^{-1}\circ[{\cal L}]^{-1}_{(2,m)}, $
(with $ [{\cal L}]^{-1}_{(2,m)} $ defined at the end of \S 2.2). $
\left[{\cal L}_2 \right]^{-1}_m $ acts as a \lq\lq quasi-inverse\rq\rq\ of
the Fourier-Laplace transformation $ {\cal L}_2={\cal L}\circ i $ in the sense
that: \par
\smallskip
i) $ {\cal L}_2\circ \left[{\cal L}_2 \right]^{-1}_m=\uniset , $ on
appropriate subspaces of analytic functions $ \tilde  F(\lambda) $ in
the domain $ \Bbb C^{(m)}_+. $ \par
ii) $ \left[{\cal L}_2 \right]^{-1}_m\circ{\cal L}_2=\Pi_ m, $ where $ \Pi_ m
$ denotes the projection which associates with
any function $ {\cal F} $ the truncated function:
$$ \left[\Pi_ m{\cal F} \right](z(\theta))  = {\cal F}(z(\theta))  - {\cal
P}_m[ {\bf F}](\theta) $$
\par
As a direct corollary of propositions 7, 7' and of lemma 1, we can
state: \par
\medskip
{\bf Theorem 5} $ (d=2)\ : $ \par
a) Let $ m \in  \Bbb R^+\backslash \Bbb N; $ then the transformation $
\left[{\cal L}_2 \right]^{-1}_m $ defines
a continuous mapping $ (\tilde  F \longrightarrow{\cal F}= \left[{\cal L}_2
\right]^{-1}_m \left(\tilde  F \right)) $ from each subspace $ \tilde {\cal
O}^{\mu ,0} \left(\Bbb C^{(m)}_+ \right) $ into the
corresponding subspace $ \left[{\cal V}_{\infty} \right]^m_\mu \left(X^{(c)}_1
\right). $ \par
Moreover, for each $ \tilde  F\in\tilde {\cal O}^{\mu ,0} \left(\Bbb C^{(m)}_+
\right), $ the triplet $ ({\cal F}, {\bf F} ,F) $ associated
with $ {\cal F }= \left[{\cal L}_2 \right]^{-1}_m \left(\tilde  F \right) $
can be computed via the formulae (2.52), (2.53), (2.54) of
theorem 4, in which the additional condition $ {\cal P}_m[ {\bf F}]=0 $ must
be imposed. \par
\smallskip
b) Let $ m \in  \Bbb N, $ $ m >0; $ then the transformation $ \left[{\cal L}_2
\right]^{-1}_m $ defines a
continuous mapping from each subspace $ \tilde {\cal O}^{\mu ,0}_\ast
\left(\Bbb C^{(m)}_+ \right) $ into the corresponding subspace
$ \left[{\cal V}^\ast_{ \infty} \right]^m_\mu \left(X^{(c)}_1 \right) $ (the
last statement of a) being unchanged, provided the prescriptions
given at the end of theorem 4 be now applied). \par
\bigskip
\noindent{\bf 3. HOROCYCLES AND RADON-ABEL TRANSFORMATIONS ON THE
(REAL AND COMPLEXIFIED) ONE-SHEETED HYPERBOLOIDS } \par
\medskip
{\bf 3.1 }$ {\bf H}^{( {\bf c})} ${\bf -invariant functions on }$ {\bf X}^{(
{\bf c})}_{ {\bf d-1}}; ${\bf\ invariant
perikernels.
} \par
\smallskip
For $ d\geq 3, $ we consider the space $ \Bbb C^d $ of variables $ z\equiv
\left(z^{(0)},\vec  z,z^{(d-1)} \right), $
with $ \vec  z = \left(z^{(1)},...,z^{(d-2)} \right) $ and the complexified
one-sheeted
hyperboloid $ X^{(c)}_{d-1} $ with equation:
$$ s(z)\equiv z^{(0)^2} - \vec  z^2 - z^{(d-1)^2}+1=0\ ,\ {\rm where\ \ }\vec
 z^2 = z^{(1)^2}+\cdot \cdot \cdot +z^{(d-2)^2}. {\rm \ } $$
\par
Two real submanifolds of $ X^{(c)}_{d-1} $
will play an important role below, namely the one-sheeted
hyperboloid $ X_{d-1} = \Bbb R^d\cap X^{(c)}_{d-1} $ and the sphere $ S_{d-1}=
\left(i\Bbb R \times  \Bbb R^{d-1} \right)\cap X^{(c)}_{d-1} $
$ (i\Bbb R $ referring to the coordinate $ z^{(0)} $ of $ z). $ We also
introduce the meridian sections of $ X^{(c)}_{d-1}, $ $ X_{d-1} $ and $
S_{d-1} $ in the $ (z^{(0)}, $
$ z^{(d-1)}) $-coordinate plane, namely respectively the complex and
real hyperbolae $ \hat  X^{(c)} $ and $ \hat  X $$ $with equation: $ \hat
s(z)\equiv z^{(0)^2}-z^{(d-1)^2}+1=0 $ and
the circle $ \hat  S $ (with equations: $ z^{(0)}=iy^{(0)}, $ $
z^{(d-1)}=x^{(d-1)}, $ $ y^{(0)^2}+x^{(d-1)^2}=1); $
the point $ z_0 = \left(0,\vec  0,1 \right) $ in $ (\hat  X\cap\hat  S) $ will
be called the \lq\lq base point\rq\rq . \par
If we introduce the complexified Lorentz group $ G^{(c)}\equiv {\rm
SO}(1;d-1)^{(c)}, $
considered as acting on $ \Bbb C^{(d)}, $ and the subgroup $ H^{(c)} $ of $
G^{(c)}, $
isomorphic to $ {\rm SO}(1,d-2)^{(c)}, $ which is the stabilizer of the
point $ z_0, $ $ X^{(c)}_{d-1} $  can be identified with the homogeneous space
$ G^{(c)}/H^{(c)}. $
We identify correspondingly $ X_{d-1} $ with the real homogeneous
space $ G/H, $ where $ G $ is the Lorentz group $ {\rm SO}(1;d-1) $ and $
H=H^{(c)}\cap G. $
Let $ \bar  \omega $ be the projection from $ X^{(c)}_{d-1} $ onto the $
z^{(d-1)} $-plane. For each $ z^{(d-1)}\not= \pm 1, $
$ \bar  \omega^{ -1} \left(z^{(d-1)} \right) $ is a $ (d-2) $-dimensional
complex hyperboloid which is a
generic orbit of $ H^{(c)} $ on $ X^{(c)}_{d-1}. $ In the following, we shall
be
concerned with $ H^{(c)} $-invariant functions $ {\cal F}(z) $ (i.e. such that
$ {\cal F}(z) = {\cal F}(hz) $
for all $ h\in H^{(c)}) $ analytic in appropriate domains $ D $ of $
X^{(c)}_{d-1}. $ Each
such function $ {\cal F} $ only depends on the coordinate $ z^{(d-1)} $ of $ z
$ and is
therefore represented by an analytic function $ \underline{f} \left(z^{(d-1)}
\right) $ in the
domain $ \underline{D} $ of $ \Bbb C $ such that $ D = \bar  \omega^{ -1}(
\underline{D}). $ It will also be
convenient to put $ z^{(d-1)}= {\rm cos} \ \theta $ and to represent
alternatively $ {\cal F}(z) = \underline{f} \left(z^{(d-1)} \right) $
by an {\sl even\/} $ 2\pi $-periodic function $ f(\theta) = \underline{f}(
{\rm cos} \ \theta) $ analytic in the
domain $ \dot {\cal J}(D) = \{ \theta \in  \dot {\Bbb C};\ {\rm cos} \ \theta
\in \underline{D}\} . $ \par
We call $ {\cal D} $
the class of $ H^{(c)} $-invariant domains $ D $ which satisfy the
following conditions: \par
i) $ D $ contains the sphere $ S_{d-1}. $ \par
ii) the corresponding domain $ \underline{D} $
is star-shaped (in $ \Bbb C) $ with respect to the point $ z^{(d-1)}=1. $ \par
Condition i) implies in particular that $ [-1,+1] \subset \underline{D} . $
\par
\smallskip
We shall be concerned more specifically by the subclass
of domains
$$ D_\mu = \left\{ z\in X^{(c)}_{d-1}\ ;\ z^{(d-1)}= {\rm cos} \ \theta  \in
\underline{D}_\mu \right\} \ \ {\rm with} \ \ \underline{D}_\mu =\Bbb
C\backslash[ {\rm cosh} \ \mu ,+\infty[ \ (\mu \geq 0) $$
also represented in the variable $ \theta $ by the cut-plane $ \dot {\cal
J}^{(\mu)}  = \dot {\cal J} \left(D_\mu \right), $
introduced in \S 2.2. \par
\smallskip
We shall call $ D_\mu $ a \lq\lq cut-domain\rq\rq , since it is bordered by
the set $ X^\mu  = X^{(c)}_{d-1}\backslash D_\mu  = \left\{ z\in X^{(c)}_{d-1}
\right.; $ $ \left.z^{(d-1)}_{ }= {\rm cosh} \ v;\ \vert v\vert  \geq  \mu
\right\} $
represented by the cut $ \{ \theta  = iv+2k\pi , $ $ \vert v\vert  \geq  \mu ,
$ $ k\in \Bbb Z\} , $ which is
the boundary of $ \dot {\cal J}^{(\mu)} . $ The set $ X^\mu $ admits two {\sl
real\/} connected
components, which are the following regions of $ X, $ (see Fig.4):
$$ \matrix{\displaystyle X^\mu_ + & \displaystyle = \left\{ z\in X_{d-1}\ ;\
z^{(0)} > 0\ ,\ z^{(d-1)} \geq  {\rm cosh} \ \mu \right\} \hfill&
\displaystyle {\rm and} \cr\displaystyle X^\mu_ - & \displaystyle = \left\{
z\in X_{d-1}\ ;\ z^{(0)} < 0\ ,\ z^{(d-1)} \geq  {\rm cosh} \ \mu \right\} \ .
\hfill& \displaystyle \cr} $$
\par
The region $ X_+ \left(\equiv X^0_+ \right)= \{ z\in X; $ $ z^{(0)}\geq 0, $ $
\left.z^{(d-1)}\geq 1 \right\} $ can be
simply characterized in terms of the \lq\lq causal ordering relation\rq\rq\
on $ X_{d-1} $ induced by the Minkowskian quadratic form in $ \Bbb C^d: $ $
z^2 = z^{(0)^2}-z^{(1)^2}-\cdot \cdot \cdot z^{(d-1)^2}; $
two points $ z, $ $ z^{\prime} $ in $ X_{d-1} $ satisfy the ordering relation
$ z \geq  z^{\prime} $ if and
only if $ \left(z-z^{\prime} \right)^2 \geq  0 $ and $ z^{(0)} \geq
z^{\prime( 0)}. $ We then see that $ X_+= \left\{ z\in X_{d-1}; \right. $ $
\left.z\geq z_0 \right\} $
(see Fig.4). \par
\smallskip
With every $ H^{(c)} $-invariant function $ {\cal F}(z) $ analytic in a
domain $ D_\mu , $ we associate the \lq\lq discontinuity function\rq\rq\ $
\Delta{\cal F }= i \left({\cal F}_+-{\cal F}_- \right), $
where $ {\cal F}_+ $ and $ {\cal F}_- $ denote the boundary values of $ {\cal
F} $ on the analytic
hypersurface $ \left\{ z\in X^{(c)}_{d-1} \right.; $ $ \left. {\rm Im}
\left(z-z_0 \right)^2=0 \right\} $ (containing the \lq\lq cut\rq\rq\ $ \Xi_
\mu ) $ from
the respective sides $ {\rm Im} \left(z-z_0 \right)^2 < 0 $ and $ {\rm Im}
\left(z-z_0 \right)^2>0. $ \par
The function $ \Delta{\cal F} $ is an $ H^{(c)} $-invariant function with
support $ X^\mu $ and is entirely characterized by the function $ F(z) = Y
\left(z^{(0)} \right)\Delta{\cal F}_{ \left\vert X_{d-1} \right.}(z) $
$ (Y $ denoting the Heaviside function), which is an $ H $-invariant
function on $ X_{d-1} $ with support contained in $ X^\mu_ + $ (in fact, for
every
$ z\in X^\mu , $ one can write $ z=hz^{\prime} , $ with $ h\in H^{(c)} $ and $
z^{\prime} \in X^\mu_ +, $ so that: $ \Delta{\cal F}(z) = \Delta{\cal F}
\left(z^{\prime} \right)=F \left(z^{\prime} \right)). $
The function $ f(\theta) $ (resp. $ \underline{f} \left(z^{(d-1)} \right)) $
which represents $ {\cal F} $ is analytic in the cut-plane $ \dot {\cal
J}^{(\mu)} $ (resp. in $ \underline{D}_\mu ), $
and one checks that $ F $ is represented by the
corresponding jump function $ {\rm f}(v) = i \dlowlim{ {\rm lim}}{
\matrix{\displaystyle \eta \longrightarrow 0 \cr\displaystyle \eta >0 \cr}}
[f(iv+\eta) -f(iv-\eta)] $
with support $ \{ v;v\geq \mu\} $ (resp. by the jump function $
\underline{{\rm f}}( {\rm cosh} \ v) = i \dlowlim{ {\rm lim}}{
\matrix{\displaystyle \varepsilon \longrightarrow 0 \cr\displaystyle
\varepsilon >0 \cr}}[ \underline{f}( {\rm cosh} \ v-i\varepsilon) -
\underline{f}( {\rm cosh} \ v+i\varepsilon)] ). $ \par
As it was explained in $\lbrack$B.V-1$\rbrack$, every $ H^{(c)} $-invariant
function $ {\cal F} $ analytic in a domain $ D_\mu $ represents a $ G^{(c)}
$-invariant perikernel $ {\cal K} \left(z,z^{\prime} \right) $ $ (\forall
g\in G^{(c)}, $ $ {\cal K} \left(z,z^{\prime} \right) = {\cal K}
\left(gz,gz^{\prime} \right)) $
analytic in the corresponding domain $ \Delta_ \mu  = X^{(c)}_{d-1}\times
X^{(c)}_{d-1}\backslash \Sigma^{( c)}_\mu , $ with $ \Sigma^{( c)}_\mu =
\left\{ \left(z;z^{\prime} \right)\in X^{(c)}_{d-1}\times X^{(c)}_{d-1};
\right. $
$ \left. \left(z-z^{\prime} \right)^2\geq 2( {\rm cosh} \ \mu -1) \right\} ; $
this follows from the identification
relation $ {\cal F}(z) = {\cal K} \left(z,z_0 \right) $ and from the fact that
$ D_\mu  = \left\{ z\in X^{(c)}_{d-1}; \right. $ $ \left. \left(z,z_0
\right)\in \Delta_ \mu \right\} . $
\par
Similarly, the jump function $ F $ associated with $ {\cal F} $
represents a $ G $-{\sl invariant\/} Volterra kernel $ K \left(z,z^{\prime}
\right) $ on $ X_{d-1}, $ with
support contained in the set $ \Sigma^ +_\mu  = \left\{ \left(z,z^{\prime}
\right); \right. $ $ \left(z-z^{\prime} \right)^2\geq  2( {\rm cosh} \ \mu
-1), $ $ \left.z^{(0)}>z^{\prime( 0)} \right\} , $
namely one has $ F(z) = K \left(z,z_0 \right), $ $ K $ being also the jump
function $ \Delta_ +{\cal K} $
of the perikernel $ {\cal K} $ represented by $ {\cal F} $ (see
$\lbrack$B,V-1$\rbrack$ Eq.(34)). The
Volterra kernels form an algebra $ V \left(X_{d-1} \right) $ for the
composition
product $ K_1\diamond K_2 $ (see $\lbrack$Fa,V$\rbrack$ and
$\lbrack$B,V-1$\rbrack$), while the {\sl invariant\/}
Volterra kernels form a subalgebra $ V \left(X_{d-1} \right)^{\natural} $ of
the latter, which
will be considered in \S 3.3. \par
\medskip
{\bf 3.2 Horosphere fibration of the real and complex
one-sheeted hyperboloids} \par
\smallskip
Let $ {\cal P} $ be the family of hyperplanes $ P_\tau $ with equation: $
z^{(0)}+z^{(d-1)}= {\rm e}^{-i\tau} , $
for $ \tau \in \dot {\Bbb C} = \Bbb C/2\pi  \Bbb Z. $ Each hyperplane $
P_\tau $
intersects the hyperbola $ \hat  X^{(c)} $ at the (unique) point $ z_\tau $ $
\left(z^{(0)}_\tau  = -i {\rm sin} \ \tau , \right. $
$ \left.z^{(d-1)}_\tau  = {\rm cos} \ \tau \right) $ and this defines a
bijection from the set $ {\cal P}= \left\{ P_\tau ;\tau \in \dot {\Bbb C}
\right\} $
onto $ \hat  X^{(c)} = \left\{ z_\tau ;\tau \in \dot {\Bbb C} \right\} . $ For
each $ \tau $ in $ \dot {\Bbb C}, $ $ P_\tau  \cap  X^{(c)}_{d-1} $
is a complex $ (d-2) $-dimensional paraboloid $ \Pi_ \tau $ which we call a
{\sl complex horosphere.\/} We notice that the set of horospheres $ \left\{
\Pi_ \tau ;\tau \in \dot {\Bbb C} \right\} $
defines a fibration with basis $ \hat  X^{(c)} $ on the (dense) domain $
X^{\prime( c)}_{d-1}= \left\{ z\in X^{(c)}_{d-1}; \right. $
$ \left.z^{(0)}+z^{(d-1)}\not= 0 \right\} $ of $ X^{(c)}_{d-1}. $ (In
fact, the section of $ X^{(c)}_{d-1} $ by the hyperplane $ z^{(0)}+z^{(d-1)}=0
$ yields
a cylinder instead of a parabolic horosphere). We denote by $ h $ the
projection associated with this fibration: $ \forall z \in  X^{(c)}_{d-1}, $ $
z_\tau  = h(z) $
is the intersection of $ \hat  X^{(c)} $ with the (unique) horosphere $ \Pi_
\tau $ which
contains $ z $ and $ \Pi_ \tau  = h^{-1} \left(z_\tau \right). $ \par
We now associate with this fibration the
following parametric representation of $ X^{\prime( c)}_{d-1}: $
$$ \eqalignnol{  &  z^{(0)} & = -i {\rm sin} \ \tau  + {1 \over 2} \vec
\zeta^ 2\ {\rm e}  ^{-i\tau}   &  &  & (3.1a) \cr z=z \left(\vec  \zeta ,\tau
\right)\ \ \  & \vec  z & \equiv  \left(z^{(1)},...,z^{(d-2)} \right) = \vec
\zeta \ {\rm e}^{-i\tau} &  &  &  (3.1b) \cr  &  z^{(d-1)} & = {\rm cos} \
\tau  - {1 \over 2} \vec  \zeta^ 2\ {\rm e}^{-i\tau} \ , &  &  & (3.1c) \cr}
$$
with $ \tau  \in  \dot {\Bbb C} $ and $ \vec  \zeta  = \left(\zeta^{(
1)},...,\zeta^{( d-2)} \right) \in  \Bbb C^{d-2}; $
in Eqs.(3.1a) and (3.1c), we make use of the notation: $ \vec  \zeta^ 2 =
\zeta^{( 1)^2}+\cdot \cdot \cdot +\zeta^{( d-2)^2}. $
For each fixed number $ \tau , $ Eqs.(3.1) clearly represent the complex
horosphere $ \Pi_ \tau $ (i.e. $ z^{(0)}+z^{(d-1)} = {\rm e}^{-i\tau} ; $ $
{\rm e}^{-i\tau} \left(z^{(0)}-z^{(d+1)} \right) - \vec  z^2 + 1 = 0), $
parametrized by the complex vector $ \vec  \zeta . $ $ (\vec  \zeta ,\tau ) $
will be called
the {\sl horocyclic coordinates\/} of the point $ z $ varying on $ X^{\prime(
c)}_{d-1}. $ \par
An alternative group-theoretical presentation of the
horocyclic coordinates relies on the introduction of the
following abelian groups of $ (d\times d) $-matrices:
$$ \eqalignnol{ A^{(c)} & = \left\{ a(\tau)  = \left[ \matrix{\displaystyle
{\rm cos} \ \tau & \displaystyle 0 & \displaystyle -i {\rm sin} \ \tau
\cr\displaystyle 0 & \displaystyle \Bbb I_{(d-2)} & \displaystyle 0
\cr\displaystyle -i {\rm sin} \ \tau & \displaystyle 0 & \displaystyle {\rm
cos} \ \tau \cr} \right]\ ;\ \tau  \in  \dot {\Bbb C} \right\} &  &  &  &
(3.2) \cr N^{(c)} & = \left\{ n \left(\vec  \zeta \right) = \left[
\matrix{\displaystyle 1+{1 \over 2}\vec  \zeta^{ \ 2} & \displaystyle\vec
\zeta^{ (t)} & \displaystyle{ 1 \over 2}\vec  \zeta^{ \ 2}
\cr\displaystyle\vec  \zeta & \displaystyle \Bbb I_{(d-2)} &
\displaystyle\vec  \zeta \cr\displaystyle -{1 \over 2}\vec  \zeta^{ \ 2} &
\displaystyle -\vec  \zeta^{ (t)} & \displaystyle 1-{1 \over 2}\vec  \zeta^{ \
2} \cr} \right]\ ;\ \vec  \zeta  \in  \Bbb C^{d-2} \right\} &  &  &  &  (3.3)
\cr} $$
(where $ \vec  \zeta^{ (t)} $ denotes the \lq\lq transposed\rq\rq\ of the
vector $ \vec  \zeta ). $ \par
One then easily checks that for every couple of horocyclic coordinates
$ (\vec  \zeta ,\tau ), $ one can rewrite Eqs.(3.1) (in view of
Eqs.(3.2)(3.3)) as follows:
$$ z \left(\vec  \zeta ,\tau \right) = n \left(\vec  \zeta \right) z_\tau  = n
\left(\vec  \zeta \right)\cdot a(\tau)  z_0 \eqno (3.4) $$
\par
Since $ H^{(c)} $ is the stabilizer of the base-point $ z_0, $ we see that
formula
(3.4) corresponds to the following Iwasawa-type decomposition
$\lbrack$H$\rbrack$ of $ G^{(c)}: $
$$ G^{(c)} = N^{(c)}\cdot A^{(c)}\cdot H^{(c)} \eqno (3.5) $$
\par
\smallskip
{\sl The invariant volume form on \/}$ X^{_{(c)}}_{d-1} $ {\sl in horocyclic
coordinates\/}: \par
\smallskip
The invariant volume form $ 2(-1)^{d-2} \left.{ {\rm d} z^{(0)}\wedge
..\wedge {\rm d} z^{(d-1)} \over {\rm d} s(z)} \right\vert_{ X^{(c)}_{d-1}},
$
admits the following representative $ \left(z^{(0)}+z^{(d-1)} \right)^{-1}
{\rm d} \left(z^{(0)}+z^{(d-1)} \right)\wedge {\rm d} z^{(1)}... {\rm d}
z^{(d-2)}; $
now, by computing $ \tau  = \tau( z) $ from Eqs.(3.1a), (3.1c), we obtain:
$$ {\rm d} \tau( z) = {i \over z^{(0)}+z^{(d-1)}} {\rm d}
\left(z^{(0)}+z^{(d-1)} \right)\ ; \eqno (3.6) $$
in view of Eqs.(3.1) and (3.6), we can thus write:
$$ \left.2(-1)^{d-2}{ {\rm d} z^{(0)}\wedge ...\wedge {\rm d} z^{(d-1)} \over
{\rm d} s(z)} \right\vert_{ X^{(c)}_{d-1}} = -i\ {\rm e}^{-i(d-2)\tau} {\rm d}
\tau \wedge {\rm d}\vec  \zeta \eqno (3.7) $$
where $ {\rm d}\vec  \zeta $ is the volume form defined by the Lebesgue
measure on $ \Bbb R^{d-2}. $ \par
\smallskip
{\sl Real horospheres\/}: for $ \tau  = iw $ (resp. $ \tau  = iw+\pi ), $ with
$ w \in  \Bbb R, $ the
point $ z_\tau $ varies on the branch $ z^{(d-1)}>0 $ (resp. $ z^{(d-1)}<0) $
of the real
hyperbola $ \hat  X; $ for each of these values of $ \tau , $ the
corresponding horosphere $ \Pi_{ iw} $
(resp. $ \Pi_{ iw+\pi} ) $ contains a real paraboloid (or {\sl real
horosphere\/}) $ \Pi^ +_w $ (resp. $ \Pi^ -_w) $
embedded in $ X_{d-1} $ and represented by the parametrization (deduced from
Eqs.(3.1)):
$$ \matrix{\displaystyle z^{(0)} \hfill& \displaystyle = \pm  ( {\rm sinh} \ w
+ {1 \over 2} \vec  \zeta^{ \ 2} {\rm e}^w)\ , \hfill& \displaystyle\vec  z =
\pm  \vec  \zeta \ {\rm e}^w\ , \hfill \cr\displaystyle z^{(d-1)} \hfill&
\displaystyle = \pm  ( {\rm cosh} \ w - {1 \over 2} \vec  \zeta^{ \ 2} {\rm
e}^w)\ , \hfill& \displaystyle {\rm with} \ \vec  \zeta  \in  \Bbb R^{d-2}\ .
\hfill \cr} \eqno (3.8) $$
By putting $ \vec  \zeta =0 $ in these equations, we obtain the point $ z_{iw}
$ (resp. $ z_{iw+\pi} ) $ which
is the apex of the corresponding paraboloid $ \Pi^ +_w $ (resp. $ \Pi^ -_w) $
(see Fig.5a). \par
{\sl Real horocycles \/}$ h^+_w: $ In each real horosphere $ \Pi^ +_\omega $
of $ X_{d-1}, $
with $ w \geq  0, $ we introduce the following $ (d-2) $-dimensional subset:
$$ \underline{h}^+_w=\Pi^ +_w\cap X_+= \left\{ z\in X_{d-1};\ z\geq z_0;\ z=n
\left(\vec  \zeta \right)a(iw)z_0,\ \vec  \zeta \in \Bbb R^{d-2} \right\} \ ,
\eqno (3.9) $$
which can be parametrized (in view of Eq.(3.1c)) by the ball
$$ B_{iw}= \left\{\vec  \zeta \in \Bbb R^{d-2};\ \vec  \zeta =\rho\vec
\alpha ,\ \vec  \alpha \in \Bbb S_{d-3},\ 0\leq \rho \leq \rho_ 0(w),\rho_
0(w)= \left[2( {\rm cosh} \ w-1) {\rm e}^{-w} \right]^{1/2} \right\} \ . $$
This set, equipped with the natural orientation of the ball $ B_{iw}, $
defines a
{\sl horocycle\/} $ h^+_w; $ the latter is a representative of a relative
homology
class in $ H_{d-2} \left(\Pi^ +_w,\partial X_+\cap \Pi^ +_w \right), $ $
\partial X_+ $ denoting the boundary of $ X_+, $ namely the
cone:
$$ \left\{ z\in X_{d-1};\ \left(z-z_0 \right)^2=0,\ z^{(0)}\geq 0 \right\} =
\left\{ z\in X_{d-1};\ z^{(0)}\geq 0,\ z^{(d-1)}=1 \right\} \ . $$
As illustrated on Fig.5b, one checks that $ X_+ = \bigcup^{ }_{ w\geq 0}
\underline{h}^+_w. $ \par
\medskip
{\bf 3.3 Transformation }$ \goth R_ {\bf d} ${\bf\ and convolution for
invariant
Volterra kernels on }$ {\bf X}_{ {\bf d-1}} $ \par
\medskip
{\sl The transformation \/}$ \goth R_d: $ \par
We shall now define an integral transformation $ F \charlvmidup{ \goth R_d}{
\rightarrowfill}  \hat  F $ of Radon-type
in $ X_{d-1}, $ in which the integration is taken on the horocycles $ h^+_w $
(the latter
playing the same role as the hyperplanes in the Radon transformation). The
transformation $ \goth R_d $ will apply to all invariant Volterra kernels $ K
\left(z,z_0 \right) \equiv  F(z) $
$ \left(K\in V \left(X_{d-1} \right)^{\natural} \right) $ and the
corresponding transforms $ \hat  F \left(z_{iw} \right) \equiv  \hat {\rm
f}(w) $ will be
functions defined on the semi-branch of hyperbola $ \hat  X_+ = \left\{
z_\tau ; \right.\allowbreak \tau =iw, $ $ w\geq 0\} , $
representing invariant Volterra kernels $ \hat  K $ on $ \hat  X $ $ (\hat  K
\in  V \left(\hat  X \right)^{\natural} ). $ \par
We put:
$$ \hat  F \left(z_{iw} \right) =2(-1)^{d-2} {\rm e}^{-(d-2)w} \int^{ }_{
h^+_w}F(z) \left.{ {\rm d} z^{(0)}\wedge ...\wedge {\rm d} z^{(d-1)} \over
{\rm d} s(z)\wedge {\rm d} w} \right\vert_{ \Pi^ +_w} \eqno (3.10) $$
or (in view of Eqs.(3.4) and (3.7)):
$$ \hat {\rm f}(w) = \hat  F \left(z_{iw} \right) = \int^{ }_{ B_{iw}}F
\left(n \left(\vec  \zeta \right)a(iw)z_0 \right) {\rm d}\vec  \zeta \eqno
(3.11) $$
\par
\smallskip
{\sl The algebras\/} $ V \left(X_{d-1} \right)^{\natural} $ {\sl and\/} $ V
\left(\hat  X \right)^{\natural} : $ \par
The composition product $ \diamond $ of the algebra $ V \left(X_{d-1} \right)
$ was defined in
$\lbrack$B,V-1$\rbrack$ as follows (see formulae\footnote {$ ^\ast $}{formula
(25) of $\lbrack$B,V.1$\rbrack$ should be corrected
by putting $ i $ at the nominator instead of the denominator} (25), (32) of
$\lbrack$B,V-1$\rbrack$):
$$ K \left(z,z^{\prime} \right) = \left(K_1\diamond K_2 \right)
\left(z,z^{\prime} \right) \equiv  \int^{ }_{ \diamond \left(z,z^{\prime}
\right)}K_1 \left(z,z^{\prime\prime} \right)K_2 \left(z^{\prime\prime}
,z^{\prime} \right) {\rm d} \sigma_{ d-1} \left(z^{\prime\prime} \right)\ ,
\eqno (3.12) $$
where:
$$ {\rm d} \sigma_{ d-1}(z) = 2 { {\rm d} z^{(0)}\wedge {\rm d} z^{(1)}\wedge
\cdot \cdot \cdot {\rm d} z^{(d-1)} \over {\rm d} s(z)} \left\vert_{ X_{d-1}}
\right. \eqno (3.12^\prime ) $$
and \par
\noindent$ \diamond \left(z,z^{\prime} \right) = \left\{ z^{\prime\prime} \in
X_{d-1}\ ;\ z \geq  z^{\prime\prime}  \geq  z^{\prime} \right\} . $ \par
\smallskip
In the subalgebra $ V \left(X_{d-1} \right)^{\natural} $ of the invariant
Volterra kernels, or
$ H $-invariant functions $ F(z) = K \left(z,z_0 \right) $ $ (\forall g\in G,
$ $ K \left(z,z^{\prime} \right) = K \left(gz,gz^{\prime} \right), $ formula
(3.12) can be seen to reduce to the following convolution formula on $ G/H $
(thanks to the transitivity property of $ X_{d-1} $ for the group $ G \equiv
{\rm SO}(1,d-1)): $
$$ F \left(gz_0 \right) = \int^{ }_{ g^{\prime} z_0\in \diamond \left(gz_0,z_0
\right)}F_1 \left(g^{\prime -1}g\ z_0 \right)F_2 \left(g^{\prime} z_0 \right)
{\rm d} \sigma_{ d-1} \left(g^{\prime} z_0 \right) \eqno (3.13) $$
(Note that in this presentation the functions $ F, $ $ F_i $ $ (i=1,2) $
appear as
$ H $-bi-invariant functions on the group $ G, $ since: $ F_i \left(h_1gh_2z_0
\right) = F_i \left(h_1g\ z_0 \right) = F_i \left(g\ z_0 \right)). $ \par
Similarly, the composition product $ \hat  \diamond $ of the algebra $ V
\left(\hat  X \right) $ is defined
as follows:
$$ \left(\hat  K_1\hat  \diamond\hat  K_2 \right) \left(z,z^{\prime} \right) =
\int^{ }_{ \diamond \left(z,z^{\prime} \right)\cap\hat  X}\hat  K_1
\left(z,z^{\prime\prime} \right)\hat  K_2 \left(z^{\prime\prime} ,z^{\prime}
\right) {\rm d} \sigma_ 1 \left(z^{\prime\prime} \right) \eqno (3.14) $$
with
$$ {\rm d} \sigma_ 1(z) = {2\ {\rm d} z^{(0)}\wedge {\rm d} z^{(d-1)} \over
{\rm d}\hat  s(z)} \left\vert_{\hat  X} \right. $$
(the points $ z, $ $ z^{\prime} $ in $ \hat  X $ being considered as points of
$ X_{d-1}, $ with the
embedding $ \hat  X\subset X_{d-1}). $ \par
In the subalgebra $ V \left(\hat  X \right)^{\natural} , $ formula (3.14)
reduces to the following
convolution formula on $ {\rm SO}(1,1)/J: $
$$ \left(\hat  K_1\hat  \diamond\hat  K_2 \right) \left(z_{iw},z_0 \right)
\equiv  \int^{ }_{ w>w^{\prime} >0}\hat  F_1 \left(z_{i \left(w-w^{\prime}
\right)} \right) \hat  F_2 \left(z_{iw^{\prime}} \right) {\rm d} \sigma_ 1
\left(z_{iw^{\prime}} \right) \eqno (3.15) $$
or equivalently (as a convolution on $ \Bbb R^+): $
$$ \left(\hat  K_1\hat  \diamond\hat  K_2 \right) \left(z_{iw},z_0 \right) =
\left(\hat {\rm f}_1\diamond\hat {\rm f}_2 \right)(w) \equiv  \int^ w_0\hat
{\rm f}_1 \left(w-w^{\prime} \right)\hat {\rm f}_2 \left(w^{\prime} \right)
{\rm d} w^{\prime} \eqno (3.15^\prime ) $$
(since : $ {\rm d} \sigma_ 1 \left(z_{iw} \right) = { {\rm d}
\left(z^{(0)}+z^{(d-1)} \right) \over z^{(0)}+z^{(d-1)}} \left\vert_{\hat  X}
\right. = {\rm d} w). $ \par
\smallskip
We will now show that the convolution product $ \diamond $ in $
V(X)^{\natural} $ is
transmuted by the transformation $ \goth R_d $ into the corresponding
convolution
product $ \hat  \diamond $ in $ V \left(\hat  X \right)^{\natural} . $ \par
\medskip
{\bf Proposition 8}: \par
For every couple of invariant kernels $ K_i \left(z,z_0 \right)\equiv F_i(z) $
in $ V \left(X_{d-1} \right)^{\natural} , $ $ i=1,2 $ with respective images $
\hat  K_i \left(z_{iw},z_0 \right)\equiv\hat  F_i \left(z_{iw} \right) $ $
\left(\hat  F_i=\goth R_d \left(F_i \right) \right) $ in
$ V \left(\hat  X \right)^{\natural} , $ the transform $ \hat  F = \goth
R_d(F) $ of the convolution product $ F(z) = \left(K_1\diamond K_2 \right)
\left(z,z_0 \right) $
is the corresponding convolution product $ \hat  F \left(z_{iw} \right) =
\left(\hat  K_1\hat  \diamond\hat  K_2 \right) \left(z_{iw},z_0 \right)\ .
$ \par
\smallskip
{\bf Proof}: In view of Eqs.(3.10), (3.11), (3.12), (3.12') and (3.7), we
can write:
$$ \goth R_d(F) \left(z_{iw} \right)= \int^{ }_{ h^+_w} {\rm d}\vec  \zeta( z)
\int^{ }_{ \diamond \left(z,z_0 \right)}K_1 \left(z,z^{\prime} \right)K_2
\left(z^{\prime} ,z_0 \right)\ {\rm e}^{(d-2)w^{\prime}} \left( {\rm d}\vec
\zeta^{ \prime} \wedge {\rm d} w^{\prime} \right) \left(z^{\prime} \right)
\eqno (3.16) $$
in which the integration variables $ z, $ $ z^{\prime} $ $ (z\in
\underline{h}^+_w, $ $ z^{\prime} \in \diamond \left(z,z_0 \right)) $ have
been
parametrized (in view of Eq.(3.4)) by:
$$ z=g\ z_0 = n \left(\vec  \zeta \right)a(iw)z_0\ ,\ \ z^{\prime}
=g^{\prime} z_0 = n \left(\vec  \zeta^{ \prime} \right)a \left(iw^{\prime}
\right)z_0\ . $$
\par
By using the invariance property of the kernels $ K_1, $ $ K_2 $ as in
Eq.(3.13), we can rewrite Eq.(3.16) as follows:
\smallskip
$$ \matrix{\displaystyle \goth R_d(F) \left(z_{iw} \right) & \displaystyle =
\int^{ }_{ h^+_w} {\rm d}\vec  \zeta( z) \int^{ }_{ \diamond \left(z,z_0
\right)}F_1 \left(a \left(iw^{\prime} \right)^{-1}n \left(\vec  \zeta^{
\prime} \right)^{-1}n \left(\vec  \zeta \right)a(iw)z_0 \right) \hfill
\cr\displaystyle  & \displaystyle \times F_2 \left(n \left(\vec  \zeta^{
\prime} \right)a \left(iw^{\prime} \right)z_0 \right) {\rm
e}^{(d-2)w^{\prime}} \left( {\rm d}\vec  \zeta^{ \prime} \wedge {\rm d}
w^{\prime} \right) \left(z^{\prime} \right) \hfill \cr} \eqno (3.17) $$
\par
We shall now take into account the following relations for the elements $ n
\left(\vec  \zeta \right), $ $ n \left(\vec  \zeta^{ \prime} \right) $
and $ a(iw), $ $ a \left(iw^{\prime} \right) $ of the subgroups $ N^{(c)} $
and $ A^{(c)}: $
$$ \matrix{\displaystyle  & \displaystyle n \left(\vec  \zeta^{ \prime}
\right)^{-1}n \left(\vec  \zeta \right) \hfill& \displaystyle = n \left(\vec
\zeta -\vec  \zeta^{ \prime} \right) \hfill \cr\displaystyle  & \displaystyle
a \left(iw^{\prime} \right)^{-1}a(iw) \hfill& \displaystyle = a \left(i
\left(w-w^{\prime} \right) \right) \hfill \cr\displaystyle {\rm and} &
\displaystyle a(iw)n \left(\vec  \zeta \right) \hfill& \displaystyle = n
\left(\vec  \zeta \ {\rm e}^{-w} \right)a(iw) \hfill \cr} $$
In view of these relations, we can then write:
$$ F_1 \left(a \left(iw^{\prime} \right)^{-1}n \left(\vec  \zeta^{ \prime}
\right)^{-1}n \left(\vec  \zeta \right)a(iw)z_0 \right) = F_1 \left(n
\left(\vec  \zeta^{ \prime\prime} \right)a \left(iw^{\prime\prime} \right)z_0
\right)\ , $$
where we have put:
$$ \vec  \zeta^{ \prime\prime}  = \left(\vec  \zeta -\vec  \zeta^{ \prime}
\right) {\rm e}^{w^{\prime}} \ \ \ \ {\rm and} \ \ \ \ \ w^{\prime\prime}
=w-w^{\prime} \eqno (3.18) $$
Eq.(3.17) then yields:
$$ \goth R_d(F) \left(z_{iw} \right) = \int^{ }_{ \Sigma_ w}F_1 \left(n
\left(\vec  \zeta^{ \prime\prime} \right)a \left(iw^{\prime\prime} \right)z_0
\right)F_2 \left(n \left(\vec  \zeta^{ \prime} \right)a \left(iw^{\prime}
\right)z_0 \right) {\rm d}\vec  \zeta^{ \prime\prime} \wedge {\rm d}\vec
\zeta^{ \prime} \wedge {\rm d} w^{\prime} \eqno (3.19) $$
where (in view of Eqs.(3.17), (3.18)) the integration cycle $ \Sigma_ w $ can
be specified as follows:
$$ \matrix{\displaystyle {\rm supp} .\ \Sigma_ w = \left\{ \left(\vec
\zeta^{ \prime\prime} ,\vec  \zeta^{ \prime} ,w^{\prime} \right)\ ;\ z=n
\left(\vec  \zeta^{ \prime} +\vec  \zeta^{ \prime\prime} {\rm
e}^{-w^{\prime}} \right)a(iw)z_0 \in  \underline{h}^+_w, \right.
\cr\displaystyle \left.z^{\prime}  = n \left(\vec  \zeta^{ \prime} \right)a
\left(iw^{\prime} \right)z_0\ ,\ \ \ z \geq  z^{\prime^{ \ }} \geq  z_0
\right\} \cr} $$
\par
Let us now show the auxiliary \par
\smallskip
{\bf Lemma\nobreak\ 2\nobreak\ }: $ {\rm Supp.} \ \Sigma_ w $ can be
alternatively described as the following subset:
$$ \matrix{\displaystyle \underline{\Sigma}_ w= \left\{ \left(\vec  \zeta^{
\prime\prime} ,\vec  \zeta^{ \prime} ,w^{\prime} \right)\ ;\ z^{\prime}  = n
\left(\vec  \zeta^{ \prime} \right) a \left(iw^{\prime} \right) z_0 \in
\underline{h}^+_{w^{\prime}} \ , \right. \hfill \cr\displaystyle \left.\ \ \ \
\ \ \ z^{\prime\prime}  = n \left(\vec  \zeta^{ \prime\prime} \right)a \left(i
\left(w-w^{\prime} \right) \right)z_0 \in  \underline{h}^+_{
\left(w-w^{\prime} \right)}\ ,\ 0 \leq  w^{\prime}  \leq  w \right\} \hfill
\cr} $$
\par
\smallskip
{\bf Proof }: \par
\smallskip
a) Let $ \left(\vec  \zeta^{ \prime\prime} ,\vec  \zeta^{ \prime} ,w^{\prime}
\right)\in {\rm supp} \ \Sigma_ w; $ since $ z^{\prime}  = n \left(\vec
\zeta^{ \prime} \right)a \left(iw^{\prime} \right)z_0 \geq  z_0 $ (i.e. $
z^{\prime}  \in  X_+), $ we have $ w^{\prime}  \geq  0 $
and $ z^{\prime}  \in  \underline{h}^+_{w^{\prime}} . $ Moreover, since $ z=n
\left(\vec  \zeta^{ \prime} +\vec  \zeta^{ \prime\prime} \ {\rm
e}^{-w^{\prime}} \right)a(iw)z_0 \geq  z^{\prime} , $ we have (in view of the
invariance of the cone $ \left\{ z \in  \Bbb R^d; \right. $ $ z \geq  0\} $
under the group $ {\rm SO}(1,d-1)): $
$$ \left[n \left(\vec  \zeta^{ \prime} \right)a \left(iw^{\prime} \right)
\right]^{-1} z \geq  \left[n \left(\vec  \zeta^{ \prime} \right)a
\left(iw^{\prime} \right) \right]^{-1}z^{\prime} \ , $$
i.e. $ n \left(\vec  \zeta^{ \prime\prime} \right) a \left(i
\left(w-w^{\prime} \right) \right)z_0 \geq  z_0, $ which implies that: $
\left(w-w^{\prime} \right) \geq  0 $ and $ z^{\prime\prime}  = n \left(\vec
\zeta^{ \prime\prime} \right)a(i(w-\allowbreak \left. \left.w^{\prime} \right)
\right)z_0\in \underline{h}^+_{ \left(w-w^{\prime} \right)}; $
and therefore that $ \left(\vec  \zeta^{ \prime\prime} ,\vec  \zeta^{ \prime}
,w^{\prime} \right) \in  \underline{\Sigma}_ w. $ \par
b) Conversely, let $ \left(\vec  \zeta^{ \prime\prime} ,\vec  \zeta^{ \prime}
,w^{\prime} \right) \in  \underline{\Sigma}_ w; $ then the conditions $
w^{\prime}  \geq  0 $ and $ z^{\prime}  \in  \underline{h}^+_{w^{\prime}} $
imply
that $ z^{\prime}  \geq  z_0; $ similarly the conditions $ w-w^{\prime}  \geq
0 $ and $ z^{\prime\prime}  \in  \underline{h}^+_{ \left(w-w^{\prime} \right)}
$ imply that $ z^{\prime\prime}  \geq  z_0, $ which
entails:
$$ n \left(\vec  \zeta^{ \prime} \right)a \left(iw^{\prime}
\right)z^{\prime\prime}  \geq  n \left(\vec  \zeta^{ \prime} \right)a
\left(iw^{\prime} \right)z_0,\ \ {\rm i.e.} $$
$$ z = n \left(\vec  \zeta^{ \prime} +\vec  \zeta^{ \prime\prime} {\rm
e}^{-w^{\prime}} \right) a(iw) z_0 \geq  z^{\prime}  \geq  z_0. $$
But since $ z \geq  z_0, $ this expression of $ z $ shows that $ z \in
\underline{h}^+_w, $ and therefore that $ \left(\vec  \zeta^{ \prime\prime}
,\vec  \zeta^{ \prime} ,w^{\prime} \right)\in {\rm supp} \ \Sigma_ w. $ \par
\smallskip
In view of this auxiliary lemma, Eq.(3.19) can thus be rewritten as follows:
$$ \goth R_d(F) \left(z_{iw} \right)= \int^ w_0 {\rm d} w^{\prime} \int^{ }_{
B_{i \left(w-w^{\prime} \right)}}F_1 \left(n \left(\zeta^{ \prime\prime}
\right)a \left(i \left(w-w^{\prime} \right) \right)z_0 \right) {\rm d}\vec
\zeta^{ \prime\prime} \int^{ }_{ B_{iw^{\prime}}} F_2 \left(n \left(\vec
\zeta^{ \prime} \right)a \left(iw^{\prime} \right)z_0 \right) {\rm d}\vec
\zeta^{ \prime} $$
and therefore, in view of Eqs.(3.11) and (3.15):
$$ \goth R_d({\cal F}) \left(z_{iw} \right) = \int^ w_0 {\rm d} w^{\prime} \
\hat  F_1 \left(z_{i \left(w-w^{\prime} \right)} \right)\hat  F_2
\left(z_{iw^{\prime}} \right) = \left(\hat  K_1 \hat  \diamond\hat  K_2
\right) \left(z_{iw},z_0 \right) \eqno  $$
which ends the proof of proposition 8. \par
\medskip
{\sl An alternative form for\/} $ \goth R_d; $ {\sl the transformation\/} $
{\cal R}_d: $ \par
By starting from Eqs.(3.11), one can obtain the following Abel-type integral
representation for $ \hat {\rm f}(w) $ in terms of $ {\rm f}(v) =
\underline{{\rm f}} \left(z^{(d-1)} \right)_{ \left\vert z^{(d-1)}= {\rm cosh}
\ v \right.}\equiv F(z); $ this
alternative form of the transformation $ \goth R_d $ (equivalent to Eq.(3.10))
will be called $ {\cal R}_d, $
namely:
$$ \hat {\rm f}(w) = \left[{\cal R}_d(F) \right](w) = \omega_{ d-2} {\rm e}^{-
\left({d-2 \over 2} \right)w} \int^ w_0 {\rm f}(v)[2( {\rm cosh} \ w- {\rm
cosh} \ v)]^{{d-4 \over 2}} {\rm sinh} \ v\ {\rm d} v\ , \eqno (3.20) $$
where the constant $ \omega_{ d-2} $ denotes the area of the sphere $ \Bbb
S^{d-3}: $ $ \omega_{ d-2}=2\cdot{ \pi^{{ d-2 \over 2}} \over \Gamma
\left({d-2 \over 2} \right)} $ $ (d\geq 3). $ \par
\smallskip
The derivation of the latter from Eq.(3.11) results directly from Eq.(3.1c)
which
yields (for $ \tau  = iw) $ the equalities: $ z^{(d-1)} = {\rm cosh} \ v =
{\rm cosh\ } w - {1 \over 2} \vec  \zeta^{ \ 2} {\rm e}^w\ , $
and therefore leads one to the parameter substitution $ (\vec  \zeta^{ \ 2}
\longrightarrow v) $ in Eq.(3.11).
In \S 3.4, we shall present in more detail a complex extension (on $
X^{(c)}_{d-1}) $ of formulae
(3.10) and (3.20) and of this equivalence (namely the transformations $ \goth
R^{(c)}_d $ and $ {\cal R}^{(c)}_d). $ \par
\medskip
{\bf Invariant Volterra kernels of moderate growth on} $ X_{d-1} $ \par
\smallskip
{\bf Definition 2:} For each couple $ (m,\mu) $ with $ \mu \geq 0, $ we call $
V^{(m)}_\mu \left(X_{d-1} \right) $ the subspace of
all kernels $ K $ in $ V \left(X_{d-1} \right)^{\natural} , $ whose
representative $ {\rm f}(v) $ (such that $ {\rm f}(v) = \underline{{\rm f}}(
{\rm cosh} \ v) $ and $ \underline{{\rm f}} \left(z^{(d-1)} \right)=F(z)=K
\left(z,z_0 \right)) $
satisfies the following conditions: \par
\smallskip
i) $ {\rm supp} \ {\rm f}  \subset  [\mu ,+\infty[ $ (i.e. $ {\rm supp} \ F
\subset  X^\mu_ + $ or $ {\rm supp} \ K\subset \Sigma^ \mu_ +) $ \par
ii) the function $ {\rm g}_m[ {\rm f}](v) = {\rm e}^{-mv}\vert {\rm
f}(v)\vert $ belongs to $ L^1 \left(\Bbb R^+ \right). $ \par
\smallskip
We also put:
$$ V^{(m)} \left(X_{d-1} \right)^{\natural}  = V^{(m)}_0 \left(X_{d-1} \right)
= \bigcup^{ }_{ \mu \geq 0}V^{(m)}_\mu \left(X_{d-1} \right) $$
\par
{\bf Remark :} \par
This definition also applies to the case $ d=2, $ for which the following
statement is
easily obtained: \par
{\it For each real number} $ m, $ $ V^{(m)} \left(X_1 \right)^{\natural} $
{\it is a subalgebra of} $ V \left(X_1 \right)^{\natural} $ {\it (and
if} $ K_i \in  V^{(m)}_{\mu_ i} \left(X_1 \right), $
$ i=1,2, $ {\it then} $ K = K_1\diamond K_2 \in  V^{(m)}_{\mu_ 1+\mu_ 2}
\left(X_1 \right)). $ \par
\smallskip
This property (strongly related to proposition 3, since, for the jump function
$ {\rm f}  = \Delta f, $
it appears as a by-product of this proposition) is immediately derived in
terms of the
representatives $ f, $ $ f_i $ of the kernels $ K, $ $ K_i $ $ (i=1,2). $ In
fact, since $ {\rm f}(v)= \int^ v_0 {\rm f}_1 \left(v-v^{\prime} \right) {\rm
f}_2 \left(v^{\prime} \right) {\rm d} v^{\prime} $
(see Eq.(3.15')), one has:
$$ \matrix{\displaystyle \left\Vert {\rm g}_m[ {\rm f}] \right\Vert_ 1 =
\int^{ \infty}_ 0 {\rm e}^{-mv}\vert {\rm f}(v)\vert {\rm d} v \leq  \int^{
\infty}_ 0 {\rm d} v^{\prime} \int^{ \infty}_ 0 {\rm d} v^{\prime\prime} {\rm
e}^{-m \left(v^{\prime} +v^{\prime\prime} \right)} \left\vert {\rm f}_2
\left(v^{\prime} \right) \right\vert \left\vert {\rm f}_1
\left(v^{\prime\prime} \right) \right\vert \cr\displaystyle = \left\Vert {\rm
g}_m \left[ {\rm f}_1 \right] \right\Vert_ 1 \left\Vert {\rm g}_m \left[ {\rm
f}_2 \right] \right\Vert_ 1 \cr} $$
\par
{\it We will now extend this statement to the case of kernels in} $ V^{(m)}
\left(X_{d-1} \right), $ {\it
for all} $ d\geq 3, $
{\it and also show that the rate of moderate growth} $ m $ {\it is preserved
by the
transformation} $ {\cal R}_d, $ {\it provided} $ m>-1. $ \par
\medskip
{\bf Proposition 9\nobreak\ :} \par
For each dimension $ d, $ $ d\geq 3, $ and for each real number $ m $ such
that $ m>-1, $ the
following properties hold: \par
\smallskip
i) $ V^{(m)} \left(X_{d-1} \right)^{\natural} $ is a subalgebra of $ V
\left(X_{d-1} \right)^{\natural} . $ More precisely, if $ K_i\in
V^{(m)}_{\mu_ i} \left(X_{d-1} \right), $
$ i=1,2, $ then $ K_1\diamond K_2\in V^{(m)}_{\mu_ 1+\mu_ 2} \left(X_{d-1}
\right). $ \par
\smallskip
ii) The transformation $ \goth R_d $ defines a morphism of
algebras between $ V^{(m)} \left(X_{d-1} \right)^{\natural} $ and $ V^{(m)}
\left(X_1 \right)^{\natural} $ (with $ X_1 $ identified to the meridian
section $ \hat  X $
of $ X_{d-1}); $ moreover, for each kernel $ K $ in $ V^{(m)}_\mu
\left(X_{d-1} \right) $ the image $ \hat {\rm f} = {\cal R}_d(F) $ of $ F(z)
\equiv  K \left(z,z_0 \right) $
satisfies the following properties: \par
\smallskip
a) $ {\rm supp} \ \hat {\rm f} \subset  [\mu ,+\infty[ $ \par
b) $ \dlowlim{ {\rm max}}{0\leq r\leq p(d)} \left\vert\hat {\rm f}^{(r)}(v)
\right\vert  \leq  {\rm Cst\ e}^{mv}\hat {\rm g}(v), $ with $ \hat {\rm g}\in
L^1 \left(\Bbb R^+ \right); $ \par
\noindent$ p(d) = {d-2 \over 2} $ if $ d $ is even and $ p(d) = {d-3 \over 2}
$ if $ d $ is odd. \par
\medskip
{\bf Proof :} \par
\smallskip
i) Being given $ K_1, $ $ K_2 $ in $ V^{(m)} \left(X_{d-1} \right), $ we will
estimate the norm $ \left\Vert {\rm g}_m[ {\rm f}] \right\Vert_ 1 $ associated
with the kernel $ K = K_1\diamond K_2 $ $ ( {\rm f}(v) = {\rm f}( {\rm cosh} \
v), $ $ \underline{{\rm f}} \left(z^{(d-1)} \right) = F(z) = K \left(z,z_0
\right)) $ by
rewriting the $ \diamond $-product formula (3.12) as follows (for $ K
\left(z,z_0 \right) = {\rm f}(v), $ in terms of the
representatives $ {\rm f}_1, $ $ {\rm f}_2 $ of $ K_1, $ $ K_2): $
$$ {\rm f}(v) = {\omega_{ d-2} \over( {\rm sinh} \ v)^{d-3}} \int^{ }_{
\dlowlim{ v^{\prime} \geq 0,v^{\prime\prime} \geq 0}{v\geq v^{\prime}
+v^{\prime\prime}}} {\rm sinh} \ v^{\prime} {\rm d} v^{\prime} {\rm sinh} \
v^{\prime\prime} {\rm d} v^{\prime\prime} \left[\Lambda \left(v,v^{\prime}
,v^{\prime\prime} \right) \right]^{{d-4 \over 2}} {\rm f}_1
\left(v^{\prime\prime} \right) {\rm f}_2 \left(v^{\prime} \right) \eqno
(3.20^\prime ) $$
where:
$$ \left. \matrix{\displaystyle \Lambda \left(v,v^{\prime} ,v^{\prime\prime}
\right) & \displaystyle = {\rm cosh}^2v+ {\rm cosh}^2v^{\prime} + {\rm
cosh}^2v^{\prime\prime} -2\ {\rm cosh} \ v\ {\rm cosh} \ v^{\prime} {\rm cosh}
\ v^{\prime\prime} -1 \hfill \cr\displaystyle  & \displaystyle = \left[ {\rm
cosh} \ v- {\rm cosh} \left(v^{\prime} +v^{\prime\prime} \right) \right]\cdot
\left[ {\rm cosh} \ v- {\rm cosh} \left(v^{\prime} -v^{\prime\prime} \right)
\right] \hfill \cr} \right\} \eqno (3.20^{\prime\prime} ) $$
\par
These formulae can be derived from Eq.(59) of $\lbrack$B.V-1$\rbrack$ which
gave an alternative
expression of the $ \diamond $-product as an integral over the argument $
z^{\prime} $ of $ F_2 $ (such that $ F_2 \left(z^{\prime} \right)=K_2
\left(z^{\prime} ,z_0 \right)), $
written in terms of the Cartan-type parametric representation $ (\Pi) : $ $
z^{\prime} =z \left(iv^{\prime} ,\varphi^{ \prime} ,\alpha^{ \prime} \right) $
(see
Eqs.(19) of $\lbrack$B.V-1$\rbrack$), the external variable $ z $ being chosen
in the special situation:
$$ z=z(iv,0,0) = ( {\rm sinh} \ v,0, {\rm cosh} \ v). $$
\par
In fact, the kernels $ F_2 $ and $ F_1 $ only depend respectively on the
variables $ {\rm cosh} \ v^{\prime}  = { \left(z^{\prime} -z_0 \right)^2
\over 2} +1 $
and $ {\rm cosh} \ v^{\prime\prime}  = { \left(z-z^{\prime} \right)^2 \over 2}
+1 = {\rm cosh} \ v\ {\rm cosh} \ v^{\prime}  - {\rm sinh} \ v\ {\rm sinh} \
v^{\prime} \ {\rm cosh} \ \varphi^{ \prime} . $ \par
By expressing the measure $ {\rm d} \sigma_{ d-1} \left(z^{\prime} \right) =
\omega \left(\vec  \alpha^{ \prime} \right) \left( {\rm sinh} \ v^{\prime}
\right)^{d-2} \left( {\rm sinh} \ \varphi^{ \prime} \right)^{d-3} {\rm d}
v^{\prime} {\rm d} \varphi^{ \prime} $ in terms
of the variables $ \left(\vec  \alpha^{ \prime} ,v^{\prime} ,v^{\prime\prime}
\right) $ (the integration on $ \vec  \alpha^{ \prime} $ being trivial), one
then obtains
formula (3.20'), the integration range $ \diamond \left(z,z_0 \right) $ (see
Eqs.(3.12)) being simply described by
the conditions $ v^{\prime} \geq 0, $ $ v^{\prime\prime} \geq 0, $ $
v^{\prime} +v^{\prime\prime} \leq v; $ we also note that the support
conditions $ v^{\prime\prime} \geq \mu_ 1, $ $ v^{\prime} \geq \mu_ 2 $
for $ K_1, $ $ K_2 $ readily imply the support condition $ v\geq \mu_ 1+\mu_ 2
$ for $ K $ (i.e. the property specified
in the second part of i)). \par
{}From Eq.(3.20'), we now deduce the following majorizations:
$$ \matrix{\displaystyle \left\Vert {\rm g}_m[ {\rm f}] \right\Vert_ 1 =
\int^{ \infty}_ 0 {\rm e}^{-mv}\vert {\rm f}(v)\vert {\rm d} v \leq  \omega_{
d-2} \times \hfill \cr\displaystyle \times \int^{ \infty}_ 0 {\rm e}^{-mw}
{\rm d} w \int^{ \infty}_ 0 {\rm sinh} v^{\prime} {\rm d} v^{\prime} \int^{
\infty}_ 0 {\rm sinh} v^{\prime\prime} {\rm d} v^{\prime\prime} {\rm g}_m
\left[ {\rm f}_2 \right] \left(v^{\prime} \right) {\rm g}_m \left[ {\rm f}_1
\right] \left(v^{\prime\prime} \right) \left.{ \left[\Lambda
\left(v,v^{\prime} ,v^{\prime\prime} \right) \right]^{{d-4 \over 2}} \over(
{\rm sinh} \ v)^{d-3}} \right\vert_{ v=w+v^{\prime} +v^{\prime\prime}} ;
\hfill \cr} $$
the r.h.s. of the latter can then be estimated with the help of the following
two
inequalities (valid for all $ v, $ $ v^{\prime} , $ $ v^{\prime\prime} $ such
that $ v^{\prime}  \geq  0, $ $ v^{\prime\prime}  \geq  0, $ $ v\geq
v^{\prime} +v^{\prime\prime} $ and easily deduced
from Eq.(3.20")): \par
\smallskip
a) $ \Lambda \left(v,v^{\prime} ,v^{\prime\prime} \right) \leq  ( {\rm cosh} \
v-1)^2, $ which yields:
$$ \forall d,\ \ \ d \geq  3,\ \ \ \ \ \ { \left[\Lambda \left(v,v^{\prime}
,v^{\prime\prime} \right) \right]^{{d-3 \over 2}} \over( {\rm sinh} \
v)^{d-3}} \leq  1 $$
\par
b) $ \Lambda \left(v,v^{\prime} ,v^{\prime\prime} \right) \geq 4 \left( {\rm
sinh} v^{\prime} {\rm sinh} v^{\prime\prime} \right)^2 {\rm e}^w \left( {\rm
e}^w-1 \right).
$ \par
By taking into account these inequalities, the previous majorization of $
\left\Vert {\rm g}_m[ {\rm f}] \right\Vert_ 1 $
can be replaced by:
$$ \left\Vert {\rm g}_m[ {\rm f}] \right\Vert_ 1 \leq  C_m\omega_{ d-2}
\left\Vert {\rm g}_m \left[ {\rm f}_1 \right] \right\Vert_ 1 \left\Vert {\rm
g}_m \left[ {\rm f}_2 \right] \right\Vert_ 1\ , $$
where the constant $ C_m $ (only defined for $ m>-1) $ is the value of the
integral $ \int^{ \infty}_ 0{ {\rm e}^{-(m+1)w} \over 2 \left(1- {\rm e}^{-w}
\right)^{1/2}} {\rm d} w. $ \par
\smallskip
Property i) is therefore established, with the following complementary result:
\par
\noindent {\it the} $ \diamond ${\it -product defines a continuous mapping
from} $ V^m \left(X_{d-1} \right)^{\natural} \times V^m \left(X_{d-1}
\right)^{\natural} $ {\it into} $ V^m \left(X_{d-1} \right)^{\natural} , $
{\it considered as a normed space equipped with the norm} $ \left\Vert {\rm
g}_m[ {\rm f}] \right\Vert_ 1. $ \par
\smallskip
ii) The proof of the properties of $ \hat {\rm f} $ (for $ K \in  V^{(m)}_\mu
\left(X_{d-1} \right)) $ relies directly on
Eq.(3.20). While property a) is immediate (since $ {\rm supp} \ {\rm f}
\subset  [\mu ,+\infty[ ), $ property b) requires estimates which (for
conciseness) we prefer to present later as a by-product of similar estimates
obtained in
the study of the complex transformation $ {\cal R}^{(c)}_d: $ this will be
done in the corollary of
proposition 13 (see \S 3.6). Properties a) and b) imply in particular that the
kernel $ \hat  K $
represented by $ \hat {\rm f} = {\cal R}_d(F) $ belongs to the algebra $
V^{(m)} \left(\hat  X \right)^{\natural} . $ Our first statement in ii) is
then a direct consequence of proposition 8 together with the result of i).
\par
\medskip
{\bf 3.4 Radon-type integrals on complex horocycles and
Abel-type transformation in star-shaped domains} \par
\smallskip
{\it In this subsection, our purpose is to define an integral
transformation} $ {\cal F } \charlvmidup{ \goth R^{(c)}_d}{ \rightarrowfill}
\hat {\cal F} $ {\it of Radon-type in} $ X^{(c)}_{d-1}, $ {\it in which
the integration is taken on appropriate cycles} $ h_\tau $ {\it of the
horospheres} $ \Pi_ \tau $ {\it introduced in \S 3.2; these cycles, called
complex
horocycles will play the same role in} $ X^{(c)}_{d-1}, $ {\it as the
horocycles} $ h^+_w $
{\it of} $ X_{d-1} $ {\it used above (see \S 3.3): in fact we shall
define} $ h_\tau $ {\it in such a
way that, for each} $ \tau  = iw, $ $ w \geq  0,, $ $ h_{iw} = h^+_w. $ \par
{\it The transformation} $ \goth R^{(c)}_d $ {\it will apply to
all} $ H^{(c)} ${\it -invariant
functions} $ {\cal F}, $ {\it analytic in a given domain} $ D $ {\it of} $
X^{(c)}_{d-1} $
{\it in the class} $ {\cal D} $
{\it defined in \S 3.1; the transforms} $ \hat {\cal F} $ {\it of these
functions will be
shown to be analytic in corresponding domains} $ \hat  D $ {\it of the
hyperbola} $ \hat  X^{(c)}. $ \par
\medskip
{\sl Complex horocycles:\/} \par
Let $ {\cal C} $ be the intersection of $ X^{(c)}_{d-1} $ with the hyperplane
with equation $ z^{(d-1)}=1. $ $ {\cal C} $ is a cone whose intersection with
each
horosphere $ \Pi_ \tau $ is a $ (d-3) $-dimensional complex sphere $ {\cal
C}_\tau $
represented in the parametrization $ z=z \left(\vec  \zeta ,\tau \right) $
(see Eqs.(3.1)) by:
$ \vec  \zeta^{ \ 2} = 2\ {\rm e}^{i\tau}( {\rm cos} \ \tau -1)\ . $ \par
\smallskip
We now introduce in each $ \Pi_ \tau $ the horocycle $ h_\tau $ as a special
representative of an element of $ H_{d-2} \left(\Pi_ \tau ,{\cal C}_\tau
\right) $ whose support is the
following set:
$$ \matrix{\displaystyle \underline{h}_\tau  = \left\{ z\in \Pi_ \tau ;\ z = z
\left(\vec  \zeta ,\tau \right),\ \vec  \zeta  = \rho\vec  \alpha ,\ \vec
\alpha  \in  \Bbb S^{d-3}\ , \right. \cr\displaystyle \left.\rho  = \left[2
{\rm e}^{i\tau}( 1-\lambda)( {\rm cos} \ \tau -1) \right]^{1/2^{\ }},\ 0 \leq
\lambda  \leq  1 \right\} \ ; \cr} \eqno (3.21) $$
$ \underline{h}_\tau $ is homeomorphic to a $ (d-2) $-ball $ B_\tau ; $ its
boundary is the {\sl real\/}
sphere described by
$ \vec  \zeta  = \lbrack 2\ {\rm e}^{i\tau}( {\rm cos} \ \tau -1)\rbrack^{
1/2}\vec  \alpha , $ $ \vec  \alpha  \in  \Bbb S^{d-3} $ (corresponding to $
\lambda =0 $
in Eq.(3.21)), and it contains the apex of the paraboloid $ \Pi_ \tau , $
namely
the point $ z_\tau , $ corresponding to the center $ \vec  \zeta  = 0 $ of the
ball
(obtained for $ \lambda =1 $ in Eq.(3.21)). \par
The orientation of $ h_\tau $ is obtained (for $ \tau  \in  \dot {\Bbb C}) $
by
continuity from the orientation of the real horocycles $ h^+_w = h_{iw} $ (see
\S 3.2). We shall now prove: \par
\medskip
{\bf Lemma 3 : }For every domain $ D $ in the class $ {\cal D}, $ the
following properties hold:
$$ \eqalignno{ {\rm a)} & \ \ \bar  \omega \left( \bigcup^{ }_{ \tau \in\dot
{\cal J}(D)} \underline{h}_\tau \right) = \underline{D} &  \cr {\rm b)} & \ \
\bigcup^{ }_{ \tau \in\dot {\cal J}(D)} \underline{h}_\tau \subset D &  \cr}
$$
\par
{\bf Proof :} \par
b) is immediately implied by a), since $ D = \bar  \omega^{ -1}(
\underline{D}). $ In order
to prove a), we notice that for each $ \tau $ in $ \dot {\cal J}(D), $ $ \bar
 \omega \left( \underline{h}_\tau \right) = \left\{ z^{(d-1)} = {\rm cos} \
\theta \right. $
$ =1+\lambda( {\rm cos} \ \tau -1); $ $ 0 \leq  \lambda  \leq  1\} , $ as it
results from Eqs.(3.1c) and (3.21); since $ \underline{D} $
is star-shaped (with respect to $ z^{(d-1)}=1), $ we then have: $ \bar
\omega \left( \underline{h}_\tau \right) \subset  \underline{D} $ and
therefore $ \bar  \omega \left( \bigcup^{ }_{ \tau \in\dot {\cal J}(D)}
\underline{h}_\tau \right)\subset  \underline{D} . $
Conversely, each point in $ \underline{D} $ is of the form $ z^{(d-1)} = {\rm
cos} \ \tau , $ with $ \tau  \in  \dot {\cal J}(D), $
and is the image by the projection $ \bar  \omega $$ $\nobreak\ of the apex of
$ \Pi_ \tau , $ namely the
point $ z_\tau  = (-i {\rm sin} \ \tau ,\ {\rm cos} \ \tau) $ of $ \hat
X^{(c)}_{d-1}, $ which belongs to $ \underline{h}_\tau . $ \par
\medskip
{\sl The integral transformation\/} $ \goth R^{(c)}_d: $ \par
Let $ {\cal F}(z) $ be an $ H^{(c)} $-invariant function analytic in a domain
$ D $
of $ X^{(c)}_{d-1} $ with $ D \in  {\cal D}. $ We then define the Radon-type
transformation $ {\cal F} \longrightarrow\hat {\cal F} = \goth
R^{(c)}_d({\cal F}) $
by the following formula (in which $ {\rm d} \tau $ is expressed by Eq.(3.6)):
$$ \hat {\cal F} \left(z_\tau \right) =2(-1)^{d-2}i {\rm e}^{i(d-2)\tau}
\int^{ }_{ h_\tau}{\cal F}(z) \left.{ {\rm d} z^{(0)}\wedge ...\wedge {\rm d}
z^{(d-1)} \over {\rm d} s(z)\wedge {\rm d} \tau} \right\vert_{ \Pi_ \tau} \ .
\eqno (3.22) $$
\par
In view of lemmas 3, the integrand at the r.h.s. of this
formula is analytic in the (compact) integration region $ \underline{h}_\tau ,
$ for all $ \tau $
in the domain $ \dot {\cal J}(D). $ Since $ h_\tau $ is a relative cycle whose
boundary
belongs to the analytic set $ {\cal C}_\tau , $ and since it \lq\lq varies
continuously\rq\rq\ with $ \tau $
(see $\lbrack$Ph$\rbrack$), the integral defines a function $ \hat {\cal F}
\left(z_\tau \right), $ analytic in the
domain $ \hat  D = \left\{ z_\tau  \in  \hat  X^{(c)}; \right. $ $ \left.\tau
\in  \dot {\cal J}(D) \right\} $ of the complex hyperbola $ \hat  X^{(c)}. $
\par
Now, since $ {\cal F}(z) = \underline{f} \left(z^{(d-1)} \right), $ we can
rewrite Eq.(3.22) (in
view of Eqs.(3.1c) and (3.7)) as follows:
$$ \hat {\cal F} \left(z_\tau \right) = \int^{ }_{ B_\tau} \underline{f}
\left( {\rm cos} \ \tau  - {1 \over 2} \vec  \zeta^{ \ 2} {\rm e}^{-i\tau}
\right) {\rm d}\vec  \zeta \ . \eqno (3.23) $$
\par
In the integral at the r.h.s. of Eq.(3.23), we now use the
following holomorphic parametrization of $ \zeta $-space, obtained by putting
$ z^{(d-1)}= {\rm cos} \ \theta $ in Eq.(3.1c):
$$ \matrix{\displaystyle\vec  \zeta  = \left[2\ {\rm e}^{i\tau}( {\rm cos} \
\tau  - {\rm cos} \ \theta) \right]^{1/2}\vec  \alpha \ , \cr\displaystyle
{\rm with} \ \ \ \ \vec  \alpha  \in  \Bbb C^{d-2}\ \ ,\ \ \ \vec  \alpha^
2=1\ ,\ \ \ \ \ \theta  \in  \Bbb C \cr} \eqno (3.24) $$
\par
In view of Eq.(3.21), the ball $ B_\tau $ is represented in $ \left(\theta
,\vec  \alpha \right)
$-space by the set $ \left\{ \left(\theta ,\vec  \alpha \right); \right. $ $
\theta  \in  \underline{\gamma}_ \tau ; $ $ \left.\vec  \alpha  \in  \Bbb
S^{d-3} \right\} $ where $ \underline{\gamma}_ \tau  = \{ \theta =\theta(
\lambda) ; $ $ {\rm cos} \ \theta( \lambda) -1=\lambda( {\rm cos} \ \tau -1);
$
$ 0\leq \lambda \leq 1, $ $ \theta( 0)=1, $ $ \theta( 1)=\tau\} , $ and
Eq.(3.24) yields:
$$ {\rm d}\vec  \zeta  = {\rm e}^{i \left({d-2 \over 2} \right)\tau} \cdot[ 2(
{\rm cos} \ \tau - {\rm cos} \ \theta)]^{{ d-4 \over 2}} {\rm sin} \ \theta \
{\rm d} \theta \ {\rm d} {\bf \sigmabf}_{ d-3} \left(\vec  \alpha \right)\ ,
\eqno (3.24^\prime ) $$
where $ {\rm d} {\bf \sigmabf}_{ d-3} \left(\vec  \alpha \right) $ denotes the
invariant measure on $ \Bbb S^{d-3}. $ \par
We can therefore replace Eq.(3.23) by the following
expression of $ \hat  f(\tau)  = \hat {\cal F} \left(z_\tau \right) $ in terms
of $ f(\theta)  = \underline{f} \left(z^{(d-1)} \right)_{ \left\vert
z^{(d-1)}= {\rm cos} \ \theta \right.}\equiv{\cal F}(z); $
this new form of the transformation $ \goth R^{(c)}_d $ (equivalent to
Eq.(3.22))
will be called $ {\cal R}^{(c)}_d, $ namely:
$$ \hat  f(\tau) = \left[{\cal R}^{(c)}_d({\cal F}) \right](\tau)  =
-\omega_{ d-2} {\rm e}^{i \left({d-2 \over 2} \right)\tau} \int^{ }_{ \gamma_
\tau} f(\theta)[ 2( {\rm cos} \ \tau - {\rm cos} \ \theta)]^{{ d-4 \over 2}}
{\rm sin} \ \theta \ {\rm d} \theta \ ; \eqno (3.25) $$
in the latter, $ \gamma_ \tau $ denotes the $ ^{\prime\prime} {\rm
ray}^{\prime\prime} $ $ \underline{\gamma}_ \tau $ oriented from 0 to $ \tau .
$
Moreover, the relevant branch of the function $ [2( {\rm cos} \ \tau - {\rm
cos} \ \theta)]^{{ d-4 \over 2}} $
is specified (for $ d $ odd) by the condition that for $ \tau =iw $ and $
\theta =iv $ (with
$ w > v), $ it takes the value $ [2( {\rm cosh} \ w- {\rm cosh} \ v)]^{{d-4
\over 2}}\geq 0. $ In fact, when $ \tau =iw $
$ (w>0), $ the horocycle $ h_\tau =h_{iw} $ is real and carried by the
hyperboloid $ X_{d-1}; $
formula (3.22) is then identical with the defining formula (3.10) of
the transformation $ \goth R_d, $ the function $ F $ being replaced here by $
{\cal F}_{ \left\vert X^+ \right.}. $
Accordingly, formula (3.25) can then be rewritten as follows:
$$ \hat {\cal F} \left(z_{iw} \right) = \hat  f(iw) = \omega_{ d-2} {\rm e}^{-
\left({d-2 \over 2} \right)w} \int^ w_0f(iv)[2( {\rm cosh} \ w- {\rm cosh} \
v)]^{{d-4 \over 2}} {\rm sinh} \ v\ {\rm d} v, \eqno (3.25^\prime ) $$
and thus coincides with the expression (3.20) of $ \hat {\rm f} ={\cal R}_d(F)
$ (see \S 3.3).
In the case when $ D $ is a cut-domain $ D_\mu , $ considered below in \S 3.6,
this
transformation will be applied to the boundary values $ {\cal F}_{\pm} $ of $
{\cal F} $ on $ X_+ $
and in particular to the corresponding discontinuity function $ \Delta{\cal F
}= i \left({\cal F}_+-{\cal F}_- \right) $
(see Eq.(3.32)). \par
We notice that the integral at the r.h.s. of Eq.(3.25)
considered with respect to the integration variable $ {\rm cos} \ \theta , $
defines an
(extended) Abel transform of the function $ \underline{f}( {\rm cos} \
\theta) = f(\theta) , $ taken along
the rays of the star-shaped domain $ \underline{D} ; $ this leads us to call $
\goth R^{(c)}_d $
a \lq\lq Radon-Abel transformation\rq\rq\ and to summarize its main properties
in
the following \par
\medskip
{\bf Proposition 10} : \par
Every $ H^{(c)} $-invariant function $ {\cal F}(z) $ analytic in a domain $ D
$ of $ X^{(c)}_{d-1} $
$ (d\geq 3) $ belonging to the class $ {\cal D} $ (i.e. such that $ {\cal
F}(z) \equiv  f(\theta) $ is analytic
and even in $ \dot {\cal J}(D)) $
admits a Radon-Abel transform $ \hat {\cal F} = \goth R^{(c)}_d({\cal F}), $
analytic in the
corresponding domain $ \hat  D = \left\{ z_\tau ;\tau \in\dot {\cal J}(D)
\right\} $ of $ \hat  X^{(c)} $ and satisfying the
following symmetry relation:
$$ \hat {\cal F} \left(z_{-\tau} \right) = (-1)^d {\rm e}^{-i(d-2)\tau}\hat
{\cal F} \left(z_\tau \right)\ . \eqno (3.26) $$
This transform, which is
defined in Eq.(3.22) by a Radon-type integral on the horocycle $ h_\tau , $ is
represented equivalently in the parameter space by an Abel-type
transform $ \hat  f(\tau)  \equiv  \hat {\cal F} \left(z_\tau \right) $ of $
f(\theta) $ (defined by Eq.(3.25)), analytic in the
domain $ \dot {\cal J}(D) $ of $ \dot {\Bbb C} $ and enjoying the following
structure:
$$ \hat  f(\tau)  = {\rm e}^{i \left({d-2 \over 2} \right)\tau} \left( {\rm
sin}  {\tau \over 2} \right)^{d-2} a( {\rm cos} \ \tau) , \eqno (3.26^\prime )
$$
the function $ a $ being analytic in $ \underline{D} . $ \par
\medskip
The proof of proposition 10 will be easily completed by deducing from
Eq.(3.25) the structural form (3.26') of $ \hat  f(\tau) . $ In fact, by using
the
parametrization $ {\rm cos} \ \theta  = 1+\lambda( {\rm cos} \ \tau -1) $ of
each ray $ \gamma_ \tau $ (with $ {\rm cos} \ \tau $ in the star-shaped
domain $ \underline{D} ), $ one can rewrite Eq.(3.25) as follows:
$$ \hat  f(\tau)  = {\omega_{ d-2} \over 2} {\rm e}^{i \left({d-2 \over 2}
\right)\tau}[ 2( {\rm cos} \ \tau -1)]^{{d-2 \over 2}} \int^ 1_0
\underline{f}(1+\lambda( {\rm cos} \ \tau -1))(1-\lambda)^{{ d-4 \over 2}}
{\rm d} \lambda \ ; \eqno (3.27) $$
since the integral at the r.h.s. of this equation is an analytic function of $
{\rm cos} \ \tau , $
defined (as $ f) $ in the domain $ \underline{D} $ (for all $ d \geq  3), $
the form (3.26') of $ \hat  f(\tau) $ follows.
{}From Eq.(3.26'), one recovers the fact that $ \hat  f(\tau) $ is $ 2\pi
$-periodic and therefore
analytic in $ \dot {\cal J}(D), $ and one checks that it satisfies the
symmetry relation
$$ \left(S_d \right)\ \ \ \ \ \ \hat  f(\tau)  = (-1)^d {\rm
e}^{i(d-2)\tau}\hat  f(-\tau) \ , $$
equivalent to Eq.(3.26). \par
\medskip
{\bf 3.5 The transformation} $ { \underline{\goth R}}_d $ {\bf for invariant
kernels on the
sphere }$ {\bf S}_{ {\bf d-1}} $ \par
\smallskip
Let us consider the case of kernels $ {\bf K} \left(z,z^{\prime} \right) $ on
the
sphere $ S_{d-1} $ satisfying the following invariance property: $ {\bf K}
\left(z,z^{\prime} \right)= {\bf K} \left(gz,gz^{\prime} \right), $
$ \forall g \in  {\rm SO}(d), $ where $ {\rm SO}(d) $ denotes here the
subgroup of $ G^{(c)} $  which
leaves the sphere $ S_{d-1} $ of $ X^{(c)} $ invariant. If $ {\bf H}  = {\rm
SO}(d) \cap  H^{(c)} $
denotes the stabilizer of $ z_0 $ in $ {\rm SO}(d), $ we can say that each
invariant kernel $ {\bf K} $ on $ S_{d-1} $ can be identified with an $ {\bf
H}
$-invariant function $ {\bf F}(z) = {\bf K} \left(z,z_0 \right) $ on $ S_{d-1}
$
(since, $ \forall g \in  {\bf H} , $ $ {\bf K} \left(z,z_0 \right) = {\bf K}
\left(gz,z_0 \right)). $ Such a
function $ {\bf F}(z) $ only depends on the coordinate $ z^{(d-1)}= {\rm cos}
\ u, $
and it can obviously be considered as well as a function defined
on the
following subset $ D_S = \left\{ z\in X^{(c)}_{d-1}; \right. $ $ \left.-1
\leq  z^{(d-1)}\leq +1 \right\} $ of $ X^{(c)}_{d-1} $ (such that
$ \bar  \omega \left(D_S \right) = \bar  \omega \left(S_{d-1} \right) =
[-1,+1]); $ we shall also put: $ {\bf f}(u)\equiv {\bf F}(z)\equiv
\underline{{\bf f}}( {\rm cos} \ u). $
The previous introduction of the transformation $ \goth R^{(c)}_d $ in \S 3.4
still applies to the present case (although the domain $ D $ is now
replaced by the \lq\lq flat\rq\rq\ set $ D_S), $ and we can now define the
transformation $ \underline{{\goth R}}_d $ which associates with $ {\bf F} $
the
function $ \hat {\bf F} \left(z_t \right) = \lbrack \underline{{\goth R}}_d(
{\bf F})\rbrack \left(z_t \right) $ (or $ \hat {\bf f}(t) = \left[
\underline{{\cal R}}_d( {\bf F}) \right](t)) $ obtained
by restricting formulae (3.22)...(3.25) to the set of {\sl real\/}
values $ (\tau =t) $ of the variable $ \tau , $ namely:
$$ \eqalignno{\hat {\bf F} \left(z_t \right) & = \int^{ }_{ B_t}
\underline{{\bf f}} \left( {\rm cos} \ t-{1 \over 2}\vec  \zeta^{ \ 2} {\rm
e}^{-it} \right) {\rm d}\vec  \zeta &  (3.28a) \cr  &  \equiv  \hat {\bf
f}(t)=-\omega_{ d-2} {\rm e}^{i \left({d-2 \over 2} \right)t} \int^ t_0 {\bf
f}(u)[2( {\rm cos} \ t- {\rm cos} \ u)]^{{d-4 \over 2}} {\rm sin} \ u\ {\rm d}
u & (3.28b) \cr} $$
\par
By taking into account the relevant branch of the factor $ [2( {\rm cos} \
\tau - {\rm cos} \ \theta)]^{{ d-4 \over 2}} $
in formula (3.25), Eq.(3.28b) can be given the following more precise forms
(involving a positive bracket):
$$ \eqalignno{ {\rm if} \ t\geq 0,\  & \hat {\bf f}(t) = (-i)^{d-2}\omega_{
d-2} {\rm e}^{i \left({d-2 \over 2} \right)t} \int^ t_0 {\bf f}(u)[2( {\rm
cos} \ u- {\rm cos} \ t)]^{{d-4 \over 2}} {\rm sin} \ u\ {\rm d} u & (3.29a)
\cr {\rm if} \ t\leq 0,\  & \hat {\bf f}(t) = (i)^{d-2}\omega_{ d-2} {\rm
e}^{i \left({d-2 \over 2} \right)t} \int^ t_0 {\bf f}(u)[2( {\rm cos} \ u-
{\rm cos} \ t)]^{{d-4 \over 2}} {\rm sin} \ u\ {\rm d} u & (3.29b) \cr} $$
\par
\vfill\eject
{\bf Remark} {\bf :} \par
In contrast with the case of the transformation $ \goth R_d $ of Volterra
kernels
on $ X_{d-1} $ (see Eqs.(3.10), (3.11)), where the integration ball $ B_{iw} $
represents a real
horocycle $ h^+_w $ of $ X_{d-1}, $ the integration ball $ B_t $ of formula
(3.28a) represents
a complex horocycle $ h_t $ whose support is {\sl not contained\/} in $
S_{d-1}, $ but has only in common with $ S_{d-1} $
the point $ z_t=(-i\ {\rm sin} \ t,\allowbreak\vec  0,\allowbreak {\rm cos} \
t) $ (represented by the center $ \left(\vec  \zeta =0 \right) $ of $ B_t). $
\par
The following property holds: \par
\smallskip
{\bf Proposition 11} : \par
For every invariant kernel $ {\bf K} $ on $ S_{d-1} $ with $ {\bf K}
\left(z,z_0 \right)\equiv {\bf F}(z), $
the corresponding transform $ \hat {\bf F}= \underline{{\goth R}}_d( {\bf F})
$ defined by
Eq.(3.28a) is a function on the meridian circle $ \hat  S $ of $ S_{d-1} $ in
the plane $ \left\{\vec  z=0 \right\} $ $ (\hat  S = \left\{ z_t;\
z^{(0)}_t=-i\ {\rm sin} \ t, \right. $ $ \left.z^{(d-1)}_t= {\rm cos} \ t;\
t\in \Bbb R \right\} ), $ which satisfies the following property:
$$ \hat {\bf F} \left(z_{-t} \right) = (-1)^d {\rm e}^{-i(d-2)t}\hat {\bf F}
\left(z_t \right) \eqno (3.30) $$
\par
Moreover, this transformation is equivalently represented by the
expression of the function $ \hat {\bf f}(t) = \hat {\bf F} \left(z_t \right)
$ in terms of $ {\bf f}(u)\equiv {\bf F}(z) $ given by
Eqs.(3.29); $ \hat {\bf f}(t) $ is a $ 2\pi $-periodic function of the
following form:
$$ \hat {\bf f}(t) = {\rm e}^{i \left({d-2 \over 2} \right)t} \left( {\rm sin}
 {t \over 2} \right)^{d-2} {\bf a}( {\rm cos} \ t) \eqno (3.30^\prime ) $$
\par
The proof of Eq.(3.30') is completely similar to that of Eq.(3.26') in
proposition 10 (see \S 3.4), the representation (3.27) being now derived and
used
for $ \tau =t $ real; Eq.(3.30) follows directly from Eq.(3.30'). \par
\medskip
{\bf 3.6 Radon-Abel transformation }$ \goth R^{( {\bf c})}_ {\bf d} ${\bf\ for
algebras of
invariant perikernels of moderate growth on }$ {\bf X}^{( {\bf c})}_{ {\bf
d-1}} $ \par
\medskip
$ H^{(c)} $-{\sl invariant triplets\/} $ ({\cal F}, {\bf F} ,F) $ {\sl on\/} $
X^{(c)}_{d-1} $ \par
We shall now consider exclusively the case of $ H^{(c)}
$-invariant functions $ {\cal F} $ analytic in \lq\lq cut-domains\rq\rq\ $
D_\mu = \left\{ z\in X^{(c)}_{d-1}; \right. $ $ \left.z^{(d-1)}\in
\underline{D}_\mu \right\} , $
where $ \underline{D}_\mu  = \Bbb C\backslash[ {\rm cosh} \ \mu ,+\infty[ . $
Each such function $ {\cal F} $ can be
represented equivalently by the function $ \underline{f} \left(z^{(d-1)}
\right)\equiv{\cal F}(z), $ analytic
in $ \underline{D}_\mu , $ or by the {\sl even\/} periodic function $
f(\theta)  = \underline{f}( {\rm cos} \ \theta) , $
analytic in $ \dot {\cal J}^{(\mu)} . $ \par
It is convenient (as in \S 2.3 for the case $ d=2) $ to
consider $ {\cal F} $ as participating in a triplet $ ({\cal F}, {\bf F} ,F),
$ where: \par
\smallskip
a) $ {\bf F}  = {\cal F}_{ \left\vert S_{d-1} \right.}, $ with the
corresponding equivalent
representations: $ {\bf F}(z) = \underline{{\bf f}} \left(z^{(d-1)} \right) $
$ (z_{d-1}\in[ -1,+1]), $ and $ \underline{{\bf f}}( {\rm cos} \ \theta) =
{\bf f}(\theta) $
$ (\theta \in  \dot {\Bbb R}). $ \par
b) $ F(z) = Y \left(z^{(0)} \right)\Delta{\cal F}_{ \left\vert X_{d-1}
\right.}(z) $ (see \S 3.1), with the
corresponding equivalent representations: $ F(z)= \underline{{\rm f}}
\left(z^{(d-1)} \right) $
(with $ {\rm supp} \ \underline{{\rm f}}  \subset[ {\rm cosh} \ \mu ,+\infty[
) $ and $ \underline{{\rm f}}( {\rm cosh} \ v) = {\rm f}(v)=\Delta f(v) $
(with $ {\rm supp} \ f\subset[ \mu ,+\infty[ ). $ \par
\smallskip
As stated in our introduction (section 1), any such
triplet $ ({\cal F}, {\bf F} ,F) $ represents an invariant triplet $ ({\cal
K}, {\bf K} ,K): $ $ {\cal K} \left(z,z^{\prime} \right) $
is a $ G^{(c)} $-invariant perikernel analytic in a domain $ \Delta_ \mu $ of
$ X^{(c)}_{d-1} $
(see \S 3.1) such that $ {\cal K} \left(z,z_0 \right)={\cal F}(z); $
$ {\bf K} $ is the invariant kernel on the sphere $ S_{d-1} $ defined by $
{\cal K}_{ \left\vert S_{d-1}\times S_{d-1} \right.}, $
and $ K $ is the invariant Volterra kernel on $ X_{d-1} $ defined by $ K
\left(z,z^{\prime} \right) = Y \left(z^{(0)}-z^{\prime( 0)} \right) $
$ \times \Delta{\cal K}_{ \left\vert X_{d-1}\times X_{d-1} \right.}
\left(z,z^{\prime} \right) $(see \S 3.1) such that: $ {\bf F}(z) = {\bf K}
\left(z,z_0 \right) $ and $ F(z) = K \left(z,z_0 \right). $ \par
\medskip
{\sl Action of the transformation\/} $ {\cal R}^{(c)}_d: $ \par
Let $ ({\cal F}, {\bf F} ,F) $ be a given triplet, with $ {\cal F} $ analytic
in the
domain $ D_\mu . $ In view of proposition 10, the Radon-Abel
transform $ \hat {\cal F} = {\cal R}^{(c)}_d({\cal F}) $ of such a function $
{\cal F} $ is then analytic in
the corresponding domain $ \hat  D_\mu  = \left\{ z_\tau ;\tau \in\dot {\cal
J}^{(\mu)} \right\} $ and satisfies there
Eqs.(3.36), (3.26'). \par
If the boundary values $ {\cal F}_{\pm} $ and the discontinuity of $ {\cal F}
$
are sufficiently regular (e.g. continuous functions), it is clear
that formula (3.22) yields the corresponding boundary values $ \hat {\cal
F}_{\pm} $
and the discontinuity function $ \hat  F = Y \left(z^{(0)} \right)\Delta\hat
{\cal F}_{ \left\vert\hat  X \right.}(z) $ of $ \hat {\cal F}. $ $ \hat  F $
is defined on
the \lq\lq right-hand branch\rq\rq\ of the hyperbola $ \hat  X, $ namely $
\left\{ z_\tau \in\hat  X;\tau =iw,w\in \Bbb R \right\} $
and its support is contained in the set $ \hat  X^\mu_ +, $ where: $ \hat
X^\mu_ + = X^\mu_ +\cap\hat  X. $
We can thus write (in view of Eq.(3.22)):
$$ \forall w,\ \ \ \ \hat  F \left(z_{iw} \right) =2(-1)^{d-2} {\rm
e}^{-(d-2)w} \int^{ }_{ h_{iw}}F(z) \left.{ {\rm d} z^{(0)}\wedge ...\wedge
{\rm d} z^{(d-1)} \over {\rm d} s(z)\wedge {\rm d} w} \right\vert_{ \Pi^ +_w}\
, \eqno (3.31) $$
or correspondingly, for $ \hat {\rm f}(w) = \hat  F \left(z_{iw} \right) $ (in
view of Eq.(3.25)):
$$ \hat {\rm f}(w) = \omega_{ d-2}\ {\rm e}^{- \left({d-2 \over 2} \right)w}
\int^ w_0 {\rm f}(v)[2( {\rm cosh} \ w- {\rm cosh} \ v)]^{{d-4 \over 2}} {\rm
sinh} \ v\ {\rm d} v\ . \eqno (3.32) $$
\par
\smallskip
{\bf Remark} : $ \hat {\cal F}_+ $ and $ \hat {\cal F}_- $ being defined (in a
way similar
to $ {\cal F}_+ $ and $ {\cal F}_-) $ as the boundary values of $ \hat {\cal
F} $ on $ \hat  X $ from the
respective sides $ {\rm Im} \left(z_{iw}-z_0 \right)^2 < 0 $ and $ {\rm Im}
\left(z_{iw}-z_0 \right)^2 > 0, $ one would
check that the discontinuity function $ \hat  F \left(z_{iw} \right)=\hat
{\rm f}(w) $ is defined (as
$ {\rm f}(v)) $ by the formula:
$$ \matrix{\displaystyle \forall w \geq 0,\ \ \ \ \ \hat {\rm f}(w) = i
\dlowlim{ {\rm lim}}{ \matrix{\displaystyle \eta \longrightarrow 0
\cr\displaystyle \eta >0 \cr}}  \left[\hat  f(iw+\eta) -\hat  f(iw-\eta)
\right] \cr} $$
\par
By comparing Eq.(3.31) with Eq.(3.10), we also see that $ \hat  F $
is nothing else than the function $ \hat  F = \goth R_d(F) $
representing an invariant Volterra kernel $ \hat  K $ in $ V \left(\hat  X
\right)^{\natural} , $ such that $ \hat  K \left(z,z_0 \right)=\hat  F(z) $
(see \S 3.3); the corresponding Abel-type representation (3.32) of
$ \hat {\rm f} = {\cal R}_d(F) $ coincides with Eq.(3.20).  \par
\smallskip
Finally, by taking Eq.(3.22) for $ \tau =t $ real (as in
proposition 11), we can deduce from the relation $ \hat {\cal F} = \goth
R^{(c)}_d({\cal F}) $
the following one: $ \hat {\bf F} = \hat {\cal F}_{ \left\vert\hat  S \right.}
= \underline{{\goth R}}_d( {\bf F}). $ \par
\medskip
{\it In this subsection, our main results will concern
algebras of invariant perikernels with moderate growth on} $ X^{(c)}_{d-1} $
{\it (in proposition 12) and the action of the Radon-Abel
transformation} $ {\goth R}^{(c)}_d $ {\it on the latter (in proposition
13): these results will parallel to some extent those of \S 3.3
concerning the algebras of invariant Volterra kernels on $ X_{d-1} $ and
the action of the transformation} $ {\goth R}_d $ {\it on the latter (cf.
proposition 9)}. \par
\medskip
{\bf Invariant perikernels of moderate growth on }$ {\bf X}^{( {\bf c})}_{
{\bf d-1}} $ \par
Let us first introduce
appropriate normed spaces of $ H^{(c)} $-invariant analytic functions
of moderate growth in the cut-domains $ D_\mu $ $ (\mu \geq 0). $ \par
\medskip
{\bf Definition 3} : For each $ H^{(c)} $-invariant function $ {\cal F}, $
analytic in a given domain $ D_\mu , $ one defines: \par
\noindent$ \forall m, $ $ m\in \Bbb R; $ $ \forall \rho , $ $ \rho \geq 1, $
\smallskip
$$ {\cal G}_m[{\cal F}](\rho)  = \rho^{ -m} \dlowlim{ {\rm sup}}{ \left\{
z;z^{(d-1)}\in E_\rho \right\}}  \vert{\cal F}(z)\vert =\rho^{ -m} \dlowlim{
{\rm sup}}{z^{(d-1)}\in E_\rho} \left\vert \underline{f} \left(z^{(d-1)}
\right) \right\vert ,\   $$
$ E_\rho $ being the ellipse specified in definition 1 (Eq.(2.47)). \par
\smallskip
a) We call $ \left[{\cal V}_1 \right]^m_\mu \left(X^{(c)}_{d-1} \right) $ the
space of such functions $ {\cal F} $ for
which:
$$ \Vert{\cal F}\Vert^{( 1)}_m =  {\cal G}_m[{\cal F}](1) + \int^{ \infty}_
1{\cal G}_m[{\cal F}](\rho)  { {\rm d} \rho \over \sqrt{ \rho^ 2-1}} < \infty
$$
\par
b) We call $ \left[{\cal V}_{\infty} \right]^m_\mu \left(X^{(c)}_{d-1} \right)
$ (resp.$ \left[{\cal V}^\ast_{ \infty} \right]^m_\mu \left(X^{(c)}_{d-1}
\right) $ in the case $ m\in {\bf N} ) $
the space of such functions $ {\cal F} $ for which $ \Vert{\cal F}\Vert^{(
\infty)}_ m<\infty $ (resp. $ \Vert{\cal F}\Vert^{( \infty) \ast}_ m < \infty
), $
$ \Vert{\cal F}\Vert^{( \infty)}_ m $ and $ \Vert{\cal F}\Vert^{( \infty)
\ast}_ m $ being defined by the same formulae (2.49) and
(2.50) as for the case $ d=2 $ (definition 1). \par
\medskip
By using the identification relations: $ {\cal F}(z) = \underline{f}
\left(z^{(d-1)} \right) $
and $ \underline{f}( {\rm cos} \ \theta)  = f(\theta) , $ one readily obtains
the following statement,
analogous to lemma 1 and whose proof is similar. \par
\medskip
{\bf Lemma 4\nobreak\ :} \par
\smallskip
a) the following inequalities hold (with suitable
constants and $ g_m[f] $ defined by Eq.(2.1)):
$$ \matrix{\displaystyle \hfill C^{\prime} g_m[f](v) & \displaystyle \leq
{\cal G}_m[{\cal F}]( {\rm cosh} \ v) \leq  C\ g_m[f](v), \hfill
\cr\displaystyle \hfill {\rm Cst}\Vert{\cal F}\Vert^{( 1)}_m & \displaystyle
\leq  \left\Vert g_m[f] \right\Vert_ 1 + g_m[f](0) \leq  {\rm
Cst}^{\prime}\Vert{\cal F}\Vert^{( 1)}_m \hfill \cr} \eqno (3.33) $$
\par
b) the following homeomorphisms hold:
$$ \left[{\cal V}_{\infty} \right]^m_\mu \left(X^{(c)}_{d-1} \right) \approx
\underline{{\cal O}}^{(2)}_m \left(\dot {\cal J}^{(\mu)} \right)\ ,\ \
\left[{\cal V}^\ast_{ \infty} \right]^m_\mu \left(X^{(c)}_{d-1} \right)
\approx  \underline{{\cal O}}^{(2)\ast}_ m \left(\dot {\cal J}^{(\mu)}
\right)\ . $$
\par
\smallskip
{\bf The algebras} $ {\cal W}^{(m)} \left(X^{(c)}_{d-1} \right)^{\natural} $
\par
We denote by $ {\cal W}^{(m)} \left(X^{(c)}_{d-1} \right)^{\natural} $the
space of invariant
perikernels $ {\cal K} \left(z,z^{\prime} \right) $ whose associated $ H^{(c)}
$-invariant function $ {\cal F}(z) = {\cal K} \left(z,z_0 \right) $
belongs to the space $ \left[{\cal V}_1 \right]^m_0 \left(X^{(c)}_{d-1}
\right) = \bigcup^{ }_{ \mu \geq 0} \left[{\cal V}_1 \right]^m_\mu
\left(X^{(c)}_{d-1} \right); $ $ {\cal W}^{(m)} \left(X^{(c)}_{d-1}
\right)^{\natural} $ is
considered as a normed space equipped with the norm $ \mid \mid \mid{\cal
K}\mid \mid \mid_ m=\Vert{\cal F}\Vert^{( 1)}_m. $ \par
We then have the following result (similar to proposition
9) whose proof will be given elsewhere. \par
\medskip
{\bf Proposition 12 :} \par
For each $ m, $ with $ m>-1, $ $ {\cal W}^{(m)} \left(X^{(c)}_{d-1}
\right)^{\natural} $ is a subalgebra of the algebra $ {\cal W}
\left(X^{(c)}_{d-1} \right)^{\natural} $ of
all invariant perikernels. More precisely, let $ {\cal K}_1, $ $ {\cal K}_2 $
be two
perikernels such that $ {\cal F}_i(z) = {\cal K}_i \left(z,z_0 \right) \in
\left[{\cal V}_1 \right]^m_{\mu_ i} \left(X^{(c)}_{d-1} \right) $ $ (i=1,2). $
Then the perikernel $ {\cal K }= {\cal K}_1\ast^{( c)}{\cal K}_2 $ is such
that: $ {\cal F}(z) = {\cal K} \left(z,z_0 \right) \in  \left[{\cal V}_1
\right]^m_{\mu_ 1+\mu_ 2} \left(X^{(c)}_{d-1} \right), $
and one has:
$$ \mid \mid \mid{\cal K}\mid \mid \mid_ m \leq  {\rm Cst} \mid \mid
\mid{\cal K}_1\mid \mid \mid_ m\mid \mid \mid{\cal K}_2\mid \mid \mid_ m
\eqno (3.34) $$
\par
\smallskip
{\bf Conservation of the rate of moderate growth }$ m $ {\bf by the
action of the transformation} $ \goth R^{(c)}_d $ \par
We shall now show that the Radon-Abel transformation
essentially preserves the behaviour at infinity (characterized by
$ \vert {\rm cos}(u+iv)\vert^ m \approx  {\rm e}^{m\vert v\vert} ); $ more
precisely, we shall prove the
following result which involves the spaces of analytic functions
$ {\cal O}^{(d)}_{m,p} $ introduced in \S 2.2 (with $ p $ integer, $ p\geq 0).
$ \par
\medskip
\vfill\eject
{\bf Proposition 13} : \par
\smallskip
The transformation $ \goth R^{(c)}_d $ defines a continuous
mapping from each subspace $ \left[{\cal V}_1 \right]^m_\mu
\left(X^{(c)}_{d-1} \right), $ with $ m>-1, $ into a corresponding
subspace $ {\cal O}^{(d)}_{m,p(d)} \left(\dot {\cal J}^{(\mu)} \right) $ and a
continuous mapping from each
subspace $ \left[{\cal V}_{\infty} \right]^m_\mu \left(X^{(c)}_{d-1} \right) $
into a corresponding subspace $ \underline{{\cal O}}^{(d)}_{m,p(d)}
\left(\dot {\cal J}^{(\mu)} \right), $
where $ p(d) $ is specified as follows: \par
\smallskip
a) for $ d $ even, $ d \geq  4, $ $ p(d) = {d-2 \over 2} $ \par
b) for $ d $ odd, \nobreak\ $ d \geq  3, $ $ p(d) = {d-3 \over 2} $ \par
\medskip
{\bf Proof} : We shall make use of Eq.(3.25) which defines the
analytic function $ \hat  f(\tau) $ in the domain $ \dot {\cal J}^{(\mu)} $
and which implies the
functional identity (3.26') (as stated in proposition 10). In order
to prove that $ \hat  f\in{\cal O}^{(d)}_{m,p(d)} \left(\dot {\cal J}^{(\mu)}
\right) $ (see \S 2.2, formula (2.33)), we now have to make relevant
estimates for all the functions
$$ g_m \left[\hat  f_r \right](w) = {\rm e}^{-mw} \dlowlim{ {\rm sup}}{0\leq
t\leq 2\pi}  \left\vert\hat  f^{(r)}(t+iw) \right\vert \ ,\ \ \ {\rm with} \ \
\ 0 \leq  r \leq  p(d) \eqno (3.35) $$
$ \hat  f^{(r)} $ denoting the derivative of order $ r $ of the function $
\hat  f(t+iw) $ $ \left(\hat  f\equiv\hat  f^{(0)} \right). $ \par
We shall first establish estimates for $ g_m \left[\hat  f \right] $ and will
treat separately the general case $ d\geq 4 $ and the special case $ d=3. $
\par
\smallskip
i) $ d\geq  4\ : $ we shall use the following majorization in the
integrand at the r.h.s. of Eq.(3.25) (valid for $ \theta =u+iv, $ and $ \tau
= t+iw, $
$ 0\leq v\leq w): $
$$ \matrix{\displaystyle \hfill\vert {\rm sin} \ \theta\vert  < {\rm e}^v\ ,\
\ \vert f(u+iv)\vert & \displaystyle \leq  {\rm e}^{mv}g_m[f](v) \hfill
\cr\displaystyle \hfill[ 2( {\rm cos} \ \tau  - {\rm cos} \ \theta)]^{{ d-4
\over 2}} & \displaystyle \leq  2^{d-4} {\rm e}^{ \left({d-4 \over 2}
\right)w}\ , \hfill \cr} $$
\par
Since the integrand is analytic in $ \dot {\cal J}^{(\mu)} , $ we can replace
the integration path $ \left(\gamma_ \tau \right) $ in (3.25) by $ [0,t] \cup
[t,t+iw] $ and then
obtain the following estimate:
$$ g_m \left[\hat  f \right](w) \leq  {\rm Cst\ e}^{-(m+1)w} \left[g_m[f](0)+
\int^ w_0g_m[f](v) {\rm e}^{(m+1)v} {\rm d} v \right] \eqno (3.35^\prime ) $$
\par
{}From the latter inequality and in view of lemma 4, we
then deduce (provided $ m>-1) $ the following continuity inequalities:
$$ \eqalignno{ \left\Vert g_m \left[\hat  f \right] \right\Vert_{ \infty} &
\leq  {\rm Cst} \left\Vert g_m[f] \right\Vert_{ \infty}  \leq  C\Vert{\cal
F}\Vert^{ \infty}_ m & (3.36) \cr \left\Vert g_m \left[\hat  f \right]
\right\Vert_ 1 & \leq  {\rm Cst}^{\prime} \left[g_m[f](0)+ \left\Vert g_m[f]
\right\Vert_ 1 \right] \leq  C_1\Vert{\cal F}\Vert^{( 1)}_m & (3.36^\prime )
\cr} $$
\par
ii) $ d=3: $ By using (as in i)) the integration path $ \{ \theta \in[
0,t]\cup[ t,t+iw]\} , $
we can rewrite Eq.(3.25) as follows:
$$ \matrix{\displaystyle\hat  f(t+iw) & \displaystyle = -2\ {\rm e}^{{i(t+iw)
\over 2}} \left[ \int^ t_0{f(u) {\rm sin} \ u\ {\rm d} u \over[ 2( {\rm
cos}(t+iw)- {\rm cos\ } u)]^{1/2}} + ... \right. \hfill \cr\displaystyle  &
\displaystyle \left.i \int^ w_0{f(t+iv) {\rm sin} \ t\ {\rm cosh\ } v\ {\rm d}
v \over[ 2( {\rm cos}(t+iw)- {\rm cos}(t+iv))]^{1/2}} - \int^ w_0{f(t+iw)
{\rm cos} \ t\ {\rm sinh} \ v\ {\rm d} v \over[ 2( {\rm cos}(t+iw)- {\rm
cos}(t+iv))]^{1/2}} \right] \hfill \cr} \eqno (3.37) $$
Majorizations on these three integrals are obtained by applying
respectively the following inequalities:
$$ \matrix{\displaystyle\vert {\rm cos}(t+iv)- {\rm cos} \ u\vert  \geq
\vert {\rm cos} \ u\cdot {\rm cosh} \ w- {\rm cos} \ t\vert ,\vert {\rm
cos}(t+iw)- {\rm cos}(t+iv)\vert \geq \left\vert {\rm sin} \ t \left({e^w-e^v
\over 2} \right) \right\vert , \hfill \cr\displaystyle\vert {\rm cos}(t+iw)-
{\rm cos} \ (t+iv)\vert  \geq  \vert {\rm cosh} \ w- {\rm cosh} \ v\vert \ ,
\hfill \cr} $$
which yields the following majorization for $ g_m \left[\hat  f \right](w) $
(by taking
Eq.(2.1) and definition 3 into account):
\smallskip
$$ \matrix{\displaystyle g_m \left[\hat  f \right](w) & \displaystyle \leq  2\
{\rm e}^{- \left(m+{1 \over 2} \right)w} \left[C_1 {g_m[f](0) \over( {\rm
cosh} \ w)^{1/2}} + \int^ w_0{g_m[f](v) {\rm e}^{(m+1)v} {\rm d} v \over
\left(e^w-e^v \right)^{1/2}} \right. \hfill \cr\displaystyle  & \displaystyle
\left.+ \int^{ {\rm cosh} \ w}_1{{\cal G}_m[{\cal F}](\rho) \rho^ m {\rm d}
\rho \over[ 2( {\rm cosh} \ w-\rho)]^{ 1/2}} \right]\ ; \hfill \cr} \eqno
(3.38) $$
In the latter, we have put $ C_1 = \dlowlim{ {\rm sup}}{0\leq t\leq \pi}
\int^ 1_{ {\rm cos} \ t}{ {\rm d} \lambda \over \left\vert \lambda  - { {\rm
cos} \ t \over {\rm cosh} \ w} \right\vert^{ 1/2}}\ . $ \par
We can make use of the inequality (3.38) in two ways: \par
\smallskip
$\alpha$) {\sl Use of \/}$ \Vert \ \ \Vert_{ \infty} ${\sl -norms\/}:
$$ \matrix{\displaystyle \forall w\ ,\ \ w\geq 0\ \ ,\ \ g_m \left[\hat  f
\right](w) \leq  {\rm Cst} \left\{ {\rm e}^{-(m+1)w}g_m[f](0) +... \matrix{^{
} \cr \cr \cr} \right. \hfill \cr\displaystyle \left. \left\Vert g_m[f]
\right\Vert_{ \infty} \left[ {\rm e}^{- \left(m+{1 \over 2} \right)w} \int^{
e^w}_1{y^m {\rm d} y \over \left(e^w-y \right)^{1/2}} \right]+\Vert{\cal
F}\Vert^{( \infty)}_ m \left[( {\rm cosh} w)^{- \left(m+{1 \over 2} \right)}
\int^{ {\rm cosh} w}_1{\rho^ m {\rm d} \rho \over( {\rm cosh} w-\rho)^{ 1/2}}
\right] \right\} \hfill \cr} \eqno (3.39) $$
\par
We now apply the following majorization, {\sl only valid for\/} $ m>-1: $
$$ \forall x,\ \ x\geq 1\ \ ,\ \ \int^ x_1{y^m {\rm d} y \over( x-y)^{1/2}}
\leq  a_mx^{m+{1 \over 2}}\ \ ,\ \ {\rm with} \ \ a_m= \int^ 1_0{t^m {\rm d} t
\over( 1-t)^{1/2}} $$
\par
In view of the latter (and of lemma 4), the inequality
(3.39) then implies the following norm inequality, valid for $ m>-1 $
(with a suitable $ m $-dependent constant):
$$ \left\Vert g_m \left[\hat  f \right] \right\Vert_{ \infty}  \leq
C\Vert{\cal F}\Vert^{( \infty)}_ m \eqno (3.40) $$
\par
$\beta$) {\sl use of \/}$ \Vert \ \ \Vert_ 1 ${\sl -norms\/}: \par
By integrating both sides of the inequality (3.38) with
respect to $ w $ (the equivalence of $ {\rm cosh} \ w $ and $ e^w $ being
taken into
account) and by making the change of variables $ x=e^w, $ $ (y=e^v), $ $ X =
{\rm cosh} \ w $
(respectively in the second and third integrals at the r.h.s. of
(3.38)), we readily obtain:
$$ \matrix{\displaystyle \left\Vert g_m \left[\hat  f \right] \right\Vert_
1\leq  {\rm Cst} \left\{ g_m[f](0) \int^{ \infty}_ 0 {\rm e}^{-(m+1)w} {\rm d}
w + ... \matrix{^{ } \cr \cr \cr} \right. \hfill \cr\displaystyle \left.\ \ \
\int^{ \infty}_ 1x^{- \left(m+{3 \over 2} \right)} {\rm d} x \int^ x_1{g_m[f](
{\rm \ell n} \ y)y^m \over( x-y)^{1/2}} {\rm d} y + \int^{ \infty}_ 1X^{-
\left(m+{1 \over 2} \right)}{ {\rm d} X \over \sqrt{ X^2-1}} \int^ X_1{{\cal
G}_m[{\cal F}](\rho) \rho^ m \over( X-\rho)^{ 1/2}} {\rm d} \rho \right\}
\hfill \cr} \eqno (3.41) $$
By making use of the following majorization, {\sl only valid for\/} $ m>-1: $
$$ \forall y,\ \ y\geq 1\ , \int^{ +\infty}_ y{x^{- \left(m+{3 \over 2}
\right)} {\rm d} x \over( x-y)^{1/2}} \leq  b_m\ y^{-(m+1)}\ \ ,\ \ {\rm with}
\ \ b_m= \int^{ \infty}_ 1{t^{- \left(m+{3 \over 2} \right)} {\rm d} t \over(
t-1)^{1/2}}\ , $$
and applying the Fubini theorem to the second and third integrals
at the r.h.s. of (3.41), we obtain (after taking into account the
inequality $ \left(X^2-1 \right)^{-1/2} \leq  \left(\rho^ 2-1 \right)^{-1/2}
\times  {\rho \over X}): $
$$ \left\Vert g_m \left[\hat  f \right] \right\Vert_ 1 \leq  {\rm Cst}
\left\{{ 1 \over m+1} g_m[f](0) + b_m \left[ \left\Vert g_m[f] \right\Vert_ 1
+ \int^{ \infty}_ 1{{\cal G}_m[{\cal F}](\rho) \over \sqrt{ \rho^ 2-1}} {\rm
d} \rho \right] \right\} $$
\par
In view of definition 3 and lemma 4, the latter
inequality readily yields the following norm inequality, valid
for $ m>-1 $ (with a suitable $ m $-dependent constant):
$$ \left\Vert g_m \left[\hat  f \right] \right\Vert_ 1 \leq  C_1\Vert{\cal
F}\Vert^{( 1)}_m \eqno (3.42) $$
\par
In order to obtain similar inequalities for the
derivatives $ \hat  f^{(r)} $ of $ \hat  f, $ we shall rely on algebraic
identities
which express the functions $ \hat  f^{(r)} $ in terms of the Radon-Abel
transforms $ \goth R^{(c)}_{d^{\prime}}({\cal F}) $ of {\sl lower orders\/} $
d^{\prime}  \leq  d. $ We thus need to use the
more precise notation $ \hat  f_d(\tau)  \equiv  \hat  f(\tau) , $ $ \hat
f^{(r)}_d(\tau)  \equiv  \hat  f^{(r)}(\tau) , $ which
specifies that the transformation $ {\cal F } \longrightarrow  \left({\cal
R}^{(c)}_d{\cal F} \right) \left(z_\tau \right)=\hat  f_d(\tau) $
(defined by Eq.(3.25)) operates on the $ (d-1) ${\sl -dimensional\/}
hyperboloid $ X^{(c)}_{d-1}. $ \par
Then we can state: \par
\medskip
{\bf Lemma 5} : For every $ d, $ with $ d\geq 4, $ and for every derivation
order $ r, $ with $ r \leq  {d-2 \over 2} $ for $ d $ even and $ r \leq  {d-3
\over 2} $ for $ d $ odd, the
corresponding function $ \hat  f^{(r)}_d(\tau) $ admits an expression of the
following form:
$$ \hat  f^{(r)}_d(\tau)  = \sum^ r_{r^{\prime} =0}a_{r^{\prime}}( \tau)
\hat  f^{(0)}_{d-2r^{\prime}}( \tau) \ , \eqno (3.43) $$
where all functions $ a_{r^{\prime}}( \tau) $ are polynomials of the variable
$ {\rm e}^{i\tau} , $
and the notation $ \hat  f^{(0)}_{d^{\prime}}  \equiv  \hat  f_{d^{\prime}} $
must be understood for $ d^{\prime} =2 $ in the
following sense: $ \hat  f^{(0)}_2(\tau)  = f(\tau) . $ \par
\medskip
{\bf Proof} : Eqs.(3.43) are obtained by taking successive
derivatives of both sides of Eq.(3.25) and applying a
straightforward recursive argument. \par
Let us now consider, for $ d \geq  4 $ and $ r \leq  p(d), $ (with $ p(d) $
defined in the statement of proposition 13) the functions
$ g_m \left[\hat  f^{(r)}_d \right](w) = {\rm e}^{-mw} \dlowlim{ {\rm
sup}}{0\leq t\leq 2\pi}  \left\vert\hat  f^{(r)}_d(t+iw) \right\vert . $ \par
\smallskip
By noticing that, for $ \tau $ varying in $ {\cal J}^{(\mu)}_ +, $ all the
functions $ a_{r^{\prime}}( \tau) $ are bounded by polynomials of $ {\rm
e}^{-w} $ and therefore
uniformly bounded, we immediately deduce from Eqs.(3.43):
$$ \left\Vert g_m \left[\hat  f^{(r)}_d \right] \right\Vert_{ \infty}  \leq  C
\sum^ r_{r^{\prime} =0} \left\Vert g_m \left[\hat  f^{(0)}_{d-2r^{\prime}}
\right] \right\Vert_{ \infty} \eqno (3.44) $$
and
$$ \left\Vert g_m \left[\hat  f^{(r)}_d \right] \right\Vert_ 1 \leq  C \sum^
r_{r^{\prime} =0} \left\Vert g_m \left[\hat  f^{(0)}_{d-2r^{\prime}} \right]
\right\Vert_ 1 \eqno (3.45) $$
Since $ d\geq 4 $ and $ r^{\prime}  \leq  p(d) $ imply $ d-2r^{\prime} \geq 2,
$ one can
majorize the r.h.s. of the inequalities (3.44), (3.45) by making use of
lemma 4 and of the
inequalities (3.36), (3.36') (for $ d-2r^{\prime} \geq 4) $ or (3.40), (3.42)
(for
$ d-2r^{\prime} =3). $ One thus obtains (for suitable constants $ C_r, $ $
C^{\prime}_ r): $
$$ \forall r \leq p(d)\ , \left\Vert g_m \left[\hat  f^{(r)}_d \right]
\right\Vert_{ \infty}  \leq  C_r\Vert{\cal F}\Vert^{ \infty}_ m,\ {\rm and} \
\ \left\Vert g_m \left[\hat  f^{(r)}_d \right] \right\Vert_ 1\leq C^{\prime}_
r\Vert{\cal F}\Vert^ 1_m, $$
and therefore, in view of the defining equations (2.34), (2.33):
$$ \mid \mid \mid\hat  f_d\mid \mid \mid^{ \infty}_{ m,p(d)} \leq  C
\Vert{\cal F}\Vert^{ \infty}_ m,\ \ {\rm and} \ \ \mid \mid \mid\hat  f_d\mid
\mid \mid^ 1_{m,p(d)} \leq  C^{\prime}\Vert{\cal F}\Vert^ 1_m\ , $$
which express the continuity properties stated in proposition 13. \par
\smallskip
As a by-product of the previous estimates, we also
obtain: \par
\smallskip
{\bf Corollary :} Let $ {\rm f}(v) $ be such that $ {\rm g}_m[ {\rm f}](v) =
{\rm e}^{-mv} {\rm f}(v) \in L^1 \left(\Bbb R^+ \right), $
and let $ \hat {\rm f}(w) $ be defined in terms of $ {\rm f}(v) $ by
Eq.(3.20). \par
Then the functions $ {\rm g}_m \left[\hat {\rm f}^{(r)} \right](w) = {\rm
e}^{-mw}\hat {\rm f}^{(r)}(w) $ satisfy
majorizations of the following type:
$$ \dlowlim{ {\rm max}}{0\leq r\leq p(d)} \left\vert {\rm g}_m \left[\hat
{\rm f}^{(r)} \right](w) \right\vert  \leq  \hat {\rm g}(w)\ , $$
where $ \hat {\rm g}\in L^1 \left(\Bbb R^+ \right) $ and $ \left\Vert\hat
{\rm g} \right\Vert_ 1 \leq  {\rm Cst} \left\Vert {\rm g}_m[ {\rm f}]
\right\Vert_ 1. $ \par
\medskip
{\bf Proof\nobreak\ :} \par
The majorizations for $ {\rm g}_m \left[\hat {\rm f} \right](w) $ are given by
the r.h.s. of
the inequalities (3.35') (for $ d\geq 4) $ and (3.38) (for $ d=3), $ from
which
one only retains the last term inside the bracket: in the
corresponding integral,
$ g_m[f](v) $ or $ {\cal G}_m[{\cal F}]( {\rm cosh} \ v) $ must be replaced by
$ {\rm g}_m[ {\rm f}](v). $
The rest of the argument is identical, including the application
of lemma 5 to the treatment of the derivatives $ \hat {\rm f}^{(r)} $ of $
\hat {\rm f}. $ \par
\medskip
{\bf 3.7 Inversion of the Radon-Abel transformation} \par
\smallskip
The inversion of the transformation $ \goth R^{(c)}_d $ (in any
analyticity domain $ D $ of the class $ {\cal D} $ introduced in \S 3.1) will
be
performed by a transformation $ \goth X^{(c)}_d $ whose definition and
properties are given in the following proposition. In the latter,
$ \hat  D $ will denote any domain of the complex hyperbola $ X^{(c)}_1 $
whose
projection $ \underline{D}  = \varpi_ 1 \left(\hat  D \right) = \{ \zeta =
{\rm cos} \ \tau  \in  \Bbb C; $ $ \left.\zeta_ \tau =(\pm i\ {\rm sin} \
\tau , {\rm cos} \ \tau) \in\hat  D \right\} $
is star-shaped with respect to the point $ \zeta =1 $ and contains the
interval $ [-1,+1]. $ \par
\medskip
{\bf Proposition 14\nobreak\ :} \par
For each dimension $ d $ $ (d\geq 3), $ one defines the following
transformation $ \goth X^{(c)}_d $ (or $ {\cal X}^{(c)}_d): $ with every
function $ \hat {\cal F} \left(z_\tau \right) $
(represented in the $ \tau $-plane by $ \hat  f(\tau)  \equiv  \hat {\cal F}
\left(z_\tau \right)) $ analytic in $ \hat  D $ and
satisfying the symmetry relation (3.26), it associates an
$ H^{(c)}_d $-invariant function $ {\cal F }= \goth X^{(c)}_d \left(\hat
{\cal F} \right)={\cal X}^{(c)}_d \left(\hat  f \right), $ defined and
analytic in the
domain $ D = \varpi^{ -1}( \underline{D}) $ of $ X^{(c)}_{d-1}; $ the function
$ f(\theta) , $ even and
analytic in $ \dot {\cal J}(D), $ which represents $ {\cal F}(z) $ $
(z^{(d-1)}= {\rm cos} \ \theta ) $ is expressed in terms of $ \hat  f(\tau) $
by the following formulae: \par
\smallskip
a) $ d $ {\sl even\/}\nobreak\ :
$$ f(\theta)  = \left(- {1 \over 2\pi}  {1 \over {\rm sin} \ \theta}  { {\rm
d} \over {\rm d} \theta} \right)^{{d-2 \over 2}} \left( {\rm e}^{-i \left({d-2
\over 2} \right)\theta}\hat  f(\theta) \right) \eqno (3.46) $$
\par
b) $ d $ {\sl odd\/}\nobreak\ :
$$ f(\theta)  = -2 \left(-{1 \over 2\pi}  {1 \over {\rm sin} \ \theta}  {
{\rm d} \over {\rm d} \theta} \right)^{{d-1 \over 2}} \int^{ }_{ \gamma_
\theta} \left[ {\rm e}^{-i \left({d-2 \over 2} \right)\tau}\hat  f(\tau)
\right][2( {\rm cos} \ \theta - {\rm cos} \ \tau)]^{ -1/2} {\rm sin} \ \tau \
{\rm d} \tau \eqno (3.47) $$
\par
For each $ d, $ the kernel $ \Xi^{( d)} $ of the transformation $ {\cal
X}^{(c)}_d $ is
the set of all functions $ \hat  f $ of the following form:
$$ \hat  f(\tau)  = {\rm e}^{i \left({d-2 \over 2} \right)\tau} P^{(d)}(\tau)
\ , \eqno (3.48) $$
where $ P^{(d)}(\tau) $ denotes a trigonometric polynomial of the following
form:
$$ P^{(d)}(\tau)  = \sum^{ }_{ \left\{ \ell \in \Bbb Z;-{d-2 \over 2}\leq
\ell <0 \right\}} b_{\ell} {\rm cos} \left\{ \left(\ell +{d-2 \over 2}
\right)\tau - \left({d-2 \over 2} \right)\pi \right\} \eqno (3.48^\prime ) $$
\par
\smallskip
{\bf Proof:} \par
In order to justify this definition of $ \goth X^{(c)}_d, $ we must
check that if Eq.(3.26) is satisfied by $ \hat {\cal F} \left(z_\tau \right),
$ or equivalently if
$ \hat  f(\tau) $ satisfies the following condition:
$$ \left(S_d \right)\ \ \ \ \ \ \ \ \ \hat  f(\tau)  = (-1)^d\ {\rm
e}^{i(d-2)\tau}\hat  f(-\tau) $$
(already considered in \S 2.2: see Eq.(2.22)), then Eqs.(3.46) and
(3.47) necessarily define an even analytic function $ f(\theta) $ in the
domain $ \dot {\cal J}(D) = \{ \theta \in  \dot {\Bbb C};\ {\rm cos} \ \theta
\in \underline{D}\} . $ Condition $ \left(S_d \right) $ leads
us to distinguish the two cases: \par
\smallskip
a) $ d $ {\sl even\/} : since $ \hat  f $ is analytic in $ \dot {\cal J}(D), $
it can be seen (in
view of $ \left(S_d \right)) $ to be of the form $ \hat  f(\tau)  = {\rm e}^{i
\left({d-2 \over 2} \right)\tau} b( {\rm cos} \ \tau) , $ with $ b $
analytic in $ \underline{D} . $ It then immediately follows from Eq.(3.46)
that $ f(\theta) $
is even and analytic in $ \dot {\cal J}(D), $ and that the only functions $
\hat  f $
whose image is $ f=0 $ are those whose associated function $ b( {\rm cos} \
\tau) $
is a polynomial of degree $ {d-4 \over 2} $ of $ {\rm cos} \ \tau , $ or
equivalently a
trigonometric polynomial $ P^{(d)}(\tau) $ of the form (3.48'); the form
(3.48) of the elements of $ \Xi^{( d)} $ is thus established for this case.
\par
\smallskip
b) $ d $ {\sl odd\/}\nobreak\ : condition $ \left(S_d \right) $ (together with
the analyticity
of $ \hat  f $ in $ \dot {\cal J}(D)) $ now implies that $ \hat  f $ is of the
form $ \hat  f(\tau)  = {\rm e}^{i \left({d-2 \over 2} \right)\tau}( 1- {\rm
cos} \ \tau)^{ 1/2}b( {\rm cos} \ \tau) , $
with $ b $ analytic in $ \underline{D} {\bf .} $ By using the same
parametrization of $ \gamma_ \theta $ as
in the proof of proposition 10, namely (with the present
notations) $ {\rm cos} \ \tau  = 1+\lambda( {\rm cos} \ \theta -1), $ with $
0\leq \lambda \leq 1, $ we can then rewrite
Eq.(3.47) as follows:
$$ f(\theta)  = i \sqrt{ 2} \left({1 \over 2\pi}  { {\rm d} \over {\rm d\ cos}
\ \theta} \right)^{{d-1 \over 2}} \left[( {\rm cos} \ \theta -1)\times \int^
1_0b[1+\lambda( {\rm cos} \ \theta -1)]\lambda^{ 1/2}(1-\lambda)^{ -1/2} {\rm
d} \lambda \right] \eqno (3.49) $$
which clearly exhibits that $ f(\theta) $ is even and analytic in $ \dot
{\cal J}(D). $
Eq.(3.49) also allows one to characterize the kernel of $ {\cal X}^{(c)}_d: $
in
fact, the relation $ f(\theta)  = 0 $ holds if and only if the integral $
\int^ 1_0b[1+\lambda( {\rm cos} \ \theta -1)]\lambda^{ 1/2}(1-\lambda)^{ -1/2}
{\rm d} \lambda $
is equal to a polynomial of degree $ {d-5 \over 2} $ of $ ( {\rm cos} \
\theta -1); $ but this
condition can be seen (e.g. by plugging the Taylor expansion of $ b( {\rm cos}
\ \tau) $
at $ {\rm cos} \ \tau =1 $ in the latter integral) to be equivalent to the
condition that $ b( {\rm cos} \ \tau) $ itself is a polynomial of degree $
{d-5 \over 2} $ of $ ( {\rm cos} \ \tau -1); $
correspondingly, the function $ (1- {\rm cos} \ \tau)^{ 1/2}b( {\rm cos} \
\tau) $ can then be
written under the trigonometric form (3.48') in such a case; the
form (3.48) of all the elements of $ \Xi^{( d)} $ is therefore
established. \par
\medskip
{\bf Remark :} \par
By restricting formulae (3.46) and (3.47) to real values
of the parameters $ (\theta =u, $ $ \tau =t) $ and by replacing the function $
\hat  f(\tau) $
and $ f(\theta) $ respectively by $ \hat {\bf f}(t) $ and $ {\bf f}(u) $ in
these formulae, we
define similarly a transformation $ \underline{{\goth X}}^{(c)}_d: $ with
every
function $ \hat {\bf F} \left(z_t \right)\equiv\hat {\bf f}(t) $ on $ S_1 $
satisfying the condition (3.30) (i.e.
the condition $ \left( {\bf S}_d \right) $ given in (2.21)), it associates an
$ {\bf H} $-invariant
function $ {\bf F}  = \underline{{\goth X}}_d \left(\hat {\bf F} \right) $ on
$ S_{d-1} $ such that $ {\bf F}(z) = {\bf f}(u) $ $ ( {\rm cos} \
u=z^{(d-1)}). $ \par
\smallskip
We shall now prove: \par
\smallskip
{\bf Proposition 15\nobreak\ :} \par
\smallskip
a) Let $ \hat {\cal F}_d = \goth R^{(c)}_d({\cal F}) $ be the Radon-Abel
transform of
an $ H^{(c)} $-invariant function $ {\cal F}, $ analytic in a domain $ D $ $
(D\in{\cal D}) $ of $ X^{(c)}_{d-1}. $
Then one has:
$$ {\cal F }= \goth X^{(c)}_d \left(\hat {\cal F}_d \right) $$
\par
b) Conversely, if $ {\cal F }= \goth X^{(c)}_d \left(\hat {\cal F} \right), $
$ \hat {\cal F} $ can be expressed
in general as follows: $ \hat {\cal F}=\goth R^{(c)}_d({\cal F})+\hat {\cal
F}^{(0)}, $ where $ \hat {\cal F}^{(0)} \left(z_\tau \right)\equiv\hat
f^{(0)}(\tau) $ is
such that $ \hat  f^{(0)}\in \Xi^{( d)}, $ $ \Xi^{( d)} $ being characterized
in proposition 14. \par
The relation $ \hat {\cal F}=\goth R^{(c)}_d({\cal F}) $ holds if and only if
$ \hat {\cal F} \left(z_\tau \right)\equiv\hat  f(\tau) $ is
of the form (3.26'), namely if
$ \hat  f(\tau) = {\rm e}^{i \left({d-2 \over 2} \right)\tau} \left( {\rm sin}
 {\tau \over 2} \right)^{d-2}b( {\rm cos} \ \tau) , $ with $ b( {\rm cos} \
\tau) $ analytic in $ \underline{D} . $ \par
\smallskip
c) Let $ {\bf F} $ be any $ {\bf H} $-invariant function on the sphere $
S_{d-1} $
and let $ \hat {\bf F}_d $ be the function on $ S_1 $ $ \left(\approx  \hat  S
\right) $ defined by $ \hat {\bf F}_d = \underline{{\goth R}}_d( {\bf F}) $
(see proposition 11); then one has:
$$ {\bf F}  = \underline{{\goth X}}_d \left(\hat {\bf F}_d \right) $$
\par
\smallskip
{\bf Proof\nobreak\ :} \par
For proving a), we will show that if $ \hat  f_d(\tau) $ is given in
terms of $ f(\theta) $ by formula (3.25) (for $ \theta $ and $ \tau $ varying
in $ \dot {\cal J}(D)), $
then $ f $ can be reobtained in terms of $ \hat  f_d $ via Eqs.(3.46) (for $ d
$
even) and (3.47) (for $ d $ odd). We can rewrite Eq.(3.25) in the
following form by putting $ {\rm cos} \ \tau  = 1+\rho \ {\rm e}^{i\beta} $ $
(\rho \geq 0, $ $ -\pi \leq \beta \leq \pi ) $ and $ \gamma_ \tau = \left\{
\tau^{ \prime} ; \right. $
$ {\rm cos} \ \tau^{ \prime} =1+\rho^{ \prime} {\rm e}^{i\beta} ; $ $
\left.0\leq \rho^{ \prime} \leq \rho \right\} : $
$$ \hat  f_d(\tau)  = \omega_{ d-2} {\rm e}^{i \left({d-2 \over 2}
\right)(\tau +\beta)} \int^ \rho_ 0 \underline{f} \left(1+\rho^{ \prime} {\rm
e}^{i\beta} \right) \left[2 \left(\rho -\rho^{ \prime} \right) \right]^{{d-4
\over 2}} {\rm d} \rho^{ \prime} \eqno (3.50) $$
\par
We now notice that the transformation $ \underline{f} \longrightarrow  \hat
f_d $ expressed by Eq.(3.50) is
conveniently represented in terms of a Riemann-Liouville integral operator
(see e.g. $\lbrack$E-1$\rbrack$
p.181-182). \par
For $ \alpha  > 0, $ let
$$ \left[I_\alpha \varphi \right](\rho)  = {1 \over \Gamma( \alpha)}  \int^
\rho_ 0\varphi \left(\rho^{ \prime} \right) \left(\rho -\rho^{ \prime}
\right)^{\alpha -1} {\rm d} \rho^{ \prime} \ ; \eqno (3.51) $$
then Eq.(3.50) can be rewritten as follows
$$ \hat  f_d(\tau)  = (2\pi)^{{ d-2 \over 2}} {\rm e}^{i \left({d-2 \over 2}
\right)(\tau +\beta)} \left[I_{{d-2 \over 2}}\varphi \right](\rho) \eqno
(3.52) $$
with:
$$ \varphi \left(\rho^{ \prime} \right) = \underline{f} \left(1+\rho^{
\prime} {\rm e}^{i\beta} \right) = f \left(\tau^{ \prime} \right) \eqno (3.53)
$$
\par
We shall now make use of the following properties of Riemann-Liouville
operators
considered here as acting on functions $ \varphi( \rho) $ analytic on some
interval $ \left[0,\rho_ 0 \right[ $ $ (\rho_ 0=\rho_ 0(\beta) $
being such that $ \left[0,\rho_ 0(\beta) \right[ = \{ \rho ; $ $ \rho  \geq
0, $ $ \left. {\rm cos} \ \tau  = 1+\rho \ {\rm e}^{i\beta} \in \underline{D}
\right\} ) $
$$ \forall \alpha ,\beta ,\ (\alpha ,\beta >0),\ \ \ \ \ I_\alpha  \circ
I_\beta  = I_{\alpha +\beta} \eqno (3.54) $$
and
$$ \forall \alpha \in \Bbb N^\ast ,\ \ \ \ \varphi( \rho)  = \left({ {\rm d}
\over {\rm d} \rho} \right)^\alpha \left[I_\alpha \varphi \right](\rho) \eqno
(3.55) $$
\par
1) For $ d $ even, let us apply formula (3.55) with $ \alpha  = {d-2 \over 2}
$ and $ \varphi $ given by
Eq.(3.53); it yields (in view of Eq.(3.52)):
$$ f(\tau)  \equiv  \underline{f} \left(1+\rho \ {\rm e}^{i\beta} \right) =
\left({1 \over 2\pi} \right)^{{d-2 \over 2}} \left[{ {\rm d} \over {\rm d\
cos} \ \tau} \right]^{{d-2 \over 2}} \left[ {\rm e}^{-i \left({d-2 \over 2}
\right)\tau}\hat  f_d(\tau) \right] $$
which coincides with Eq.(3.46). \par
2) For $ d $ odd, we shall apply the following identity which immediately
results from
Eqs.(3.51), (3.54) and (3.55):
$$ \left({ {\rm d} \over {\rm d} \rho} \right)^{{d-1 \over 2}} \left[I_{1/2}
\left(I_{{d-2 \over 2}}\varphi \right) \right] = {1 \over \sqrt{ \pi}}
\left({ {\rm d} \over {\rm d} \rho} \right)^{{d-1 \over 2}} \int^ \rho_ 0
\left(I_{{d-2 \over 2}}\varphi \right) \left(\rho^{ \prime} \right)
\left(\rho -\rho^{ \prime} \right)^{-1/2} {\rm d} \rho^{ \prime}  = \varphi(
\rho) \eqno (3.56) $$
\par
By reexpressing $ I_{{d-2 \over 2}}\varphi $ in terms of $ \hat  f_d(\tau) $
via Eq.(3.52) and rewriting the integral
in Eq.(3.56) as a path integral on $ \gamma_ \tau , $ we then obtain (in view
of Eq.(3.53), and after
cancellation of $ \beta $-dependent exponents):
$$ f(\tau)  = - {1 \over \sqrt{ \pi}}  \left({1 \over 2\pi} \right)^{{d-2
\over 2}} \left({ {\rm d} \over {\rm d\ cos} \ \tau} \right)^{{d-1 \over 2}}
\int^{ }_{ \gamma_ \tau} \left[ {\rm e}^{-i \left({d-2 \over 2} \right)\tau^{
\prime}}\hat  f_d \left(\tau^{ \prime} \right) \right] \left( {\rm cos} \
\tau - {\rm cos} \ \tau^{ \prime} \right)^{-{1 \over 2}} {\rm sin} \ \tau^{
\prime} \ {\rm d} \tau^{ \prime} \ , $$
which coincides with Eq.(3.47). \par
We notice that, by restricting the previous computations to the special case
when
$ \tau $ is real (i.e. $ \beta  = \pm \pi ), $ one obtains the proof of c) as
a by-product. \par
The first part of b) readily follows from a), since $ \goth X^{(c)}_d
\left(\hat {\cal F}-\goth R^{(c)}_d({\cal F}) \right)=0. $ Let
us now show the second part of b): in view of proposition 10, $ \hat {\cal F}
$ is of the form (3.26') if
and only if $ \hat {\cal F}^{(0)}=\hat {\cal F}-\goth R^{(c)}_d({\cal F}) $ is
itself of the form (3.26'), i.e.
$ \hat {\cal F}^{(0)} \left(z_\theta \right) = {\rm e}^{i \left({d-2 \over 2}
\right)\theta} \ell_ 0(\theta) , $ where:
$$ \ell_ 0(\theta)  = (1- {\rm cos} \ \theta)^{{ d-2 \over 2}}a_0( {\rm cos} \
\theta) \ , \eqno (3.57) $$
with $ a_0(=b-a) $ analytic in $ \underline{D} . $ But we will now check that,
since $ \hat {\cal F}^{(0)}\in \Xi^{( d)}, $ this is equivalent to saying that
$ \hat {\cal F}^{(0)}=0; $
in fact, as a result of
proposition 14 (Eq.(3.48)), it follows that $ \ell_ 0(\theta) $ is of the form
$ P^{(d)}(\theta) $ given by Eq.(3.48'). When $ d $
is even (resp. odd) the function $ \ell_ 0(\theta) $ must then be equal to a
polynomial of degree $ {d-4 \over 2} $ of $ {\rm cos} \ \theta $ (resp. to
$ (1- {\rm cos} \ \theta)^{ 1/2} $ times a polynomial of degree $ {d-5 \over
2} $ of $ {\rm cos} \ \theta ); $ in both cases, the only function $ \ell_
0(\theta) $ which also
satisfies Eq.(3.57) is $ \ell_ 0=0, $ which ends the proof of proposition
15. \par
In order to apply the transformation $ \goth X^{(c)}_d $ to analytic functions
$ \hat {\cal F} \left(z_\tau \right)\equiv\hat  f_d(\tau) $ in appropriate
subspaces $ \underline{{\cal O}}^{(d)}_{m,p} \left(\dot {\cal J}^{(\mu)}
\right), $ we now need to establish the following properties. \par
\medskip
{\bf Lemma 6} : Identities of the following type hold, in which all functions
$ \ell_ r, $ $ \ell^{ \prime}_ r, $ $ \ell^{ \prime\prime}_ r $ are {\sl
uniformly bounded\/}
in any region of the form $ \{ \theta \in  \dot {\Bbb C}; $ $ \left. {\rm Im}
\ \theta  \geq  v_0 \right\} $ (with $ v_0 > 0): $ \par
\smallskip
a) for $ d $ even
$$ \left({ {\rm d} \over {\rm d\ cos} \ \theta} \right)^{{d-2 \over 2}}[b(
{\rm cos} \ \theta)]  = \sum^{{ d-2 \over 2}}_{r=0}\ell_ r(\theta)  \hat
f^{(r)}_d(\theta) \eqno (3.58) $$
where
$$ b( {\rm cos} \ \theta)  = {\rm e}^{-i \left({d-2 \over 2}
\right)\theta}\hat  f_d(\theta) \eqno (3.59) $$
\par
\smallskip
b) for $ d $ odd,
$$ \eqalignno{ \left({ {\rm d} \over {\rm d\ cos} \ \theta} \right)^{{d-1
\over 2}}[b( {\rm cos} \ \theta)] &  = {\rm e}^{i\theta} \sum^{{ d-1 \over
2}}_{r=0}\ell^{ \prime}_ r(\theta)  \hat  f^{(r)}_d(\theta) \ , &
(3.58^\prime ) \cr \left({ {\rm d} \over {\rm d\ cos} \ \theta} \right)^{{d-3
\over 2}}[b( {\rm cos} \ \theta)] &  = \sum^{{ d-3 \over 2}}_{r=0}\ell^{
\prime\prime}_ r(\theta)  \hat  f^{(r)}_d(\theta) \ , & (3.58^{\prime\prime} )
\cr} $$
where
$$ b( {\rm cos} \ \theta)  = {1 \over {\rm sin}  {\theta \over 2}} {\rm e}^{-i
\left({d-2 \over 2} \right)\theta}\hat  f_d(\theta) \eqno (3.59^\prime ) $$
\par
\smallskip
{\bf Proof} : \par
a) In view of Eq.(3.59), the $\ell$.h.s. of Eq.(3.58) is a sum of terms of the
form
$$ \left({1 \over {\rm sin} \ \theta} \right)^{q_0} \left[ \prod^ m_{i=1}
\left({ {\rm d} \over {\rm d} \theta} \right)^{r_i} \left({1 \over {\rm sin} \
\theta} \right)^{q_i} \right] \left({ {\rm d} \over {\rm d} \theta}
\right)^{r_0} \left[ \left( {\rm e}^{-i \left({d-2 \over 2} \right)\theta}
\right)\hat  f^{(r)}_d \right], $$
such that $ q_0+ \sum^ m_{i=1}q_i = {d-2 \over 2}, $ $ r_i>0 $ $ (1\leq i\leq
m), $  $ r_0\geq 0, $ and $ 0\leq  r \leq  {d-2 \over 2}. $ Each factor $
\left({ {\rm d} \over {\rm d} \theta} \right)^{r_i} \left({1 \over {\rm sin} \
\theta} \right)^{q_i} $
is a homogeneous rational function of $ ( {\rm cos} \ \theta , $ $ {\rm sin} \
\theta ) $ with degree $ -q_i; $ therefore each term of the previous
form can be written as follows:
$$ H^{(r)}_{\alpha ,\beta}  = {\rm e}^{-i \left({d-2 \over 2} \right)\theta}
\times{ P_\alpha \over Q_\beta}  \hat  f^{(r)}_d, \eqno (3.60) $$
where $ Q_\beta =( {\rm sin} \ \theta)^ \beta $ and $ P_\alpha $ is a
homogeneous polynomial of degree $ \alpha $ with respect to $ ( {\rm cos} \
\theta , $ $ {\rm sin} \ \theta ), $ such
that $ \beta -\alpha  = {d-2 \over 2}. $ We can thus write:
$$ {P_\alpha \over Q_\beta}  = {\rm e}^{i\theta( \beta -\alpha)} \times{
P^{(1)}_\alpha \left( {\rm e}^{i\theta} \right) \over Q^{(1)}_\beta \left(
{\rm e}^{i\theta} \right)}\ , $$
where $ Q^{(1)}_\beta( z)= \left(1-z^2 \right)^\beta $ and $ P^{(1)}_\alpha(
z) $ is a certain polynomial; it follows that
each term $ H^{(r)}_{\alpha \beta}( \theta) $ is of the form $ \ell^{(
r)}_{\alpha \beta}( \theta)\hat  f^{(r)}_d(\theta) , $ where the
function $ \ell^{( r)}_{\alpha \beta}( \theta)  = {P^{(1)}_\alpha \left( {\rm
e}^{i\theta} \right) \over Q^{(1)}_\beta \left( {\rm e}^{i\theta} \right)} $
is uniformly bounded for $ {\rm Im} \ \theta  \geq  v_0 $ $ \left(v_0>0
\right); $ we have thus justified the expressions at
the r.h.s. of Eq.(3.58). \par
\smallskip
b) In view of Eq.(3.59'), the $\ell$.h.s. of Eq.(3.58') (resp.(3.58")) are
sums of terms which can be
analysed as above, up to the following differences: $ \beta -\alpha  = {d-1
\over 2} $ (resp. $ {d-3 \over 2}) $ and additional factors of the
form $ {\rm e}^{i\theta /2}\times P^{\prime} \left( {\rm e}^{i\theta} \right)/
\left(1- {\rm e}^{i\theta} \right)^{\beta^{ \prime}} $ must be inserted in the
expressions (3.60) of the functions $ H^{(r)}_{\alpha ,\beta} ; $ thus, one
still obtains bounded functions of $ \theta , $ of the form $ {P^{(1)}_\alpha
\over Q^{(1)}_\beta}  \times  {P^{\prime} \over \left(1- {\rm e}^{i\theta}
\right)^{\beta^{ \prime}}} $ whose sums are the multipliers $ \ell^{ \prime}_
r $
(resp. $ \ell^{ \prime\prime}_ r) $ of the functions $ \hat  f^{(r)}_d, $ as
expected in formulae (3.58') (resp. (3.58")). \par
We are now in a position to establish the following statement concerning the
inverse Radon-Abel
transformation in cut-domains $ D_\mu $ of $ X^{(c)}_{d-1}. $ \par
\medskip
{\bf Proposition 16} : \par
Let $ {\bf F}(z) $ be an $ {\bf H}^{(c)} $-invariant function on $ S_{d-1} $
and let $ \hat {\bf F} = \underline{{\goth R}}^{(c)}_d( {\bf F}) $
be its Radon-Abel transform, defined (according to Eqs.(3.28), (3.29) and
proposition 11) on the circle $ \hat  S \left(\equiv S_1 \right) $ of $
S_{d-1}. $ \par
The following properties hold: \par
i) If $ \hat {\bf F} $ admits an analytic continuation $ \hat {\cal F} $ in
the domain $ \hat  D_\mu = \left\{ z_\tau ; \tau \in\dot {\cal J}^{(\mu)}
\right\} , $ $ (\mu \geq 0), $ then $ {\bf F} $
admits an $ (H^{(c)} $-invariant) analytic continuation $ {\cal F} $ in the
corresponding domain $ D_\mu $ of $ X^{(c)}_{d-1}, $ and one has $ \hat {\cal
F} = \goth R^{(c)}_d({\cal F}). $ \par
ii) The inverse transformation $ \hat {\cal F} \longrightarrow  {\cal F}=
\goth X^{(c)}_d \left(\hat {\cal F} \right) = \left(\goth R^{(c)}_d
\right)^{-1} \left(\hat {\cal F} \right) $ is given in
$ \hat  D_\mu $ respectively by Eq.(3.46) for the case $ \gq\gq d $
even\rq\rq , and by Eq.(3.47) for the case $ \gq\gq d $ odd\rq\rq . \par
\smallskip
Correspondingly, the following inversion formulae express the discontinuity $
F(z) = \underline{{\rm f}} \left(z^{(d-1)} \right) $
of $ {\cal F} $ (on $ X_{d-1}) $ in terms of the discontinuity $ \hat  F =
\goth R_d(F) $ of $ \hat {\cal F} $ (on $ \hat  X \approx  X_1), $ namely $
\hat  F \left(z_{iv} \right) \equiv  \hat {\rm f}(v): $ \par
\smallskip
a) $ d $ {\sl even\nobreak\ :\/}
$$ \underline{{\rm f}}( {\rm cosh} \ v) \equiv  {\rm f}(v) = \left({1 \over
2\pi \ {\rm sinh} \ v} { {\rm d} \over {\rm d} v} \right)^{{d-2 \over 2}}
\left[ {\rm e}^{{d-2 \over 2}v}\hat {\rm f}(v) \right] \eqno (3.64) $$
\par
b{\bf )} $ d ${\sl\ odd\/}\nobreak\ :{\bf\ }
$$ \underline{{\rm f}}( {\rm cosh} \ v) \equiv  {\rm f}(v) = 2 \left({1 \over
2\pi \ {\rm sinh} \ v} { {\rm d} \over {\rm d} v} \right)^{{d-1 \over 2}}
\int^ v_0 \left[ {\rm e}^{{d-2 \over 2}w}\hat {\rm f}(w) \right][2( {\rm cosh}
\ v- {\rm cosh} \ w)]^{-1/2} {\rm sinh} \ w\ {\rm d} w \eqno (3.65) $$
\par
iii) The transformation $ {\cal X}^{(c)}_d $ defines a continuous mapping from
each subspace \par
\noindent$ \underline{{\cal O}}^{(d)}_{m,\hat  p(d)} \left(\dot {\cal
J}^{(\mu)} \right) $ (resp.
$ \underline{{\cal O}}^{(d)\ast}_{ m,\hat  p(d)} \left(\dot {\cal J}^{(\mu)}
\right)), $ with subindex $ \hat  p(d) $ specified below, into the
corresponding subspace $ \left[{\cal V}_{\infty} \right]^m_\mu
\left(X^{(c)}_{d-1} \right) $ (resp. $ \left[{\cal V}^\ast_{ \infty}
\right]^m_\mu \left(X^{(c)}_{d-1} \right)): $ \par
\smallskip
a) for $ d $ even, $ m \in  \Bbb R $ and $ \hat  p(d) = {d-2 \over 2} $ \par
b) for $ d $ odd, $ m > - {d \over 2} $ and $ \hat  p(d) = {d-1 \over 2}\ . $
\par
\medskip
{\bf Proof\nobreak\ :} \par
In view of proposition 11, the function $ \hat {\bf f}(t) \equiv  \hat {\bf F}
\left(z_t \right) $ is of the form given by Eq.(3.30'), and
therefore its analytic continuation $ \hat {\cal F} $ is necessarily of the
form $ \hat  f(\tau)  \equiv  \hat {\cal F} \left(z_\tau \right) = {\rm e}^{i
\left({d-2 \over 2} \right)\tau} \left( {\rm sin}  {\tau \over 2}
\right)^{d-2}a( {\rm cos} \ \tau) $
(i.e. (3.26')) with $ a( {\rm cos} \ \tau) $ analytic in the cut-place $
\underline{D}_\mu . $ We can then associate with $ \hat {\cal F} $ the $
H^{(c)} $-invariant
function $ {\cal F }= \goth X^{(c)}_d \left(\hat {\cal F} \right) $ which, in
view of proposition 13, is analytic in the corresponding domain $ D_\mu $ of $
X^{(c)}_{d-1}. $
The restriction of $ {\cal F} $ to $ S_{d-1} $ coincides with the given
function $ {\bf F} , $ since $ {\cal F}_{ \left\vert S_{d-1} \right.} =
\underline{{\goth X}}_d \left(\hat {\bf F} \right) $ which, by
proposition 15c), is equal to $ {\bf F} . $ Moreover, the conditions of
proposition 15b) being satisfied by $ \hat {\cal F}, $ it
follows that $ \hat {\cal F} = \goth R^{(c)}_d({\cal F}) $ which completes
property i). In view of proposition 14, the inversion
formulae are then given by the expressions (3.46) and (3.47) of the
transformation $ \goth X^{(c)}_d; $ these expressions
also allow one to compute the boundary values $ f_{\pm}( v) $ of $ f $ on the
semi-axis $ \{ \theta =iv; $ $ v\geq 0\} $ in terms of the
corresponding boundary values $ \hat  f_{\pm}( w) $ of $ \hat  f $ (with $
\underline{\gamma}_{ iv} = \{ \tau =iw; $ $ 0\leq w\leq v\} ). $ \par
By subtraction, one thus obtains the expressions (3.64) and (3.65) (since $
{\rm f}  = i \left(f_+-f_- \right) $ and $ \hat {\rm f}=i \left(\hat
f_+-\hat  f_- \right)); $
the latter perform the inversion of the transformation $ \goth R_d, $ which
completes property ii). \par
For the proof of property iii), we distinguish two cases \par
\smallskip
a) $ d $ {\sl even\/}: \par
In view of lemma 6 (formulae (3.58), (3.59)), we can rewrite Eq.(3.46) as
follows:
$$ f(\theta)  = \left({1 \over 2\pi} \right)^{{d-2 \over 2}} \left({ {\rm d}
\over {\rm d\ cos} \ \theta} \right)^{{d-2 \over 2}}[b( {\rm cos} \ \theta)]
= \left({1 \over 2\pi} \right)^{{d-2 \over 2}} \sum^{{ d-2 \over
2}}_{r=0}\ell_ r(\theta)  \hat  f^{(r)}_d(\theta) \ , \eqno (3.66) $$
which implies (in view of Eq.(2.1)):
$$ \forall v,\ \ {\rm with} \ \ v\geq v_0,\ \ g_m[f](v) \leq  {\rm Cst} \
\dlowlim{ {\rm max}}{0\leq r\leq{ d-2 \over 2}} g_m \left[\hat  f^{(r)}_d
\right](v) \eqno (3.67) $$
\par
On the other hand, since the function $ b( {\rm cos} \ \theta) $ is analytic
in $ \underline{D}_\mu $ (see the proof of proposition 13)
and in particular (provided $ v_0 < \mu ) $ in the ellipse $ E_{ {\rm cosh} \
v_0}, $ the Cauchy inequalities applied to the
derivatives of $ b( {\rm cos} \ \theta) $ yield (by making use of Eqs.(3.66)
and (3.59)):
$$ \forall v,\ {\rm with} \ 0\leq v\leq v_0,\ g_m[f](v)\leq {\rm Cst} \
\dlowlim{ {\rm sup}}{0\leq u\leq 2\pi} \left\vert b \left( {\rm cos}
\left(u+iv_0 \right) \right) \right\vert  \leq  {\rm Cst}^{\prime} g_m
\left[\hat  f_d \right] \left(v_0 \right) \eqno (3.68) $$
\par
In view of Eqs.(2.34) and (2.38), the previous inequalities (3.67) and (3.68)
then imply the
following continuity inequalities:
$$ \matrix{\displaystyle \left\Vert g_m[f] \right\Vert_{ \infty} &
\displaystyle \leq  {\rm Cst} \mid \mid \mid\hat  f_d\mid \mid \mid^{(
\infty)}_{ m,{d-2 \over 2}} \hfill \cr\displaystyle \left\Vert g^\ast_ m[f]
\right\Vert_{ \infty} & \displaystyle \leq  {\rm Cst} \mid \mid \mid\hat
f_d\mid \mid \mid^{( \infty) \ast}_{ m,{d-2 \over 2}}\ , \hfill \cr} $$
which (via the identifications $ {\cal F}(z) = \underline{f} \left(z^{(d-1)}
\right), $ $ \underline{f}( {\rm cos} \ \theta)  = f(\theta) $ and lemma 4)
are equivalent to the
statement of property iii) (even case). \par
\smallskip
b) $ d $ {\sl odd\/} : we shall make use of the form (3.49) of the
transformation $ \goth X^{(c)}_d $ (see the proof of
proposition 13), which we can obviously rewrite as follows:
$$ \matrix{\displaystyle f(\theta) & \displaystyle = {\rm Cst} \left\{( {\rm
cos} \ \theta -1) \int^ 1_0 \left({ {\rm d} \over {\rm d\ cos} \ \theta}
\right)^{{d-1 \over 2}}[b(1+\lambda( {\rm cos} \ \theta -1))] \lambda^{
1/2}(1-\lambda)^{ -1/2} {\rm d} \lambda \right. \hfill \cr\displaystyle  &
\displaystyle \left.+ {d-1 \over 2} \int^ 1_0 \left({ {\rm d} \over {\rm d\
cos} \ \theta} \right)^{{d-3 \over 2}}[b(1+\lambda( {\rm cos} \ \theta
-1))]\lambda^{ 1/2}(1-\lambda)^{ -1/2} {\rm d} \lambda \right\} \hfill \cr}
\eqno (3.69) $$
\par
We call $ \tau  = \tau( u,v,\lambda) , $ and in particular $ w = {\rm Im} \
\tau  = w(u,v,\lambda) , $ the complex angle defined by the
relation $ {\rm cos} \ \tau  = 1+\lambda( {\rm cos} \ \theta -1), $ with $
\theta =u+iv. $ By applying lemma 6b and Eqs.(2.1), (2.34), we can
then write:
$$ \matrix{\displaystyle \left\vert \left({ {\rm d} \over {\rm d\ cos} \
\theta} \right)^{{d-1 \over 2}}[b(1+\lambda( {\rm cos} \ \theta -1))]
\right\vert & \displaystyle = \lambda^{{ d-1 \over 2}} \left\vert \left({
{\rm d} \over {\rm d\ cos} \ \tau} \right)^{{d-1 \over 2}}[b( {\rm cos} \
\tau)] \right\vert \hfill \cr\displaystyle  & \displaystyle \leq  \lambda^{{
d-1 \over 2}} {\rm e}^{(m-1)w}\mid \mid \mid\hat  f_d\mid \mid \mid^{(
\infty)}_{ m,{d-1 \over 2}} \hfill \cr} $$
and
$$ \matrix{\displaystyle \left\vert \left({ {\rm d} \over {\rm d\ cos} \
\theta} \right)^{{d-3 \over 2}}[b(1+\lambda( {\rm cos} \ \theta -1))]
\right\vert & \displaystyle = \lambda^{{ d-3 \over 2}} \left\vert \left({
{\rm d} \over {\rm d\ cos} \ \tau} \right)^{{d-3 \over 2}}[b( {\rm cos} \
\tau)] \right\vert \hfill \cr\displaystyle  & \displaystyle \leq  \lambda^{{
d-3 \over 2}} {\rm e}^{mw}\mid \mid \mid\hat  f_d\mid \mid \mid^{( \infty)}_{
m,{d-1 \over 2}}\ \ ; \hfill \cr} $$
in the r.h.s. of these inequalities, $ w $ stands for $ w(u,v,\lambda) . $ By
taking these majorizations into account in
Eq.(3.69), we then obtain (in view of Eq.(2.1)):
$$ \matrix{\displaystyle g_m[f](v) & \displaystyle \leq  {\rm Cst} \mid \mid
\mid\hat  f_d\mid \mid \mid^{( \infty)}_{ m,{d-1 \over 2}} \left[ {\rm
e}^{(1-m)v} \int^ 1_0\lambda^{{ d \over 2}}(1-\lambda)^{ -1/2} \dlowlim{ {\rm
sup}}{u} \left[ {\rm e}^{(m-1)w(u,v,\lambda)} \right] {\rm d} \lambda \right.
\hfill \cr\displaystyle  & \displaystyle \left.+ {\rm e}^{-mv} \int^
1_0\lambda^{{ d-2 \over 2}}(1-\lambda)^{ -1/2} \dlowlim{ {\rm sup}}{u} \left[
{\rm e}^{mw(u,v,\lambda)} \right] {\rm d} \lambda \right] \hfill \cr} \eqno
(3.70) $$
\par
In order to show that these majorizations are meaningful, we now need to check
that both
integrals at the r.h.s. of this inequality are finite. We make use of the fact
that: \par
\smallskip
\noindent \item {\nobreak\ \nobreak\ i)} $ \dlowlim{ {\rm sup}}{u}
w(u,v,\lambda)  = w_{ {\rm max}} $ is defined by $ {\rm cosh\ } w_{ {\rm max}}
= 1+\lambda( {\rm cosh} \ v-1) $ \par
\noindent \item {\nobreak\ ii)} $ \dlowlim{ {\rm sup}}{u} w(u,v,\lambda)  =
w_{ {\rm min}} $ is defined by $ {\rm cosh} \ w_{ {\rm min}} = -1+\lambda(
{\rm cosh} \ v+1) $ \par
\noindent which implies:
$$ \forall \lambda ,\ \ \ 0\leq \lambda \leq 1,\ \ {\rm e}^{w_{ {\rm max}}}
\leq  4\ {\rm e}^v\ \ ,\ \ {\rm and} \ \ \ {\rm e}^{w_{ {\rm min}}} \geq
{\lambda {\rm e}^v \over 2} - 1 \eqno (3.71) $$
\par
{}From (3.71), it follows that if $ m\geq 1 $ (resp. $ m\geq 0) $ the first
(resp. second) term at the r.h.s. of
(3.70) is uniformly bounded (in $ v) $ by $ {\rm Cst} \ \int^
1_0(1-\lambda)^{ -1/2} {\rm d} \lambda  < \infty . $ \par
Let us now study this first term in the case $ m<1; $ we can then write (in
view of (3.71))
$$ \matrix{\displaystyle \dlowlim{ {\rm sup}}{v} \left[ {\rm e}^{(1-m)v}
\int^ 1_0\lambda^{{ d \over 2}}(1-\lambda)^{ -1/2} \dlowlim{ {\rm sup}}{u}
\left[ {\rm e}^{(m-1)w(u,v,\lambda)} \right] {\rm d} \lambda \right] \leq
\hfill \cr\displaystyle \dlowlim{ {\rm sup}}{v} \left[ {\rm e}^{(1-m)v}
\int^{ 4 {\rm e}^{-v}}_0\lambda^{{ d \over 2}} {\rm d} \lambda  + {\rm
e}^{(1-m)v} \int^ 1_{4 {\rm e}^{-v}} \left({\lambda \ {\rm e}^v \over 4}
\right)^{m-1} \lambda^{{ d \over 2}}(1-\lambda)^{ -1/2} {\rm d} \lambda
\right] \hfill \cr} \eqno (3.72) $$
\par
We then easily see that the quantity at the r.h.s. of the inequality (3.72) is
finite if and only
if $ m>-{d \over 2}. $ \par
A similar analysis (implying the same condition on $ m) $ could be done for
the second term at the
r.h.s. of (3.70). Therefore we have established that, provided $ m > - {d
\over 2}, $ the following continuity
inequality holds:
$$ \left\Vert g_m[f] \right\Vert_{ \infty}  \leq  {\rm Cst} \ \mid \mid
\mid\hat  f_d\mid \mid \mid^{( \infty)}_{ m,{d-1 \over 2}} \eqno (3.73) $$
By making use of Eq.(2.38) instead of (2.34), we would be led (by the same
argument) to the following
inequality in place of (3.70):
$$ \matrix{\displaystyle g^\ast_ m[f](v) & \displaystyle \leq  {\rm Cst} \
\mid \mid \mid\hat  f_d\mid \mid \mid^{( \infty) \ast}_{ m,{d-1 \over 2}}
\left[ {\rm e}^{(1-m)v}(1+v) \int^ 1_0\lambda^{{ d \over 2}}(1-\lambda)^{
-1/2} { \dlowlim{ {\rm sup}}{u} \left[ {\rm e}^{(m-1)w(u,v,\lambda)} \right]
\over 1+w[u,v,\lambda]}  {\rm d} \lambda \right. \hfill \cr\displaystyle  &
\displaystyle \left.+ {\rm e}^{-mv}(1+v) \int^ 1_0\lambda^{{ d-2 \over
2}}(1-\lambda)^{ -1/2} { \dlowlim{ {\rm sup}}{u} \left( {\rm
e}^{(m-1)w(u,v,\lambda)} \right) \over 1+w(u,v,\lambda)}  {\rm d} \lambda
\right] \hfill \cr} $$
Then, by establishing similarly the convergence of these two integrals under
the condition $ m>-{d \over 2}, $ one would
derive the following inequality:
$$ \left\Vert g^\ast_ m[f] \right\Vert_{ \infty}  \leq  {\rm Cst} \mid \mid
\mid\hat  f_d\mid \mid \mid^{( \infty) \ast}_{ m,{d-1 \over 2}} \eqno (3.74)
$$
The inequalities (3.73) and (3.74) are equivalent to the statement of property
iii) (odd case). \par
\medskip
{\bf 3.8 A family of elementary perikernels} \par
\smallskip
The transformation $ \goth X^{(c)}_d $ $ (d\geq 3) $ studied in \S 3.7 will
allow
us to introduce a family of elementary perikernels in the
following way. \par
We consider the family of functions $ \left\{\hat  f^{(d)}_\lambda( \tau) ;\
\lambda \in \Bbb C \right\} $
defined by:
$$ \hat  f^{(d)}_\lambda( \tau)  = {\rm e}^{i \left({d-2 \over 2}
\right)(\tau -\pi)} {\rm cos} \left[ \left(\lambda  + {d-2 \over 2}
\right)(\tau -\pi) \right]\ , \eqno (3.75) $$
for $ \tau $ varying in the strip $ \{ \tau ; $ $ 0\leq {\rm Re} \ \tau  \leq
2\pi\} , $ and such that $ \hat  f^{(d)}_\lambda( \tau +2\pi n)=\hat
f^{(d)}_\lambda( \tau) $
$ (\forall n \in  \Bbb Z). $ For each $ \lambda , $ this defines $ \hat
f^{(d)}_\lambda $ as an analytic
function (with continuous boundary values) in the $ 2\pi $-periodic cut-plane
$ \dot {\cal J}^{(0)}, $ and therefore
correspondingly $ \hat {\cal F}^{(d)}_\lambda \left(z_\tau \right) = \hat
f^{(d)}_\lambda( \tau) $ as an analytic function in the
domain $ \hat  D_0 = \left\{ z_\tau ;\ \tau  \in  \dot {\cal J}^{(0)}
\right\} $ of $ X^{(c)}_1. $ \par
Moreover these functions satisfy the following properties
(which can be checked directly from Eq.(3.75) together with the
$ 2\pi $-periodicity condition): \par
i) $ \hat  f^{(d)}_\lambda( \tau) $ satisfies the symmetry condition $ S_d $
(see
Eq.(2.22)), or correspondingly $ \hat {\cal F}^{(d)}_\lambda $ satisfies the
relation (3.26). \par
\smallskip
ii) the discontinuity $ \hat {\rm f}^{(d)}_\lambda( w) = i \left(\hat
f^{(d)}_{\lambda ,+}(w)-\hat  f^{(d)}_{\lambda ,-}(w) \right) $ of $ \hat
f^{(d)}_\lambda $
across the cut $ \{ \tau ;\ {\rm Re} \ \tau  = 0\} $ is given by the following
expressions: \par
\smallskip
a) {\sl for d even\/},
$$ \hat {\rm f}^{(d)}_\lambda( w) = -2\ {\rm e}^{- \left({d-2 \over 2}
\right)w} {\rm sinh} \left(\lambda +{d-2 \over 2} \right)w\cdot {\rm sin} \
\pi \lambda \eqno (3.76) $$
\par
b) {\sl for d odd,\/}
$$ \hat {\rm f}^{(d)}_\lambda( w) = -2\ {\rm e}^{- \left({d-2 \over 2}
\right)w} {\rm cosh} \left(\lambda +{d-2 \over 2} \right)w\cdot {\rm sin} \
\pi \lambda \eqno (3.77) $$
\par
Then, property i) allows us to apply proposition 14 to
all functions of the family $ \left\{\hat {\cal F}^{(d)}_\lambda \left(z_\tau
\right);\ \lambda \in \Bbb C \right\} $ and to obtain
correspondingly a family of $ H^{(c)}_d $-invariant functions $ \left\{{\cal
F}^{(d)}_\lambda =\goth X^{(c)}_d \left(\hat {\cal F}^{(d)}_\lambda \right);\
\lambda \in \Bbb C \right\} , $
which are given (in view of Eqs.(3.46), (3.47) and (3.75)) by the
following formulae, in which we introduce the special notation $ \Psi^{(
d)}_\lambda $
(instead of $ \underline{f}^{(d)}_\lambda ): $
$$ \Psi^{( d)}_\lambda( {\rm cos} \ \theta)  = {\cal F}^{(d)}_\lambda( z)_{
\left\vert z^{(d-1)}= {\rm cos} \ \theta \right.} $$
\par
a) {\sl for d even\/} : $ \forall \theta , $ with $ 0 < {\rm Re} \ \theta  <
2\pi \ : $
$$ \Psi^{( d)}_\lambda( {\rm cos} \ \theta)  = \left({1 \over 2\pi}  {1 \over
{\rm sin} \ \theta}  { {\rm d} \over {\rm d} \theta} \right)^{{d-2 \over 2}}
\left[ {\rm cos} \left(\lambda +{d-2 \over 2} \right)(\theta -\pi) \right]
\eqno (3.78) $$
\par
b) {\sl for d odd\/} : $ \forall \theta , $ with $ 0 < {\rm Re} \ \theta  <
2\pi \ : $
$$ \Psi^{( d)}_\lambda( {\rm cos} \ \theta)  = -2i \left({1 \over 2\pi}  {1
\over {\rm sin} \ \theta}  { {\rm d} \over {\rm d} \theta} \right)^{{d-1
\over 2}} \int^{ }_{ \gamma_ \theta}{ {\rm cos} \left(\lambda +{d-2 \over 2}
\right)(\tau -\pi) \over[ 2( {\rm cos} \ \theta - {\rm cos} \ \tau)]^{ 1/2}}
{\rm sin} \ \tau \ {\rm d} \tau \eqno (3.79) $$
\par
According to proposition 14, the functions $ \Psi^{( d)}_\lambda( {\rm cos} \
\theta) $ thus defined
are, in both cases, analytic in the cut-plane $ \underline{D}_0=\Bbb
C\backslash[ 1,+\infty[ . $ Moreover,
each function $ \hat  f^{(d)}_\lambda( \tau) $ such that $ {\rm Re} \ \lambda
\geq  -{d-2 \over 2} $ can be seen to belong to the
corresponding functional space $ \underline{{\cal O}}^{(d)}_{ {\rm Re} \
\lambda ,\hat  p(d)} \left(\dot {\cal J}^{(0)} \right) $ (see Eq.(2.34)) and
therefore the following property follows from proposition 16\nobreak\ iii).
\par
\medskip
{\bf Proposition 17 :} \par
For every $ \lambda =m+i\nu $ such that $ m \geq  - {d-2 \over 2} $ the
function $ \Psi^{( d)}_\lambda( {\rm cos} \ \theta) $
defines an invariant perikernel on $ X^{(c)}_{d-1}, $ which belongs to the
subspace $ \left[{\cal V}_{\infty} \right]^m_0 \left(X^{(c)}_{d-1} \right). $
\par
\smallskip
In section 4, it will appear, as a by-product of the Fourier-Laplace
transformation of perikernels, that the family $ \left\{ \Psi^{( d)}_\lambda
;\ \lambda \in \Bbb C \right\} $ provides
a {\sl basic set of elementary\/} {\sl perikernels\/}, entering in appropriate
integral
representations of arbitrary invariant perikernels of moderate growth. In the
derivation of such representations, it will be necessary to make use of
asymptotic properties of the family $ \left\{ \Psi^{( d)}_\lambda( {\rm cos} \
\theta) ;\ \lambda \in \Bbb C \right\} $ with respect to $ \lambda $
(when $ \lambda \longrightarrow \infty , $ at any fixed value of $ {\rm cos} \
\theta ); $ we will therefore establish the
following property. \par
\medskip
{\bf Proposition 18 :} \par
For every $ \theta $ such that $ 0 < {\rm Re} \ \theta  < 2\pi , $ the mapping
$ \lambda  \longrightarrow  \Psi^{( d)}_\lambda( {\rm cos} \ \theta) $
defines an entire function of $ \lambda $ which satisfies the following
properties: \par
a)
$$ \Psi^{( d)}_\lambda( {\rm cos} \ \theta)  = \Psi^{( d)}_{-\lambda -d+2}(
{\rm cos} \ \theta) \eqno (3.80) $$
\par
b) $ \forall \lambda  \in  \Bbb C, $ $ \lambda  = m+i\nu , $
$$ \left\vert \Psi^{( d)}_\lambda( {\rm cos} \ \theta) \right\vert  \leq  C_d(
{\rm cos} \ \theta)  {\rm e}^{ {\rm max}(m,-m-d+2)\vert {\rm Im} \
\theta\vert} \times ( {\rm cosh} \ \pi \nu )(1+\vert \lambda\vert)^{\hat
p(d)} \eqno (3.80^\prime ) $$
$ C_d( {\rm cos} \ \theta) $ denoting a suitable function of $ {\rm cos} \
\theta , $ bounded at infinity in $ \Bbb C\backslash[ 1,+\infty[ . $ \par
\medskip
{\bf Proof :} \par
Eq.(3.80) is an immediate consequence of the definition of $ \Psi^{(
d)}_\lambda $ (see
Eqs.(3.78), (3.79)). We shall now prove the majorization (3.80'). \par
For the case $ d $ odd, it will
be convenient to use the parametrization $ {\rm cos} \ \tau -1=\lambda( {\rm
cos} \ \theta -1) $
at the r.h.s. of Eq.(3.79) (as in Eq.(3.49)).
We can then rewrite Eq.(3.79) in the following form:
$$ \Psi^{( d)}_\lambda( \zeta)  = -i \left(-{1 \over 2\pi}  { {\rm d} \over
{\rm d} \zeta} \right)^{{d-1 \over 2}} \left[{(\zeta -1)^{1/2} \over 2} \int^
1_0{ {\rm d} \rho \over( 1-\rho)^{ 1/2}} \left(Z(\rho ,\zeta)^ \alpha {\rm
e}^{-i\pi \alpha} +Z(\rho ,\zeta)^{ -\alpha} {\rm e}^{i\pi \alpha} \right)
\right] \eqno (3.81) $$
where we have put: $ \zeta  = {\rm cos} \ \theta , $ $ z=z(\rho ,\zeta)  =
{\rm cos} \ \tau  = 1+\rho( \zeta -1), $ $ Z=Z(\rho ,\zeta) = {\rm e}^{i\tau}
=z+ \left(z^2-1 \right)^{1/2}, $ $ \alpha =\lambda +{d-2 \over 2}. $ \par
The action of the derivative operator at the r.h.s. of Eq.(3.81) can be
checked to result in the following alternative expression:
$$ \Psi^{( d)}_\lambda( \zeta)  = \int^ 1_0{ {\rm d} \rho \over( 1-\rho)^{
1/2}} \left[A_d(\alpha ,\zeta ,\rho) Z^\alpha {\rm e}^{-i\pi \alpha}
+A_d(-\alpha ,\zeta ,\rho) Z^{-\alpha} {\rm e}^{i\pi \alpha} \right] \eqno
(3.82) $$
where the function $ A_d $ denotes a polynomial of degree $ {d-1 \over 2} $ in
$ \alpha , $ whose all coefficients
are rational functions of $ \sqrt{ \rho} , $ $ \sqrt{ \zeta -1}, $ $ \sqrt{
z(\rho ,\zeta)^ 2-1} $ bounded by $ {\rm Cst}\vert \zeta\vert^{ -{d-2 \over
2}}, $ provided $ \zeta \in \underline{D}_0\backslash] -\infty ,-1]. $
On the other hand, one can use the following majorization:
$$ {\rm max} \left( \left\vert Z^\alpha {\rm e}^{-i\pi \alpha} \right\vert ,
\left\vert Z^{-\alpha} {\rm e}^{i\pi \alpha} \right\vert \right)\leq {\rm Cst}
\ {\rm cosh} \ \pi \nu \ {\rm e}^{\vert {\rm Re\ \alpha}\vert\vert {\rm Im} \
\theta\vert} $$
\par
Since $ {\rm Re} \ \alpha  = m+{d-2 \over 2} $ and $ \vert \zeta\vert  \leq
{\rm e}^{\vert {\rm Im} \ \theta\vert} , $ the previous analysis allows one to
majorize the integral of Eq.(3.82) by the r.h.s. of (3.80), at least for $
\zeta \in \underline{D}_0\backslash] -\infty ,-1]. $
For $ \zeta <-1, $ a slight distortion of the integration path $ \{ \rho \in[
0,1]\} $ (in Eq.(3.82)) is
required, in order to avoid the singularity of the factor $ [z(\rho ,\zeta)
+1]^{-{d-2 \over 2}} $ in the
integrand and to extend therefore the majorization (3.80) in $
\underline{D}_0\backslash\{ -1\} , $ and finally (via
the maximum modulus principle) in $ \underline{D}_0. $ \par
In view of Eq.(3.78), the case $ d $ even is treated in a simpler way: the
integration over $ \lambda $ is essentially replaced by putting $ \lambda =1 $
in the previous argument,
whose all other points are similar. \par
\medskip
{\bf Connection with the Legendre functions in dimension }$ d $ \par
We shall introduce the Legendre functions in dimension $ d $ by their
Abel-type
integral representations (see e.g. $\lbrack$E-2$\rbrack$ and
$\lbrack$Fa-1$\rbrack$). \par
\smallskip
{\sl The first-kind Legendre functions\/} : We put:
$$ P^{(d)}_\lambda( {\rm cos} \ \theta)  = 2{\omega_{ d-2} \over \omega_{
d-1}} ( {\rm sin} \ \theta)^{ -(d-3)} \int^ \theta_ 0 {\rm cos} \left[
\left(\lambda +{d-2 \over 2} \right)\tau \right][2( {\rm cos} \ \tau - {\rm
cos} \ \theta)]^{{ d-4 \over 2}} {\rm d} \tau \eqno (3.83) $$
For $ d\geq 3, $ this formula defines $ P^{(d)}_\lambda( {\rm cos} \ \theta) $
as a ramified analytic function, whose
fundamental sheet is specified by the following conditions: $ {\rm cos} \
\theta \in \Bbb C\backslash] -\infty ,-1], $ with
$ -\pi  < {\rm Re} \ \theta  < -\pi ; $ it also shows that for each $ {\rm
cos} \ \theta $ in this domain, the function $ \lambda \longrightarrow
P^{(d)}_\lambda( {\rm cos} \ \theta) $
is an entire function such that: $ P^{(d)}_\lambda( {\rm cos} \ \theta)  =
P^{(d)}_{-\lambda -d+2}( {\rm cos} \ \theta) . $ \par
Let us now introduce the generalized dimensional factor $ h_d(\lambda) $
(connected for $ \lambda =\ell \in \Bbb N $
with the $ {\rm SO}(d) $-representation theory), which coincides up to a
multiplicative constant
with the function $ h^{(0)}_d(\lambda) $ of the Appendix (cf. Eq.(A.3)):
$$ h_d(\lambda)  = {2 \over( d-2)!} h^{(0)}_d(\lambda)  = {(2\lambda +d-2)
\over( d-2)!} \cdot  {\Gamma( \lambda +d-2) \over \Gamma( \lambda +1)} \eqno
(3.84) $$
which is such that:
$$ h_d(\lambda)  = (-1)^dh_d(-\lambda -d+2) \eqno (3.84^\prime ) $$
\par
By applying the algebraic identities (A.32) and (A.33) (see proposition A.2),
we then obtain the
following expressions for the products $ h_d(\lambda)  P^{(d)}_\lambda , $
after having taken into account the
identities:
$$ \omega_ d = 2 {\pi^{ d/2} \over \Gamma \left({d \over 2} \right)} \ \ (
{\rm for} \ d,\ d-1,\ d-2)\ \ {\rm and} \ \ (d-2)!\omega_ d\omega_{
d-1}=2(2\pi)^{ d-1} $$
\par
a) {\sl d even\/} :
$$ h_d(\lambda) \cdot P^{(d)}_\lambda( {\rm cos} \ \theta)  = {2\omega_ d
\over( 2\pi)^{ d/2}} \left({-1 \over {\rm sin} \ \theta}  { {\rm d} \over
{\rm d} \theta} \right)^{{d-2 \over 2}} \left[ {\rm cos} \left(\lambda +{d-2
\over 2} \right)\theta \right] \eqno (3.85) $$
\par
b) {\sl d odd\/} :
$$ h_d(\lambda) \cdot P^{(d)}_\lambda( {\rm cos} \ \theta)  = {4i\omega_ d
\over( 2\pi)^{{ d+1 \over 2}}} \left({-1 \over {\rm sin} \ \theta}  { {\rm d}
\over {\rm d} \theta} \right)^{{d-1 \over 2}} \int^ \theta_ 0{ {\rm sin}
\left(\lambda +{d-2 \over 2} \right)\tau \ {\rm sin} \ \tau \over[ 2( {\rm
cos} \ \theta - {\rm cos} \ \tau)]^{ 1/2}} {\rm d} \tau \eqno (3.86) $$
\par
\smallskip
{\sl The second-kind Legendre functions\/}\nobreak\ : \par
We put :
$$ Q^{(d)}_\lambda( {\rm cos} \ \theta)  =-i {\omega_{ d-2} \over \omega_{
d-1}} (-i\ {\rm sin} \ \theta)^{ -(d-3)} \int^{ +i\infty}_ \theta {\rm e}^{i
\left(\lambda +{d-2 \over 2} \right)\tau}[ 2( {\rm cos} \ \tau - {\rm cos} \
\theta)]^{{ d-4 \over 2}} {\rm d} \tau . \eqno (3.87) $$
\par
For $ d\geq 3, $ this formula defines $ Q^{(d)}_\lambda( {\rm cos} \ \theta) $
as a ramified analytic function, whose
fundamental sheet is specified by the following conditions: $ {\rm cos} \
\theta \in \Bbb C\backslash] -\infty ,+1], $ with $ -\pi < {\rm Re} \ \theta
<\pi $
and $ {\rm Im} \ \theta  > 0 $ (in particular: for $ \lambda $ real, $
\lambda  > -1, $ $ Q^{(d)}_\lambda( {\rm cosh} \ v) > 0); $ it also shows
that,
for each $ {\rm cos} \ \theta $ in this domain, the function $ \lambda
\longrightarrow Q^{(d)}_\lambda( {\rm cos} \ \theta) $ is analytic in the
half-plane
$ \Bbb C^{(-1)}_+ = \{ \lambda ;\ {\rm Re} \ \lambda  > -1\} . $ \par
By applying the algebraic identities (A.28) and (A.29) (see proposition A.1),
we also obtain the
following expressions for the products $ h_d(\lambda) Q^{(d)}_\lambda , $
similar to Eqs.(3.85) and (3.86): \par
\smallskip
a) {\sl d even\/}\nobreak\ :
$$ h_d(\lambda) Q^{(d)}_\lambda( {\rm cos} \ \theta)  = {\omega_ d \over(
2\pi)^{{ d \over 2}}} \left({1 \over {\rm sin} \ \theta}  { {\rm d} \over
{\rm d} \theta} \right)^{ \left({d-2 \over 2} \right)} \left[ {\rm e}^{i
\left(\lambda +{d-2 \over 2} \right)\theta} \right] \eqno (3.88) $$
Eq.(3.88) performs a meromorphic continuation of the function $ \lambda
\longrightarrow Q^{(d)}_\lambda( {\rm cos} \ \theta) $ (for each $ {\rm cos} \
\theta ) $
in the whole complex $ \lambda $-plane, whose poles are localized (in view of
Eq.(3.84)) at the
integers $ \lambda =-1,-2,...,-(d-3). $ \par
\smallskip
b) {\sl d odd\/}\nobreak\ :
$$ h_d(\lambda) Q^{(d)}_\lambda( {\rm cos} \ \theta)  = - {2\omega_ d \over(
2\pi)^{{ d+1 \over 2}}} \left({1 \over {\rm sin} \ \theta}  { {\rm d} \over
{\rm d} \theta} \right)^{{d-1 \over 2}} \int^{ \theta +i\infty}_ \theta{ {\rm
e}^{i \left(\lambda +{d-2 \over 2} \right)\tau} {\rm sin} \ \tau \over[ 2(
{\rm cos} \ \tau - {\rm cos} \ \theta)]^{ 1/2}} {\rm d} \tau \eqno (3.89) $$
\par
Eq.(3.89) performs a meromorphic continuation of the function $ \lambda
\longrightarrow  Q^{(d)}_\lambda( {\rm cos} \ \theta) $
(for each $ {\rm cos} \ \theta ) $ in the half-plane $ \Bbb C^{ \left(-{d-3
\over 2} \right)}_+, $ whose poles are localized (in view of
Eq.(3.84)) at the integers $ -1,...,-{d-5 \over 2}. $ \par
\medskip
{\sl Expressions of the elementary perikernels\/} $ \Psi^{( d)}_\lambda : $
\par
\smallskip
a) {\sl d even\nobreak\ :\/} \par
In view of Eqs.(3.78) and (3.85), we immediately obtain the following
relations
$$ \Psi^{( d)}_\lambda( {\rm cos} \ \theta)  = {\pi \over \omega_ d}
h_d(\lambda)  P^{(d)}_\lambda( {\rm cos}(\theta -\pi)) \eqno (3.90) $$
(These identities of analytic functions are valid for $ {\rm cos} \ \theta
\in \Bbb C\backslash[ 1,+\infty[ , $ with $ 0< {\rm Re} \ \theta  < 2\pi , $
and for all $ \lambda \in \Bbb C). $ \par
\smallskip
b) {\sl d odd\nobreak\ :\/} \par
We first consider the case $ \lambda =\ell \in \Bbb Z $ which is specially
simple. In this case,
since $ {\rm cos} \left(\ell +{d-2 \over 2} \right)(\tau -\pi)  = (-1)^{\ell
+{d-3 \over 2}} {\rm sin} \left(\ell +{d-2 \over 2} \right)\tau , $ the r.h.s.
of Eqs.(3.79) and (3.86) become
identical up to a constant factor, and yield the following relations:
$$ \forall \ell \in \Bbb Z\ ,\ \Psi^{( d)}_{\ell}( {\rm cos} \ \theta)  =
(-1)^{\ell}{ \pi \over \omega_ d} h_d(\ell) P^{(d)}_{\ell}( {\rm cos} \
\theta) \eqno (3.91) $$
\par
{\sl Remarks\/}\nobreak\ : \par
i) For $ \lambda  = \ell $ integer, the expression (3.90) of $ \Psi^{(
d)}_\lambda $ in the {\sl even\/} case also
reduces to the form (3.91). \par
ii) For $ -d+3\leq \ell \leq -1, $ one has $ h_d(\ell) =0, $ and therefore the
corresponding functions $ \Psi^{( d)}_{\ell} $
vanish identically; this also results from proposition 14, since the
corresponding
functions $ \hat  f^{(d)}_{\ell} $ belong to the kernel of $ {\cal X}^{(c)}_d,
$ i.e.
$$ {\cal F}^{(d)}_{\ell}  = {\cal X}^{(c)}_d \left(\hat  f^{(d)}_{\ell}
\right) = 0\ . $$
\par
In the general case $ \lambda \in \Bbb C, $ we shall show that $ \Psi^{(
d)}_\lambda $ can be expressed in terms
of the corresponding second-kind Legendre function $ Q^{(d)}_\lambda $ and of
a residual analytic
function $ R^{(d)}_\lambda $ which enjoys special regularity properties. \par
In view of Eq.(3.79), we can write:
$$ \Psi^{( d)}_\lambda( {\rm cos} \ \theta)  = \left({1 \over 2\pi}
\right)^{{d-1 \over 2}} \left[{\cal J}^{(d)}_{\lambda +{d-2 \over 2}}( {\rm
cos} \ \theta) +{\cal J}^{(d)}_{- \left(\lambda +{d-2 \over 2} \right)}( {\rm
cos} \ \theta) \right] \eqno (3.92) $$
where $ {\cal J}^{(d)}_\alpha( {\rm cos} \ \theta) $ is defined, for all $
\alpha \in \Bbb C, $ by the following formulae:
$$ {\cal J}^{(d)}_\alpha( {\rm cos} \ \theta)  = {\rm e}^{-i\pi \alpha}
\left({1 \over {\rm sin} \ \theta}  { {\rm d} \over {\rm d} \theta}
\right)^{{d-1 \over 2}}{\cal J}_\alpha( {\rm cos} \ \theta) \eqno (3.93) $$
and
$$ {\cal J}_\alpha( {\rm cos} \ \theta)  = -i \int^ \theta_ 0{ {\rm
e}^{i\alpha \tau} {\rm sin} \ \tau \over[ 2( {\rm cos} \ \theta - {\rm cos} \
\tau)]^{ 1/2}} {\rm d} \tau \eqno (3.94) $$
\par
For $ {\rm Re} \ \alpha  > {1 \over 2}, $ one can easily show that the Cauchy
formula $ \int^{ }_{ \Gamma( \theta)}{ {\rm e}^{i\alpha \tau} {\rm sin} \
\tau \over[ 2( {\rm cos} \ \theta - {\rm cos} \ \tau)]^{ 1/2}} {\rm d} \tau  =
0 $
holds for a \lq\lq polygonal-type contour\rq\rq\ $ \Gamma( \theta) $ with
successive vertices $ (0, $ $ \theta , $ $ \theta +i\infty , $ $ \theta -2\pi
+i\infty , $
$ \theta -2\pi , $ $ -2\pi , $ $ 0) $ in the complex $ \tau $-plane (see
Fig.6). By regrouping together the
contributions of the paths $ (0,\theta) $ and $ (\theta -2\pi ,\ -2\pi) $ $ (
{\rm resp}(\theta ,\ \theta +i\infty) $ and $ (\theta -2\pi +i\infty ,\
\theta -2\pi) $ in
the latter, one then obtains the following alternative expression for the
function $ {\cal J}_\alpha : $
$$ \matrix{\displaystyle{\cal J}_\alpha( {\rm cos} \ \theta) & \displaystyle
=+i\ {\rm cotg} \ \pi \alpha \int^{ \theta +i\infty}_ \theta{ {\rm
e}^{i\alpha \tau} {\rm sin} \ \tau \over[ 2( {\rm cos} \ \tau - {\rm cos} \
\theta)]^{ 1/2}} {\rm d} \tau  \hfill \cr\displaystyle  & \displaystyle + {
{\rm e}^{i\pi \alpha} \over 2\ {\rm sin} \ \pi \alpha}  \int^ 0_{-2\pi}{ {\rm
e}^{i\alpha t} {\rm sin} \ t \over[ 2( {\rm cos} \ \theta - {\rm cos} \
t)]^{1/2}} {\rm d} t \hfill \cr} \eqno (3.95) $$
\par
In view of Eq.(3.93), we immediately obtain from Eq.(3.95) a similar
expression
for the function $ {\cal J}^{(d)}_\alpha $ in which the first term can be
interpreted (by taking Eq.(3.89)
into account) in terms of the function $ Q^{(d)}_{\alpha -{d-2 \over 2}}(
{\rm cos}(\theta -\pi)) =- {\rm e}^{-i\pi \left(\alpha -{d-2 \over 2}
\right)}Q^{(d)}_{\alpha -{d-2 \over 2}}( {\rm cos} \ \theta) : $ \par
\noindent$ \forall \theta , $ $ 0 < {\rm Re} \ \theta  < 2\pi , $ $ {\rm Im} \
\theta  \geq  0, $
$$ \matrix{\displaystyle{\cal J}^{(d)}_\alpha( {\rm cos} \ \theta) &
\displaystyle = {(-2\pi)^{{ d+1 \over 2}} \over 2\omega_ d} {\rm cotg} \ \pi
\alpha \ h_d{\scriptstyle \left(\alpha -{d-2 \over 2} \right)}
Q^{(d)}_{\alpha -{d-2 \over 2}}( {\rm cos}(\theta -\pi)) \hfill
\cr\displaystyle  & \displaystyle + { {\rm e}^{-i\pi \alpha} \over 2\ {\rm
sin} \ \pi \alpha}  R^{(d)}_{\alpha -{d-2 \over 2}}( {\rm cos} \ \theta) \ ,
\hfill \cr} \eqno (3.96) $$
where we have put:
$$ R^{(d)}_\lambda( {\rm cos} \ \theta)  = \left(- { {\rm d} \over {\rm d\
cos} \ \theta} \right)^{{d-1 \over 2}} \int^{ +\pi}_{ -\pi}{ {\rm e}^{i
\left(\lambda +{d-2 \over 2} \right)t} {\rm sin} \ t \over[ 2( {\rm cos} \
\theta + {\rm cos} \ t)]^{1/2}} {\rm d} t \eqno (3.97) $$
\par
{\bf Remark\nobreak\ :} \par
The function $ R^{(d)}_\lambda $ enjoys the following properties: \par
i) It is analytic and bounded by $ {\rm Cst}\vert {\rm cos} \ \theta\vert^{
-d/2} $ in the cut-plane $ \Bbb C\backslash[ -1,+\infty[ ; $
moreover its jump across the cut $ [-1,+1] $ cancels the corresponding jump
due to the
function $ Q^{(d)}_\lambda $ at the r.h.s. of Eq.(3.96). \par
ii) For each value of $ {\rm cos} \ \theta , $ the function $ \lambda
\longrightarrow R^{(d)}_\lambda( {\rm cos} \ \theta) $ is an entire function
which satisfies the following uniform majorization:
$$ \forall \lambda \in \Bbb C,\ \ \lambda  = \sigma +i\nu ,\ \ \left\vert
R^{(d)}_\lambda( {\rm cos} \ \theta) \right\vert  \leq  C( {\rm cos} \
\theta) {\rm cosh} \ \pi \nu \ , \eqno (3.98) $$
and the following identity:
$$ R^{(d)}_\lambda( {\rm cos} \ \theta)  = -R_{-\lambda -d+2}( {\rm cos} \
\theta) \eqno (3.98^\prime ) $$
\par
Since the function $ \alpha \longrightarrow{\cal J}^{(d)}_\alpha( {\rm cos} \
\theta) $ is also (in view of Eqs.(3.93) and (3.94)) an
entire function, it follows that Eq.(3.96) can be analytically continued with
respect to
$ \alpha $ as an identity between meromorphic functions in the whole complex $
\alpha $-plane, defining
thereby the function $ \lambda  \longrightarrow  {\rm sin} \ \pi \lambda \
h_d(\lambda) Q^{(d)}_\lambda( {\rm cos} \ \theta) $ as an entire function of $
\lambda $ for each $ {\rm cos} \ \theta $
in $ \Bbb C\backslash] -\infty ,-1], $ $ (-\pi < {\rm Re} \ \theta  < \pi) . $
\par
\medskip
{}From Eqs.(3.92) and (3.96), and by taking Eqs.(3.84') and (3.98') into
account,
we then obtain the following expression for $ \Psi^{( d)}_\lambda : $
$$ \eqalignno{ \Psi^{( d)}_\lambda( {\rm cos} \ \theta) &  = (-1)^{{d-1 \over
2}} {\pi \over \omega_ d} {\rm tg} \ \pi \lambda \ h_d(\lambda)
\left[Q^{(d)}_\lambda( {\rm cos}(\theta -\pi)) +Q^{(d)}_{-\lambda -d+2}( {\rm
cos}(\theta -\pi)) \right] &  \cr  &  - \left({1 \over 2\pi} \right)^{{d-1
\over 2}} {\rm tg} \ \pi \lambda \ R^{(d)}_\lambda( {\rm cos} \ \theta) &
(3.99) \cr} $$
\par
We now notice that, in view of Eqs.(3.86), (3.93) and (3.94), the function $
h_d(\lambda) P^{(d)}_\lambda( {\rm cos} \ \theta) $
can also be expressed in terms of the functions$ $ $ {\cal J}^{(d)}_\alpha(
{\rm cos} \ \theta) $ by the following relation:
$$ h_d(\lambda) P^{(d)}_\lambda( {\rm cos} \theta)  = i(-1)^{{d-1 \over
2}}{2\omega_ d \over( 2\pi)^{{ d+1 \over 2}}} \left[ {\rm e}^{i\pi
\left(\lambda +{d-2 \over 2} \right)} {\cal J}^{(d)}_{\lambda +{d-2 \over 2}}(
{\rm cos} \theta) - {\rm e}^{-i\pi \left(\lambda +{d-2 \over 2} \right)}{\cal
J}^{(d)}_{- \left(\lambda +{d-2 \over 2} \right)}( {\rm cos} \theta) \right]
$$
Then, by making use of Eqs.(3.96) and (3.98'), we obtain the following
functional
relation:
$$ P^{(d)}_\lambda( {\rm cos} \ \theta)  = (-1)^{{d+1 \over 2}} {\rm tg} \
\pi \lambda \left[Q^{(d)}_\lambda( {\rm cos} \ \theta) -Q^{(d)}_{-\lambda
-d+2}( {\rm cos} \ \theta) \right] \eqno (3.100) $$
which yields in particular the simple pole structure of the function $
\lambda \longrightarrow Q_{-\lambda -d+2}( {\rm cos} \ \theta) $
at all $ \lambda  = \ell  \in  \Bbb N, $ namely:
$$ P^{(d)}_{\ell}( {\rm cos} \ \theta)  = (-1)^{{d-1 \over 2}} \left[ {\rm tg}
\ \pi \lambda \ Q^{(d)}_{-\lambda -d+2}( {\rm cos} \ \theta) \right]_{\vert
\lambda =\ell} \eqno (3.101) $$
In view of Eq.(3.100), the expression (3.99) of $ \Psi^{( d)}_\lambda $ can
also be rewritten (for all $ \lambda  \in  \Bbb C) $
under the following form, whose first term is identical to the expression of $
\Psi^{( d)}_\lambda $ for
the case $ \gq\gq d $ even\rq\rq\ (see Eq.(3.90)):
$$ \eqalignno{ \Psi^{( d)}_\lambda( {\rm cos} \ \theta) &  = {\pi \over
\omega_ d} h_d(\lambda) P^{(d)}_\lambda( {\rm cos}(\theta -\pi)) &  (3.102)
\cr  &  + {\rm tg} \ \pi \lambda \left[(-1)^{{d-1 \over 2}}{2\pi \over
\omega_ d} h_d(\lambda) Q^{(d)}_\lambda( {\rm cos}(\theta -\pi)) - \left({1
\over 2\pi} \right)^{{d-1 \over 2}}R^{(d)}_\lambda( {\rm cos} \ \theta)
\right] &  \cr} $$
for $ \lambda  = \ell \in \Bbb N, $ the latter reduces to the previously
obtained expression (3.91) (since
$ P^{(d)}_{\ell}( {\rm cos}(\theta -\pi))  = (-1)^{\ell} P^{(d)}_{\ell}( {\rm
cos} \ \theta) ). $ \par
\medskip
{\bf Discontinuity formulae\nobreak\ :} \par
\smallskip
We shall now compute the discontinuity $ \Delta \Psi^{( d)}_\lambda $ of $
\Psi^{( d)}_\lambda $ across the cut $ \{ {\rm cos} \ \theta  \in
[1,+\infty[\} $
by applying the inverse Abel transformation formulae (3.64) and (3.65) to the
corresponding discontinuity $ \hat {\rm f}^{(d)}_\lambda $ of $ \hat
f^{(d)}_\lambda , $ given by Eqs.(3.76) and (3.77). \par
We thus obtain the following representations for $ \Delta \Psi^{(
d)}_\lambda( {\rm cosh} \ v) $ and the
corresponding expressions in terms of the Legendre functions $
P^{(d)}_\lambda $ and $ Q^{(d)}_\lambda . $ \par
\smallskip
a) {\sl d even\/} :
$$ \eqalignno{ \Delta \Psi^{( d)}_\lambda( {\rm cosh} \ v) & = -2\ {\rm sin} \
\pi \lambda \left({1 \over 2\pi}  {1 \over {\rm sinh} \ v} { {\rm d} \over
{\rm d} v} \right)^{{d-2 \over 2}} \left[ {\rm sinh} \left(\lambda +{d-2
\over 2} \right)v \right] &  \cr  &  = - {2\pi \over \omega_ d} {\rm sin} \
\pi \lambda \ h_d(\lambda) \left[P^{(d)}_\lambda( {\rm cosh} \ v)+2(-1)^{{d
\over 2}}Q^{(d)}_\lambda( {\rm cosh} \ v) \right] & (3.103) \cr} $$
(the r.h.s. of the latter follows from formulae (3.85) and (3.88)). \par
\medskip
{\bf Remark\nobreak\ :} \par
Computing $ \Delta \Psi^{( d)}_\lambda $ directly from Eq.(3.90) and comparing
with Eq.(3.103) yields the
following discontinuity formula for the function $ P^{(d)}_\lambda : $
$$ i \left[P^{(d)}_\lambda( {\rm cos}(iv-\pi)) -P^{(d)}_\lambda( {\rm
cos}(iv+\pi)) \right]=-2\ {\rm sin} \ \pi \lambda \left[P^{(d)}_\lambda( {\rm
cosh} \ v)+2(-1)^{{d \over 2}}Q^{(d)}_\lambda( {\rm cosh} \ v) \right] \eqno
(3.104) $$
\par
b) {\sl d odd\/} :
$$ \eqalignno{ \Delta \Psi^{( d)}_\lambda( {\rm cosh} \ v) & = -4\ {\rm sin} \
\pi \lambda \left({1 \over 2\pi}  {1 \over {\rm sinh} \ v} { {\rm d} \over
{\rm d} v} \right)^{{d-1 \over 2}} \int^ v_0{ {\rm cosh} \left(\lambda +{d-2
\over 2} \right)w\ {\rm sinh} \ w \over[ 2( {\rm cosh} \ v- {\rm cosh} \
w)]^{1/2}} {\rm d} w &  \cr  &  = -2\ {\rm sin} \ \pi \lambda \left[(-1)^{{d-1
\over 2}}{\pi \over \omega_ d} h_d(\lambda) {\rm tg} \ \pi \lambda
\left(Q^{(d)}_\lambda( {\rm cosh} \ v)+Q^{(d)}_{-\lambda -d+2}( {\rm cosh} \
v) \right) \right. &  \cr  & \left.+ \left({1 \over 2\pi} \right)^{{d-1 \over
2}} {iR^{(d)}_\lambda( {\rm cosh} \ v) \over {\rm cos} \ \pi \lambda} \right]
& (3.105) \cr} $$
\par
The r.h.s. of the latter can be obtained by recognizing the expression $
{\cal J}^{(d)}_{\lambda +{d-2 \over 2}}( {\rm cosh} \ w)+{\cal J}^{(d)}_{-
\left(\lambda +{d-2 \over 2} \right)}( {\rm cosh} \ w) $
in the previous integral and by applying Eq.(3.96). It can also
be derived directly from the defining equation
$ \Delta \Psi^{( d)}_\lambda( {\rm cosh} \ v) = i \left[\Psi^{( d)}_\lambda(
{\rm cos}(iv))-\Psi^{( d)}_\lambda( {\rm cos}(iv+2\pi)) \right] $ by applying
formula (3.99) and the identity
(deduced from Eq.(3.87)):
$$ Q^{(d)}_\lambda( {\rm cos}(iv-\pi))  - Q^{(d)}_\lambda( {\rm cos}(iv+\pi))
= 2i\ {\rm sin} \ \pi \lambda \ Q^{(d)}_\lambda( {\rm cosh} \ v) $$
\par
{\bf Remark :} \par
{}From the latter equation and Eq.(3.100), we also deduce the following useful
identity:
$$ i\left[P^{(d)}_\lambda( {\rm cos}(iv-\pi))-P^{(d)}_\lambda( {\rm
cos}(iv+\pi)) \right] = - 2 {\rm sin} \ \pi \lambda \ P^{(d)}_\lambda( {\rm
cosh} \ v) \eqno (3.105') $$
and thereby (in view of Eq.(3.102)) the following alternative form of $ \Delta
\Psi^{(d)}_\lambda : $
$$ \eqalignno{ \Delta \Psi^{( d)}_\lambda( {\rm cosh\ } \ v) & = -2\ {\rm sin}
\ \pi \lambda \left\{{ \pi \over \omega_ d}h_d(\lambda)
\left[P^{(d)}_\lambda( {\rm cosh} \ v)+2(-1)^{{d-1 \over 2}} {\rm tg} \ \pi
\lambda \ Q^{(d)}_\lambda( {\rm cosh} \ v) \right] \right. &  \cr  & \left.+
\left({1 \over 2\pi} \right)^{{d-1 \over 2}}{iR^{(d)}_\lambda( {\rm cosh} \ v)
\over {\rm cos} \ \pi \lambda} \right\} &  (3.106) \cr} $$
\par
In connection with the expressions (3.102) (resp.(3.106)) of $ \Psi^{(
d)}_\lambda $ (resp. $ \Delta \Psi^{( d)}_\lambda ), $ it
will be useful to complement the results of proposition 18 by similar
majorizations for the
functions $ P^{(d)}_\lambda $ and $ Q^{(d)}_\lambda $ (for $ d $ odd). We
shall prove: \par
\medskip
{\bf Proposition 19 :} \par
\smallskip
For every $ \theta $ such that $ -\pi < {\rm Re} \ \theta <\pi , $ the
following majorizations hold, in which $ C^{\prime}_ d $ and $
C^{\prime\prime}_ d $
denote suitable functions of $ {\rm cos} \ \theta , $ bounded at infinity in $
\Bbb C\backslash] -\infty ,+1]: $ \par
i) $ \forall \lambda  \in  \Bbb C, $ $ \lambda  = m+i\nu , $
$$ \left\vert h_d(\lambda) P^{(d)}_\lambda( {\rm cos} \ \theta) \right\vert
\leq  C^{\prime}_ d( {\rm cos} \ \theta)  {\rm e}^{ {\rm max}(m,-m-d+2)\vert
{\rm Im} \ \theta\vert} {\rm cosh} \ \pi \nu( 1+\vert \lambda\vert)^{{ d-1
\over 2}} \eqno (3.107) $$
\par
ii) $ \forall \lambda  \in  \Bbb C, $ $ \lambda  =m+i\nu , $ $ m>-{d-3 \over
2}, $ $ {\rm Im} \ \theta  > 0, $
$$ \left\vert h_d(\lambda) Q^{(d)}_\lambda( {\rm cos} \ \theta) \right\vert
\leq  C^{\prime\prime}_ d( {\rm cos} \ \theta)  {\rm e}^{(-m-d+2) {\rm Im} \
\theta} {\rm cosh}[\nu \varphi( {\rm Re} \ \theta)] \cdot( 1+\vert
\lambda\vert)^{{ d-1 \over 2}} \eqno (3.108) $$
where
$$ \varphi( u)=\vert u\vert \ \ {\rm for} \ \ 0\leq\vert u\vert \leq{ \pi
\over 2}\ \ {\rm and} \ \ \varphi( u) = \pi \ \ {\rm for} \ \ {\pi \over
2}<\vert u\vert \leq \pi \eqno (3.108^\prime ) $$
\par
{\bf Proof} : \par
The proof of i) is completely similar to that of proposition 18 by making use
of the
same parametrization in the expression (3.86) of $ h_d(\lambda)
P^{(d)}_\lambda $ one obtains:
$$ h_d(\lambda) P^{(d)}_\lambda( \zeta)  = {\rm Cst} \int^ 1_0{ {\rm d} \rho
\over( 1-\rho)^{ 1/2}} \left[A_d(\alpha ,\zeta ,\rho) Z^\alpha -A_d(-\alpha
,\zeta ,\rho) Z^{-\alpha} \right], $$
which can be majorized as the r.h.s. of Eq.(3.82) and yields the inequality
(3.107). \par
The proof of ii) is obtained by treating Eq.(3.89) via a similar
parametrization, namely
(with the same notations): $ z = \rho \zeta , $ with $ \rho \in[ 1,+\infty[ ,
$ which yields:
$$ h_d(\lambda) Q^{(d)}_\lambda( \zeta)  = {\rm Cst} \int^{ +\infty}_ 1{ {\rm
d} \rho \over( \rho -1)^{1/2}} B_d(\alpha ,\zeta ,\rho) Z^\alpha \ , \eqno
(3.109) $$
with $ \alpha =\lambda +{d-2 \over 2}. $ \par
In the latter, the function $ B_d $ denotes a polynomial of degree $ {d-1
\over 2} $ in $ \alpha , $ whose all
coefficients are rational functions of $ \sqrt{ \rho} , $ $ \sqrt{ \zeta} $
and $ \sqrt{ \rho^ 2\zeta^ 2-1} $ bounded by $ {\rm Cst.}\vert \zeta\vert^{
-{d-2 \over 2}}, $ and
therefore by $ {\rm Cst}^{\prime} \ {\rm e}^{- \left({d-2 \over 2} \right)
{\rm Im} \ \theta} . $ On the other hand, one has $ \left\vert Z^\alpha
\right\vert = \left\vert {\rm e}^{i\alpha \tau} \right\vert \leq {\rm e}^{-
\left(m+{d-2 \over 2} \right) {\rm Im} \ \tau} {\rm cosh}(\nu \ {\rm Re} \
\tau) , $
and, since (for $ {\rm cos} \ \tau  = \rho \ {\rm cos} \ \theta , $ $ \rho
\geq 1) $ $ \vert {\rm Re} \ \tau\vert  \leq  \varphi( {\rm Re} \ \theta) $
and $ \rho \ {\rm e}^{ {\rm Im} \ \theta} \leq {\rm Cst} \ {\rm e}^{ {\rm Im}
{\bf \ } \tau} $ (for $ {\rm Im} \ \theta  > v_0 > 0): $
$$ \vert Z\vert^ \alpha  \leq  {\rm Cst\ \rho}^{ - \left(m+{d-2 \over 2}
\right)} {\rm e}^{- \left(m+{d-2 \over 2} \right) {\rm Im} \ \theta} {\rm
cosh}[\nu \varphi( {\rm Re} \ \theta)] $$
Under the convergence condition $ m>-{d-3 \over 2} $ for the integral over $
\rho $ in Eq.(3.109), the previous
majorizations imply the inequality (3.108). \par
As a direct application of the majorizations (3.108) and (3.98) on the
functions $ Q^{(d)}_\lambda $
and $ R^{(d)}_\lambda , $ to the expressions (3.102) and (3.106) of $ \Psi^{(
d)}_\lambda $ and $ \Delta \Psi^{( d)}_\lambda , $ we can then state the
following: \par
\medskip
{\bf Corollary} : The following majorizations hold: $ \forall \lambda =m+i\nu
, $ with $ m > - {d-3 \over 2} $ \par
\smallskip
a) $ \forall \theta , $ $ 0< {\rm Re} \ \theta <2\pi , $ $ {\rm Im} \ \theta
> 0, $
$$ { \left\vert \Psi^{( d)}_\lambda( {\rm cos} \ \theta) -{\pi \over \omega_
d} h_d(\lambda) P^{(d)}_\lambda {\rm cos}(\theta -\pi) \right\vert \over\vert
{\rm sin} \ \pi \lambda\vert}  \leq  C_1( {\rm cos} \ \theta)  {\rm max}
\left(1,(1+\vert \lambda\vert)^{{ d-1 \over 2}} {\rm e}^{-(m+d-2) {\rm Im} \
\theta} \right) \eqno (3.110) $$
\par
b) $ \forall \theta =iv, $ $ v\geq 0 $
$$ \left\vert{ \Delta \Psi^{( d)}_\lambda( {\rm cosh} \ v) \over {\rm sin} \
\pi \lambda}  + 2{\pi \over \omega_ d} h_d(\lambda) P^{(d)}_\lambda( {\rm
cosh} \ v) \right\vert  \leq  C_2( {\rm cosh} \ v) {\rm max} \left(1,(1+\vert
\lambda\vert)^{{ d-1 \over 2}} {\rm e}^{-(m+d-2)v} \right) \eqno (3.111) $$
\par
We notice that in the derivation of these formulae, the poles of the functions
$ {\rm tg} \ \pi \lambda $ and
$ {1 \over {\rm cos} \ \pi \lambda} $ which appear in the two last terms of
Eq.(3.106) are in fact not present in the sum of
these terms (since $ \Psi^{( d)}_\lambda $ and $ P^{(d)}_\lambda $ are entire
functions of $ \lambda ), $ so that the uniform majorizations
(3.110) and (3.111) are valid in the whole half-plane $ \Bbb C^{ \left(-{d-3
\over 2} \right)}_+. $ \par
\medskip
\noindent{\bf 4. HARMONIC ANALYSIS FOR INVARIANT PERIKERNELS ON THE
COMPLEX HYPERBOLOID }$ {\bf X}^{( {\bf c})}_{ {\bf d-1}} $ \par
\smallskip
{\bf 4.1 Laplace transformation for invariant Volterra
kernels on the real hyperboloid }$ {\bf X}_{ {\bf d-1}} $ \par
\smallskip
Invariant Volterra kernels $ K \left(z,z^{\prime} \right) $ $ (K \in  V
\left(X_{d-1} \right)^{\natural} ) $ will be
represented
equivalently by $ F(z) \equiv  K \left(z,z_0 \right) $ or by $ \underline{{\rm
f}}( {\rm cosh\ } v)= {\rm f}(v) $
(with $ {\rm cosh} \ v = z^{(d-1)}). $ \par
\smallskip
For every invariant kernel in a class $ V^{(m)}_\mu \left(X_{d-1} \right) $
(see
\S 3.3, definition 2), such that $ m>-1, $ we shall
define the corresponding Laplace transform $ \tilde  F = L_d(F) $ by
the composition of the two mappings $ F \longrightarrow\hat {\rm f} = {\cal
R}_d(F), $ (i.e. $ \hat  F={\goth R}_d(F), $ with $ \hat  F \left(z_{iw}
\right) \equiv\hat {\rm f}(w), $
see \S 3.3) and $ \hat {\rm f} \longrightarrow\tilde  F = L \left(\hat {\rm f}
\right) $ (see \S 2.1). We shall say indifferently
that $ \tilde  F $ is the Laplace transform of $ F $ or of the invariant
Volterra kernel $ K. $ \par
This transformation $ L_d = L\circ  {\cal R}_d $ enjoys the following
properties: \par
\medskip
{\bf Proposition 20} : \par
For every kernel $ K $ in $ V^{(m)}_\mu \left(X_{d-1} \right), $ $ m>-1, $
with $ K \left(z,z_0 \right) = F(z), $ the
corresponding Laplace transform $ \tilde  F(\lambda)  \equiv  \left[
\left(L\circ{\cal R}_d \right)(F) \right](\lambda) $ is
analytic in the half-plane $ \Bbb C^{(m)}_+ = \{ \lambda \in \Bbb C; $
$ {\rm Re} \ \lambda  > m\} , $ continuous in the closed half-plane $ \bar
{\Bbb C}^{(m)}_+ $
and satisfies uniform bounds of the following form:
$$ \left\vert\tilde  F(\lambda) \right\vert  \leq  { {\rm Cst.} \over(
1+\vert \lambda\vert)^{ p(d)}} {\rm e}^{-[ {\rm Re} \ \lambda -m]\mu} \eqno
(4.1) $$
where $ p(d)= {d-2 \over 2} $ if $ d $ is even and $ p(d)= {d-3 \over 2} $ if
$ d $ is odd. \par
\medskip
{\bf Proof} : In view of proposition 9 ii) and of definition 2
(see \S 3.3), the function
$ \hat {\rm f}={\cal R}_d(F) $ satisfies the following conditions:
$$ \eqalignno{ {\rm i)} & \ \ \ {\rm supp} .\ \hat {\rm f}\subset[ \mu
,+\infty[ &  (4.2) \cr {\rm ii)} & \ \ \ \dlowlim{ {\rm max}}{0\leq r\leq
p(d)} \left\vert\hat {\rm f}^{(r)}(w) \right\vert  \leq  {\rm e}^{mw}.\hat
g(w)\ ,\ \ {\rm with} \ \ \hat  g\subset L^1 \left(\Bbb R^+ \right)\ \ \ \ \ \
\  & (4.3) \cr} $$
\par
\noindent The corresponding properties of the function
$ \tilde  F(\lambda)  = \left[L \left(\hat {\rm f} \right) \right](\lambda)  =
\int^{ \infty}_ 0 {\rm e}^{-\lambda w}\hat {\rm f}(w) {\rm d} w $
are then obtained by applying the argument of proposition 1i) (see
\S 2.1) to the functions $ \hat {\rm f}^{(r)} $ $ (0\leq r\leq p(d)). $ \par
\smallskip
{\bf Remark} : If condition ii) of the definition of $ V^{(m)}_\mu
\left(X_{d-1} \right) $ is
replaced by $ \vert {\rm f}(v)\vert  \leq  {\rm Cst.} \ {\rm e}^{mv}, $ the
corresponding condition on $ \left\{\hat {\rm f}^{(r)}(w) \right\} $
is obtained by replacing $ \hat  g(w) $ by a constant at the r.h.s. of
Eq.(4.3); it follows that the Laplace transform $ \tilde  F $ of $ F, $
still analytic in $ \Bbb C^{(m)}_+, $ only admits a boundary value in the
sense of distributions (namely a distribution of order 1) on the
boundary line $ {\rm Re} \ \lambda  = m $ of $ \Bbb C^{(m)}_+. $ \par
\medskip
{\bf Proposition 21 :} \par
Let $ K_1 \in  V^{(m)}_{\mu_ 1} \left(X_{d-1} \right), $ $ K_2 \in
V^{(m)}_{\mu_ 2} \left(X_{d-1} \right) $ be two kernels with respective
Laplace transforms $ \tilde  F_1 $ and $ \tilde  F_2. $ Then the Laplace
transform of the Volterra convolution product $ K = K_1\diamond K_2 $ is the
function $ \tilde  F(\lambda)  = \tilde  F_1(\lambda) \cdot\tilde
F_2(\lambda) . $ \par
\smallskip
{\bf Proof} : \par
Let $ F(z) = K \left(z,z_0 \right); $ proposition 8 implies that
$$ \hat  f(w) = {\cal R}_d(F)(w) = \left(\hat  K_1\diamond\hat  K_2 \right)
\left(z_{iw},z_0 \right)\ , $$
or (in view of Eq.(3.15')):
$$ \hat {\rm f}(w) = \int^ w_0\hat {\rm f}_1 \left(w-w^{\prime} \right)\hat
{\rm f}_2 \left(w^{\prime} \right) {\rm d} w^{\prime} \ \ ,\ \ {\rm with} \ \
\hat {\rm f}_i={\cal R}_d \left(F_i \right),\ \ i=1,2\ ; $$
the proposition follows from the standard
property of the Laplace transform $ L $ of a convolution product on
$ \Bbb R^+. $ \par
We shall now compute the integral representation of the transformation $ L_d $
and show that it coincides (up to a constant factor) with
the \lq\lq spherical Laplace transformation\rq\rq\ of invariant Volterra
kernels
in the sense of J. Faraut (see $\lbrack$Fa1$\rbrack$ $\lbrack$Fa,V$\rbrack$
for the case $ d=3 $ and
$\lbrack$Fa2$\rbrack$ for a more general presentation). \par
In view of the Abel-type expression (3.32) of the
transformation $ {\cal R}_d, $ we can write for any kernel $ K $ in $
V^{(m)}_\mu \left(X_{d-1} \right) $ (with the
notations: $ K \left(z,z_0 \right) \equiv  F(z) \equiv  \underline{{\rm f}}(
{\rm cosh} \ v) \equiv  {\rm f}(v), $ $ {\rm cosh} \ v = z^{(d-1)}): $
$$ \tilde  F(\lambda)  = \left[ \left(L\circ{\cal R}_d \right)(F)
\right](\lambda)  = \omega_{ d-2} \int^{ \infty}_ 0 {\rm e}^{- \left(\lambda
+{d-2 \over 2} \right)w} {\rm d} w \int^ w_0 {\rm f}(v)[2( {\rm cosh} \ w-
{\rm cosh} \ v)]^{{d-4 \over 2}} {\rm sinh} \ v\ {\rm d} v \eqno (4.4) $$
\par
In view of proposition 9, which puts a uniform bound of the form $ \hat {\rm
g}(w) {\rm e}^{ \left(m+{d-2 \over 2} \right)w}, $ with
$ \hat {\rm g}\in L_1, $ on the subintegral over $ v $ in formula (4.4), it
follows that the double integral is
uniformly convergent for all $ \lambda $ in $ \Bbb C^{(m)}_+ $ (if $ m > -1).
$ \par
By inverting the integration order at the r.h.s. of Eq.(4.4), we then obtain
the
following result: \par
\medskip
{\bf Proposition 22} : \par
For every kernel $ K \left(z,z_0 \right) \equiv  {\rm f}(v) $ in $
V^{(m)}_\mu \left(X_{d-1} \right), $ $ m>-1, $ (with $ z^{(d-1)}= {\rm cosh} \
v), $ the corresponding
Laplace transform $ \tilde  F(\lambda) $ can be represented in its analyticity
domain $ \Bbb C^{(m)}_+ $ by the
following formula:
$$ \tilde  F(\lambda)  = \omega_{ d-1} \int^{ \infty}_ 0 {\rm f}(v)
Q^{(d)}_\lambda( {\rm cosh} \ v)( {\rm sinh} \ v)^{d-2} {\rm d} v\ , \eqno
(4.5) $$
where $ Q^{(d)}_\lambda $ is the second-kind Legendre function in dimension $
d, $ defined, for $ {\rm Re} \ \lambda  > -1, $ by the integral
representation (see in this connection \S 3.8, Eq.(3.87)):
$$ Q^{(d)}_\lambda( {\rm cosh} \ v) = {\omega_{ d-2} \over \omega_{ d-1}} {1
\over( {\rm sinh} \ v)^{d-3}} \int^{ \infty}_ v {\rm e}^{- \left(\lambda +{d-2
\over 2} \right)w}[2( {\rm cosh} \ w- {\rm cosh} \ v)]^{{d-4 \over 2}} {\rm d}
w \eqno (4.6) $$
\par
{\sl Identification with the spherical Laplace transform of J. Faraut\/}: \par
Let us now introduce, as in $\lbrack$Fa1$\rbrack$, $\lbrack$Fa,V$\rbrack$, the
Poisson function on $ X: $
$$ p^{(d)}_\lambda( z) = \left(z^{(0)}+z^{(d-1)} \right)^{-(\lambda +d-2)} =
{\rm e}^{-(\lambda +d-2)w} \eqno (4.7) $$
and consider the associated \lq\lq spherical function\rq\rq
$$ \Phi^{( d)}_\lambda( {\rm cosh} \ v) = \int^{ }_ Hp^{(d)}_\lambda \left(h
\left(\varphi ,\vec  \alpha \right)a(iv)z_0 \right) {\rm d} h \eqno (4.8) $$
where we have used the Cartan-type decomposition $ z = h(\varphi ,\alpha)
a(iv)z_0 $ and where integration is
taken on the stabilizer $ H $ of $ z_0 $ in $ {\rm SO}(1,d-1). $ The measure $
{\rm d} h $ is such that:
$$ {\rm d} h\wedge ( {\rm sinh} \ v)^{d-2} {\rm d} v = 2(-1)^{d-2} \left.{
{\rm d} z^{(0)}\wedge ...\wedge {\rm d} z^{(d-1)} \over {\rm d} s(z)}
\right\vert_{ X_{d-1}} $$
But the latter can also be expressed as follows (in view of Eqs.(3.7) and
(3.24'), written
for $ \tau =iw $ and $ \theta =iv): $
$$ ( {\rm d} v\wedge {\rm d} h)( {\rm sinh} \ v)^{d-2}= {\rm e}^{(d-2)w} {\rm
d} w\wedge {\rm d}\vec  \zeta  = - {\rm e}^{{d-2 \over 2}w}[2( {\rm cosh} w-
{\rm cosh} v)]^{{d-4 \over 2}} {\rm d} w\wedge( {\rm sinh} v {\rm d}
v)_{\wedge} {\rm d} {\bf \sigmabf}_{ d-3} \left(\vec  \alpha \right), $$
so that we can choose to represent $ {\rm d} h $ by:
$$ {\rm d} h = {\rm e}^{{d-2 \over 2}w}[2( {\rm cosh} \ w- {\rm cosh} \
v)]^{{d-4 \over 2}} ( {\rm sinh} \ v)^{-(d-3)} {\rm d} w\wedge {\rm d} {\bf
\sigmabf}_{ d-3} \left(\vec  \alpha \right) \eqno (4.9) $$
\par
By taking Eqs.(4.7) and (4.9) into account, one then checks that Eq.(4.8)
yields the
identification: $ Q^{(d)}_\lambda( {\rm cosh} \ v) = {1 \over \omega_{ d-1}}
\Phi^{( d)}_\lambda( {\rm cosh} \ v), $ and therefore Eq.(4.5) coincides with
the
definition of the spherical Laplace transform (see $\lbrack$Fa1$\rbrack$,
$\lbrack$Fa,V$\rbrack$, $\lbrack$Fa2$\rbrack$). \par
\medskip
{\bf 4.2 Fourier-Legendre expansion for invariant kernels on the sphere }$
{\bf S}_{ {\bf d-1}} {\bf :} ${\bf\
interpretation as a one-dimensional Fourier series} \par
\smallskip
By replacing the hyperboloid $ X_{d-1} $ by the sphere $ S_{d-1} $ and the
invariant Volterra kernels
on $ X_{d-1} $ by the invariant kernels on $ S_{d-1}, $ we shall now give
results which parallel those of
\S 4.1. \par
Being given an invariant kernel $ {\bf K} \left(z,z_0 \right) \equiv  {\bf
F}(z) $ on the sphere $ S_{d-1}, $ the
corresponding transform $ \hat {\bf F}= \underline{{\goth R}}_d( {\bf F}), $
defined on the circle $ \hat  S $
(see proposition 11) admits a Fourier series with coefficients $ \left[\hat
{\bf f} \right]_{\ell} $ normalized as follows:
$$ \left[\hat {\bf f} \right]_{\ell}  = \int^ \pi_{ -\pi} {\rm e}^{i\ell t}
\hat {\bf f}(t) {\rm d} t\ \ \ \ ( {\rm with} \ \hat {\bf f}(t)\equiv  \hat
{\bf F} \left(z_t \right)) \eqno (4.10) $$
Since $ \hat {\bf F} $ satisfies Eq.(3.30), the set $ \left\{ \left[\hat {\bf
f} \right]_{\ell} ;\ \ell \in \Bbb Z \right\} $ satisfies the equivalent
set of relations:
$$ \forall \ell  \in  \Bbb Z\ ,\ \ \ \left[\hat {\bf f} \right]_{\ell}  =
(-1)^d \left[\hat {\bf f} \right]_{-(\ell +d-2)} \eqno (4.11) $$
We then consider more specially the subset of coefficients $ \left[\hat {\bf
f} \right]_{\ell} $ such that $ \ell  \geq  0 $ and define
the mapping $ {\bf L} $ by $ \hat {\bf F} \longrightarrow  {\bf L} \left(\hat
{\bf F} \right) = \left\{ \left[\hat {\bf f} \right]_{\ell} ;\ \ell \in \Bbb N
\right\} . $ \par
On the other hand, we define the mapping $ {\bf L}_d $ which associates with $
{\bf F} $ the
\lq\lq Fourier-Legendre\rq\rq\ sequence $ \left\{[ {\bf f}]_{\ell} ;\ell \in
\Bbb N \right\}  = {\bf L}_d( {\bf F}) $ whose elements $ [ {\bf f}]_{\ell} $
are
defined by the following formulae (where: $ \underline{{\bf f}}( {\rm cos} \
u) = \underline{{\bf f}} \left(z^{[d-1]} \right)= {\bf F}(z)): $
$$ [ {\bf f}]_{\ell} =\omega_{ d-1} \int^ \pi_ 0 \underline{{\bf f}}( {\rm
cos} \ u) P^{(d)}_{\ell}( {\rm cos} \ u)( {\rm sin} \ u)^{d-2} {\rm d} u\ ;
\eqno (4.12) $$
in the latter, $ P^{(d)}_{\ell} $ denotes the ultraspherical Legendre
polynomial defined by the
following integral representation (see e.g.$\lbrack$Vi$\rbrack$ and
$\lbrack$Fa-3$\rbrack$):
$$ P^{(d)}_{\ell}( {\rm cos} \ u) = {\Gamma \left({d-1 \over 2} \right) \over
\sqrt{ \pi} \ \Gamma \left({d-2 \over 2} \right)} \int^ \pi_ 0( {\rm cos} \
u+i\ {\rm sin} \ u\ {\rm cos} \ \eta)^{ \ell}( {\rm sin} \ \eta)^{( d-3)}
{\rm d} \eta \ , \eqno (4.13) $$
where the normalization constant $ {\Gamma \left({d-1 \over 2} \right) \over
\sqrt{ \pi} \ \Gamma \left({d-2 \over 2} \right)} = {\omega_{ d-2} \over
\omega_{ d-1}} $ is the one adopted in $\lbrack$Fa1$\rbrack$. \par
We shall now prove: \par
\smallskip
{\bf Proposition 23} : \par
For every invariant kernel $ {\bf K} \left(z,z_0 \right)\equiv {\bf F}(z) $ on
$ S_{d-1}, $ the two associated sequences
of Fourier coefficients $ \left\{ \left[\hat {\bf f} \right]_{\ell} ;\ell \in
\Bbb N \right\} $ and $ \left\{[ {\bf f}]_{\ell} ;\ell \in \Bbb N \right\} $
coincide, i.e. one has $ {\bf L}_d = {\bf L}  \circ \underline{{\goth R}}_d. $
\par
\medskip
{\bf Proof :} \par
In view of formulae (3.29a), (3.29b) and (4.10) (with $ {\bf f}(u)\equiv
\underline{{\bf f}}( {\rm cos} \ u)), $
we can write:
$$ \matrix{\displaystyle{ i^{d-2} \left[\hat {\bf f} \right]_{\ell} \over
\omega_{ d-2}} & \displaystyle = \int^ \pi_ 0 {\rm e}^{i \left(\ell +{d-2
\over 2} \right)t} {\rm d} t \int^ t_0 \underline{{\bf f}}( {\rm cos} \ u)[2(
{\rm cos} \ u- {\rm cos} \ t)]^{{d-4 \over 2}} {\rm sin} \ u\ {\rm d} u
\hfill \cr\displaystyle  & \displaystyle + {\rm e}^{i\pi( d-2)} \int^
0_{-\pi} {\rm e}^{i \left(\ell +{d-2 \over 2} \right)t} {\rm d} t \int^{ -t}_0
\underline{{\bf f}}( {\rm cos} \ u)[2( {\rm cos} \ u- {\rm cos} \ t)]^{{d-4
\over 2}} {\rm sin} \ u\ {\rm d} u \hfill \cr} \eqno (4.14) $$
Then, by making the change of variable $ (t,u) \longrightarrow  (t+2\pi ,u) $
in the second term at the r.h.s.
of this equation, regrouping the two terms and inverting the integrations over
$ t $ and $ u, $ we
can rewrite (Eq.(4.14) as follows:
$$ \left[\hat {\bf f} \right]_{\ell}  = \int^ \pi_ 0 \underline{{\bf f}}(
{\rm cos\ } u) \Pi^{( d)}_{\ell}( {\rm cos} \ u)( {\rm sin} \ u)^{d-2} {\rm d}
u\ , \eqno (4.15) $$
with:
$$ \Pi^{( d)}_{\ell}( {\rm cos} \ u) = (-i)^{d-2}\omega_{ d-2}( {\rm sin} \
u)^{-(d-3)} \int^{ 2\pi -u}_u {\rm e}^{i \left(\ell +{d-2 \over 2}
\right)t}[2( {\rm cos} \ u- {\rm cos} \ t)]^{{d-4 \over 2}} {\rm d} t \eqno
(4.16) $$
In order to show that $ \left[\hat {\bf f} \right]_{\ell}  = [ {\bf
f}]_{\ell} , $ it remains to check (according to Eq.(4.12)) that one has: $
\Pi^{( d)}_{\ell} =\omega_{ d-1}P^{(d)}_{\ell} . $ \par
In the integral representation (4.13) of $ P^{(d)}_{\ell} , $ let us
substitute to $ \eta $ the
complex integration variable $ \tau $ defined by
$$ {\rm e}^{i\tau}  = {\rm cos} \ u + i\ {\rm sin} \ u\ {\rm cos} \ \eta
\eqno (4.17) $$
We then check that:
$$ 2 {\rm e}^{i\tau}( {\rm cos} \ \tau  - {\rm cos} \ u) = \left( {\rm
e}^{i\tau} - {\rm e}^{iu} \right) \left( {\rm e}^{i\tau} - {\rm e}^{-iu}
\right) = {\rm sin}^2u\ {\rm sin}^2\eta \ , \eqno (4.18) $$
and since $ {\rm e}^{i\tau} {\rm d} \tau  = - {\rm sin} \ u\ {\rm sin} \ \eta
\ {\rm d} \eta , $ the integrand at the r.h.s. of formula (4.13) can be
rewriten as follows (in view of Eqs.(4.17) and (4.18)):
$$ \matrix{\displaystyle - {\rm sin} \ u^{-(d-3)} {\rm e}^{i(\ell +1)\tau}
\left[ \left( {\rm e}^{i\tau} - {\rm e}^{iu} \right) \left( {\rm e}^{i\tau} -
{\rm e}^{-iu} \right) \right]^{{d-4 \over 2}} {\rm d} \tau  = \hfill
\cr\displaystyle \ \ \ \ \ \ \ = -( {\rm sin} \ u)^{-(d-3)} {\rm e}^{i
\left(\ell +{d-2 \over 2} \right)\tau}[ 2( {\rm cos} \ \tau - {\rm cos} \
u)]^{{d-4 \over 2}} {\rm d} \tau \hfill \cr} \eqno (4.19) $$
\par
In order to determine the integration path, let us consider an intermediate
step
where $ {\rm e}^{i\tau} $ is chosen as the integration variable; the original
path (corresponding to $ \eta  \in  [0,\pi] ) $
is the (oriented) linear segment $ \delta_ 0(u) $ starting at $ {\rm e}^{iu} $
and ending at $ {\rm e}^{-iu}. $ \par
Since (as exhibited by the $\ell$.h.s. of Eq.(4.19)), the integrand is an
analytic
function of $ {\rm e}^{i\tau} $ in the disk $ \left\vert {\rm e}^{i\tau}
\right\vert  < 1 $ (since $ \ell  \in  \Bbb N), $ the integration path $
\delta_ 0(u) $ can
be replaced by the circular path $ \delta_ +(u) = \left\{ {\rm e}^{i\tau} ;
\right. $ $ \tau =t, $ $ u \leq  t \leq  2\pi -u\} $ (see Fig.7). Moreover,
by using the fact that $ \left[ \left( {\rm e}^{i\tau} - {\rm e}^{iu} \right)
\left( {\rm e}^{i\tau} - {\rm e}^{-iu} \right) \right]^{{d-4 \over 2}} $ is
positive for $ {\rm e}^{i\tau} \in \delta_ 0(u)\cup  \Bbb R, $ and therefore
at
$ {\rm e}^{i\tau} = {\rm e}^{i\pi} , $ we conclude from the left equality
(4.18) that in the r.h.s. of Eq.(4.19), the
following specification holds (for $ \tau =t; $ $ u \leq  t \leq  2\pi -u): $
$$ [2( {\rm cos} \ t- {\rm cos} \ u)]^{{d-4 \over 2}}=(-i)^{d-4}[2( {\rm cos}
\ u- {\rm cos} \ t)]^{{d-4 \over 2}} $$
\par
By taking the latter into account, the integral representation (4.13) can be
replaced by the following one:
$$ P^{(d)}_{\ell}( {\rm cos} \ u) = (-i)^{d-2}{\omega_{ d-2} \over \omega_{
d-1}}( {\rm sin} \ u)^{-(d-3)} \int^{ 2\pi -u}_u {\rm e}^{i \left(\ell +{d-2
\over 2} \right)t}[2( {\rm cos} \ u- {\rm cos} \ t)]^{{d-4 \over 2}} {\rm d} t
\eqno (4.20) $$
\par
By comparing with Eq.(4.16), we conclude that $ \Pi^{( d)}_{\ell} =\omega_{
d-1}P^{(d)}_{\ell} , $ which ends the proof
of proposition 23. \par
\medskip
{\bf Remark} {\bf 1.} In the argument of the previous proof, one could have
chosen
as well the integration path $ \delta_ -(u) = \left\{ {\rm e}^{i\tau} ;
\right. $
$ \tau =t, $ $ -u\leq t\leq u\} , $ (with the orientation shown on Fig.7),
thus obtaining the alternative integral
representation of the hyperspherical Legendre polynomials:
$$ P^{(d)}_{\ell}( {\rm cos} \ u) = +{\omega_{ d-2} \over \omega_{ d-1}} (
{\rm sin} \ u)^{-(d-3)} \int^ u_{-u} {\rm e}^{i \left(\ell +{d-2 \over 2}
\right)t}[2( {\rm cos} \ t- {\rm cos} \ u)]^{{d-4 \over 2}} {\rm d} t \eqno
(4.20^\prime ) $$
\par
The latter extends in a natural way to a representation of the first-kind
Legendre
functions (by replacing $ \ell \in \Bbb N $ by $ \lambda \in \Bbb C) $ which
is obviously equivalent to the one
introduced in \S 3.8 (see Eq.(3.83)). \par
\medskip
{\bf Remark 2.} Since $ {\bf F} $ is completely characterized by its
Fourier-Legendre sequence $ \left\{[ {\bf f}]_{\ell} ; \right. $
$ \ell \in \Bbb N\} , $ the Fourier sequence $ \left\{ \left[\hat {\bf f}
\right]_{\ell} ;\ \ell \in \Bbb Z \right\} $ of $ \hat {\bf F} =
\underline{{\goth R}}_d( {\bf F}) $ must also be completely
determined by the former. The relations $ \left\{[ {\bf f}]_{\ell} =
\left[\hat {\bf f} \right]_{\ell} =(-1)^d \left[\hat {\bf f} \right]_{-(\ell
+d-2)},\ \ell \in \Bbb N \right\} $
of proposition 23 could in fact be completed by $ E \left({d-2 \over 2}
\right) $ relations expressing the remaining
coefficients $ \left[\hat {\bf f} \right]_{\ell} , $ $ -(d-3) \leq  \ell
\leq  -1, $ as (linear) functionals of $ \left\{[ {\bf f}]_{\ell} ;\ \ell \in
\Bbb N \right\} ; $ these relations
can be derived from the special form (3.30') of the function $ \hat {\bf F}
\left(z_t \right)\equiv\hat {\bf f}(t), $ established in
proposition 11. \par
\medskip
{\bf 4.3 The \lq\lq theorem F.G.\rq\rq , part I: Fourier-Laplace
transformation for invariant perikernels on the complex
hyperboloid }$ X^{(c)}_{d-1} $ \par
We shall now establish our main result which concerns the Fourier-Laplace
transformation for invariant perikernels of moderate growth on $ X^{(c)}_{d-1}
$ and corresponds to a
refined version of the part (I) of the theorem F.G., presented in our section
1. \par
As in \S 3.6, we shall consider $ H^{(c)} $-invariant functions $ {\cal F} $
in subspaces $ \left[{\cal V}_1 \right]^m_\mu \left(X^{(c)}_{d-1} \right) $
with $ m>-1, $ $ \mu \geq 0 $ (see definition 3) and the associated invariant
triplets $ ({\cal F}, {\bf F} ,F) $ (each
triplet of functions $ ({\cal F}, {\bf F} ,F) $ representing a triplet $
({\cal K}, {\bf K} ,K) $ associated with the invariant perikernel $ {\cal K}
\left(z,z^{\prime} \right) $
such that $ {\cal F}(z) = {\cal K} \left(z,z_0 \right)). $ \par
We can then state: \par
\smallskip
{\bf Theorem 1 :} With every $ H^{(c)} $-invariant triplet $ ({\cal F}, {\bf
F} ,F) $ on $ X^{(c)}_{d-1}, $ such that $ {\cal F}\in \left[{\cal V}_1
\right]^m_\mu \left(X^{(c)}_{d-1} \right), $
$ m>-1, $ one can associate: \par
\smallskip
a) the transform $ \tilde  F = {\cal L}_d({\cal F}), $ where $ {\cal L}_d =
{\cal L}\circ{\cal R}^{(c)}_d, $ which can be expressed as follows
in terms of the function $ \underline{{\rm f}} \left(z^{(1)} \right) = F(z) $
$ (z\in X^\mu_ +): $
$$ \tilde  F(\lambda)  = \omega_{ d-1} \int^{ \infty}_{ {\rm cosh} \ \mu}
\underline{{\rm f}}(\zeta)  Q^{(d)}_\lambda( \zeta) \left(\zeta^ 2-1
\right)^{{d-3 \over 2}} {\rm d} \zeta \eqno (4.21) $$
\par
b) the Fourier-Legendre sequence $ {\bf L}_d( {\bf F}) = \left\{[ {\bf
f}]_{\ell} ; \right. $ $ \ell \in \Bbb N\} , $ where:
$$ [ {\bf f}]_{\ell}  = \omega_{ d-1} \int^{ +1}_{-1} \underline{{\bf
f}}(\zeta) P^{(d)}_{\ell}( \zeta) \left(1-\zeta^ 2 \right)^{{d-3 \over 2}}
{\rm d} \zeta \ , \eqno (4.22) $$
with $ \underline{{\bf f}} \left(z^{(d-1)} \right) = {\bf F}(z) $ $ (z\in \Bbb
S_{d-1}). $ \par
\smallskip
Then the following properties hold: \par
\smallskip
i) $ \tilde  F \in  \tilde { \underline{{\cal O}} }^{\mu ,p(d)} \left(\Bbb
C^{(m)}_+ \right) $ (with $ p(d) $ specified in propositions 13
and 20) and the transformation $ {\cal L}_d $ defines a
continuous mapping from $ \left[{\cal V}_1 \right]^m_\mu \left(X^{(c)}_{d-1}
\right) $ into $ \tilde { \underline{{\cal O}} }^{\mu ,p(d)} \left(\Bbb
C^{(m)}_+ \right). $ \par
\smallskip
ii) $ \forall \ell  \in  \Bbb N, $ $ \ell  \geq  m, $ $ [ {\bf f}]_{\ell}  =
\tilde  F(\ell) . $ \par
\medskip
{\bf Proof\nobreak\ :} \par
It is straightforward to check that the transform $ \tilde  F={\cal
L}_d({\cal F}) $ coincides with the
Laplace transform $ L_d(F) $ of the jump function $ F $ of $ {\cal F} $ (with
support $ X^\mu_ + $ on $ X_{d-1}). $ In fact, we
can write: $ \tilde  F = {\cal L} \left(\hat  f \right) $ with $ \hat
f={\cal R}^{(c)}_d({\cal F}), $ and therefore (in view of propositions 1 and
13) $ \tilde  F=L \left(\Delta\hat  f \right), $
with $ \Delta\hat  f = \hat {\rm f} = {\cal R}_d(F), $ which implies that $
\tilde  F = \left[L \circ{\cal R}_d \right](F) = L_d(F). $ Moreover, since the
assumption $ {\cal F }\in  \left[{\cal V}_1 \right]^m_\mu \left(X^{(c)}_{d-1}
\right) $ obviously implies that $ F \in  V^m_\mu \left(X_{d-1} \right) $
(compare definitions 2 and
3), the transform $ \tilde  F $ of $ {\cal F} $ enjoys all the properties
established in \S 4.1. In particular, the
expression (4.21) of $ \tilde  F $ is equivalent to the formula (4.5) obtained
in proposition 22. \par
Let us now prove the properties described under i) and ii). We first notice
that
the relation $ \tilde  F \in  \tilde { \underline{{\cal O}} }^{\mu ,p(d)}
\left(\Bbb C^{(m)}_+ \right) $ is just the content of the inequality
(4.1) of proposition 20. Moreover, since (for each $ m>-1) $ $ {\cal
R}^{(c)}_d $ defines a continuous mapping
from $ \left[{\cal V}_1 \right]^m_\mu \left(X^{(c)}_{d-1} \right) $ into $
{\cal O}^{(d)}_{m,p(d)} \left(\dot {\cal J}^{(\mu)} \right) $ (see proposition
13) and since $ {\cal L} $ defines a continuous
mapping from $ {\cal O}^{(d)}_{m,p(d)} \left(\dot {\cal J}^{(\mu)} \right) $
into $ \tilde { \underline{{\cal O}} }^{\mu ,p(d)} \left(\Bbb C^{(m)}_+
\right) $ (see proposition 6), the
transformation $ {\cal L}_d = {\cal L }\circ  {\cal R}_d $ enjoys the
continuity properties described under i). \par
Let us now consider the function $ \hat {\bf f} = \hat  f_{\vert \dot {\Bbb
R}}, $ which is by definition equal
to $ \underline{{\cal R}}_d( {\bf F}) $ (see \S 3.5). In view of proposition
20, we then have:
$$ \left\{[ {\bf f}]_{\ell} \ ;\ \ell \in \Bbb N \right\}  = {\bf L}_d( {\bf
F}) = {\bf L} \left(\hat {\bf f} \right)= \left\{ \left[\hat {\bf f}
\right]_{\ell} \ ;\ \ell \in \Bbb N \right\} \ . $$
Then, by applying proposition 4 to the function $ \hat  f = {\cal
R}^{(c)}_d({\cal F}) $ (which belongs to $ {\cal O}^{(d)}_m \left(\dot {\cal
J}^{(\mu)} \right) $
in view of proposition 13), we obtain (in view of Eq.(2.25)): \par
$ \forall \ell  \in  \Bbb N, $ with $ \ell  \geq  m, $
$$ \left[\hat {\bf f} \right]_{\ell}  = \left[{\cal L} \left(\hat  f \right)
\right](\ell)  = \left[{\cal L}_d({\cal F}) \right](\ell)  = \tilde  F(\ell) ,
$$
which proves property ii). \par
The previous statement can be supplemented by a theorem on convolution
products.
According to proposition 12, the invariant perikernels $ {\cal K} $ with
moderate
growth $ m $ (for $ m>-1) $ form a (normed) algebra $ {\cal W}^{(m)}
\left(X^{(c)}_{d-1} \right)^{\natural} $ for the $ \ast^{( c)} $-convolution
product,
namely a subalgebra of the general algebra of perikernels introduced in
$\lbrack$B.V-1$\rbrack$ (i.e. $ {\cal W} \left(X^{(c)}_{d-1} \right) =
\bigcup^{ }_{ \mu \geq 0}{\cal W}_\mu ). $
Let us now call $ \tilde { \underline{{\cal O}} }^{[p]} \left(\Bbb C^{(m)}_+
\right) $ (for $ p\in \Bbb N) $ the (normed) algebra of
functions $ \tilde  F(\lambda) $ which are analytic in $ \Bbb C^{(m)}_+, $ and
such that $ (1+\vert \lambda\vert)^ p \left\vert\tilde  F(\lambda)
\right\vert $ belongs to $ {\cal C}^0 \left(\bar {\Bbb C}^{(m)}_+ \right) $
(i.e. $ \tilde { \underline{{\cal O}} }^{[p]} \left(\Bbb C^{(m)}_+ \right) =
\bigcup^{ }_{ \mu \geq 0}\tilde { \underline{{\cal O}} }^{\mu ,p} \left(\Bbb
C^{(m)}_+ \right)). $ We then have: \par
\medskip
{\bf Theorem 2\nobreak\ :} For each dimension $ d, $ $ d\geq 3, $ and for each
real number $ m $ such that $ m>-1, $
the transformation $ {\cal L}_d $ defines a morphism of (normed) algebras
between $ {\cal W}^{(m)} \left(X^{(c)}_{d-1} \right)^{\natural} $ and $ \tilde
{ \underline{{\cal O}} }^{[p(d)]} \left(\Bbb C^{(m)}_+ \right). $ \par
More precisely let $ \left({\cal K}_i, {\bf K}_i,K_i \right), $ $ i=1,2, $ and
$ ({\cal K}, {\bf K} ,K) $ be three invariant triplets on $ X^{(c)}_{d-1} $
such that:
\smallskip
$$ \matrix{\displaystyle {\rm i)} & \displaystyle  & \displaystyle{\cal K }=
{\cal K}_1\ast^{( c)}{\cal K}_2\ \ \ \ \ \ \ ( {\bf K}  = {\bf K}_1\ast {\bf
K}_2\ \ ,\ \ \ K = K_1 \diamond  K_2) \hfill \cr\displaystyle  &
\displaystyle  & \displaystyle \hfill \cr\displaystyle  & \displaystyle  &
\displaystyle \cr\displaystyle {\rm ii)} & \displaystyle  &
\displaystyle{\cal F}_i(z) = {\cal K}_i \left(z,z_0 \right)\in \left[{\cal
V}_1 \right]^m_{\mu_ i} \left(X^{(c)}_{d-1} \right)\ ,\ \ {\cal F}(z)={\cal K}
\left(z,z_0 \right)\in \left[{\cal V}_1 \right]^m_\mu \left(X^{(c)}_{d-1}
\right)\ . \hfill \cr\displaystyle  & \displaystyle  & \displaystyle
\left(\mu =\mu_ 1+\mu_ 2 \right) \hfill \cr} $$
\par
Then, the corresponding transforms $ \tilde  F_i = {\cal L}_d \left({\cal F}_i
\right), $ $ \tilde  F = {\cal L}_d({\cal F}) $ (such that $ \tilde  F_i\in
\tilde { \underline{{\cal O}} }^{\mu_ i,p(d)} \left(\Bbb C^{(m)}_+ \right), $
$ \tilde  F\in \tilde { \underline{{\cal O}} }^{\mu ,p(d)} \left(\Bbb
C^{(m)}_+ \right)) $ satisfy the relation
$ \tilde  F(\lambda)  = \tilde  F_1(\lambda) \cdot\tilde  F_2(\lambda) , $
which interpolates in $ \bar {\Bbb C}^{(m)}_+ $ the corresponding set of
relations $ [ {\bf f}]_{\ell}  = \left[ {\bf f}_1 \right]_{\ell} \cdot \left[
{\bf f}_2 \right]_{\ell} $ for the Fourier coefficients of $ {\bf F} , $ $
{\bf F}_1, $ $ {\bf F}_2 $ $ ( {\bf F}_i(z)= {\bf K}_i \left(z,z_0 \right), $
$ {\bf F}(z) = {\bf K} \left(z,z_0 \right)). $ \par
\medskip
{\bf Proof :} \par
Since $ {\cal L}_d \left({\cal F}_i \right) = L_d \left(F_i \right) $ and $
{\cal L}_d({\cal F}) = L_d(F) $ (with $ F_i(z) = K_i \left(z,z_0 \right) $ and
$ F(z) = K \left(z,z_0 \right)), $
and since the assumptions $ {\cal F}_i\in \left[{\cal V}_1 \right]^m_{\mu_ i}
\left(X^{(c)}_{d-1} \right), $ $ {\cal F}\in \left[{\cal V}_1 \right]^m_\mu
\left(X^{(c)}_{d-1} \right) $ imply that $ K_i\in V^{(m)}_{\mu_ i}
\left(X_{d-1} \right) $ and $ K\in V^{(m)}_\mu \left(X_{d-1} \right), $
the factorization property $ \tilde  F(\lambda)  = \tilde  F_1(\lambda)
\cdot\tilde  F_2(\lambda) $ (in $ \bar {\Bbb C}^{(m)}_+) $ follows directly
from
proposition 21. The fact that $ {\cal L}_d $ defines a morphism of normed
algebras is then an
immediate consequence of theorem 1 (property i)), while the
interpolation property for the relations $ [ {\bf f}]_{\ell}  = \left[ {\bf
f}_1 \right]_{\ell} \cdot \left[ {\bf f}_2 \right]_{\ell} $ results from
property ii) of
theorem 1.
\par
\smallskip
{\bf 4.4 The \lq\lq theorem F.G.\rq\rq , part II: the inverse problem } \par
\smallskip
We shall now establish the following properties which can be regarded as the
converse of theorem 1. \par
\medskip
{\sl Invariant kernels on the sphere whose coefficients of
the Fourier-Legendre expansion admit an analytic interpolation in a half plane
\/} \par
\smallskip
{\bf Theorem 3} {\bf :
}Let $ {\bf F}(z) $ be an $ {\bf H} $-invariant function on the sphere $
S_{d-1} $ whose set of Fourier-Legendre coefficients $ [ {\bf f}]_{\ell} $
(see Eq.(4.12)) admit an analytic interpolation $ \tilde  F(\lambda) $ in $
\Bbb C^{(m)}_+ $ such that: \par
a) $ m\in \Bbb R\backslash \Bbb Z, $ $ m > E \left[-{d-2 \over 2} \right], $
and $ \tilde  F \in  \tilde {\cal O}^{\mu ,\hat  p(d)} \left(\Bbb C^{(m)}_+
\right) $ (with $ \hat  p(d) $ specified
in proposition 16). \par
\smallskip
b) $ \forall \ell , $ $ \ell \in \Bbb N, $ $ \ell >m, $ $ \tilde  F(\ell)  = [
{\bf f}]_{\ell} . $ \par
Then $ {\bf F} $ admits an analytic continuation $ {\cal F} $ in the domain $
D_\mu $ such that: \par
\smallskip
i) $ {\cal F }\in  \left[{\cal V}_{\infty} \right]^m_\mu \left(X^{(c)}_{d-1}
\right) $ \par
ii) $ \tilde  F_{ \left\vert \Bbb C^{ \left(m^{\prime} \right)}_+ \right.} =
{\cal L}_d({\cal F}), $ where $ m^{\prime} = {\rm max}(m,-1). $ \par
\smallskip
For completeness, this theorem can be supplemented by the following one: \par
\smallskip
{\bf Theorem 3' :} If the assumptions of theorem 3 are satisfied, with the
following
conditions (replacing a) and b)): \par
\smallskip
a') $ m \in  \Bbb Z, $ $ m > -{d-2 \over 2} $ and $ \tilde  F \in  \tilde
{\cal O}^{\mu ,\hat  p(d)}_\ast \left(\Bbb C^{(m)}_+ \right) $ \par
b') $ \forall \ell , $ $ \ell \in \Bbb N, $ $ \ell \geq m, $ $ \tilde
F(\ell)  = [ {\bf f}]_{\ell} , $ \par
\smallskip
\noindent then the conclusions of theorem 3 remain valid, with the additional
specification: \par
\smallskip
i') $ {\cal F }\in  \left[{\cal V}^\ast_{ \infty} \right]^m_\mu
\left(X^{(c)}_{d-1} \right). $ \par
\smallskip
{\bf Proof of theorem 3} : Let us introduce, as in \S 3.5, the transform $
\hat {\bf F}= = \underline{{\goth R}}^{(c)}_d( {\bf F}) $
of the function $ {\bf F}(z) \equiv  {\bf K} \left(z,z_0 \right); $ according
to proposition 11, the function $ \hat {\bf F} \left(z_t \right)\equiv\hat
{\bf f}(t) $
satisfies the functional relation (3.30), already considered in \S 2.2 under
the name of
\lq\lq symmetry condition $ \left( {\bf S}_d \right)\dq\dq $ (see Eq.(2.21));
moreover, in view of proposition 23, the
coefficients $ \left[\hat {\bf f} \right]_{\ell} $ $ (\ell \in \Bbb Z) $ of
the Fourier series of $ \hat {\bf f}(t) $ coincide, for $ \ell \geq 0, $ with
the corresponding coefficients $ [ {\bf f}]_{\ell} $ of the Fourier-Legendre
expansion of $ {\bf F} . $
Therefore, the function $ \tilde  F_{ \left\vert \Bbb C^{ \left(m^{\prime}
\right)}_+ \right.} $ can be seen as an analytic interpolation in $ \Bbb C^{
\left(m^{\prime} \right)}_+ $
of the set of Fourier coefficients $ \left[\hat {\bf f} \right]_{\ell} $ of $
\hat {\bf f}(t) $ satisfying the conditions of proposition 5
(since one has $ G^{(\mu)}_ m \left[\tilde  F \right]\in L^1(\Bbb R) $ and
therefore $ G^{(\mu)}_{ m^{\prime}} \left[\tilde  F \right]\in L^1(\Bbb R)); $
according to the
latter, there exists an analytic function $ \hat  f(\tau) $ $ (\equiv\hat
{\cal F} \left(z_\tau \right)) $ in $ \underline{{\cal O}}^{(d)}_{m^{\prime}}
\left(\dot {\cal J}^{(\mu)} \right) $ such that $ \hat  f_{\vert \dot {\Bbb
R}}=\hat {\bf f} $
and $ \tilde  F_{ \left\vert \Bbb C^{ \left(m^{\prime} \right)}_+ \right.} =
{\cal L} \left(\hat  f \right). $ It follows that the function $ \hat {\bf F}=
\underline{{\goth R}}^{(c)}_d( {\bf F}) $ satisfies the
conditions of proposition 16 i), which implies that $ {\bf F} $ admits an
analytic continuation $ {\cal F} $ in
the domain $ D_\mu $ of $ X^{(c)}_{d-1} $ and that $ \hat  f={\cal
R}^{(c)}_d({\cal F}). $ By construction, one then has:
$$ \tilde  F_{ \left\vert \Bbb C^{ \left(m^{\prime} \right)}_+ \right.}=
\left[{\cal L}\circ{\cal R}^{(c)}_d \right]({\cal F}) = {\cal L}_d({\cal F})\
. $$
\par
In order to show that $ {\cal F}\in \left[{\cal V}_{\infty} \right]^m_\mu
\left(X^{(c)}_{d-1} \right), $ we introduce the following auxiliary analytic
function $ \hat  f^{\prime}( \tau) $ (different from $ \hat  f $ only if $
m<-1): $
$$ \hat  f^{\prime}( \tau)  = \hat  f(\tau)  + {1 \over \pi} {\rm e}^{i
\left({d-2 \over 2} \right)(\tau -\pi)} \sum^{ }_{ m<\ell \leq -1}
\left(\tilde  F(\ell) - \left[\hat {\bf f} \right]_{\ell} \right) {\rm cos}
\left[ \left(\ell +{d-2 \over 2} \right)\tau - \left({d-2 \over 2} \right)\pi
\right]\ . \eqno (4.22^\prime ) $$
\par
It is clear that $ \hat {\bf f}^{\prime} \equiv\hat  f^{\prime}_{\vert \Bbb R}
$ satisfies like $ \hat {\bf f} $ the symmetry condition $ \left( {\bf S}_d
\right) $ and that
its Fourier coefficients $ \left[\hat {\bf f}^{\prime} \right]_{\ell} $ (such
that $ \left[\hat {\bf f}^{\prime} \right]_{\ell}  = \left[\hat {\bf f}
\right]_{\ell} $ for $ \ell \geq 0) $ satisfy the set of
relations:
$$ \forall \ell ,\ \ \ \ \ell >m,\ \ \ \ \ \left[\hat {\bf f}^{\prime}
\right]_{\ell}  = \tilde  F(\ell) \ . $$
\par
It follows that proposition 5 can be applied to the function $ \hat {\bf
f}^{\prime} , $ with the result
that $ \tilde  F={\cal L} \left(\hat  f^{\prime} \right) $ (considered in its
full analyticity domain $ \Bbb C^{(m)}_+), $ and that $ \hat  f^{\prime} \in
\underline{{\cal O}}^{(d)}_m \left(\dot {\cal J}^{(\mu)} \right). $
Moreover, by applying similarly proposition 2 to the successive
derivatives $ \left(D_\alpha\hat {\bf f}^{\prime} \right)(t) $ with $ \alpha
\leq\hat  p(d), $ we can even say that the function $ \hat  f^{\prime}( \tau)
$ belongs to the
subspace $ \underline{{\cal O}}^{(d)}_{m,\hat  p(d)} \left(\dot {\cal
J}^{(\mu)} \right) $ (since the Laplace transforms $ (-i\lambda)^
\alpha\tilde  F(\lambda) $ of the respective
functions $ \left(D_\alpha\hat  f^{\prime} \right)(\tau) $ satisfy the
conditions $ G^{(\mu)}_ m \left[\lambda^ \alpha\tilde  F \right]\in L^1(\Bbb
R)). $ \par
We can also assert (in view of proposition 14) that the function $ \hat
f^{\prime}( \tau) -\hat  f(\tau) $ belongs
to $ \Xi^{( d)}, $ so that one has (in view of proposition 15): $ {\cal
F}={\cal X}^{(c)}_d \left(\hat  f \right)={\cal X}^{(c)}_d \left(\hat
f^{\prime} \right). $ We are thus in a
position to apply proposition 16 iii), which implies that $ {\cal F} $ belongs
to $ \left[{\cal V}_{\infty} \right]^m_\mu \left(X^{(c)}_{d-1} \right). $ \par
The proof of theorem 3' is completely similar, up to the following changes: in
view of the additional conditions a') and b') on the function $ \tilde  F, $
proposition 5' (instead
of proposition 5) can be applied and implies that $ \hat  f^{\prime} \in
\underline{{\cal O}}^{(d)\ast}_ m \left(\dot {\cal J}^{(\mu)} \right); $
moreover,
proposition 2' (instead of proposition 2) can be applied to the functions $
D_\alpha\hat {\bf f}^{\prime}( t), $ $ 0 \leq \alpha \leq\hat  p(d) $
(since the corresponding functions $ \lambda^ \alpha\tilde  F(\lambda) $
satisfy conditions $\alpha$), $\beta$), $\gamma$) of proposition 2'), which
implies (via formula (2.8')) that $ \hat  f^{\prime} \in \underline{{\cal
O}}^{(d)\ast}_{ m,\hat  p(d)} \left(\dot {\cal J}^{(\mu)} \right). $ Finally,
in view of proposition
16 iii), we conclude that $ {\cal F }= {\cal X}^{(c)}_d \left(\hat
f^{\prime} \right) $ belongs to the subspace $ \left[{\cal V}^\ast_{ \infty}
\right]^m_\mu \left(X^{(c)}_{d-1} \right). $ \par
\medskip
{\sl Inversion formulae\/} \par
We are now going to derive two types of inversion formulae for the
transformation $ {\cal L}_d $ which involve either the elementary perikernels
$ \Psi^{( d)}_\lambda( {\rm cos} \ \theta) $ introduced
in \S 3.8 or the first-kind Legendre functions $ P^{(d)}_\lambda( {\rm
cos}(\theta -\pi)) $ (see also \S 3.8 and our
geometrical interpretation below in \S 4.5). \par
\medskip
{\bf Theorem 4 :} Let $ ({\cal F}, {\bf F} ,F) $ be an $ H^{(c)} $-invariant
triplet, $ {\cal F} $ being of moderate growth,
whose transforms $ \tilde  F={\cal L}_d({\cal F}) $ and
$ \left\{[ {\bf f}]_{\ell} ; \right. $ $ \ell  \in  \Bbb N\}  = {\bf L}_d(
{\bf F}) $ satisfy either the conditions a), b) of theorem 3 (for $ m $
non-integer, $ m>E \left(-{d-2 \over 2} \right)) $ or the conditions a'), b')
of theorem 3' (for $ m $ integer, $ m >-{d-2 \over 2}). $ Then, this triplet
can be expressed in terms of $ \tilde  F $ by the following (double) system of
inversion formulae (up
to the restriction of our remark iii) below): \par
\smallskip
a) $ \forall z\in X^{(c)}_{d-1}, $ $ z^{(1)} = {\rm cos} \ \theta  \in
\underline{D}_0 = \Bbb C\backslash[ 1,+\infty[ , $
$$ {\cal F}(z) = \underline{f}( {\rm cos} \ \theta)  = - {1 \over 2\pi}
\int^{ +\infty}_{ -\infty}{\tilde  F(m+i\nu) \Psi^{( d)}_{m+i\nu}( {\rm cos} \
\theta) \over {\rm sin} \ \pi( m+i\nu)}  {\rm d} \nu  + \goth P^{(d)}_m( {\bf
F})( {\rm cos} \ \theta) \eqno (4.23) $$
or:
$$ {\cal F}(z) = - {1 \over 2\omega_ d} \int^{ +\infty}_{ -\infty}{\tilde
F(m+i\nu) h_d(m+i\nu) P^{(d)}_{m+i\nu}( {\rm cos}(\theta -\pi)) \over {\rm
sin} \ \pi( m+i\nu)}  {\rm d} \nu  + \goth P^{(d)}_m( {\bf F})( {\rm cos} \
\theta) \eqno (4.23^\prime ) $$
\par
In the case $ m\in \Bbb Z, $ the integrals at the r.h.s. of formulae (4.23)
(4.23'), must
be understood in the sense of distributions with $ {1 \over {\rm sin} \ \pi(
m+i\nu)}  = \dlowlim{ {\rm lim}}{ \doublelow{ \varepsilon \longrightarrow 0
\cr \varepsilon >0 \cr}}  {1 \over {\rm sin} \ \pi( m-\varepsilon +i\nu)} . $
\par
\smallskip
b) $ \forall z \in  S_{d-1}, $ $ z^{(1)}= {\rm cos} \ u \in  [-1,+1], $
$$ {\bf F}(z) = \underline{{\bf f}}( {\rm cos} \ u) = {1 \over \pi}  \sum^{
}_{ \ell \in \Bbb N}(-1)^{\ell}[ {\bf f}]_{\ell} \Psi^{( d)}_{\ell}( {\rm cos}
\ u) \eqno (4.24) $$
or:
$$ {\bf F}(z) = {1 \over \omega_ d} \sum^{ }_{ \ell \in \Bbb N}[ {\bf
f}]_{\ell} h_d(\ell) P^{(d)}_{\ell}( {\rm cos} \ u) \eqno (4.24^\prime ) $$
\par
c) $ \forall z \in  X_{d-1}, $ $ z^{(0)} \geq  0, $ $ z^{(1)} = {\rm cosh} \ v
\geq  1, $
$$ F(z) = \underline{{\rm f}}( {\rm cosh} \ v) = {1 \over \pi}  \int^{
+\infty}_{ -\infty}\tilde  F(m+i\nu)  {\Delta \Psi^{( d)}_{m+i\nu}( {\rm cosh}
\ v) \over {\rm sin} \ \pi( m+i\nu)}  {\rm d} \nu \eqno (4.25) $$
or:
$$ F(z) = {1 \over \omega_ d} \int^{ +\infty}_{ -\infty}\tilde  F(m+i\nu)
h_d(m+i\nu) P^{(d)}_{m+i\nu}( {\rm cosh} \ v) {\rm d} \nu \eqno (4.25^\prime )
$$
where we have put:
$$ \goth P^{(d)}_m( {\bf F})( {\rm cos} \ \theta)  = {1 \over \omega_ d}
\sum^{ }_{ 0\leq \ell <m}[ {\bf f}]_{\ell} h_d(\ell) \ P^{(d)}_{\ell}( {\rm
cos} \ \theta)  \eqno (4.26) $$
\par
\smallskip
{\bf Remarks :} \par
i) For $ m $ such that $ E \left(-{d-2 \over 2} \right) < m \leq  0, $ it is
understood that $ \goth P^{(d)}_m( {\bf F})\equiv 0. $ \par
ii) for $ d $ even, the two inversion formulae (4.23) and (4.23') are
identical in view
of Eq.(3.90). \par
iii) for $ d $ odd, the inversion formulae (4.23), (4.25) will be proved under
the
precise condition $ m>-{d-1 \over 2} $ assumed in theorem 3; the equivalence
of formulae (4.23'),
(4.25') with the latter will be established under the somewhat restrictive
condition $ m> {\rm min} \left(-1,-{d-3 \over 2} \right). $ \par
\smallskip
{\bf Proof\nobreak\ :} \par
a) According to the proof of theorem 3, we can assert that $ {\cal F }= {\cal
X}^{(c)}_d \left(\hat  f^{\prime} \right), $ where the
function $ \hat  f^{\prime} , $ defined by Eq.(4.22'), is given in terms of $
\tilde  F $ by proposition 5 (or 5' if $ m \in  \Bbb Z), $
namely by Eq.(2.30) (with $ f $ replaced by $ \hat  f^{\prime} ). $ \par
The expression of the transformation $ {\cal X}^{(c)}_d $ given in proposition
14 then allows one
to compute $ {\cal F} $ in terms of $ \hat  f^{\prime} $ by Eqs.(3.46) or
(3.47) according to whether $ d $ is even or odd.
The
corresponding differential (resp. integro-differential) operator can now be
applied
directly to the factor $ {\rm cos} \left[ \left(m+{d-2 \over 2} + i\nu
\right)(\theta -\pi) \right] $ under the integration sign of Eq.(2.30); as a
result, this factor is replaced (in view of Eqs.(3.78) and (3.79)) by the
integration
kernel $ \Psi^{( d)}_{m+i\nu}( {\rm cos} \ \theta) , $ which yields the
integral at the r.h.s. of Eq.(4.23). This procedure
is fully justified since, in view of the inequality (3.80) (see proposition
18) and of the
assumption that $ \tilde  F\in\tilde {\cal O}^{\mu ,\hat  p(d)} \left(\Bbb
C^{(m)}_+ \right) $ (i.e. condition a) of theorem 3), the considered
integral is absolutely and uniformly convergent: it is in fact majorized by:
$$ C_d( {\rm cos} \ \theta) {\rm e}^{m\vert {\rm Im} \ \theta\vert} \int^{
+\infty}_{ -\infty} \left\vert\tilde  F(m+i\nu) \right\vert  { {\rm cosh} \
\pi \nu \over\vert {\rm sin} \ \pi( m+i\nu)\vert} {\rm d} \nu \leq C_d( {\rm
cos} \ \theta) {\rm e}^{m\vert {\rm Im} \ \theta\vert} \mid \mid \mid\tilde
F\mid \mid \mid^{ m\mu\hat  p(d)}_1 $$
The last term at the r.h.s. of Eq.(4.23) results similarly from the action of
the
transformation $ {\cal X}^{(c)}_d $ (via Eqs.(3.46), (3.47)) on the last term
$ {\cal P}^{(d)}_m \left[\hat {\bf F} \right](\theta) $ at the r.h.s. of
Eq.(2.30). In fact, in the expression (2.31) of $ {\cal P}^{(d)}_m \left[\hat
{\bf f} \right](\theta) , $ the only part of the sum which
survives to the action of $ {\cal X}^{(c)}_d $ is
$$ \matrix{\displaystyle 2{\cal X}^{(c)}_d \left[ \sum^{ }_{ 0\leq \ell
<m}(-1)^{\ell} \left[\hat {\bf f} \right]_{\ell} {\rm cos} \left[ \left(\ell
+{d-2 \over 2} \right)(\theta -\pi) \right] {\rm e}^{-i \left({d-2 \over 2}
\right)(\theta -\pi)} \right] \cr\displaystyle = {1 \over \pi}  \sum^{ }_{
0\leq \ell <m}(-1)^{\ell} \left[\hat {\bf f} \right]_{\ell} \Psi^{(
d)}_{\ell}( {\rm cos} \ \theta) \cr} $$
since $ \Psi^{( d)}_{\ell} \equiv 0 $ for $ -d-3 \leq  \ell  \leq  -1 $ (see
in \S 3.8 our remark ii) after Eq.(3.91)). Now, in view
of the equalities $ \left\{ \left[\hat {\bf f} \right]_{\ell} = \left[ {\bf
f}_{\ell} \right];\ \ell \in \Bbb N \right\} $ (or $ {\bf L} \left(\hat {\bf
f} \right)= {\bf L}_d( {\bf f}), $ see proposition 23) and of the
expression of $ \Psi^{( d)}_{\ell} $ $ (\ell \in \Bbb N) $ given in \S 3.8
(see Eqs.(3.90), (3.91)), we obtain the last
term $ \goth P^{(d)}_m[ {\bf F}]( {\rm cos} \ \theta) $ of Eq.(4.23), as it is
specified in Eq.(4.26). \par
We shall now prove that (for $ d $ odd), the formula (4.23') is equivalent to
Eq.(4.23), although $ \Psi^{( d)}_\lambda( {\rm cos} \ \theta) $ is (according
to Eq.(3.102)) a linear combination of $ P^{(d)}_\lambda( {\rm cos}(\theta
-\pi)) , $
$ Q^{(d)}_\lambda( {\rm cos}(\theta -\pi)) $ and $ R^{(d)}_\lambda( {\rm cos}
\ \theta) . $ Firstly, it follows from proposition 19 i) that the
integral at the r.h.s. of Eq.(4.23') is absolutely convergent; on the other
hand, we can
assert (provided $ m> {\rm min} \left(-1,-{d-3 \over 2} \right)) $ that the
residual integral
$$ \int^{ }_{ {\rm L}_m}\tilde  F(\lambda)  \left\{{ \Psi^{( d)}_\lambda(
{\rm cos} \ \theta)  - {\pi \over \omega_ d} h_d(\lambda)  P^{(d)}_\lambda(
{\rm cos}(\theta -\pi)) \over {\rm sin} \ \pi \lambda} \right\} {\rm d}
\lambda $$
vanishes in view of the Cauchy theorem, since its integrand is
a holomorphic function in $ \Bbb C^{(m)}_+ $ (namely the linear combination
of $ h_d(\lambda)  Q^{(d)}_\lambda $ and $ R^{(d)}_\lambda $ given in
Eq.(3.102)) which, in view of the
corollary of proposition 19 (see \S 3.8), can be majorized in $ \Bbb C^{(m)}_+
$
by: \par
\noindent$ {\rm Cst\ e}^{-\mu( {\rm Re} \ \lambda -m)}G^{(\mu)}_ m
\left[\lambda^{{ d-1 \over 2}}\tilde  F \right]. $ \par
\smallskip
b) Since $ {\bf f}(u) $ satisfies the symmetry condition $ \left( {\bf S}_d
\right) $ and
therefore (see proposition 4, Eq.(2.24)) admits the Fourier
expansion:
$$ \hat {\bf f}(u) {\bf  = }{1 \over 2\pi}  {\rm e}^{i \left({d-2 \over 2}
\right)(u-\pi)} \sum^{ +\infty}_{ \ell =-\infty}( -1)^{\ell} \left[\hat {\bf
f} \right]_{\ell} {\rm cos} \left[ \left(\ell +{d-2 \over 2} \right)(u-\pi)
\right]\ , $$
the last argument of a) can be extended to this series, namely we
have (in view of propositions 17, 22 and of Eqs.(3.90), (3.91))
$$ {\bf F}  = \underline{{\goth X}}^{(c)}_d \left(\hat {\bf F} \right) = {1
\over \pi}  \sum^{ +\infty}_{ \ell =0}(-1)^{\ell} \left[\hat {\bf f}
\right]_{\ell} \Psi^{( d)}_{\ell}( {\rm cos} \ u) = {1 \over \omega_ d}
\sum^{ +\infty}_{ \ell =0}[ {\bf f}]_{\ell} h_d(\ell) P^{(d)}_{\ell}( {\rm
cos} \ u) $$
\par
The absolute convergence of this series follows (at least for $ \mu >0) $
from the following facts: \par
\smallskip
i) $ \forall \ell \geq 0, $ $ \left\vert \left[\hat {\bf f} \right]_{\ell}
\right\vert  \leq  {\rm Cst\ e}^{-\ell \mu} $ (see our remark i) in \S 2.1,
before
proposition 2); \par
ii) $ \left\vert \Psi^{( d)}_{\ell}( {\rm cos} \ u) \right\vert  \leq  {\rm
Cst}(1+\ell)^{\hat  p(d)}.
$ \par
\noindent In the present framework, we have thus reobtained a direct algebraic
proof
of the Fourier-Legendre expansion (4.24) of $ {\bf F} $ (see e.g.
$\lbrack$Vi$\rbrack$ and
references therein for the standard derivation based on group
representation theory). \par
\smallskip
c) By taking the discontinuities of both sides of Eq.(4.23) across
the cut $ \{ {\rm cos} \ \theta  \in  [1,+\infty[\} , $ we readily obtain the
inversion formula (4.25)
in which $ \Delta \Psi^{( d)}_{m+i\nu}( {\rm cosh} \ v) $ can be expressed via
Eq.(3.103) when $ d $ is even,
and via Eqs.(3.105) or (3.106) when $ d $ is odd. Note that the latter could
also be obtained from Eq.(2.27) for $ d $ even (resp. from Eq.(2.28) for $ d $
odd), by applying the inverse Abel operator $ {\cal X}_d $ to the
corresponding
expressions of $ \hat {\rm f}(v) ${\bf\ }and then recognizing the
representation (3.103) (resp.
(3.105)) of $ \Delta \Psi^{( d)}_{m+i\nu}( {\rm cosh} \ v). $ \par
On the other hand, one can start from the alternative form of $ \hat  f, $
namely the one given by Eq.(2.28) for $ d $ even (resp. (2.27) for $ d $ odd),
since both are equivalent in view of proposition 5. Then applying the
inverse Abel operator $ {\cal X}_d $ yields expressions in which one
recognizes the
representation (3.85) (resp. (3.86)) of the function $ h_d(m+i\nu)
P^{(d)}_{m+i\nu}( {\rm cosh} \ v); $
this gives a first (direct) derivation of Eq.(4.25'). The alternative
derivation, which starts from Eq.(4.25) and then makes use of the
expression (3.106) of $ \Delta \Psi^{( d)}_\lambda $ in terms of $
P^{(d)}_\lambda , $ $ Q^{(d)}_\lambda , $ $ R^{(d)}_\lambda $ could be carried
out in the same way as our derivation of Eq.(4.23') from Eq.(4.23) above
(see a)), by using again the corollary of proposition 19. q.e.d. \par
\medskip
{\sl Continuity properties of the inverse transformation of \/}$ {\cal L}_d $
{\sl on
appropriate functional spaces\/} \par
Let us define, for each $ m $ such that $ m>E \left(- {d-2 \over 2} \right), $
the
transformation $ \left[{\cal L}_d \right]^{-1}_m = {\cal
X}^{(c)}_d\circ[{\cal L}]^{-1}_{(d,m)}, $ where
$ [{\cal L}]^{-1}_{(d,m)} $ is the \lq\lq inverse Fourier-Laplace
transformation\rq\rq\ (in one variable) introduced in \S 2.2. $ \left[{\cal
L}_d \right]^{-1}_m $ can be shown
to have continuity properties and to act as a
\lq\lq quasi-inverse\rq\rq\ of the Fourier-Laplace transformation $ {\cal
L}_d={\cal L}\circ{\cal R}^{(c)}_d $
in the sense specified below in theorem 5. \par
\medskip
{\bf Theorem 5 :} \par
a) For each $ m\in \Bbb R\backslash \Bbb Z, $ with $ m > E \left[-{d-2 \over
2} \right], $ the transformation
$ \left[{\cal L}_d \right]^{-1}_m $  defines a continuous mapping $ (\tilde  F
\longrightarrow{\cal F}= \left[{\cal L}_d \right]^{-1}_m \left(\tilde  F
\right)) $ from each subspace $ \tilde {\cal O}^{\mu ,\hat  p(d)} \left(\Bbb
C^{(m)}_+ \right) $
into the corresponding subspace $ \left[{\cal V}_{\infty} \right]^m_\mu
\left(X^{(c)}_{d-1} \right). $ \par
Moreover, for each $ \tilde  F\in\tilde {\cal O}^{\mu ,\hat  p(d)} \left(\Bbb
C^{(m)}_+ \right), $ the triplet $ ({\cal F}, {\bf F} ,F) $
associated with $ {\cal F}= \left[{\cal L}_d \right]^{-1}_m \left(\tilde  F
\right) $ can be computed via the formulae (4.23), (4.24),
(4.25) of theorem 4, in which the additional condition $ \goth P^{(d)}_m[
{\bf F}]=0 $ must
be imposed. \par
\smallskip
b) For each $ m\in \Bbb Z, $ with $ m > - {d-2 \over 2}, $ the transformation
$ \left[{\cal L}_d \right]^{-1}_m $
defines a continuous mapping from each subspace $ \tilde {\cal O}^{\mu ,\hat
p(d)}_\ast \left(\Bbb C^{(m)}_+ \right) $ into
the corresponding subspace $ \left[{\cal V}^\ast_{ \infty} \right]^m_\mu
\left(X^{(c)}_{d-1} \right) $ (the last statement of a) being
unchanged, provided the prescriptions given at the end of theorem 4 a) be now
applied). \par
\smallskip
c) For $ m>-1, $ let $ [{\cal V}]^m_\mu \left(X^{(c)}_{d-1} \right) $ be the
subspace of functions $ {\cal F} $ in $ \left[{\cal V}_1 \right]^m_\mu
\left(X^{(c)}_{d-1} \right)\cap \left[{\cal V}_{\infty} \right]^m_\mu
\left(X^{(c)}_{d-1} \right) $
satisfying the condition: $ \goth P^{(d)}_m({\cal F})=0. $ Then, if $ {\cal
F}\in[{\cal V}]^m_\mu \left(X^{(c)}_{d-1} \right), $ the
Laplace transform $ \tilde  F={\cal L}_d({\cal F}) $ is such that, $ \forall
r, $ $ 0 \leq  r \leq  p(d), $ $ \lambda^ r\times\tilde  F(\lambda) $ belongs
to the Hardy space $ {\bf H}^2 \left(\Bbb C^{(m)}_+ \right) $ and moreover,
one has:
$$ {\cal F }= \left[{\cal L}_d \right]^{-1}_m \left(\tilde  F \right)\ .
\eqno (4.27) $$
\par
\smallskip
{\bf Proof\nobreak\ :} \par
Since $ \left[{\cal L}_d \right]^{-1}_m = {\cal X}^{(c)}_d\circ[{\cal
L}]^{-1}_{(d,m)}, $ the continuity properties stated in
a) (resp. b)) are immediate consequences of: \par
i) the continuity of the mapping $ [{\cal L}]^{-1}_{d,m} $ from $ \tilde
{\cal O}^{\mu ,\hat  p(d)} \left(\Bbb C^{(m)}_+ \right) $
into $ \underline{{\cal O}}^{(d)}_{m,\hat  p(d)} \left(\dot {\cal J}^{(\mu)}
\right) $ (resp. from $ \tilde {\cal O}^{\mu ,\hat  p(d)}_\ast \left(\Bbb
C^{(m)}_+ \right) $ into $ \underline{{\cal O}}^{(d)\ast}_{ m,\hat  p(d)}
\left(\dot {\cal J}^{(\mu)} \right)) $
established in proposition 7 (resp. proposition 7'), \par
ii) the continuity of the mapping $ {\cal X}^{(c)}_d $ from $ \underline{{\cal
O}}^{(d)}_{m,\hat  p(d)} \left(\dot {\cal J}^{(\mu)} \right) $
into $ \left[{\cal V}_{\infty} \right]^m_\mu \left(X^{(c)}_{d-1} \right) $
(resp. from $ \underline{{\cal O}}^{(d)\ast}_{ m,\hat  p(d)} \left(\dot {\cal
J}^{(\mu)} \right) $ into $ \left[{\cal V}^\ast_{ \infty} \right]^m_\mu
\left(X^{(c)}_{d-1} \right)), $
established in proposition 16 iii). \par
\smallskip
The last statement of a) is readily obtained by noticing that the
expression of $ [{\cal L}]^{-1}_{(d,m)} \left(\tilde  F \right) $ is given by
Eq.(2.30) {\sl without\/} the term $ {\cal P}^{(d)}_m \left[\hat {\bf f}
\right](\theta) $
at the r.h.s. of the latter (see \S 2.2), which implies the absence of the
corresponding term $ \goth P^{(d)}_m[ {\bf F}] = {\cal X}^{(c)}_d \left({\cal
P}^{(d)}_m \left[\hat {\bf f} \right] \right) $ in the respective
expressions (4.23), (4.23') of the function $ {\cal F }= \left[{\cal L}_d
\right]^{-1}_m \left(\tilde  F \right). $ \par
For proving c), we first notice that if $ {\cal F}\in[{\cal V}]^m_\mu
\left(X^{(c)}_{d-1} \right), $ then the
function $ \hat  f={\cal R}^{(c)}_d({\cal F}) $ belongs (in view of
proposition 13) to $ {\cal O}^{(d)}_{m,p(d)} \left(\dot {\cal J}^{(\mu)}
\right)\cap \underline{{\cal O}}^{(d)}_{m,p(d)} \left(\dot {\cal J}^{(\mu)}
\right) $
and satisfies the additional conditions $ \left[\hat {\bf f} \right]_{\ell}
=0, $ for $ 0\leq \ell <m $ (in view of
proposition 23). \par
The corresponding \lq\lq Hardy space property\rq\rq\ of $ \tilde  F={\cal
L}_d({\cal F})={\cal L} \left(\hat  f \right) $ then
follows from our remark after proposition 6. \par
In order to prove Eq.(4.27), let us now introduce the auxiliary
function
$$ \hat  f^1(\tau)  = \hat  f(\tau)  - {1 \over \pi}  {\rm e}^{i \left({d-2
\over 2} \right)(\tau -\pi)} \sum^{ }_{ -{d-2 \over 2}\leq \ell \leq -1}
\left[\hat {\bf f} \right]_{\ell} {\rm cos} \left[ \left(\ell +{d-2 \over 2}
\right)\tau - \left({d-2 \over 2} \right)\pi \right], $$
which satisfies the conditions $ \left[\hat {\bf f}^1 \right]_{\ell} =0 $ for
all the integers $ \ell $ such that $ -m-d+2<\ell <m. $ \par
Since $ \hat  f-\hat  f^1 $ belongs to the kernel $ {\cal N}^{(d)}_m $ of $
{\cal L} $ in $ {\cal O}^{(d)}_m \left(\dot {\cal J}^{(\mu)} \right) $ (see
\S 2.2) and to the kernel $ \Xi^{( d)} $ of $ {\cal X}^{(c)}_d $ (see
proposition 14), one can then
write: \par
\smallskip
i) $ \tilde  F = {\cal L} \left(\hat  f \right) = {\cal L} \left(\hat  f^1
\right) $ which entails (since $ \hat  f^1\in  {\cal M}^{(d)}_m \left(\dot
{\cal J}^{(\mu)} \right)) $ that:
$$ \hat  f^1 = [{\cal L}]^{-1}_{(d,m)} \left(\tilde  F \right) \eqno (4.28) $$
\par
ii) (in view of proposition 15)
$$ {\cal F }= {\cal X}^{(c)}_d \left(\hat  f \right) = {\cal X}^{(c)}_d
\left(\hat  f^1 \right) \eqno (4.29) $$
Eqs.(4.28) and (4.29) then imply Eq.(4.27). \par
\medskip
{\bf 4.5 Complements on the representation of perikernels in
terms of the Legendre functions} $ P^{(d)}_\lambda $ \par
According to formula (4.23) of theorem 4, completed by
the statement of theorem 5, any general perikernel $ {\cal K} $ of moderate
growth on $ X^{(c)}_{d-1} $ whose representative $ {\cal F}(z) = {\cal K}
\left(z,z_0 \right) $ belongs to a
given class $ [{\cal V}]^m_\mu \left(X^{(c)}_{d-1} \right), $ with $ m>-1, $
is decomposable in terms of
the family of elementary perikernels $ \left\{{\cal K}^{(d)}_{m+i\nu} \right.
\left(z,z^{\prime} \right)=\Psi^{( d)}_{m+i\nu} \left(1+{ \left(z-z^{\prime}
\right)^2 \over 2} \right)\ ; $
$ \nu  \in  \Bbb R\} , $ introduced in \S 3.8 (see proposition 17). \par
It is then tempting to give a similar interpretation to
the alternative formula (4.23'), according to which any function
$ {\cal F} $ in a class $ [{\cal V}]^m_\mu \left(X^{(c)}_{d-1} \right) $ can
be expressed in terms of the
corresponding family of Legendre functions of the first kind $ \left\{
P^{(d)}_{m+i\nu} \ ;\ \nu \in \Bbb R \right\} . $ \par
In fact, each function $ P^{(d)}_{m+i\nu}( {\rm cos}(\theta -\pi)) $ is (as $
\Psi^{( d)}_{m+i\nu}( {\rm cos} \ \theta) ) $
an analytic function of $ {\rm cos} \ \theta $ in the cut-plane $
\underline{D}_0, $ which is of
moderate growth of order $ m $ (for $ m \geq  - {d-2 \over 2}) $ (see
proposition 19
and compare with proposition 18). One is thus led to introduce
the following family of invariant perikernels $ {\cal K}^{\ast( d)}_\lambda ,
$ similar to $ {\cal K}^{(d)}_\lambda $
$ (\lambda \in \Bbb C): $
$$ {\cal K}^{\ast( d)}_\lambda \left(z,z^{\prime} \right) = {\pi \over
\omega_ d} h_d(\lambda)  P^{(d)}_\lambda \left(-1 - { \left(z-z^{\prime}
\right)^2 \over 2} \right) \eqno (4.30) $$
\par
In order to complete the parallel between the
representations (4.23) and (4.23'), one would then like to give a
direct geometrical interpretation of the latter similar to that
of the former (since the derivation of Eq.(4.23') from Eq.(4.23)
given in theorem 4 purely relies on the algebraic relations of
\S 3.8 and does not shed light on the underlying geometrical
aspects). \par
While Eq.(4.23) was directly interpretable (e.g. for $ {\cal F }\in
[V]^m_\mu \left(X^{(c)}_{d-1} \right)) $
in terms of the inversion formula $ \left[{\cal L}_d \right]^{-1}_m = {\cal
X}^{(c)}_d \circ[{\cal L}]^{-1}_{d,m} $ of $ {\cal L}_d={\cal L}\circ{\cal
R}^{(c)}_d, $
Eq.(4.23') requires the introduction of an alternative complex
Radon-Abel transformation $ {\cal R}^{\ast( c)}_d $ and of its inverse $
{\cal X}^{\ast( c)}_d, $ such
that $ {\cal L}_d $ and $ \left[{\cal L}_d \right]^{-1}_m $ also admit the
following decompositions (to be
specified below):
$$ {\cal L}_d = {\cal L }\circ  {\cal R}^{\ast( c)}_d \eqno (4.31) $$
and
$$ \left[{\cal L}_d \right]^{-1}_m = {\cal X}^{\ast( c)}_d \circ[{\cal
L}]^{\ast -1}_{(d,m)} \eqno (4.32) $$
\par
\smallskip
{\sl The \/}$ {\rm \ast} ${\sl -Radon-Abel transformation \/} $ \goth
R^{\ast( c)}_d $ (or $ {\cal R}^{\ast( c)}_d) $ \par
One defines a transformation $ {\cal F} \longrightarrow\hat {\cal F}^\ast
=\goth R^{\ast( c)}_d({\cal F}) $ by a
formula similar to Eq.(3.22) in which the integration cycle $ h_\tau $ is
replaced by another cycle $ h^\ast_ \tau $ specified below (the integrand
being
the same). The function $ \hat {\cal F}^\ast \left(z_\tau \right) $ (analytic
in the same domain $ \hat  D $ of
$ \hat  X^{(c)} $ as $ \hat {\cal F} \left(z_\tau \right)) $ is represented by
a function $ \hat  f^\ast( \tau)  = {\cal R}^{\ast( c)}_d({\cal F}), $
defined as an Abel-type transform similar to $ \hat  f(\tau) $ (see
Eq.(3.25)), except that the path $ \gamma_ \tau $ is now replaced by a path $
\gamma^ \ast_ \tau $
with support:
$$ \underline{\gamma}^ \ast_ \tau =\{ \theta =\theta( \lambda) \ ;\ {\rm cos}
\ \theta( \lambda) +1 = \lambda( {\rm cos} \ \tau +1)\ ;\ 0\leq \lambda \leq
1,\ \theta( 0)=\varepsilon( \tau) \pi \ ,\ \theta( 1)=\tau\} $$
where $ \varepsilon( \tau)  = {\rm sgn} \ {\rm Re} \ \tau $ and $ \vert {\rm
Re} \ \tau\vert  \leq  2\pi . $ \par
A more detailed geometrical description (on $ X^{(c)}_{d-1}) $ can be
given as follows. Let $ {\cal C}^\ast = \left\{ z\in X^{(c)}_{d-1}\ ;\
z^{(d-1)}=-1 \right\} $ and for each $ \tau \in \Bbb C, $
let $ {\cal C}^\ast_ \tau  = {\cal C}^\ast \cap \Pi_ \tau $ $ (\Pi_ \tau $
being the complex horosphere introduced in
\S 3.2): $ {\cal C}^\ast_ \tau $ is a $ (d-3) $-dimensional complex sphere
parametrized by $ z=z \left(\vec  \zeta ,\tau \right) $
(see Eqs.(3.1)), with $ \vec  \zeta^ 2=2\ {\rm e}^{i\tau}( {\rm cos} \ \tau
+1). $ We then define the
horocycle $ h^\ast_ \tau $ as a representative of an element of $ H_{d-2} $ $
(\Pi_ \tau , $ $ {\cal C}^\ast_ \tau ) $
whose support is:
$$ \underline{h}^\ast_ \tau = \left\{ z\in \Pi_ \tau \ ;\ z=z \left(\rho\vec
\alpha ,\tau \right)\ ,\ \vec  \alpha  \in  \Bbb S^{d-3},\ \rho = \left[2\
{\rm e}^{i\tau}( 1-\lambda)( {\rm cos} \ \tau +1) \right]^{1/2},\ 0\leq
\lambda \leq 1 \right\} \ , $$
and whose orientation is obtained by continuity from that of the \lq\lq real
horocycle\rq\rq\ $ h^{\ast +}_w \equiv  h^\ast_{ iw}, $ the latter being
itself represented by the
(oriented) ball
$$ B^\ast_{ iw} = \left\{\vec  \zeta \in \Bbb R^{d-2}\ ;\ \vec  \zeta
=\rho\vec  \alpha \ ,\ \vec  \alpha \in \Bbb S^{d-3}\ ,\ 0\leq \rho \leq
\left[2( {\rm cosh} \ w+1) {\rm e}^{-w} \right]^{1/2} \right\} \ . $$
(Note that $ \underline{h}^{\ast +}_w \supset  \underline{h}^+_w\ ; $ see
Fig.5b)). \par
For all domains $ D $ of $ X^{(c)}_{d-1} $ in the class $ {\cal D} $ (see \S
3.1) whose
projection $ \underline{D}  = \varpi( D) $ in the $ z^{(d-1)} $-plane is also
{\sl star-shaped with respect
to the point\/} $ z^{(d-1)}\equiv {\rm cos} \ \theta  = -1, $ we then have the
following property,
similar to lemma 3 (see \S 3.4): $ \varpi \left( \bigcup^{ }_{ \tau \in\dot
{\cal J}(D)} \underline{h}^\ast_ \tau \right) = \underline{D} $ (since $
\varpi \left( \underline{h}^\ast_ \tau \right) $ is the linear
segment with end-points $ (-1,\ {\rm cos} \ \tau) $ in $ \Bbb C); $
this is valid in particular for all domains $ D_\mu . $ \par
\medskip
The following properties of the transformation $ {\cal R}^{\ast( c)}_d $ hold:
\par
\smallskip
{\bf Proposition 24 :} \par
i) The function $ \hat  f^\ast( \tau)  = \hat {\cal F}^\ast \left(z_\tau
\right) $ is analytic in the $ 2\pi $-periodic
cut-plane $ \dot {\cal J}^{(0)} $ and satisfies the symmetry property:
$$ \left(S^\ast_ d \right)\ \ \ \ \ \ \hat  f^\ast( \tau)  = {\rm
e}^{i(d-2)\tau}\hat  f^\ast( -\tau) \eqno (4.33) $$
\par
ii) For $ w \in  \Bbb R^+, $
$$ \Delta\hat  f^\ast( w) = \Delta\hat  f(w) = {\cal R}_d(\Delta{\cal F})
\eqno (4.34) $$
The proof of (4.34) follows from the fact that for $ \tau =iw, $ $
\Delta{\cal F}_{ \left\vert \underline{h}^\ast_{ iw} \right.} $ has its
support in the subset $ \underline{h}^+_w $ of $ \underline{h}^{\ast +}_w, $
or equivalently that $ \underline{h}^{\ast +}_w $ is represented
in the $ \theta $-plane by the two equivalent paths $ \gamma^{ \ast
\varepsilon}_{ iw}= \left( \charlvupup{ \rightarrowfill}{ \varepsilon \pi ,0}
\right)+ \left( \charlvupup{ \rightarrowfill}{ 0,iw} \right) $ (where $
\varepsilon =+ $ or
$ -, $ according to the sign of $ {\rm Re} \ \tau $ from which the limit of $
\gamma^ \ast_ \tau $ is taken),
the corresponding contributions of $ \left( \charlvupup{ \rightarrowfill}{
+\pi ,0} \right) $ and $ \left( \charlvupup{ \rightarrowfill}{ -\pi ,0}
\right) $ to the Abel
transforms $ \left[{\cal R}^{\ast( c)}_d({\cal F}) \right]_{\pm}( w) $ being
equal. \par
The inversion of the transformation $ {\cal R}^{\ast( c)}_d $ is performed via
a
transformation $ {\cal X}^{\ast( c)}_d, $ similar to $ {\cal X}^{(c)}_d $ (see
\S 3.7), with the following main
results: \par
\smallskip
a) $ d $ {\sl even\/}: $ f(\theta)  = {\cal F}(z) $ $ (z^{(d-1)}= {\rm cos} \
\theta ) $ is expressed in terms of $ \hat  f^\ast( \theta) $
again by formula (3.46) (with $ \hat  f $ replaced by $ \hat  f^\ast ); $ it
follows that $ {\rm e}^{-i \left({d-2 \over 2} \right)\theta} \left(\hat
f^\ast( \theta) -\hat  f(\theta) \right) $
is a polynomial of degree $ {d-4 \over 2} $ of $ {\rm cos} \ \theta . $ We
conclude that $ {\cal R}^{\ast( c)}_d({\cal F})-{\cal R}^{(c)}_d({\cal F}) $
belongs to the kernel $ \Xi^{( d)} $ of $ {\cal X}^{(c)}_d: $ the
transformation $ {\cal R}^{\ast( c)}_d $ reduces
trivially to $ {\cal R}^{(c)}_d({\cal F}) $ (note that the symmetry conditions
$ \left(S_d \right) $ and $ \left(S^\ast_ d \right) $ are
identical in this case). \par
\smallskip
b) $ d $ {\sl odd\/}: $ f(\theta)  = {\cal F}(z) $ is expressed in terms of $
\hat  f^\ast( \theta) $ by the
following formula (similar to Eq.(3.47), but different from it):
$$ f(\theta)  = -2 \left(-{1 \over 2\pi}  {1 \over {\rm sin} \ \theta}  {
{\rm d} \over {\rm d} \theta} \right)^{{d-1 \over 2}} \int^ \theta_{
\varepsilon( {\rm Re} \ \theta) \pi} \left[ {\rm e}^{-i \left({d-2 \over 2}
\right)\tau}\hat  f^\ast( \tau) \right][2( {\rm cos} \ \theta - {\rm cos} \
\tau)]^{ -1/2} {\rm sin} \ \tau \ {\rm d} \tau \eqno (4.35) $$
(where $ \varepsilon( {\rm Re} \ \theta) \pi $ and $ \theta $ stand for the
end-points of the integration path $ \gamma^ \ast_ \theta , $ and $ [( {\rm
cos} \ \theta - {\rm cos} \ \tau)]^{ -1/2}>0 $
for $ \theta $ real $ (\vert \theta\vert \leq\vert \tau\vert  \leq  \pi) ). $
\par
Statements similar to propositions 14 and 15 could thus be written for the
inverse
$ {\cal X}^{\ast( c)}_d $ of $ {\cal R}^{\ast( c)}_d, $ including also the
fact that the kernel $ \Xi^{ \ast( d)} $ of the transformation $ {\cal
X}^{\ast( c)}_d $
is the set of all functions of the form:
$$ \hat  f^\ast( \tau)  = {\rm e}^{i \left({d-2 \over 2} \right)\tau}
P^{\ast( d)}(\tau) \ , \eqno (4.36) $$
with:
$$ P^{\ast( d)}(\tau)  = \sum^{ }_{ \left(\ell \in \Bbb Z;-{d-2 \over 2}\leq
\ell <0 \right)}b_{\ell} {\rm sin} \left[ \left(\ell +{d-2 \over 2}
\right)\tau - \left({d-2 \over 2} \right)\pi \right] \eqno (4.36^\prime ) $$
\par
\smallskip
{\sl The Legendre functions\/} $ P^{(d)}_\lambda $ {\sl considered as
perikernels\/} \par
Let us introduce the following family of functions $ \left\{\hat  f^{\ast(
d)}_\lambda( \tau) \ ;\ \lambda \in \Bbb C \right\} , $ analytic
in $ \dot {\cal J}^{(0)}: $ \par
\smallskip
a) for $ d $ even, $ \hat  f^{\ast( d)}_\lambda =\hat  f^{(d)}_\lambda $ (see
Eq.(3.75)) \par
\smallskip
b) for $ d $ odd, $ 0 \leq  {\rm Re} \ \tau  \leq  2\pi \ : $
$$ \matrix{\displaystyle\hat  f^{\ast( d)}_\lambda( \tau)  = {\rm e}^{i
\left({d-2 \over 2} \right)(\tau -\pi)} \ {\rm sin} \left[ \left(\lambda +{d-2
\over 2} \right)(\tau -\pi) \right] \cr\displaystyle \left(\hat  f^{\ast(
d)}_\lambda( \tau +2\pi n) = \hat  f^{\ast( d)}_\lambda( \tau) \ ,\ n\in  \Bbb
Z \right)\ . \cr} $$
\par
Then, one obtains the following statement which exhibits the functions
$ \displaystyle{ \pi \over \omega_ d} h_d(\lambda) P^{(d)}_\lambda( {\rm
cos}(\theta -\pi)) $
as a second family of elementary perikernels on $ X^{(c)}_{d-1} $ via the $
\ast $-Radon-Abel inverse
transformation. \par
\medskip
{\bf Proposition 25 :} \par
i) For each integer $ d $ $ (d \geq  3), $ the functions $ \hat  f^{\ast(
d)}_\lambda $ satisfy (for all $ \lambda \in \Bbb C) $ the
symmetry property $ \left(S^\ast_ d \right). $ \par
\smallskip
ii) $ \forall \lambda  \in  \Bbb C, $ let $ {\cal X}^{\ast( c)}_d \left(\hat
f^{\ast( d)}_\lambda \right) = {\cal F}^{\ast( d)}_\lambda ; $ then one has:
$$ {\cal F}^{\ast( d)}_\lambda( z)_{ \left\vert z^{(d-1)}= {\rm cos} \ \theta
\right.} = {\pi \over \omega_ d} h_d(\lambda)  P^{(d)}_\lambda( {\rm
cos}(\theta -\pi)) \eqno (4.37) $$
\par
To show formula (4.37), one notices that formula (4.35) written for $ \hat
f^\ast =\hat  f^{\ast( d)}_\lambda $ (case
$ d $ odd), with $ 0 \leq  \theta  \leq  2\pi , $ reduces to the r.h.s. of
Eq.(3.86) after making the substitution
$ \theta \longrightarrow \theta -\pi $ (and the change of integration variable
$ \tau =\tau^{ \prime} -\pi ) $ in the latter). In the case $ d $ even,
$ {\cal X}^{(c)}_d \left(\hat  f^{(d)}_\lambda \right) $ coincides with $
{\cal X}^{\ast( c)}_d \left(\hat  f^{\ast( d)}_\lambda \right) $ and the
coincidence of the functions $ \psi^{( d)}_\lambda( {\rm cos} \ \theta) $
and $ {\pi \over \omega_ d} h_d(\lambda)  P^{(d)}_\lambda( {\rm cos}(\theta
-\pi)) $ is just confirmed (see Eq.(3.90)). \par
\medskip
{\sl The alternative form for the Laplace transformation \/}$ {\cal L}_d $
{\sl (case\/} $ d $ {\sl odd)\nobreak\ :\/} \par
The transformation $ {\cal L}_d, $ defined for $ m>-1 $ on all the
corresponding subspaces $ \left[{\cal V}_1 \right]^m_\mu \left(X^{(c)}_{d-1}
\right) $
(see theorem 1), admits the alternative expression (4.31): this follows
directly from
proposition 24-ii), since the identity $ \Delta\hat  f=\Delta\hat  f^\ast $
implies: $ \tilde  F=L \left(\Delta\hat  f \right)=L \left(\Delta\hat  f^\ast
\right) $ or $ \tilde  F={\cal L} \left(\hat  f \right)={\cal L} \left(\hat
f^\ast \right). $ \par
In order to justify the expression (4.32) of the transformation $ \left[{\cal
L}_d \right]^{-1}_m $ (namely the
inverse of $ {\cal L}_d $ on the corresponding subspaces $ [{\cal V}]^m_\mu
\left(X^{(c)}_{d-1} \right); $ see theorem 5-c), we shall
introduce (for $ m\geq -{d-2 \over 2}) $ \par
\smallskip
i) the space $ {\cal O}^{\ast( d)}_m \left(\dot {\cal J}^{(\mu)} \right) $ of
all functions $ \hat  f^\ast( \theta) $ analytic in $ \dot {\cal J}^{(\mu)} $
whose restriction
to $ \dot {\cal J}^{(\mu)}_ + $ belongs to $ {\cal O}_m \left(\dot {\cal
J}^{(\mu)}_ + \right) $ and which satisfy the symmetry condition $
\left(S^\ast_ d \right). $ These functions
have Fourier coefficients $ \left[\hat {\bf f}^\ast \right]_{\ell} $ which
satisfy the set of relations: \par
\noindent$ \forall \ell  \in  \Bbb Z\ \ ,\ \ \ \left[\hat {\bf f}^\ast
\right]_{\ell}  = \left[\hat {\bf f}^\ast \right]_{-(\ell +d-2)}. $ \par
\smallskip
ii) the kernel $ {\cal N}^{\ast( d)}_m $ of $ {\cal L} $ in $ {\cal O}^{\ast(
d)}_m \left(\dot {\cal J}^{(\mu)} \right), $ namely
$$ {\cal N}^{\ast( d)}_m = \left\{\hat  f^\ast( \theta) = {\rm e}^{i
\left({d-2 \over 2} \right)(\theta -\pi)} \times \sum^{ }_{ \left\{ \ell \in
\Bbb Z,-{d-2 \over 2}\leq \ell <m \right\}} a_{\ell} {\rm sin} \left[
\left(\ell +{d-2 \over 2} \right)\theta - \left({d-2 \over 2} \right)\pi
\right] \right\} \eqno (4.38) $$
\par
iii) the subspace $ {\cal M}^{\ast( d)}_m \left(\dot {\cal J}^{(\mu)} \right)
$ of functions $ \hat  f^\ast $ in $ {\cal O}^{\ast d}_m \left(\dot {\cal
J}^{(\mu)} \right) $ such that $ [ {\bf f}]_{\ell} =0 $ for $ -m-d+2 < \ell  <
m. $ \par
\smallskip
iv) the inverse of $ {\cal L}, $ considered as acting on $ {\cal M}^{\ast(
d)}_m \left(\dot {\cal J}^{(\mu)} \right), $ which we denote by $ [{\cal
L}]^{\ast -1}_{(d,m)}. $ \par
(Compare with the similar notions $ {\cal O}^{(d)}_m \left(\dot {\cal
J}^{(\mu)} \right), $ $ {\cal N}^{(d)}_m, $ $ {\cal M}^{(d)}_m \left(\dot
{\cal J}^{(\mu)} \right), $ $ [{\cal L}]^{-1}_{(d,m)} $
introduced in \S 2.2). \par
\smallskip
A statement similar to proposition 5 can be proved, namely: \par
\medskip
{\bf Proposition 26 :} \par
For $ d $ odd, the following inversion formula applies to the Laplace
transform $ \tilde  F={\cal L} \left(\hat  f^\ast \right) $
of any function $ \hat  f^\ast $ in $ {\cal M}^{\ast( d)}_m \left(\dot {\cal
J}^{(\mu)} \right): $
$$ \hat  f^\ast( \tau)  = \left\{[{\cal L}]^{\ast -1}_{(d,m)} \left(\tilde  F
\right) \right\}( \tau) =-{1 \over 2\pi} {\rm e}^{i \left({d-2 \over 2}
\right)(\tau -\pi)} \int^{ +\infty}_{ -\infty}{\tilde  F(m+i\nu) {\rm sin}
\left[ \left(m+{d-2 \over 2}+i\nu \right)(\tau -\pi) \right] \over {\rm sin} \
\pi( m+i\nu)} {\rm d} \nu \eqno (4.39) $$
Formula (4.32) can now be justified as an inversion formula for $ {\cal L}_d,
$ valid on the
corresponding subspaces $ [{\cal V}]^m_\mu \left(X^{(c)}_{d-1} \right), $
(with $ m>-1): $ the argument is similar to the proof of
theorem 5-c). \par
Let then $ {\cal F} $ be a perikernel in $ [{\cal V}]^m_\mu
\left(X^{(c)}_{d-1} \right), $ for which we write:
$$ {\cal F}= \left[{\cal L}_d \right]^{-1}_m \left(\tilde  F \right) = {\cal
X}^\ast_ d \left(\hat  f^\ast \right),\ \ \ {\rm with} \ \ \ \hat  f^\ast =
[{\cal L}]^{\ast -1}_{(d,m)} \left(\tilde  F \right)\ {\rm given\ by\ Eq.}
(4.39). $$
\par
As in the proof of Eq.(4.23) (concerning the action of $ {\cal X}^{(c)}_d), $
we can here justify
the action of the operator $ {\cal X}^{\ast( c)}_d $ under the integration
sign at the r.h.s. of Eq.(4.39) and
therefore write:
$$ {\cal X}^{\ast( c)}_d \left(\hat  f^\ast \right)(\tau)  = - {1 \over 2\pi}
\int^{ +\infty}_{ -\infty}{\tilde  F(m+i\nu){\cal X}^{\ast( c)}_d \left(\hat
f^{\ast( d)}_{m+i\nu} \right)(\tau) \over {\rm sin} \ \pi( m+i\nu)}  {\rm d}
\nu \ , $$
which, in view of proposition 25, is identical with the representation (4.23')
of $ {\cal F} $ (with
$ \goth P^{(d)}_m( {\bf F})=0). $ \par
\medskip
{\bf Connection with Cauchy representations for the Laplace transform:} \par
Our two families of elementary perikernels being defined by the equations: $
{\cal F}^{(d)}_{\lambda_ 0}={\cal X}^{(c)}_d \left(\hat  f^{(d)}_{\lambda_ 0}
\right) $
and $ {\cal F}^{\ast( d)}_{\lambda_ 0} = {\cal X}^{\ast( c)}_d \left(\hat
f^{\ast( d)}_{\lambda_ 0} \right), $ $ (\lambda_ 0\in \Bbb C), $ it is natural
to look for the Laplace transforms $ \tilde  F^{(d)}_{\lambda_ 0}(\lambda) $
and $ \tilde  F^{\ast( d)}_{\lambda_ 0}(\lambda) $ of these perikernels. First
of all, in view of proposition 15, we can say
that $ {\cal R}^{(c)}_d \left({\cal F}^{(d)}_{\lambda_ 0} \right)-\hat
f^{(d)}_{\lambda_ 0}\in \Xi^{( d)} $ and similarly $ {\cal R}^{\ast( c)}_d
\left({\cal F}^{\ast( d)}_{\lambda_ 0} \right)-\hat  f^{\ast( d)}_{\lambda_
0}\in \Xi^{ \ast( d)}. $ \par
Now, for any fixed number $ m $ such that $ m>-1, $ one sees that: $ \Xi^{(
d)}\subset{\cal N}^{(d)}_m $ and $ \Xi^{ \ast( d)}\subset{\cal N}^{\ast(
d)}_m; $
this implies that, {\sl for every\/} $ \lambda_ 0 $ {\sl such that:\/} $ -
{d-2 \over 2} \leq  {\rm Re} \ \lambda_ 0 < m, $ the following identities
hold, as analytic functions in $ \Bbb C^{(m)}_+ $ (since $ {\cal
F}^{(d)}_{\lambda_ 0} $ and $ {\cal F}^{\ast( d)}_{\lambda_ 0} $ then belong
to $ \left[{\cal V}_{\infty} \right]^m_0 \left(X^{(c)}_{d-1} \right)): $
$$ \tilde  F^{(d)}_{\lambda_ 0} = {\cal L}_d \left({\cal F}^{(d)}_{\lambda_ 0}
\right) = {\cal L} \left(\hat  f^{(d)}_{\lambda_ 0} \right) \eqno (4.40) $$
and similarly (from Eq.(4.31)):
$$ \tilde  F^{\ast( d)}_{\lambda_ 0} = {\cal L}_d \left({\cal F}^{\ast(
d)}_{\lambda_ 0} \right) = {\cal L} \left(\hat  f^{\ast( d)}_{\lambda_ 0}
\right)\ . \eqno (4.41) $$
The $\ell$.h.s. of Eqs.(4.40) and (4.41) are immediately computable as Laplace
transforms of
trigonometric functions, thus yielding the following rational functions (in
view of
Eqs.(2.3), (3.76), (3.77) in the case of $ \tilde  F^{(d)}_{\lambda_ 0}, $ and
of similar equations for the
discontinuity of $ \hat  f^{\ast( d)}_{\lambda_ 0} $ in the case of $ \tilde
F^{\ast( d)}_{\lambda_ 0} $ with $ d $ odd): \par
\smallskip
a) for all $ d: $
$$ \tilde  F^{\ast( d)}_{\lambda_ 0}(\lambda)  = - {\rm sin} \ \pi \lambda_ 0
\left[{1 \over \lambda -\lambda_ 0} - {1 \over \lambda +\lambda_ 0+d-2}
\right] \eqno (4.42) $$
\par
b) for $ d $ odd,
$$ \tilde  F^{(d)}_{\lambda_ 0}(\lambda)  = - {\rm sin} \ \pi \lambda_ 0
\left[{1 \over \lambda -\lambda_ 0} + {1 \over \lambda +\lambda_ 0+d-2}
\right] \eqno (4.43) $$
(for $ d $ even, $ \tilde  F^{(d)}_{\lambda_ 0}(\lambda)  \equiv  \tilde
F^{\ast( d)}_{\lambda_ 0}(\lambda) ). $ \par
Now, by making use of the explicit integral form (4.21) with $ \mu =0 $ (or
(4.5)) of the
transformation $ {\cal L}_d $ in terms of the discontinuities $
\underline{f}_{\lambda_ 0}( {\rm cosh} \ v) = \Delta \Psi^{( d)}_{\lambda_ 0}(
{\rm cosh} \ v) $ and:
$$ \matrix{\displaystyle \underline{f}^{\ast( d)}_{\lambda_ 0}( {\rm cosh} \
v) =  {\pi \over \omega_ d} h_d \left(\lambda_ 0 \right) \Delta
P^{(d)}_{\lambda_ 0}( {\rm cos}(iv-\pi))  = \hfill \cr\displaystyle \ \ =-2\
{\rm sin} \ \pi \lambda_ 0 \times  {\pi \over \omega_ d} h_d \left(\lambda_ 0
\right) \left[P^{(d)}_{\lambda_ 0}( {\rm cosh} \ v)+2\ {\rm cos}  {\pi d
\over 2} Q^{(d)}_{\lambda_ 0}( {\rm cosh} \ v) \right], \hfill \cr} $$
associated respectively with the perikernels $ {\cal K}^{(d)}_{\lambda_ 0} $
and $ {\cal K}^{\ast( d)}_{\lambda_ 0} $ (see Eqs.(3.103), (3.105'), (3.106)),
we obtain the following \par
\smallskip
{\bf Proposition 27 :} \par
To the families of perikernels $ {\cal K}^{(d)}_{\lambda_ 0}, $ $ {\cal
K}^{\ast( d)}_{\lambda_ 0} $ $ ( {\rm Re} \ \lambda_ 0 \geq  - {d-2 \over 2}),
$ there correspond respectively the
families of Fourier-Laplace transforms $ \tilde  F^{(d)}_{\lambda_ 0}, $ $
\tilde  F^{\ast( d)}_{\lambda_ 0} $ given by Eqs.(4.42), (4.43), considered as
analytic functions in the corresponding half-plane $ \Bbb C^{ \left( {\rm Re}
\ \lambda_ 0 \right)}_+. $ This property is expressed by the following
integral relations: \par
For all $ d, $ $ \lambda , $ $ \lambda_ 0, $ such that: $ {\rm Re} \ \lambda
> {\rm Re} \ \lambda_ 0 \geq  - {d-2 \over 2}\ : $
$$ \matrix{\displaystyle 2\pi  {\omega_{ d-1} \over \omega_ d} h_d
\left(\lambda_ 0 \right) \int^{ \infty}_ 0 \left[P^{(d)}_{\lambda_ 0}( {\rm
cosh} \ v) + 2\ {\rm cos}  {\pi d \over 2} Q^{(d)}_{\lambda_ 0}( {\rm cosh} \
v) \right]Q^{(d)}_\lambda( {\rm cosh} \ v)( {\rm sinh} \ v)^{d-2} {\rm d} v
\hfill \cr\displaystyle \ \ \ \ \ = {2\lambda_ 0+d-2 \over \left(\lambda
-\lambda_ 0 \right) \left(\lambda +\lambda_ 0+d-2 \right)} \hfill \cr} \eqno
(4.44) $$
\par
For $ d $ odd:
$$ \matrix{\displaystyle 2 \omega_{ d-1} \int^{
\infty}_ 0 \left\{ {\pi \over \omega_d} h_d \left(\lambda_ 0 \right)
\left[P^{(d)}_{\lambda_ 0}(
{\rm cosh} v) + 2(-1)^{{d-1 \over 2}} {\rm tg} \ \pi \lambda_ 0
Q^{(d)}_{\lambda_ 0}( {\rm cosh} v) \right]+ i \left({1 \over 2\pi}
\right)^{{d-1 \over 2}}{R^{(d)}_{\lambda_ 0}( {\rm cosh} v) \over {\rm cos} \
\pi \lambda_ 0} \right\} \hfill \cr\displaystyle \ \ \ \ \ \times
Q^{(d)}_\lambda( {\rm cosh} \ v)( {\rm sinh} \ v)^{d-2} {\rm d} v = {2\lambda
+d-2 \over \left(\lambda -\lambda_ 0 \right) \left(\lambda +\lambda_ 0+d-2
\right)} \hfill \cr} \eqno (4.45) $$
\par
\smallskip
{\sl Final Remark :\/} \par
By applying to both sides of Eqs.(4.23) and (4.23') the operator $ {\cal L}_d,
$ and by taking into account
the result of proposition 27, we see that, for $ d $ odd, the representations
(4.23) and (4.23') of an
invariant perikernel $ {\cal K}, $ such that $ {\cal F}(z) = {\cal K}
\left(z,z_0 \right) \in  [{\cal V}]^m_\mu \left(X^{(c)}_{d-1} \right) $ $ (m>
{\rm max} \left(-1,-{d-2 \over 2} \right)) $ are respectively equivalent to
two
(equally valid) Cauchy integral representations of its Fourier-Laplace
transform $ \tilde  F={\cal L}_d({\cal F}), $ namely (for $ \lambda  \in  \Bbb
C^{ \left( {\rm Re} \ \lambda_ 0 \right)}_+): $
$$ \tilde  F(\lambda)  = {1 \over 2\pi i} \int^{ }_{ {\rm L}_{ \left( {\rm Re}
\ \lambda_ 0 \right)}}\tilde  F \left(\lambda_ 0 \right) \left[{1 \over
\lambda -\lambda_ 0} \mp  {1 \over \lambda +\lambda_ 0+d-2} \right] {\rm d}
\lambda_ 0 \eqno (4.46) $$
\par
For $ d $ even, Eq.(4.23) (identical with Eq.(4.23')) is equivalent to
Eq.(4.46) with the
specification given by the minus sign in the integrand. \par
\vfill\eject
\smallskip
\centerline{{\bf APPENDIX }}
\bigskip
{\bf D\'efinitions :} for $ d\in \Bbb N, $ we introduce \par
a) the following differential operators:
$$ \eqalignno{ d\geq 2\ ,\  & \Delta_ d = { {\rm d}^2 \over {\rm d} v^2} -
\left({d-2 \over 2} \right)^2\ , & (A.1) \cr d\geq 3\ ,\  & D_d = ( {\rm sinh}
\ v)^{d-1} \left({1 \over {\rm sinh} \ v} { {\rm d} \over {\rm d} v} \right){1
\over( {\rm sinh} \ v)^{d-3}} = {\rm sinh} \ v { {\rm d} \over {\rm d} v}
-(d-3) {\rm cosh} \ v & (A.2) \cr} $$
\par
b) the following polynomial function:
$$ d\geq 3\ ,\ \ \ h^{(0)}_d(\lambda)  = \left(\lambda  + {d-2 \over 2}
\right) {\Gamma( \lambda +d-2) \over \Gamma( \lambda +1)}\ , \eqno (A.3) $$
which satisfies the identity:
$$ h^{(0)}_d(\lambda)  = (-1)^dh^{(0)}_d(-\lambda -d+2)\ , \eqno (A.3^\prime )
$$
\par
c) the operator $ H^{(0)}_d $ associated with the function $ h^{(0)}_d
\left(-{d-2 \over 2}+i\nu \right) $
through the correspondence $ i\nu  \longleftrightarrow  \displaystyle{ {\rm d}
\over {\rm d} v}, $ namely: \par
{\sl for d even\/} :
$$ d\geq 4,\ \ \ H^{(0)}_d= \left[{ {\rm d}^2 \over {\rm d} v^2} - \left({d-4
\over 2} \right)^2 \right]... \left[{ {\rm d}^2 \over {\rm d} v^2}-1 \right] {
{\rm d}^2 \over {\rm d} v^2} = \Delta_{ d-2}...\Delta_ 4\Delta_ 2 \eqno (A.4)
$$
\par
{\sl for d odd\/} :
$$ H^{(0)}_3 = {d \over dv}\ , \eqno (A.5) $$
$$ d\geq 5,\ \ \ H^{(0)}_d= \left[{ {\rm d}^2 \over {\rm d} v^2} - \left({d-4
\over 2} \right)^2 \right]... \left[{ {\rm d}^2 \over {\rm d} v^2}- \left({1
\over 2} \right)^2 \right] { {\rm d} \over {\rm d} v} = \Delta_{
d-2}...\Delta_ 5\Delta_ 3{d \over dv}\ . \eqno (A.6) $$
\par
The operators $ \Delta_ d $ and $ D_d $ are mutually related by the properties
described in the following \par
\smallskip
{\bf Lemma A.1} : Let $ {\cal I}_1(v) = [0,v] $ and $ {\cal I}_{-1}(v) =
[v,\infty[ , $ considered as
oriented intervals. \par
Then for every integer $ d, $ $ d\geq 3, $ the following identities hold: \par
$ \varepsilon =+1 $ or $ -1, $
$$ \matrix{\displaystyle \varepsilon( d-2)D_d \int^{ }_{{\cal I}_\varepsilon(
v)}g(w)[2\varepsilon( {\rm cosh} \ v- {\rm cosh} \ w)]^{{d-4 \over 2}} {\rm d}
w \hfill \cr\displaystyle \ \ \ \ \ = \int^{ }_{{\cal I}_\varepsilon(
v)}\Delta_ dg(w)[2\varepsilon( {\rm cosh} \ v- {\rm cosh} \ w)]^{{d-2 \over
2}} {\rm d} w\ , \hfill \cr} \eqno (A.7) $$
provided the function $ g(v) $ belongs to $ {\cal C}^2 \left(\Bbb R^+ \right)
$ and satisfies the
following respective conditions: \par
\smallskip
i) for $ \varepsilon  = +1, $ $ \displaystyle{ {\rm d} g \over {\rm d} v}(0) =
0 $ \par
ii) for $ \varepsilon  = -1, $ the functions $ {\rm e}^{ \left({d-2 \over 2}
\right)v}\cdot g(v) $ and $ {\rm e}^{ \left({d-2 \over 2} \right)v}\cdot
\displaystyle{ {\rm d} g \over {\rm d} v}(v) $ both
belong to $ L^1 \left(\Bbb R^+ \right) $ and tend to zero when $ v
\longrightarrow  +\infty . $ \par
\smallskip
{\bf Proof :} \par
Let $ A_\varepsilon  = 2\varepsilon( {\rm cosh} \ v- {\rm cosh} \ w); $ the
following computations are
valid for both cases $ \varepsilon =+1 $ and $ \varepsilon =-1: $ the
corresponding conditions i) and ii)
are used for the convergence properties of integrals and the cancellation of
integrated terms at the end-points $ w=0 $ (resp. $ +\infty ) $ of $ {\cal
I}_1(v) $ (resp. $ {\cal I}_{-1}(v)); $
the assumption $ d\geq 3 $ is used for similar properties at the end-point $
w=v $ of $ {\cal I}_\varepsilon( v). $ \par
\smallskip
a)
$$ \varepsilon \int^{ }_{{\cal I}_\varepsilon( v)}{ {\rm d} \over {\rm d} w}
\left[g(w)A^{{d-4 \over 2}}_\varepsilon( {\rm sinh} \ v- {\rm sinh} \ w)
\right] {\rm d} w=- \left({1+\varepsilon \over 2} \right)g(0) {\rm sinh} \
v[2( {\rm cosh} \ v-1)]^{{d-4 \over 2}}, $$
which yields:
$$ \matrix{\displaystyle \varepsilon \int^{ }_{{\cal I}_\varepsilon( v)}{
{\rm d} g \over {\rm d} w}(w)A^{{d-4 \over 2}}_\varepsilon( {\rm sinh} \ v-
{\rm sinh} \ w) {\rm d} w+ \left({1+\varepsilon \over 2} \right)g(0) {\rm
sinh} \ v[2( {\rm cosh} \ v-1)]^{{d-4 \over 2}} \hfill \cr\displaystyle =
\int^{ }_{{\cal I}_\varepsilon( v)}g(w)A^{{d-6 \over 2}}_\varepsilon[( d-4)
{\rm sinh} \ w( {\rm sinh} \ v- {\rm sinh} \ w)+2\ {\rm cosh} \ w( {\rm cosh}
\ v- {\rm cosh} \ w)] {\rm d} w \hfill \cr} \eqno (A.8) $$
\par
b) By applying the operator $ { {\rm d} \over {\rm d} v} $ to both sides of
the identity:
$$ \left[g(w)A^{{d-2 \over 2}}_\varepsilon \right]_{\partial{\cal
I}_\varepsilon( v)}= \int^{ }_{{\cal I}_\varepsilon( v)}{ {\rm d} g \over
{\rm d} w}(w) A^{{d-2 \over 2}}_\varepsilon {\rm d} w-\varepsilon( d-2)
\int^{ }_{{\cal I}_\varepsilon( v)}g(w)A^{{d-4 \over 2}}_\varepsilon {\rm
sinh} \ w\ {\rm d} w, $$
we obtain:
$$ \matrix{\displaystyle \varepsilon \ {\rm sinh} \ v \int^{ }_{{\cal
I}_\varepsilon( v)}{ {\rm d} g \over {\rm d} w}(w) A^{{d-4 \over
2}}_\varepsilon {\rm d} w + \left({1+\varepsilon \over 2} \right)g(0) {\rm
sinh} \ v[2( {\rm cosh} \ v-1)]^{{d-4 \over 2}} \cr\displaystyle =
\varepsilon  { {\rm d} \over {\rm d} v} \left[ \int^{ }_{{\cal
I}_\varepsilon( v)}g(w) A^{{d-4 \over 2}}_\varepsilon {\rm sinh} \ w\ {\rm d}
w \right] \cr} \eqno (A.9) $$
\par
c) Integration by parts yields:
$$ \varepsilon( d-2) \int^{ }_{{\cal I}_\varepsilon( v)}{ {\rm d} g \over
{\rm d} w} A^{{d-4 \over 2}}_\varepsilon {\rm sinh} \ w\ {\rm d} w = \int^{
}_{{\cal I}_\varepsilon( v)}{ {\rm d}^2g \over {\rm d} w^2}(w) A^{{d-2 \over
2}}_\varepsilon  {\rm d} w \eqno (A.10) $$
\par
d) By taking Eqs.(A.9) and (A.10) into account at the r.h.s. of Eq.(A.8),
we then obtain:
$$ \matrix{\displaystyle \varepsilon( d-2){ {\rm d} \over {\rm d} v} \left[
\int^{ }_{{\cal I}_\varepsilon( v)}g(w)A^{{d-4 \over 2}}_\varepsilon {\rm
sinh} \ w\ {\rm d} w \right]- \int^{ }_{{\cal I}_\varepsilon( v)}{ {\rm d}^2g
\over {\rm d} w^2}(w) A^{{d-2 \over 2}}_\varepsilon  {\rm d} w \hfill
\cr\displaystyle = (d-2) \int^{ }_{{\cal I}_\varepsilon( v)}g(w)A^{{d-6 \over
2}}_\varepsilon[( d-4) {\rm sinh} \ w( {\rm sinh} \ v- {\rm sinh} \ w)+2 {\rm
cosh} \ w( {\rm cosh} \ v- {\rm cosh} \ w)] {\rm d} w \hfill \cr} \eqno (A.11)
$$
\par
e) By computing in two ways the expression $ \varepsilon \displaystyle{ {\rm
d} \over {\rm d} v} \left\{ \int^{ }_{{\cal I}_\varepsilon( v)}g(w)A^{{d-4
\over 2}}_\varepsilon( {\rm sinh} \ v- {\rm sinh} \ w) {\rm d} w \right\} , $
we obtain the identity:
$$ \matrix{\displaystyle( d-4) {\rm sinh} \ v \int^{ }_{{\cal I}_\varepsilon(
v)}g(w) A^{{d-6 \over 2}}_\varepsilon( {\rm sinh} \ v- {\rm sinh} \ w) {\rm d}
w \hfill \cr\displaystyle = \varepsilon \ {\rm sinh} \ v{ {\rm d} \over {\rm
d} v} \left[ \int^{ }_{{\cal I}_\varepsilon( v)}g(w)A^{{d-4 \over
2}}_\varepsilon {\rm d} w \right]-\varepsilon{ {\rm d} \over {\rm d} v} \left[
\int^{ }_{{\cal I}_\varepsilon( v)}g(w)A^{{d-4 \over 2}}_\varepsilon {\rm
sinh} \ w\ {\rm d} w \right] \hfill \cr} \eqno (A.12) $$
\par
f) In view of Eq.(A.12), Eq.(A.11) can then be rewritten as follows:
$$ \matrix{\displaystyle \varepsilon( d-2) {\rm sinh} \ v { {\rm d} \over
{\rm d} v} \left[ \int^{ }_{{\cal I}_\varepsilon( v)}g(w) A^{{d-4 \over
2}}_\varepsilon {\rm d} w \right]- \int^{ }_{{\cal I}_\varepsilon( v)}{ {\rm
d}^2g \over {\rm d} w^2} (w) A^{{d-2 \over 2}}_\varepsilon  {\rm d} w \hfill
\cr\displaystyle = (d-2) \int^{ }_{{\cal I}_\varepsilon( v)}g(w) A^{{d-6
\over 2}}_\varepsilon \left[(d-4) \left( {\rm sinh}^2v- {\rm sinh}^2w
\right)+2\ {\rm cosh} \ w( {\rm cosh} \ v- {\rm cosh} \ w) \right] {\rm d} w
\hfill \cr\displaystyle = (d-2) \int^{ }_{{\cal I}_\varepsilon( v)}g(w)A^{{d-6
\over 2}}_\varepsilon( {\rm cosh} \ v- {\rm cosh} \ w)[(d-4)( {\rm cosh} \ v+
{\rm cosh} \ w)+2\ {\rm cosh} \ w] {\rm d} w \hfill \cr} \eqno (A.13) $$
\par
g) By using the expression of the operator $ D_d $ given at the r.h.s. of
Eq.(A.2),
and Eq.(A.13), we then obtain the following identity:
$$ \matrix{\displaystyle \varepsilon( d-2)D_d \left[ \int^{ }_{{\cal
I}_\varepsilon( v)}g(w) A^{{d-4 \over 2}}_\varepsilon {\rm d} w \right]-
\int^{ }_{{\cal I}_\varepsilon( v)}{ {\rm d}^2g \over {\rm d} w^2} (w) A^{{d-2
\over 2}}_\varepsilon  {\rm d} w \hfill \cr\displaystyle = (d-2) \int^{
}_{{\cal I}_\varepsilon( v)}g(w) A^{{d-6 \over 2}}_\varepsilon( {\rm cosh} \
v- {\rm cosh} \ w)[(d-4)( {\rm cosh} \ v+ {\rm cosh} \ w) \hfill
\cr\displaystyle +2\ {\rm cosh} \ w-2(d-3) {\rm cosh} \ v] {\rm d} w
\cr\displaystyle = -(d-2)^2 \int^{ }_{{\cal I}_\varepsilon( v)}g(w)A^{{d-6
\over 2}}_\varepsilon( {\rm cosh} \ v- {\rm cosh} \ w)^2= - {(d-2)^2 \over 4}
\int^{ }_{{\cal I}_\varepsilon( v)}g(w)A^{{d-2 \over 2}}_\varepsilon {\rm d} w
\hfill \cr} \eqno (A.14) $$
In view of the definition (A.1) of $ \Delta_ d, $ Eq.(A.14) coincides with the
announced result, i.e.
Eq.(A.7). \par
\medskip
We shall now make use of the following \lq\lq iterated form\rq\rq\ of lemma
A.1. \par
\medskip
{\bf Lemma A.2 :} Let $ \left\{ \Xi^{( \varepsilon)}_ d,\ d\geq d_0,\ d_0\geq
3 \right\} $ denote a sequence of (integro-differential) operators
on $ C^{\infty} \left(\Bbb R^+ \right) $ which satisfy the following
relations:
$$ \matrix{\displaystyle \forall d,\ d=d_0+2n,\ n\in \Bbb N \cr\displaystyle
\Xi^{( \varepsilon)}_ d = {1 \over( {\rm sinh} \ v)^{d-3}} \cdot  D_{d-2}(
{\rm sinh} \ v)^{d-5}\cdot \Xi^{( \varepsilon)}_{ d-2} \cr} \eqno (A.15) $$
Then, a sufficient condition for the following identity to hold for a given
value of $ d $ (with $ \varepsilon $
fixed, $ \varepsilon =+1 $ or $ -1): $
$$ \Xi^{( \varepsilon)}_ dg(v) = {C_d \over( {\rm sinh} \ v)^{d-3}} \int^{
}_{{\cal I}_\varepsilon( v)} \left[H^{(0)}_dg(w) \right]\{ 2\varepsilon( {\rm
cosh} \ v- {\rm cosh} \ w)\}^{{ d-4 \over 2}} {\rm d} w \eqno (A.16) $$
is that Eq.(A.16) holds for the value $ d=d_0 $ and that the function $ g(v) $
satisfies the following
additional conditions: \par
i) for $ \varepsilon  = +1, $ $ g\in{\cal C}^{\infty}([ 0,v]) $ and $ \forall
n, $  $ n \leq  {d-d_0 \over 2}, $ $ \displaystyle \left[{d^{(2n+1)} \over
dv^{(2n+1)}}H^{(0)}_{d_0}g \right](0) =0 $ \par
\smallskip
ii) for $ \varepsilon =-1, $ $ g\in{\cal C}^{\infty}([ v,+\infty[) $ and $
\forall n, $ $ n \leq  d+1, $ $ \displaystyle{ d^ng \over dv^n}(v)\cdot {\rm
e}^{ \left({d-2 \over 2} \right)v}\in L^1 \left(\Bbb R^+ \right) $ and tends
to zero
for $ v \longrightarrow \infty . $ \par
Moreover, the constants $ C_d $ are such that:
$$ C_d = (\varepsilon)^{{ d-d_0 \over 2}} \left[(d-4)\cdot( d-6)...d_0
\left(d_0-2 \right) \right]C_{d_0} \eqno (A.17) $$
\par
{\bf Proof :} \par
This is an immediate application of lemma A.1, used through the following
recursive
argument. \par
Since $ H^{(0)}_d=\Delta_{ d-2}H^{(0)}_{d-2} $ (see Eqs.(A.4), (A.6)), we can
write the r.h.s. of Eq.(A.16) as
follows (in view of Eq.(A.7)):
$$ \matrix{\displaystyle{ C_d \over( {\rm sinh} \ v)^{d-3}} \int^{ }_{{\cal
I}_\varepsilon( v)}\Delta_{ d-2} \left(H^{(0)}_{d-2}g \right)(w)\{
2\varepsilon( {\rm cosh} \ v- {\rm cosh} \ w)\}^{{ d-4 \over 2}} {\rm d} w
\hfill \cr\displaystyle = {\varepsilon( d-4)C_d \over( {\rm sinh} \ v)^{d-3}}
D_{d-2} \int^{ }_{{\cal I}_\varepsilon( v)} \left[H^{(0)}_{d-2}g(w) \right]\{
2\varepsilon( {\rm cosh} \ v- {\rm cosh} \ w)\}^{{ d-6 \over 2}} {\rm d} w\ ,
\hfill \cr} $$
which, by recursive assumption (i.e. Eq.(A.16) for $ d \longrightarrow d-2) $
is equal to the expression
$$ {\varepsilon( d-4)C_d \over( {\rm sinh} \ v)^{d-3}} D_{d-2} \left[{( {\rm
sinh} \ v)^{d-5} \over C_{d-2}} \Xi^{( \varepsilon)}_{ d-2}g(v) \right]\ , $$
and therefore (in view of Eq.(A.15)) to $ \varepsilon( d-4) {C_d \over
C_{d-2}} \Xi^{( \varepsilon)}_ dg(v)\cdot $ Eq.(A.16) is thus recursively
established, provided the constants $ C_d $ satisfy the relation $ C_d =
\varepsilon( d-4)^{-1}C_{d-2} $ and therefore
Eq.(A.17). In fact, one checks that the application of lemma 1 is justified at
each step of the
recursion (from $ d_0 $ to $ d) $ if the previous conditions i) or ii) are
satisfied by the function $ g(v). $ \par
We are now in a position to apply lemma A.2 to the case where the operators $
\Xi^{( \varepsilon)}_ d $ are
essentially \lq\lq inverse Abel operators\rq\rq\ (such as in particular the
operator $ {\cal X}_d $ of \S 3.7 in the
situation where the chosen integration set is $ {\cal I}_1(v) = [0,v]). $ The
following theorems A.1 and A.2
will correspond respectively to the cases when the integer $ d $ is even or
odd. \par
\medskip
{\bf Theorem A.1} : Let $ d $ be an even integer, $ d\geq 4. $ \par
i) For any function $ g(v) $ in $ {\cal C}^{\infty} \left(\Bbb R^+ \right) $
such that $ \forall n, $ $ 0\leq n\leq{ d-4 \over 2}, $ $ \displaystyle{ {\rm
d}^{(2n+1)}g \over {\rm d} v^{(2n+1)}}(0)=0, $ the
following identity holds:
$$ \matrix{\displaystyle \left({1 \over {\rm sinh} \ v} { {\rm d} \over {\rm
d} v} \right)^{ \left({d-2 \over 2} \right)}g(v)= \hfill \cr\displaystyle = {1
\over 2^{ \left({d-4 \over 2} \right)} \left({d-4 \over 2} \right)!} {1
\over( {\rm sinh} \ v)^{(d-3)}} \int^ v_0 \left[H^{(0)}_dg(w) \right]\{ 2(
{\rm cosh} \ v- {\rm cosh} \ w)\}^{ \left({d-4 \over 2} \right)} {\rm d} w
\hfill \cr} \eqno (A.18) $$
\par
ii) For any function $ g(v) $ in $ {\cal C}^{\infty} \left(\Bbb R^+ \right) $
such that $ \forall n, $ $ n\leq d+1, $ $ {\rm e}^{ \left({d-2 \over 2}
\right)v}\cdot \displaystyle{ {\rm d}^ng \over {\rm d} v^n}\in L^1 \left(\Bbb
R^+ \right) $
and tends to zero when $ v \longrightarrow +\infty , $ the following identity
holds:
$$ \matrix{\displaystyle \left({1 \over {\rm sinh} \ v} { {\rm d} \over {\rm
d} v} \right)^{ \left({d-2 \over 2} \right)}g(v)= \hfill \cr\displaystyle =
{(-1)^{ \left({d-2 \over 2} \right)} \over 2^{ \left({d-4 \over 2} \right)}
\left({d-4 \over 2} \right)!} {1 \over( {\rm sinh} \ v)^{(d-3)}} \cdot
\int^{ +\infty}_ v \left[H^{(0)}_dg(w) \right]\{ 2( {\rm cosh} \ w- {\rm cosh}
\ v)\}^{ \left({d-4 \over 2} \right)} {\rm d} w \hfill \cr} \eqno (A.19) $$
\par
{\bf Proof :} \par
The operators $ \Xi_ d = \left({1 \over {\rm sinh} \ v} { {\rm d} \over {\rm
d} v} \right)^{{d-2 \over 2}} $ obviously satisfy the relations (A.15). \par
If $ g $ satisfies the conditions described in i) or ii), then we can write
(with the
notations of lemma A.2):
$$ { {\rm d} g \over {\rm d} v}(v) = \varepsilon \int^{ }_{{\cal
I}_\varepsilon( v)}\Delta_ 2g(w) {\rm d} w\ , $$
or equivalently :
$$ \Xi_ 4g(v) = \varepsilon  {1 \over {\rm sinh} \ v} \int^{ }_{{\cal
I}_\varepsilon( v)}H^{(0)}_4g(w) {\rm d} w\ , $$
where $ \varepsilon =+1 $ in the case i) and $ \varepsilon  = -1 $ in the case
ii); Eq.(A.16) of lemma A.2 is therefore
satisfied for $ d=d_0=4 $ (with $ C_{d_0}=\varepsilon ). $ \par
For $ d\geq 6, $ one checks that the function $ g $ satisfies the
corresponding conditions i) or ii)
of lemma A.2; in view of the latter, Eq.(A.16) is therefore valid for the
considered value of $ d, $
with $ C_d $ given by Eq.(A.17), namely:
$$ C_d=(\varepsilon)^{{ d-2 \over 2}}[2\times 4\times ...\times(
d-4)]^{-1}=(\varepsilon)^{{ d-2 \over 2}} \left[2^{{d-4 \over 2}} \left({d-4
\over 2} \right)! \right]^{-1}\ ; $$
Eqs.(A.18) and (A.19) are therefore established. \par
\medskip
{\bf Theorem A.2 :} Let $ d $ be an odd integer, $ d\geq 3. $ \par
i) For any function $ g(v) $ in $ {\cal C}^{\infty} \left(\Bbb R^+ \right) $
such that $ \forall n, $ $ 0\leq n\leq{ d-3 \over 2}, $ $ \displaystyle{ {\rm
d}^{2n}g \over {\rm d} v^{2n}}(0) = 0, $ the
following identity holds:
$$ \matrix{\displaystyle \left({1 \over {\rm sinh} \ v} { {\rm d} \over {\rm
d} v} \right)^{ \left({d-1 \over 2} \right)} \int^ v_0{g(w) {\rm sinh} \ w
\over\{ 2( {\rm cosh} \ v- {\rm cosh} \ w)\}^{ 1/2}} {\rm d} w = \hfill
\cr\displaystyle = {2^{ \left({d-3 \over 2} \right)}\Gamma \left({d-1 \over 2}
\right) \over \Gamma( d-2)} {1 \over( {\rm sinh} \ v)^{(d-3)}} \int^ v_0
\left[H^{(0)}_dg(w) \right]\{ 2( {\rm cosh} \ v- {\rm cosh} \ w)\}^{
\left({d-4 \over 2} \right)} {\rm d} w \hfill \cr} \eqno (A.20) $$
\par
ii) For any function $ g(v) $ in $ {\cal C}^{\infty} \left(\Bbb R^+ \right) $
such that $ \forall n, $ $ 0 \leq  n \leq  d-2, $ $ {\rm e}^{ \left({d-2
\over 2} \right)v}\cdot \displaystyle{ {\rm d}^ng \over {\rm d} v^n} \in  L^1
\left(\Bbb R^+ \right) $
and tends to zero when $ v \longrightarrow  +\infty , $ the following identity
holds:
$$ \matrix{\displaystyle \left({1 \over {\rm sinh} \ v} { {\rm d} \over {\rm
d} v} \right)^{ \left({d-1 \over 2} \right)} \int^{ +\infty}_ v{g(w) {\rm
sinh} \ w \over\{ 2( {\rm cosh} \ w- {\rm cosh} \ v)\}^{ 1/2}} {\rm d} w =
\hfill \cr\displaystyle = {(-2)^{ \left({d-3 \over 2} \right)}\Gamma
\left({d-1 \over 2} \right) \over \Gamma( d-2)} {1 \over( {\rm sinh} \
v)^{(d-3)}} \int^{ +\infty}_ v \left[H^{(0)}_dg(w) \right]\{ 2( {\rm cosh} \
w- {\rm cosh} \ v)\}^{ \left({d-4 \over 2} \right)} {\rm d} w \hfill \cr}
\eqno (A.21) $$
\par
{\bf Proof :} \par
The operators $ \Xi^{( \varepsilon)}_ d $ defined by:
$$ \Xi^{( \varepsilon)}_ dg(v) = \left({1 \over {\rm sinh} \ v} { {\rm d}
\over {\rm d} v} \right)^{{d-1 \over 2}} \int^{ }_{{\cal I}_\varepsilon(
v)}{g(w) {\rm sinh} \ w \over[ 2\varepsilon( {\rm cosh} \ v- {\rm cosh} \
w)]^{1/2}} {\rm d} w \eqno (A.22) $$
obviously satisfy the relations (A.15). \par
In view of the condition $ g(0) = 0 $ (resp. $ {\rm e}^{v/2}\displaystyle{
{\rm d} g \over {\rm d} v} \in  L^1 \left(\Bbb R^+ \right)), $ we have for $
\varepsilon =+1 $ (resp. $ \varepsilon =-1): $
$$ \int^{ }_{{\cal I}_\varepsilon( v)}{ {\rm d} \over {\rm d} w} \left\{
g(w)[2\varepsilon( {\rm cosh} \ v- {\rm cosh} \ w)]^{1/2} \right\} {\rm d}
w=0\ , $$
and therefore :
$$ \int^{ }_{{\cal I}_\varepsilon( v)}{ {\rm d} g \over {\rm d}
w}(w)[2\varepsilon( {\rm cosh} \ v- {\rm cosh} \ w)]^{1/2} {\rm d} w =
\varepsilon \int^{ }_{{\cal I}_\varepsilon( v)}{g(w) {\rm sinh} \ w \over[
2\varepsilon( {\rm cosh} \ v- {\rm cosh} \ w)]^{1/2}} {\rm d} w, \eqno (A.23)
$$
On the other hand, we also have:
$$ \matrix{\displaystyle{ 1 \over {\rm sinh} \ v} { {\rm d} \over {\rm d} v}
\left[\varepsilon \int^{ }_{{\cal I}_\varepsilon( v)}{ {\rm d} g \over {\rm d}
w}(w)[2\varepsilon( {\rm cosh} \ v- {\rm cosh} \ w)]^{1/2} {\rm d} w \right] =
\cr\displaystyle = \int^{ }_{{\cal I}_\varepsilon( v)}{ {\rm d} g \over {\rm
d} w}(w)[2\varepsilon( {\rm cosh} \ v- {\rm cosh} \ w)]^{-1/2} {\rm d} w \cr}
\eqno (A.24) $$
Comparing Eqs.(A.23) and (A.24) then yields:
$$ {1 \over {\rm sinh} \ v} { {\rm d} \over {\rm d} v} \int^{ }_{{\cal
I}_\varepsilon( v)}{g(w) {\rm sinh} \ w \over[ 2\varepsilon( {\rm cosh} \ v-
{\rm cosh} \ w)]^{1/2}} {\rm d} w = \int^{ }_{{\cal I}_\varepsilon(
v)}{\displaystyle{ {\rm d} g \over {\rm d} w}(w) \over[ 2\varepsilon( {\rm
cosh} \ v- {\rm cosh} \ w)]^{1/2}} {\rm d} w \eqno (A.25) $$
or equivalently (in view of Eqs.(A.5) and (A.22)):
$$ \Xi^{( \varepsilon)}_ 3g(v) = \int^{ }_{{\cal I}_\varepsilon( v)}
\left[H^{(0)}_3g(w) \right][2\varepsilon( {\rm cosh} \ v- {\rm cosh} \
w)]^{-1/2} {\rm d} w\ , \eqno (A.25^\prime ) $$
which is Eq.(A.16) of lemma A.2, written for $ d=d_0=3 $ (with $ C_{d_0}=1). $
\par
For $ d\geq 5, $ one checks that the function $ g $ satisfies the
corresponding conditions i) or ii)
of lemma A.2; in view of the latter, Eq.(A.16) is therefore valid for the
considered value of $ d, $
with $ C_d $ given by Eq.(A.17), namely:
$$ C_d = \varepsilon^{{ d-3 \over 2}}[1\times 3
..\times( d-4)]^{-1} = (2\varepsilon)^{{ d-3 \over 2}}{\Gamma \left({d-1
\over 2} \right) \over \Gamma( d-2)}\ . $$
\par
Eqs.(A.20) and (A.21) are therefore established. \par
\medskip
{\bf Application to the exponential functions} \par
Let us introduce for $ \lambda \in \Bbb C $ the exponential functions
$$ g^-_\lambda( v) = {\rm e}^{- \left(\lambda +{d-2 \over 2} \right)v}\ \ ,\ \
\ g^+_\lambda( v) = {\rm e}^{ \left(\lambda +{d-2 \over 2} \right)v}\ ; $$
the latter are easily checked to be eigenfunctions of the corresponding
differential operators $ H^{(0)}_d, $
namely one has (in view of Eqs.(A.3), (A.3'), ...,(A.6))
$$ \eqalignno{ H^{(0)}_dg^+_\lambda( v) & = h^{(0)}_d(\lambda) g^+_\lambda( v)
& (A.26) \cr H^{(0)}_dg^-_\lambda( v) & = (-1)^dh^{(0)}_d(\lambda)
g^-_\lambda( v) & (A.27) \cr} $$
\par
For every $ \lambda $ such that $ {\rm Re} \ \lambda  > 0, $ the corresponding
function $ g^-_\lambda( v) $ obviously satisfies the
conditions ii) of theorem A.1 or A.2 (according to whether $ d $ is even or
odd). \par
Then, by taking into account Eq.(A.27), we obtain as an immediate application
of
formulae (A.19) and (A.21) of these theorems to the class of functions $
\left\{ g=g^-_\lambda \ ; {\rm Re} \ \lambda  >0 \right\} : $ \par
\smallskip
{\bf Proposition A.1 :} \par
For every $ \lambda $ such that $ {\rm Re} \ \lambda  > 0, $ the following
formulae hold: \par
\smallskip
i) $ \forall d $ even, $ d\geq 4: $
$$ \matrix{\displaystyle \left({1 \over {\rm sinh} \ v} { {\rm d} \over {\rm
d} v} \right)^{ \left({d-2 \over 2} \right)} \left( {\rm e}^{- \left(\lambda
+{d-2 \over 2} \right)v} \right)= \hfill \cr\displaystyle = {(-1)^{{d-2 \over
2}} \over 2^{{d-4 \over 2}} \left({d-4 \over 2} \right)!} \cdot
{h^{(0)}_d(\lambda) \over( {\rm sinh} \ v)^{d-3}} \cdot  \int^{ +\infty}_ v
{\rm e}^{- \left(\lambda +{d-2 \over 2} \right)w}\{ 2( {\rm cosh} \ w- {\rm
cosh} \ v)\}^{ \left({d-4 \over 2} \right)} {\rm d} w \hfill \cr} \eqno (A.28)
$$
\par
ii) $ \forall d $ odd, $ d\geq 3\ : $
$$ \matrix{\displaystyle \left({1 \over {\rm sinh} \ v} { {\rm d} \over {\rm
d} v} \right)^{ \left({d-1 \over 2} \right)} \int^{ +\infty}_ v{ {\rm e}^{-
\left(\lambda +{d-2 \over 2} \right)w} {\rm sinh} \ w \over\{ 2( {\rm cosh} \
w- {\rm cosh} \ v)\}^{ 1/2}} {\rm d} w = \hfill \cr\displaystyle = -
{(-2)^{{d-3 \over 2}} \over \Gamma( d-2)} \cdot  {h^{(0)}_d(\lambda) \over(
{\rm sinh} \ v)^{d-3}} \cdot  \int^{ +\infty}_ v {\rm e}^{- \left(\lambda
+{d-2 \over 2} \right)w}\{ 2( {\rm cosh} \ w- {\rm cosh} \ v)\}^{ \left({d-4
\over 2} \right)} {\rm d} w \hfill \cr} \eqno (A.29) $$
\par
\smallskip
{}From Eqs.(A.26), (A.27), one also obtains the following couple of relations:
$$ \eqalignno{ {\rm for\ } d {\bf \ } {\rm even\ } ,\  & H^{(0)}_d \left[
{\rm cosh} \left(\lambda +{d-2 \over 2} \right)v \right] = h^{(0)}_d(\lambda)
{\rm cosh} \left(\lambda +{d-2 \over 2} \right)v & (A.30) \cr {\rm for\ } d
{\bf \ } {\rm odd\ ,} \ {\rm \ } & H^{(0)}_d \left[ {\rm sinh} \left(\lambda
+{d-2 \over 2} \right)v \right] = h^{(0)}_d(\lambda) {\rm cosh} \left(\lambda
+{d-2 \over 2} \right)v & (A.31) \cr} $$
\par
For every $ \lambda $ in $ \Bbb C, $ the functions $ {\rm cosh} \left(\lambda
+ {d-2 \over 2} \right)v, $ with $ d $ even (resp. $ {\rm sinh} \left(\lambda
+{d-2 \over 2} \right)v, $
with $ d $ odd) satisfy the conditions i) of theorem A.1 (resp. theorem A.2).
Then, by taking into
account Eqs.(A.30), (A.31), we obtain as an immediate application of formulae
(A.18), (A.20) of
these theorems to the respective classes of functions $ \left\{ {\rm cosh}
\left(\lambda +{d-2 \over 2} \right)v\ ;\lambda \in \Bbb C,\ d\ {\rm even}
\right\} $ and $ \left\{ {\rm sinh} \left(\lambda +{d-2 \over 2} \right)v,\
\lambda \in \Bbb C,\ d\ {\rm odd} \right\} : $ \par
\smallskip
{\bf Proposition A.2 :} \par
For every $ \lambda $ in $ \Bbb C, $ the following formulae hold: \par
\smallskip
i) $ \forall d $ even, $ d\geq 4: $
$$ \matrix{\displaystyle \left({1 \over {\rm sinh} \ v} { {\rm d} \over {\rm
d} v} \right)^{ \left({d-2 \over 2} \right)} \left[ {\rm cosh} \left(\lambda
+{d-2 \over 2} \right)v \right] = \hfill \cr\displaystyle = {1 \over 2^{{d-4
\over 2}} \left({d-4 \over 2} \right)!} \cdot  {h^{(0)}_d(\lambda) \over(
{\rm sinh} \ v)^{d-3}} \int^ v_0 \left[ {\rm cosh} \left(\lambda +{d-2 \over
2} \right)w \right]\{ 2( {\rm cosh} \ v- {\rm cosh} \ w)\}^{ \left({d-4 \over
2} \right)} {\rm d} w \hfill \cr} \eqno (A.32) $$
\par
ii) $ \forall d $ odd, $ d \geq  3\ : $
$$ \matrix{\displaystyle \left({1 \over {\rm sinh} \ v} { {\rm d} \over {\rm
d} v} \right)^{ \left({d-1 \over 2} \right)} \int^ v_0{ {\rm sinh}
\left(\lambda +{d-2 \over 2} \right)w\ {\rm sinh} \ w \over\{ 2( {\rm cosh} \
v- {\rm cosh} \ w)\}^{ 1/2}} {\rm d} w = \hfill \cr\displaystyle = {2^{
\left({d-3 \over 2} \right)}\Gamma \left({d-1 \over 2} \right) \over \Gamma(
d-2)} \cdot  {h^{(0)}_d(\lambda) \over( {\rm sinh} \ v)^{d-3}} \int^ v_0
\left[ {\rm cosh} \left(\lambda +{d-2 \over 2} \right)w \right]\{ 2( {\rm
cosh} \ v- {\rm cosh} \ w)\}^{ \left({d-4 \over 2} \right)} {\rm d} w\ .
\hfill \cr} \eqno (A.33) $$
\par
\smallskip
{\bf Remark :} \par
The integrals at the r.h.s. of Eqs.(A.32) and (A.33) (resp. at the r.h.s. of
Eqs.(A.28)
and (A.29)) coincide, up to constant factors, with the first-kind (resp.
second-kind) Legendre
functions in dimension $ d, $ denoted by $ P^{(d)}_\lambda( {\rm cosh} \ v) $
(resp. $ Q^{(d)}_\lambda( {\rm cosh} \ v)), $ and the analytic
continuations of Eqs.(A.32), (A.33) (resp. Eqs.(A.28), (A.29)) provide useful
identities for
these analytic functions, which are exploited in \S 3.8. \par
\vfill\eject
\centerline{{\bf REFERENCES}}
\bigskip
\noindent \item {$\lbrack$B$\rbrack$}Bros J.: \lq\lq Complexified de Sitter
Space: Analytic Causal
Kernels and K\"allen-Lehmann-Type Representation\rq\rq\ in \lq\lq Recent
Advances in Field Theory\rq\rq , P. Bin\'etruy, G. Girardi, P. Sorba
Eds. ({\sl Nucl. Phys. B\/} (Proc. Suppl.) {\bf 18B} (1990)) p.22-28. \par
\noindent \item {$\lbrack$Bo$\rbrack$}Boas R.P.: \lq\lq Entire
Functions\rq\rq , Academic Press, New-York
(1954). \par
\noindent \item {$\lbrack$B.G.M.$\rbrack$}Bros J., Gazeau J.P. and Moschella
U.: \lq\lq Quantum Field
Theory on de Sitter Space\rq\rq , Saclay preprint (1993). \par
\noindent \item {$\lbrack$B,V-1$\rbrack$}Bros J. and Viano G.A.: \lq\lq
Connection Between the
Algebra of Kernels on the Sphere and the Volterra Algebra on the
One-Sheeted Hyperboloid: Holomorphic Perikernels\rq\rq , {\sl Bull. Soc.
Math. France\/}, vol.{\bf 120} (1992) p.169-225. \par
\noindent \item {$\lbrack$B,V-2$\rbrack$}Bros J. and Viano G.A.: \lq\lq
Complex Angular Momentum
Analysis in Axiomatic Quantum Field Theory\rq\rq\ in \lq\lq Rigorous Methods
in Particle Physics\rq\rq , S. Ciulli, F. Scheck, W. Thirring Eds.
({\sl Springer tracts in Mod. Phys.\/} vol.{\bf 119} (1990)) p.53-76. \par
\noindent \item {$\lbrack$E-1$\rbrack$}Erd\'elyi A., Magnus W., Oberhettinger
F. and Tricomi F.G.:
\lq\lq Tables of Integral Transforms\rq\rq , vol.II, McGraw-Hill, New-York
(1954). \par
\noindent \item {$\lbrack$E-2$\rbrack$}Erd\'elyi A., Magnus W., Oberhettinger
F. and Tricomi
F.G.: \lq\lq Higher Transcendental Functions\rq\rq , vol.I, McGraw-Hill,
New-York (1953). \par
\noindent \item {$\lbrack$F$\rbrack$}Froissart M.: \lq\lq Asymptotic Behaviour
and Subtractions in the
Mandelstam Representation\rq\rq , {\sl Phys. Rev.\/} vol.{\bf 123},
N\nobreak$^\circ$3 (1961)
p.1053-1057. \par
\noindent \item {$\lbrack$Fa-1$\rbrack$}Faraut J.: \lq\lq Alg\`ebre de
Volterra et Transformation de
Laplace Sph\'erique\rq\rq , S\'eminaire d'Analyse Harmonique de Tunis
(1981) (Expos\'e 29). \par
\noindent \item {$\lbrack$Fa-2$\rbrack$}Faraut J.: \lq\lq Alg\`ebres de
Volterra et Transformation de
Laplace Sph\'erique sur Certains Espaces Sym\'etriques Ordonn\'es\rq\rq ,
{\sl Symposia Mathematica\/} vol.XXIX (1984) p.183-196. \par
\noindent \item {$\lbrack$Fa-3$\rbrack$}Faraut J.: \lq\lq Analyse Harmonique
et Fonctions Sp\'eciales\rq\rq ,
Ecole d'\'et\'e d'Analyse Harmonique de Tunis (1984). \par
\noindent \item {$\lbrack$Fu$\rbrack$}Furstenberg H.: \lq\lq
Translation-Invariant Cones of Functions
on Semi-Simple Lie Groups\rq\rq , {\sl Bulletin American Math. Soc.\/}
vol.{\bf 71}
(1965) p.271-326. \par
\noindent \item {$\lbrack$Fa,V$\rbrack$}Faraut J. and Viano G.A.: \lq\lq
Volterra Algebra and the
Bethe-Salpeter Equation\rq\rq , {\sl Journal of Math. Phys.\/} {\bf 27, }3
(1986)
p.840-848. \par
\noindent \item {$\lbrack$G$\rbrack$}Gribov V.N.: \lq\lq Partial Waves with
Complex Angular Momenta
and the Asymptotic Behaviour of the Scattering Amplitude\rq\rq , {\sl J. Exp.
Theor. Phys.\/} vol.{\bf 14} (1962) p.1395. \par
\noindent \item {$\lbrack$H$\rbrack$}Helgason S.: \lq\lq Differential Geometry
and Symmetric Spaces\rq\rq ,
Academic Press, New-York (1962). \par
\noindent \item {$\lbrack$K$\rbrack$}Khuri N.M.: \lq\lq Regge Poles, Power
Series and a
Crossing-Symmetric Watson-Sommerfeld Transformation\rq\rq , {\sl Phys. Rev.\/}
vol.{\bf 132} (1963) p.914-926. \par
\noindent \item {$\lbrack$L$\rbrack$}Le Roy E.: \lq\lq Sur les S\'eries
Divergentes et les Fonctions
D\'efinies par un D\'eveloppement de Taylor\rq\rq , {\sl Ann. Facult\'e des
Sci.
Univ. de Toulouse\/}, vol.{\bf 2} (1900) p.317-430. \par
\noindent \item {$\lbrack$M$\rbrack$}Martin A.: \lq\lq Remarks on the
Prolongation of a Scattering
Amplitude in the Complex Angular Momentum Plane\rq\rq , {\sl Phys. Letters\/}
{\bf 1}, 72 (1962) p.72-75. \par
\noindent \item {$\lbrack$P$\rbrack$}Poincar\'e H.: Oeuvres, vol.{\bf 10}.
p.94-203 (section 10). \par
\noindent \item {$\lbrack$Pa$\rbrack$}Painlev\'e P.: \lq\lq Sur les Lignes
Singuli\`eres des Fonctions
Analytiques\rq\rq , {\sl Ann. Fac. Toulouse\/} {\bf 2} (1888) p.27. \par
\noindent \item {$\lbrack$Ph$\rbrack$}Pham F.: \lq\lq Introduction \`a l'Etude
Topologique des
Singularit\'es de Landau\rq\rq , {\sl M\'emorial des Sciences
Math\'ematiques\/} {\bf 164},
Gauthier-Villars, Paris (1967). \par
\noindent \item {$\lbrack$R$\rbrack$}Regge T.: \lq\lq Introduction to Complex
Orbital Momenta\rq\rq , {\sl Nuovo
Cimento\/} vol.{\bf 14} (1959) 951-976. \par
\noindent \item {$\lbrack$Ru$\rbrack$}Rudin W.: \lq\lq Fourier Analysis on
Groups\rq\rq , Interscience
Publishers, New-York (1962). \par
\noindent \item {$\lbrack$So$\rbrack$}Sommerfeld A.: \lq\lq Partial
Differential Equations in
Physics\rq\rq , New-York (1949). \par
\noindent \item {$\lbrack$S,W$\rbrack$}Stein E.M. and Wainger S.: \lq\lq
Analytic Properties of
Expansions and Some Variants of Parseval-Plancherel Formula\rq\rq ,
{\sl Arkiv f\"ur Math.\/} vol.{\bf 5} (1963) p.553-567. \par
\noindent \item {$\lbrack$Vi$\rbrack$}Vilenkin N. Ja.: \lq\lq Special
Functions and the Theory of
Group Representations\rq\rq , {\sl Amer. Math. Soc. Transl.\/} vol.{\bf 22},
Providence R.I. (1968). \par
\noindent \item {$\lbrack$Wa$\rbrack$}Watson G.N.: \lq\lq The Diffraction of
Electric Waves by the
Earth\rq\rq , {\sl Proc. Roy. Soc.,\/} London vol.{\bf 95} (1918) p.83-99.
\par
\bye
\listrefs
\draftend
\end